\documentclass{amsart}
\usepackage{amssymb, amsmath}
\usepackage{bbm}
\usepackage{calrsfs}
\usepackage{mathrsfs}

\font\fndef=cmbxti12
\font\sps=cmti5

\def\largerightarrow{   -\negthinspace\negthinspace -\negthinspace\negthinspace
                                       \negthinspace\longrightarrow }

\begin{document}
\title{\Large Regular Algebraic K-Theory for groups -- Part II}
\author{U. Haag}
\date{\today}

\maketitle
\Large
\bigskip
\par\noindent
\hfil\hfil \Large{\bf{Section V -- $K^J$-theory.}} \hfil
\bigskip\bigskip
\par\noindent
In this section we introduce the basic features of group $K$-theory (or rather $K^J$-theory ). $K^J$-theory (as opposed to regular $K$-theory) is not in general exact (in low dimensions), but there are some conditions ensuring exactness in all dimensions. These are subject of the Mapping Cone Theorem (Theorem 5) and its Corollary 5.1. We have already seen the definition of the groups $\, K^J_2 ( N , F )\, $ in section 1. There is the problem of defining higher dimensional 
$K^J$-groups which are naturally affiliated with the $K^J_2$-groups in a certain sense which needs further explanation. In particular one wants the higher groups to fit into a long exact sequence associated with a group extension which entails the existence of boundary maps decreasing the dimension by one when passing from the quotient to the kernel of the extension. The importance of  this theory may be seen from the simplicity and naturality of its definition (also facilitating its computation). The regular theory depends heavily on the results of this section and is thoroughly discussed in section 6 following this one. Let us begin by introducing some basic notions. Although the regular theory is not developed in this chapter we will provide some basic concepts below for convenience of the reader.
\par\bigskip\noindent
{\bf Definition 10. (a)} The pair $\, ( N, F )\, $ will be called {\it admissible} iff the normal subgroup $\, N\, $ admits a semiuniversal $F$-central extension satisfying part (i) of Definition 2(c), such that the central kernel of this extension $\, {\overline N}^{ _{\sps F}}\, $ is contained in $\, [\, {\overline N}^{ _{\sps F}},\, F\, ]\, $ ( this implies that it is equal to $\, K^J_2 ( N , F )\, $). The subgroup of the central kernel contained in the k-fold commutator subgroup 
$\, [\, [\cdots [\, {\overline N}^{ _{\sps F}},\, F\, ],\cdots ], F\, ]\, $ will be referred to as the 
{\it k-regular part} of $\, K^J_2 ( N , F )\, $, and the intersection of these subgroups for all k will be called the {\it regular part}. Even if $\, ( N, F )\, $ is not admissible and 
$\, {\overline N}^{ _{\sps F}}\, $ is not defined one can define the k-regular parts by
\smallskip 
$$ \bigl(\, J_{N, F}\cap \bigl[\, \bigl[\cdots \bigl[\, U_{N, F},\, U_F\, \bigr],\cdots \bigr], U_F\, \bigr]\, \bigr)\, /\, 
\bigl(\, \bigl[\, J_{N, F},\, U_F\, \bigr]\, +\, 
\bigl[\, U_{N, F},\, J_F\, \bigr]\, \bigr) . $$
\par\medskip\noindent
We will also consider elements of $\, K^J_2 ( N , F )\, $ which become ($k$-)regular in some quotient, most often in the abelianization of $\, N\, $ ( dividing $\, ( N , F )\, $ by $\, [ N , N ]\, $). These are referred to as {\it abelian k-regular elements}.
If a morphism of normal pairs $\, ( N, F ) \rightarrow ( M, G )\, $ induces an isomorphism (surjection, injection) such that the inverse image of a k-regular element contains a k-regular element for each k,
 it will be called a {\it regular isomorphism (surjection, injection)}. For an admissible pair one always gets an almost canonical map $\, U_{N , F}\buildrel\alpha\over\rightarrow {\overline N}^{ _{\sps F}}\, $ which is uniquely determined on commutators $\, [\, U_{N , F}, U_F\, ]\, $. Let $\, R_{N , F}\, $ be its kernel. Given a morphism $\, ( N , F )\buildrel\varphi\over\rightarrow ( M , G )\, $ of admissible pairs and almost canonical maps
$$\, U_{N , F}\buildrel {\alpha_N}\over\rightarrow {\overline N}^{ _{\sps F}}\, , \quad U_{M , G}\buildrel {\alpha_M}\over\rightarrow {\overline M}^{ _{\sps G}}\, ,$$ 
$\, \alpha_N\, $ (resp. $\, R_{N , F}\, $) is said to be {\fndef restricted} from $\,\alpha_M\, $ (resp. 
$\, R_{M , G}\, $) if there exists an equivariant lift 
$\, {\overline N}^{ _{\sps F}}\buildrel {\overline\varphi }\over\rightarrow {\overline M}^{ _{\sps G}}\, $ of 
$\,\varphi\, $ giving a commutative
\par\noindent
diagram
\smallskip
$$ \vbox{\halign{ #&#&#&#&#\cr
     $ R_{N , F} $ & $\largerightarrow $ & $ U_{N , F} $ & $\buildrel {\alpha }_N\over\largerightarrow $ &  
     ${\overline N}^{ _{\sps F}} $ \cr
     \hfil $\Biggm\downarrow $\hfil && $\;\;\Biggm\downarrow\negthinspace u_{\varphi }$ && 
     $\,\Biggm\downarrow {\overline\varphi }$ \cr
     $ R_{M , G} $ & $\largerightarrow $ & $ U_{M , G} $ & $\buildrel {\alpha }_M\over\largerightarrow $&
     $ {\overline M}^{ _{\sps G}} $ \cr  } }.  $$
\par\smallskip\noindent
In this case we also say that $\, {\alpha }_M\, $ (resp. $\, R_{M , G}\, $ ) is {\fndef induced } from 
$\, {\alpha }_N\, $ (resp. $\, R_{N , F}\, $ ). 
\par\medskip\noindent
{\bf (b)}\quad For any normal pair $\, ( N , F )\, $ the {\fndef cone} over $\, ( N , F )\, $ is the pair 
$\, ( U_{N , F} , U_F )\, $. The {\fndef J-suspension } of $\, ( N , F )\, $ is the pair $\, ( J_{N , F} , U_F )\, $ (so that $\, K^J_3 ( N , F )\, =\, K^J_2 ( J_{N , F}\, ,\, U_F )\, $, compare with Definition 9). Let  
$\, \varphi : ( N , F )\rightarrow ( M , G )\, $ be a morphism of normal pairs. Define the 
{\fndef mapping cone } of $\,\varphi\, $ to be the pull back 
$\, ( C_\varphi , U_\varphi ) = ( N {\times }_\varphi  U_{M , G} , F {\times }_\varphi  U_G )\, $ of the cone extension over $\, ( M , G )\, $ by $\,\varphi\, $. 
\par\bigskip\noindent
{\bf Lemma 9.}\quad Assume given an exact sequence
$$ 1\largerightarrow ( N' , F ) \largerightarrow ( N , F ) \largerightarrow ( N'' , E ) \largerightarrow 1 $$
of normal pairs ( i.e. $\, ( N'' , E ) = ( N /  N' , F /  N' )\, $). Then
$$ K^J_2 ( N' , F )\largerightarrow K^J_2 ( N , F )\largerightarrow 
K^J_2 ( N'' , E ) $$
is halfexact and the same holds for the corresponding sequence of extended $K^J_2$-functors. Further if $\, K^J_2 ( N , F )\longrightarrow K^J_2 ( N'' , E )\, $ is surjective and $\, ( N' , F )\, $ and $\, ( N'' , E )\, $ are admissible, then also $\, ( N , F )\, $ is admissible.
\par\bigskip\noindent
{\it Proof.} Assume that an element of $\, J_{N , F}\cap [ U_{N , F} , U_F ]\, $ is mapped to an element of 
$\, [ J_{N'' , E} , U_E ] + [ U_{N'' , E} , J_E ]\, $. This element can be lifted to a trivial element in 
$\, [ J_{N , F} , U_F ] + [ U_{N , F} , J_F ]\, $, so that the difference of the original element with the trivial one maps to zero in $\, J_{N'' , E}\,$ and hence comes from 
$\, J_{N , F}\cap [ U_{N , F} , U_F ]\cap U_{N' , F} = J_{N' , F}\cap [ U_{N' , F} , U_F ]\, $. The argument in the extended case is much the same. Suppose that 
$\, K^J_2 ( N , F )\rightarrow K^J_2 ( N'' , E )\, $ is surjective. Let 
$\, U_{N , F}\, /\,\sim\, $ denote the quotient of $\, U_{N , F}\, $ by $\, J_{N , F} \cap  [ U_{N , F}\, ,\, U_F ]\, $ and
$\, U_{N'' , E}\, /\,\sim\, $ the quotient of $\, U_{N'' , E}\, $ by $\, J_{N'' , E} \cap [ U_{N'' , E} , U_E ]\, $. Then one gets an exact sequence
$$ 1\largerightarrow {U_{N' , F}\, /\,\sim }\largerightarrow {U_{N ; F}\, /\,\sim }\largerightarrow 
{U_{N'' , E}\, /\,\sim }\largerightarrow 1 .$$
If $\, ( N' , F )\, $ and $\, ( N'' , E )\, $ are admissible, $\, U_{N' , F}\, /\,\sim\, $ contains a $F$-normal copy of 
$\, N'\, $ and $\, U_{N'' , E}\, /\,\sim\, $ a $E$-normal copy of $\, N''\,$. Dividing by the copy of $\, N'\,$ one gets a splitting of the quotient $\, U_{N'' , E}\, /\,\sim\, $ by the universal property of the pair 
$\, ( U_{N'' , E} , U_E )\, $ (and the fact that $\, [ U_{N' , F} , U_F ]\, $ is divided out). Then the preimage of the lifted copy of $\, N''\, $ with respect to the copy of $\, N'\, $ determines a $F$-normal copy of $\, N\,$ in $\, U_{N , F}\, /\,\sim\, $, so 
$\, ( N , F )\, $ is admissible\qed
\par\bigskip\noindent 
{\bf Remark.}\quad
The assumption of the Lemma is fulfilled if either $\, ( N' , F )\, $ or $\, ( N'' , E )\, $ is a cone and the other is admissible.
\par\bigskip\noindent
{\bf Lemma 10.}\quad
Assume that $\, ( N , F ) \buildrel\varphi\over\longrightarrow ( M , G )\, $ is a morphism with $\, ( N , F )\, $ admissible, such that the map $\, K^J_2 ( N , F )\rightarrow K^J_2 ( M , G )\, $ is surjective. Then any almost canonical map 
$\, \beta : U_{M , G} \longrightarrow \overline M\, $ has a restriction 
$\, \alpha : U_{N , F} \longrightarrow {\overline N}^{ _{\sps F}}\, $. If, on the other hand 
$\, \varphi \, $ is injective then any almost canonical map 
$\, \alpha : U_{N , F} \longrightarrow {\overline N}^{ _{\sps F}}\, $ can be induced to a map 
$\, \beta : U_{M , G} \longrightarrow \overline M\, $, where $\,\overline M\, $ denotes some minimal 
almost canonical $G$-central extension of $\, M\, $ by $\, K^J_2 ( M , G )\, $  (equal to 
$\, {\overline M}^{ _{\sps G}}\, $ in case that $\, ( M , G )\, $ is admissible).
\par\bigskip\noindent
{\it Proof.} Let $\, K^J_2 ( N , F )\rightarrow K^J_2 ( M , G )\, $ be surjective and choose some almost canonical map $\, \beta : U_{M , G}\longrightarrow \overline M\, $. Consider the map
\smallskip
$$ { U_{N , F}\over [\, J_{N , F} , U_F\, ] + [\, U_{N , F} , J_F\, ]}\buildrel {\overline u_\varphi }\over\largerightarrow { U_{M , G}\over [\, J_{M , G} , U_G\, ] + [\, U_{M , G} , J_G\, ] }   $$
\par\smallskip\noindent
induced by $\, U_{N , F}\buildrel u_\varphi\over\longrightarrow U_{M , G}\, $. Choose 
$\, {\overline N}^{ _{\sps F}}\, $ as a minimal $F$-normal subextension on the left hand side, and define
$\, \overline\varphi : {\overline N}^{ _{\sps F}}\longrightarrow \overline M\, $ by the composition
\smallskip 
$$ {\overline N}^{ _{\sps F}}\subseteq { U_{N , F}\over [\, J_{N , F} , U_F\, ] + [\, U_{N , F} , J_F\, ] }\buildrel
     \beta\circ {\overline u_\varphi }\over\largerightarrow \overline M  . $$
\par\smallskip\noindent     
By surjectivity of $\, K^J_2 ( N , F ) \rightarrow K^J_2 ( M , G )\, $ the image of 
$\, {\overline N}^{ _{\sps F}}\rightarrow \overline M\, $ is the same as the image of 
$\, U_{N , F}\buildrel\beta\circ {\overline u_\varphi }\over\longrightarrow \overline M\, $. Then by the universal property of the pair $\, ( U_{N , F} , U_F )\, $ the diagonal map 
$\,\beta\circ {\overline u_\varphi }\, $ can be lifted to  
$\, U_{N , F}\buildrel\alpha\over\longrightarrow {\overline N}^{ _{\sps F}}\, $ making a commutative diagram
\smallskip
$$ \vbox{\halign{#&#&#\cr
 $U_{N , F}$ & $\buildrel\alpha\over\longrightarrow $\hfil & \hfil ${\overline N}^{ _{\sps F}}$\hfil \cr
 \hfil $\beta\negthinspace\circ\negthinspace {\overline u_\varphi }\negthinspace\negthinspace\negthinspace $ & 
 $\searrow $ & \hfil $\Bigm\downarrow\negthinspace {\overline\varphi }$\hfil \cr
                                                                                    && \hfil $ {\overline M}^{ _{\sps G}} $\hfil \cr}} . $$
\par\smallskip\noindent 
Now let $\, ( N , F )\subseteq ( M , G )\, $ be an inclusion. Put $\, M_F = M \cap F\, $. Then the inclusion factors into $\, ( N , F ) \subseteq (M_F , F )\, $ and $\, ( M_F , F ) \subseteq ( M , G )\, $. One notes the following facts. Any almost canonical map
$\, U_{N , F}\longrightarrow {\overline N}\, $ is obtained in the following way: it must factor over the canonical $F$-central extension of $\, N\, $, so one can divide $\, U_{N , F}\, $ by the subgroup 
$\, [\, J_{N , F}\, ,\, U_F\, ]\, +\, [\, U_{N , F}\, ,\, J_F\, ]\, $. Then one may lift the free abelian group 
$\, J_{N , F}\, /\, (\, J_{N , F}\,\cap\, [\, U_{N , F}\, ,\, U_F\, ]\, )\, $ normally to the image of $\, J_{N , F}\, $ in the canonical $F$-central extension and divide by this subgroup to obtain a minimal $F$-central extension of $\, N\, $ with kernel $\, K^J_2 ( N , F )\, $. In case that $\, ( N , F )\, $ is admissible this extension is unique up to equivalence and equal to the semiuniversal $F$-central extension. Now consider the commutative diagram
\smallskip
$$\vbox{\halign{#&#&#&#&#&#&#&#&#\cr
$1$ & $\longrightarrow$ & $ {J_{N , F} \over J_{N , F}\,\cap\, \bigl[\, U_{N , F}\, ,\, U_F\, \bigr]}$ & 
$\longrightarrow$ & ${U_{N , F} \over \bigl[\, U_{N , F}\, ,\, U_F\,\bigr]}$ & $\longrightarrow$ & 
${ N \over \bigl[ N , F \bigr]}$ & $\longrightarrow$ & $1$ \cr
&& \hfil$\Bigm\downarrow$\hfil && \hfil$\Bigm\downarrow$\hfil && \hfil$\Bigm\downarrow$\hfil &&\cr 
$1$ & $\longrightarrow$ & $ {J_{M_F , F} \over J_{M_F , F}\,\cap\, \bigl[\, U_{M_F , F}\, ,\, U_F\, \bigr]}$ & 
$\longrightarrow$ & ${U_{M_F , F} \over \bigl[\, U_{M_F , F}\, ,\, U_F\,\bigr]}$ & $\longrightarrow$ & 
${ M_F \over \bigl[ M_F , F \bigr]}$ & $\longrightarrow$ & $1$ \cr
&& \hfil$\Bigm\downarrow$\hfil && \hfil$\Bigm\downarrow$\hfil && \hfil$\Bigm\downarrow$\hfil &&\cr 
$1$ & $\longrightarrow$ & ${J_{M , G}\over J_{M , G}\,\cap\,\bigl[\, U_{M , G}\, ,\, U_G\,\bigr]}$ & 
$\longrightarrow$ & ${ U_{M , G} \over \bigl[\, U_{M , G}\, ,\, U_G\,\bigr]}$ & $\longrightarrow$ &
$ { M \over \bigl[ M , G \bigr]}$ & $\longrightarrow$ & $1$  \cr }}\> . $$
\par\smallskip\noindent
Both vertical maps in the middle are split injective, the first one by the exact sequence 
\smallskip
$$ 1 \longrightarrow {U_{N , F} \over \bigl[\, U_{N , F}\, ,\, U_F\,\bigr]} \longrightarrow 
{U_{M_F , F} \over \bigl[\, U_{M_F , F}\, ,\, U_F\,\bigr]} \longrightarrow {U_{M_F / N , F / N} \over 
\bigl[\, U_{M_F / N , F / N}\, ,\, U_{F / N}\,\bigr]} \longrightarrow 1 $$
\par\smallskip\noindent
of free abelian groups, and the second since the core of $\, U_{M_F , F}\, $ for $\, U_F\, $ is naturally contained in the core of $\, U_{M , G}\, $ for $\, U_G\, $,  Then their composition is split injective and the composition of the first two vertical maps is injective. Put 
$\, B_0\, =\, {U_{N , F} \over \bigl[\, U_{N , F}\, ,\, U_F\,\bigr]}\, ,\, 
B\, =\, { U_{M , G} \over \bigl[\, U_{M , G}\, ,\, U_G\,\bigr]}$ and let $\, B_1\, $ denote a complement in 
$\, B\, $ for $\, B_0\, $. Put $\, A\, =\, {J_{M , G}\over J_{M , G}\,\cap\,\bigl[\, U_{M , G}\, ,\, U_G\,\bigr]}\, $ and $\, A_0\, =\,  {J_{N , F} \over J_{N , F}\,\cap\, \bigl[\, U_{N , F}\, ,\, U_F\, \bigr]}\, $ and let $\, A_1\, $ be the image of $\, A\, $ in $\, B_1\, $ which is a free abelian group. The kernel of the map 
$\, A \longrightarrow A_1\, $ is an extension of $\, ( N \cap [ M , G ] )\, /\, [ N , F]\, $ by $\, A_0\, $. From the assumption that $\, ( N , F )\, $ is admissible one gets that $\, N / [ N , F ]\, $ is free abelian (compare the Remark after Definition 11). Then also the quotient $\, ( N \cap [ M , G ] )\, /\, [ N , F]\, $ is free abelian and the map 
$\, A_0 \longrightarrow A\, $ splits, so that any almost canonical map 
$\, U_{N , F}\longrightarrow {\overline N}^F\, $ can be induced to 
$\, U_{M , G}\longrightarrow \overline M\, $ \qed
\par\bigskip\noindent
{\bf Remark.}\quad The assumptions of the preceding Lemma may not be the most general to ensure the existence of induction or restriction maps. However some assumptions must be made. In general, there will be no pair of almost canonical maps 
$\, \alpha : U_{N , F}\rightarrow {\overline N}^{ _{\sps F}}\, $ and
$\,\beta : U_{M , G}\rightarrow {\overline M}^{ _{\sps G}}\, $ such that one gets an induced map
$\, u_\varphi : R_{N , F} \rightarrow R_{M , G}\, $ on the kernels. Things become different restricting to
$\, u_\varphi : U_{0 . N , F} \rightarrow U_{0 , M , G}\, $ where $\, U_{0 , N , F} = [\, U_{N , F}\, ,\, U_F\, ]\, $ considering the induced map $\, R_{0 , N , F} \rightarrow R_{0 , M , G}\, $. One gets a natural commutative diagram
\smallskip
$$ \vbox{\halign{ #&#&#&#&#&#&#&#&#\cr
$1$ & $\largerightarrow $ & $R_{0 , N , F}$ & $\largerightarrow $ & $U_{0 , N , F}$ &
$\buildrel\alpha\over\largerightarrow $ & $ [\, {\overline N}^{ _{\sps F}}\, ,\, F\, ] $ & 
$\largerightarrow $ & $1$ \cr
&&\hfil $\Biggm\downarrow $\hfil && $ u_\varphi\negthinspace\negthinspace\Biggm\downarrow $\hfil && 
$\quad {\overline\varphi }\negthinspace\negthinspace\Biggm\downarrow\quad $ &&\cr
$1$ & $\largerightarrow $ & $R_{0 , M , G}$ & $\largerightarrow $ & $U_{0 , M , G}$ &
$\buildrel\beta\over\largerightarrow $ & $[\, {\overline M}^{ _{\sps G}}\, ,\, G\, ]$ &
$\largerightarrow $ & $1$  \cr  }} $$
\par\smallskip\noindent
with $\, {\overline\varphi } , \alpha , \beta\, $ uniquely determined. 
\par\bigskip\noindent
{\bf Lemma 11.}\quad Let $\, N\, $ and $\, N'\, $ be two complementary normal subgroups of a group 
$\, F\, $ and put $\, E = F / N'\, ,\, E' = F / N\, $. Then the pair $\, ( N , F )\, $ is admissible if and only if the pair $\, ( N , E )\, $ is admissible, and the natural map $\, ( N , F ) \rightarrow ( N , E )\, $ induces a regular surjection $\, K^J_2 ( N , F ) \twoheadrightarrow K^J_2 ( N , E )\, $. The kernel of this map is naturally a quotient of the kernel of 
$\, \widetilde K^J_2 ( N' , F )\,\twoheadrightarrow\, \widetilde K^J_2 ( N' , E' )\, $ and if $\, N'\, $ has a weak core for $\, E'\, $ a quotient of the kernel of 
$\, K^J_2 ( N' , F )\,\twoheadrightarrow\, K^J_2 ( N' , E' )\, $.
\par\bigskip\noindent
{\it Proof. }\quad $\, U_{N , F}\, $ differs from $\, U_{N , E}\, $ by 
$\, U_{N , F} \cap U_{N' , F} = J_{N , F}\cap J_{N' , F}\, $. Construct a splitting 
$\, \sigma : ( U_{N , E} , U_E ) \rightarrow ( U_{N , F} , U_F )\, $. The quotient map induces a surjection 
$\, J_{N , F} \cap [\, U_{N , F}\, ,\, U_F\, ] \rightarrow J_{N , E} \cap [\, U_{N , E}\, ,\, U_E\, ]\, $ because 
$\, \sigma ( J_{N , E} \cap [\, U_{N , E}\, ,\, U_E\, ] ) \subseteq J_{N , F} \cap [\, U_{N , F}\, ,\, U_F\, ]\, $. So 
$\, ( N , E )\, $ is admissible if and only if $\, ( N , F )\, $ is admissible and the map 
$\, K^J_2 ( N , F )\rightarrow K^J_2 ( N , E )\, $ is surjective. By the same reasoning $\,\sigma\, $ maps the subgroup $\, J_{N , E} \cap [\, [\cdots [\, U_{N , E}\, ,\, U_E\, ] ,\cdots ], U_E\, ]\, $ of k-fold commutators into the corresponding subgroup 
$\, J_{N , F} \cap [\, [\cdots [\, U_{N , F}\, ,\, U_F\, ],\cdots ], U_F\, ]\, $ so that the map 
$\, K^J_2 ( N , F )\rightarrow K^J_2 ( N , E )\, $ is surjective also on the k-regular parts. Suppose given an element in the kernel of this map. One readily checks that it can be represented by an element of
\smallskip 
$$ J_{N , F} \cap [\, U_{N , F}\, ,\, U_F\, ] \cap U_{N' , F} = J_{N' , F} \cap [\, U_{N' , F}\, ,\, U_F\, ] \cap 
U_{N , F}  $$
\par\smallskip\noindent 
so that the same representative also gives a representative for the kernel of  
$\, K^J_2 ( N' , F )\rightarrow K^J_2 ( N' , E' )\, $. To obtain a map on the level of 
$K^J$-theory one must ensure that trivial elements in one group correspond to trivial elements in the other. This is not always the case and the question of when it is is a subtle and difficult problem which is closely linked to the exactness properties of the $K^J$-groups. One does have the following partial results.  Assume an element is in $\, [\, J_{N' , F}\, ,\, U_F\, ] \cap J_{N , F}\, $. Then it becomes an elementary relation in 
$\, [\, J_{N' , E'}\, ,\, U_{E'}\, ]\, $ and lifting this relation by $\,\sigma\, $ one finds that the element is in
$\, [\, J_{N' , F}\, ,\, U_{N , F}\, ]\, +\, [\, J_{N , F}\, ,\, U_F\, ]\, $. If on the other hand an element is in
$\, J_{N , F} \cap [\, U_{N' , F}\, ,\, J_F\, ]\, $ it becomes elementary in 
$\, [\, U_{N' , E'}\, ,\, J_{E'}\, ] \subseteq U_F / J_{N , F}\, $, so by lifting with respect to a basis of 
$\, U_{N' , E'} + J_{E'}\, $ one gets that the element is in $\, [\, J_{N , F}\, ,\, J_F + U_{N' , F}\, ]\, $. In particular one gets a well defined (and surjective) map from the kernel of 
$\, \widetilde K^J_2 ( N' , F )\,\twoheadrightarrow\, \widetilde K^J_2 ( N' , E' )\, $ to the kernel of 
$\, K^J_2 ( N , F ) \twoheadrightarrow K^J_2 ( N , E )\, $. Now assume that $\, N'\, $ has a weak core for 
$\, E'\, $. Consider the subgroup $\, [\, U_{N' , F}\, ,\, J_F\, ]\, $ of $\, J_{N' , F}\, $ modulo 
$\, [\, J_{N' , F}\, ,\, U_F\, ]\, $. Denote $\, C_{N'}\, $ the subgroup of $\, N' / [ N' , N' ]\, $ generated by the elements $\,\{ e_i \}\, $ of Definition 8. Then the commutator subgroup above maps to a quotient of 
$\, C_{N'} \otimes J_E^{ab}\, $ which is easily checked from looking at the usual basis of 
$\, [\, U_{N'}\, ,\, \lambda (J_E)\, ]\, $ and the corresponding relations for the tensor product which lie inside $\, [\, J_{N' , F}\, ,\, U_F\, ]\, $. In fact a basis of the $B_3'$-reduced $\, U_{F / [ N' , N' ]}\, $ is given by the union of the $\,\{ u_{e_i} \}\, $, $B_2$-elements involving commutators with the $\, u_{e_i}\, $ and the $\,\{ u_z = u_{s (x)} \}\, $ and the subgroup $\, [\, J_{N' , F }\, ,\, U_F\, ]\,$ comprises modulo 
$\, [\, J_{N' , F}\, ,\, J_F\, ]\, $ of elements $\, [\, B_{1,2}\, ,\, u_z\, ]\, $ and $\, [\, B_{1 , 2}\, ,\, u_c\, ]\, $ whereas our commutator subgroup comprises of elements $\, [\, u_{e_i}\, ,\, u_x\, u_y\, u^{-1}_{xy}\, ]\, $ which serve together with the corresponding relations of the form 
$\, [\, u_x\, ,\, B^{y , e_i}_2\, ]\, $ to express an element of the form $\, B^{y , x e_i x^{-1}}\, $ corresponding to a uniquely determined element $\, [\, u_{e_i}\, ,\, u_x\, u_y\, u^{-1}_{xy}\, ]\, $ by basis elements from 
$\, \{ B^{z , e_i}_2 \}\, $. In particular no combination of such elementary commutators will intersect with 
the image of $\, [\, J_{N' , F}\, ,\, U_F\, ]\, $ unless it is trivial in the tensor product. This shows that the quotient of $\, [\, U_{N' , F}\, ,\, J_F\, ]\, $ modulo $\, [\, J_{N' , F}\, ,\, U_F\, ]\, $ is precisely equal to the tensor product. Then dividing by the image of $\, U_{N , F}\, $ it maps to 
$\, C_{N'} \otimes J_{E / N}^{ab}\, $ and the kernel of this map which is 
$\, C_{N'} \otimes J_{N , E}^{ab}\, $ clearly consists of (images of) trivial elements for 
$\, K^J_2 ( N , F )\, $ so that putting things together elements which are trivial for 
$\, K^J_2 ( N' , F )\, $ are also trivial for $\, K^J_2 ( N , F )\, $ and the result follows     \qed
\par\bigskip\noindent
{\bf Lemma 12.}\quad In the situation of Lemma 11 assume that $\, \overline N = N / [ N , N ]\, $ is free abelian and the action of $\, F\, $ on $\, \overline N' = N' / [ N' , N' ]\, $ reduces to the action of a quotient 
$\, U\, $ containing $\, N\, $ which has a lift to $\, F\, $ such that $\, N\, $ has a weak core for $\, U\, $. Then the kernel of the map 
$\, K^J_2 ( N , F ) \rightarrow K^J_2 ( N , E )\, $ is equal to the tensor product of the core of 
$\,\overline N\, $ for $\, U\, $ with $\,\overline N'\, $. In particular it is free abelian whenever 
$\, \overline N'\, $ is free abelian.
\par\bigskip\noindent
{\it Proof.}\quad Put $\,\overline F = F / ([ N , N ] + [ N' , N' ])\, ,\, \overline E = E / [ N , N ]\, ,$ Since any element of the kernel of $\, K^J_2 ( N , F ) \rightarrow K^J_2 ( N , E )\, $ is represented by an element in 
$\, [\, U_N\, ,\, U_{N'}\, ]\, $ we may reduce $\, U_F\, $ to 
$\, U_{\overline F}\, $ modulo trivial elements and then to the subgroup generated by 
$\, U_{\overline N}\, ,\, U_{\overline N'}\, $ and a semicanonical lift of $\, U_{E / N}\, $. Denote elements in $\, \overline N\, $ by letters 
$\, c , d\, $ etc., elements of $\, \overline N'\, $ by greek letters $\,\gamma ,\delta\, $ etc. and elements in the image of a section $\, E / N \nearrow F\, $ by letters $\, x, y, z\, $ etc. One easily finds a basis of 
$\, J_{\overline N , \overline F} \cap J_{\overline N' , \overline F} \cap 
[\, U_{\overline F}\, , U_{\overline F}\, ]\, $ modulo 
$\, [\, J_{\overline N , \overline F} \cap J_{\overline N' , \overline F}\, ,\, U_{\overline N , \overline F}\, ]\, $  consisting of conjugates of elementary commutators 
$\, [\, u_{f_{k , x}}\, ,\, u_{\gamma }\, ]\, $ by arbitrary elements of $\, U_{\overline N'} + U_{E / N}\, $ where 
$\,\{\, f_{k , x}\,\}\, $ is a basis of $\, N / [ N , N ]\, $ as in Definition 8. Using the interdependence of the trivial relations given in the proof of Theorem 4 one finds that reducing this group modulo 
$\, [\, J_{N , F}\, ,\, U_F\, ]\, $ amounts to reducing the conjugacy coefficients to 1 and in addition replacing 
the elementary commutators by commutators of the form $\, [\, u_{e_i}\, ,\, u_{\gamma }\, ]\, $ where 
$\,\{\, e_i\,\}\, $ is a basis of the core of $\, \overline N\, $ for $\, U\, $ (at least after dividing by 
$\, [\, U_{N , F}\, ,\, J_F\, ]\, $ which also reduces the right side of these elementary commutators to their image in $\,\overline N'\, $, i.e. one gets the relations 
\smallskip
$$ \bigl[\, u_{e_i}\, ,\, u_{\gamma\delta }\, \bigr]\, =\, \bigl[\, u_{e_i}\, ,\, u_{\gamma }\, \bigr]\, +\, 
\bigl[\, u_{e_i}\, ,\, u_{\delta }\, \bigr]\, =\, \bigl[\, u_{e_i}\, ,\, u_{\delta\gamma }\, \bigr]\> . $$
\par\bigskip\noindent
One checks that there are no other relations, so the result follows from the definition of the tensor product by generators and relations \qed 
\par\bigskip\noindent 
{\bf Lemma 13.} (Cone Lemma)\quad The cone $\, (\, U_{N , F}\, ,\, U_F\, )\, $ is regular for any normal pair $\, ( N , F )\, $ and its $K^J$-groups (resp. $K$-groups) are trivial in all dimensions 
$\,\geq 2\, $.
\par\bigskip\noindent
{\it Proof.}\quad Clearly $\, ( U_{N , F} , U_F )\, $ is admissible and 
$\, {\overline U_{N , F}}\, =\, U_{N , F}\, $ as there is always a compatible splitting 
$\, (\, U_{N , F}\, ,\, U_F\, )\,\buildrel\sigma\over\longrightarrow\, (\, V\, ,\, W\, )\, $ to any projection in the reverse direction, so the pair is 0-regular and 
$\, R_{U_{N , F}\, ,\, U_F}\, =\, J_{U_{N , F}\, ,\, U_F}\, $. In particular 
$\, (\, J_{U_{N , F}\, ,\, U_F}\, ,\, U^2_F\, )\, $ is admissible by the proof of Theorem 2 (or by Lemma 14 below), but more is true. We will construct a strict splitting $\, J \longrightarrow U\, $ where 
$\, J\, =\, J_{U_{N , F}\, ,\, U_F}\, $ and 
$\, U\, =\, U_{J_{U_{N , F}\, ,\, U_F}\, ,\, U^2_F}\, $, in the sense of Definition 7 of section 3. This will imply that each $\, (\, J^n_{U_{N , F}\, ,\, U_F}\, ,\, U^{n+1}_F\, )\, $ 
is admissible and  
\smallskip
$$   J^n_{U_{N , F}\, ,\, U_F}\> =\> R_{J^{n-1}_{U_{N , F}\, ,\, U_F}\, ,\, U^n_F}\> =\> 
R^n_{U_{N , F}\, ,\, U_F}\> ,  $$
\smallskip
so that it is regular with trivial $K^J$-groups (and trivial $K$-groups). Let $\, I\, =\, J_{U_F}\, $ and 
$\, {\overline I}\, =\, J_{U_{F/N}}\, $, so that $\, I\,\simeq\, J\rtimes\, {\overline I}\, $. Start by considering the standard core 
\smallskip
$$  C_{I\, ,\, U}\> =\> \bigl\{\, u^{\pm}_{u_{x_1}}\cdots u^{\pm}_{u_{x_n}}\> 
u^{-1}_{u^{\pm}_{x_1}\cdots u^{\pm}_{x_n}}\,\bigr\}   $$
\smallskip
for the canonical lift $\ U_F\,\rightarrow\, U^2_F\, $, so that $\, C_{I\, ,\, U}\, $ is complemented by 
$\, C_{U/I}\, =\, \{\, u_{u_x}\,\}\, $ to give a full basis of $\, U^2_F\, $. Then 
$\, \langle\, C_{I\, ,\, U}\,\rangle\, $ is equally well expressed by a basis 
$\, C_{J\, ,\, U}\, +\, C_{{\overline I}\, ,\, {\overline U}}\, $ with $\, C_{J\, ,\, U}\, $ a core for the semicanonical lift of a core $\, C_{{\overline I}\, ,\, {\overline U}}\, $ of $\, {\overline I}\, $ with respect to the canonical lift $\, U_{F/N}\, \rightarrow\, U^2_{F/N}\, $. It is clear that the standard core $\, C_{I\, ,\, U}\, $ maps into the standard core $\, C_{{\overline I}\, ,\, {\overline U}}\, $ modulo $\, J\, $, so that on choosing a section $\, s\, :\, F/N\,\nearrow\, F\, $ it will lift $\, C_{{\overline I}\, ,\, {\overline U}}\, $ into 
$\, C_{I\, ,\, U}\, $, thus giving a core $\, C_{J\, ,\, U}\, $ for this lift. $\, C_{J\, ,\, U}\, $ is completed to a core 
$\, C_{J\, ,\, I}\, $ of $\, J\, $ for $\, I\, $ by putting
\smallskip
$$ C_{J\, ,\, I}\; =\; C_{U/I}\> C_{J\, ,\, U}\; C^{-1}_{U/I}\; +\; C_{U/I}\; 
\bigl[\, C_{{\overline I}}\, ,\, C_0\,\bigr]\; C^{-1}_{U/I}   $$
where $\, C_0\, =\, C_{U_{N,F}\, ,\, U_F}\,\subseteq\, C_{U/I}\, $. Then a full basis of $\, J\, $ is given by the set 
\smallskip
$$ C_{\overline I}\; C_{U / I}\; C_{J\, ,\, U}\; C^{-1}_{U / I}\; C^{-1}_{\overline I }\; +\;
C_{\overline I }\; C_{U / I}\; \bigl[\, C_{\overline I }\, ,\, C_0\,\bigr]\> C^{-1}_{U / I}\; C^{-1}_{\overline I }  $$
\smallskip
which means that writing $\, U^2_F / J\, $ as a semidirect product 
$\, {\overline I }\,\rtimes\, U_F\, =:\, {\tilde U}\, $ and applying the section which is the semicanonical lift of $\, \overline I\, $ to $\, I\, $ as above and the copy of $\, U_F\, $ is given by the canonical lift 
$\, U_F\,\rightarrow \, U^2_F\, $ one gets a basis element for each expression 
$\, x\, f\, x^{-1}\, $ with $\, x\,\in\, {\tilde U }\, $ and $\, f\,\in\, C_{J\, ,\, U}\, $ or 
$\, f\,\in\, [\, C_{\overline I }\, ,\, C_0\, ]\, $. It is then immediate that lifting $\, f\, $ to $\, u_f\, $ and $\, x\, $ to 
$\, U^3_F\, $ by the canonical lift $\, U^2_F\,\rightarrow U^3_F\, $ determines a lift $\, J\,\rightarrow\, U\, $ 
which is semicanonical and normal modulo $\, [\, U\, ,\, J_{U^2_F}\, ]\, $ and that it satisfies the strictness condition of Definition 7 \qed 
\par\bigskip\noindent
An important technical condition in $K^J$-theory is (relative or partial) admissibility of a pair 
$\, ( N , F )\, $.  
\par\bigskip\noindent
{\bf Definition 11.}\quad Let $\, N_0\, $ be a normal subgroup of $\, F\, $ such that 
$\, [ N , F ] \subseteq N_0 \subseteq N\, $.  The pair $\, ( N , F )\, $ will be called {\it partially admissible on $\, N_0\, $}  if the preimage of $\, N_0\, $ in the canonical $F$-central extension of 
$\, N\, $ contains an $F$-normal subextension of $\, N_0\, $ by $\, K^J_2 ( N , F )\, $ which then is uniquely determined up to equivalence. 
\par\smallskip\noindent
Considering an inclusion $\, ( N , F ) \subseteq ( M , G )\, $ the pair $\, ( N , F )\, $ is called 
{\it relatively admissible for $\, ( M , G )\, $} if the canonical $G$-central extension of $\, M\, $ contains an 
$F$-normal subextension of $\, N\, $ by $\, K^J_2 ( N , F )\, $. (This is equivalent to the content of Definition 6 in case of a relative semiuniversal $F$-central extension).
\par\bigskip\noindent
{\bf Remark.}\quad The condition of partial admissibility might seem rather artificial and far-fetched at first, but it is needed as a technical bridge in a certain context to deal with the (first) suspension of a normal pair and thus arises naturally. Note that any pair $\, ( N , F )\, $ is partially admissible on 
$\, [ N , F ]\, $. If $\, ( N , F ) \, $ is partially admissible on $\, N_0\, $ and if
$\, {\overline N_0}^{ _{\sps {( N , F )}}}\, $ denotes the subextension in the canonical $F$-central extension of 
$\, N\, $ as above, then there exists an "almost canonical" map from the preimage of 
$\, N_0\, $ in $\, U_{N , F}\, $ to $\, {\overline N_0}^{ _{\sps {( N , F )}}}\, $ , where "almost canonical" simply means an equivariant projection over the evaluation map. To see this construct a surjection of 
$\, U_{N , F}\, $ to a certain $F$-central extension $\, {\overline N}\, $ of $\, N\, $ with kernel 
$\, K^J_2 ( N , F )\, $ as follows. Consider the commutative diagram
\smallskip
$$ \vbox{\halign{#&#&#&#&#\cr
$ J_{N,F} \cap \bigl[ U_{N,F} , U_F \bigr] $ & $\largerightarrow $ & $\bigl[\, U_{N,F}\, ,\, U_F\,\bigr] $ &
$\largerightarrow $ & \hfil $\bigl[\, N\, ,\, F\,\bigr]\quad $\hfil \cr
\hfil $\qquad\Biggm\downarrow $\hfil && \hfil $\Biggm\downarrow $\hfil && \hfil $\Biggm\downarrow\quad $\hfil \cr
\hfil $\qquad J_{N,F} $\hfil & $\largerightarrow $ & \hfil $ U_{N,F} $\hfil & $\largerightarrow $ & 
\hfil $ N\quad $\hfil \cr
\hfil $\qquad\Biggm\downarrow $\hfil && \hfil $\Biggm\downarrow $\hfil && 
\hfil $\Biggm\downarrow\quad $\hfil \cr
\hfil $\qquad J_0 $\hfil & $\largerightarrow $ & 
\hfil $ {U_{N;F} \over \bigl[\, U_{N,F}\, ,\, U_F\, \bigr] }$\hfil & 
$\largerightarrow $ & \hfil $ {N \over \bigl[\, N\, ,\, F\,\bigr] }\quad$\hfil \cr}}    .         $$
\smallskip
Divide the diagram by
\smallskip 
$$ \bigl[\, J_{N,F}\, ,\, U_F\, \bigr]\, +\, \bigl[\, U_{N;F}\, ,\, J_F\, \bigr]\, \subseteq\, 
     J_{N;F}\,\cap\, \bigl[\, U_{N,F}\, ,\, U_F\, \bigr]\>   .   $$
\smallskip
Then as $\, J_0\, $ is free abelian (as a subgroup of $\, U_F\, /\, [\, U_F\, ,\, U_F\, ]\, $) it has a lift to
$\, J_{N,F}\, /\, (\, [\, J_{N,F}\, ,\, U_F\, ] + [\, U_{N,F}\, ,\, J_F\, ]\, )\, $ and dividing the image of $\, J_0\, $ gives a projection of $\, U_{N , F}\, $ onto an $\, F$-central extension $\, {\overline N}\, $ of $\, N\, $ by 
$\, K^J_2 ( N , F )\, $, hence an isomorphism of $\, {\overline N_0}^{ _{\sps {( N , F )}}}\, $ with the preimage of $\, N_0\, $ in $\, {\overline N}\, $ which is the identity on $\, K^J_2 ( N , F )\, $ 
since the image of $\, [\, U_{N , F}\, ,\, U_F\, ]\, $, which is $\, [\, {\overline N}\, ,\, F\, ]\, $, is canonical, so that the isomorphism is canonical on 
$\, [\, \overline N\, ,\, F\, ]\,\cap\, {\overline N_0}^{ _{\sps {( N , F )}}}\,\supset\, K^J_2 ( N , F )\, $. Then 
$\, {\overline N_0}^{ _{\sps {( N , F )}}}\, $ has the following "semiuniversal" property. Given any $F$-central extension of $\, N\, $ there exists a map from $\, {\overline N_0}^{ _{\sps {( N , F )}}}\, $ to the preimage of $\, N_0\, $ in this extension which is canonical on the subgroup 
$\, [\, \overline N\, ,\, F\, ]\,\cap\, {\overline N_0}^{ _{\sps {( N , F )}}}\, $. Lemmas 9 and 11 can now be extended to the notions almost admissible and relatively admissible in the following way. In Lemma 9 assume that the quotient map induces a surjection on $\, K^J_2\, $, that the kernel is admissible and the quotient is partially admissible. Then the extension is partially admissible (on the appropriate subgroup). Lemma 11 is also generalized in a straightforward manner. We have already seen some examples of admissible pairs, e.g. the full subgroups of arbitrary groups $\, F\, $ (by the proof of Theorem 1). An example of a pair which is not admissible is given by a central copy of 
$\, {\mathbb Z}\,/\, n{\mathbb Z}\, $ in an arbitrary group $\, F\, $. Then the extension
\smallskip
$$ 1\,\largerightarrow\, {\mathbb Z}\,\buildrel n\over\largerightarrow\, {\mathbb Z}\,\largerightarrow\, {\mathbb Z}\, /\, n{\mathbb Z}\,\largerightarrow\, 1  $$
\smallskip
(with trivial $\, F$-action) is a minimal $\, F$-central extension, but its kernel $\, ( n{\mathbb Z} )\, $ is not contained in $\, [\, {\mathbb Z}\, ,\, F\, ]\> =\> 1\, $. This result is easily generalized to finite central torsion groups, so one finds that a finitely generated central (abelian) group is admissible if and only if it is free abelian. In fact the preceding argument exhibits a surprisingly simple characterization of admissibility in terms of the group $\, K_1 ( N , F )\, $, which is defined just below. Namely the pair $\, ( N , F )\, $ is admissible if and only if $\, K_1 ( N , F )\, $ is free abelian. \par\bigskip\noindent
{\bf Definition 12.}\quad For any normal pair $\, ( N , F )\, $ define
\smallskip
$$ K_1\,\bigl(\, N ,\, F\, \bigr)\; =\; {N\over \bigl[\, N\, ,\, F\,\bigr] }    $$
\smallskip
Also if $\, (\, C_\varphi\, ,\, U_\varphi\, )\, $ denotes the mapping cone of an injective morphism 
$\, ( N , F )\,\subseteq\, ( M , G )\, $, let $\, {\widetilde C}_\varphi\, $ be the image of $\, C_\varphi\, $ by an "almost canonical" map 
\smallskip
$$ C_\varphi\,\subseteq\, U_{M , G}\,\longrightarrow\, {\overline M}^{ _{\sps G}} $$
\smallskip
where $\, {\overline M}^{ _{\sps G}}\, $ is the semiuniversal (or canonical if $\, ( M , G )\, $ is not admissible) $G$-central extension of 
$\, M\, $ and put 
\smallskip
$$ {\widetilde K}_1\, \bigl(\, C_\varphi\, ,\, U_\varphi\, \bigr)\; =\; 
{{\widetilde C}_\varphi\,\cap\, \bigl[\, {\overline M}^{ _{\sps G}}\, ,\, G\,\bigr]\over
\bigl[\, {\widetilde C}_\varphi\, ,\, F\,\bigr] }\>    .   $$
\smallskip
Using the semiuniversal $G$-central extension $\, {\widetilde C}_\varphi\, $ is well defined only up to isomorphism over the identity of$\, N\, $ fixing the kernel but the group 
$\, {\widetilde K}_1\, (\, C_\varphi\, ,\, U_\varphi\, )\, $ is well defined in any case and functorial from mapping cones to abelian groups.
\par\bigskip\noindent
A useful result is the following.
\par\bigskip\noindent 
{\bf Lemma 14.}\quad (Suspension Lemma) Assume that $\, ( R , U )\, $ is a normal pair such that 
$\, R\, $ is a free group. Then the pair $\, ( J_{R , U}\, ,\, U_U )\, $ is admissible. In particular for any normal pair 
$\, ( N , F )\, $ the second $J$-suspension $\, (\, J^2_{N , F}\, ,\, U^2_F\, )\, $ is admissible. 
\par\bigskip\noindent
{\it Proof.}\quad By the proof of Corollary 5.1 below one gets for any extension of normal pairs 
\smallskip
$$ 1 \longrightarrow ( N , F ) \longrightarrow ( M , F ) \longrightarrow ( P , E ) \longrightarrow 1 $$
\par\smallskip\noindent 
the following exact sequence of $K^J$-groups
\smallskip
$$ K^J_2 ( N , F ) \longrightarrow K^J_2 ( M , F ) \longrightarrow K^J_2 ( P , E ) \longrightarrow\quad\quad $$
$$ \quad\quad K_1 ( N , F ) \longrightarrow K_1 ( M , F ) \longrightarrow K_1 ( P , E ) \longrightarrow 1 $$
\par\smallskip\noindent
without further assumptions. Then consider the extension of normal pairs
\smallskip
$$ 1 \longrightarrow ( J_{R , U}\, ,\, U_U ) \longrightarrow ( U_{R , U}\, ,\, U_U ) \longrightarrow 
( R \, ,\, U_U / J_{R , U} ) \longrightarrow 1 $$
\par\smallskip\noindent
where $\, K^J_2 ( U_{R , U}\, ,\, U_U )\, =\, 0\, $ and $\, K_1 ( U_{R , U}\, ,\, U_U )\, $ is free abelian since $\, ( U_{R , U}\, ,\, U_U )\, $ is a cone. Then $\, K_1 ( J_{R , U}\, ,\, U_U )\, $ is an extension of some free abelian group by 
$\, K^J_2 ( R\, ,\, U_U / J_{R , U} )\, $. The quotient of $\, U_U / J_{R , U}\, $ by 
$\, R\, $ is the free group $\, U_{U / R}\, $. It now follows easily from the proof of the Excision Theorem that $\, K^J_2 ( N , F )\, $ is a free abelian group whenever $\, N\, $ and $\, F / N\, $ are free, because dividing $\, J_{N , F}\, $ by the trivial relations $\, [\, J_{N , F}\, ,\, U_F\, ]\, +\, 
[\, U_{N , F}\, ,\, J_F\, ]\, $ using $B_3$-reduction and the usual basis of $\, J_{N , F}\, $ given by the core 
$\,\{\, u_z\, B_{1 , 2}\, u^{-1}_z\,\}\, $ the difference between the relations 
$\, [\, u_c\, ,\, u_x\, u_y\, u^{-1}_{xy}\, ]\, $ and the corresponding elements $\, [\, B_2\, ,\, u_x\, ]\, $ can be used to reduce the $B_2$-elements to the generating set $\, [ u_{x_j}\, ,\, u_c\, ]\, u_c\, u^{-1}_{c_{x_j}}\, $ 
where $\,\{\, x_j\,\}\, $ is a basis of $\, F / N\, $, whereas the difference between the relations 
$\, [\, B_2\, ,\, u_c\, ]\, $ and the corresponding relations 
$\, [\, B_1\, ,\, u_z\, ]\, $ can be used to reduce the elements $\, B_2\, $ further to the generating set 
$\, [\, u_{x_j}\, ,\, u_{e_k}\, ]\, $ with $\,\{ e_k\,\}\, $ a basis of $\, N\, $. Then as there are no more relations one is left with a free abelian group and the intersection with the image of 
$\, J_{N , F}\,\cap\, [\, U_{N , F}\, ,\, U_F\, ]\, $ is also free abelian. Returning to the situation above one gets that $\, K_1 ( J_{R , U}\, ,\, U_U )\, $ is the extension of a free abelian group by another free abelian group, hence free abelian. By the observation preceding Definition 10 the pair 
$\, ( J_{R , U}\, ,\, U_U )\, $ is admissible
\qed 
\par\bigskip\noindent
The notion of (partial) admissibility given by Definition 11 is one ingredient in the proof of exactness of the Mapping Cone sequence but not a really essential one. There exists a more fundamental obstruction which is captured in the following definition.
\par\bigskip\noindent
{\bf Definition 13.} Let $\, ( P , E )\, $ be a normal pair. Then consider 
$\, \overline F := U_E / J_{P , E} \simeq P \rtimes U_{E / P}\, $ and let $\,\overline N\, $ denote the copy of 
$\, J_{E / P}\, $ complementary to $\, P\, $ which is the image of 
$\, J_E\, $. Put $\, \overline M := J_{E / P} \times P\, $ and consider 
$\, U_{\overline F} / J_{\overline M , \overline F}\,\simeq\, ( J_{E / P} \times P ) \rtimes U_{E / P}\, $ which contains another copy of $\, J_{E / P}\, $ complementary to $\, \overline M\, $ (which is the image of 
$\, J_{\overline F}\, $). It will be denoted $\, J^0_{E / P}\, $. Then $\, ( P , E )\, $ is said to be an
{\it exact pair} or {\it exact of first order} if one of the following equivalent conditions hold:
\par\smallskip\noindent
(i) the following map is injective 
\smallskip
$$ K^J_2 ( J_{E / P}\, ,\, \overline F )\,\rightarrowtail\, 
K^J_2 ( J_{E / P} \times P\, ,\, \overline F )\> , $$
\par\medskip\noindent
(ii) the following map is injective
\smallskip
$$ K^J_2 ( J_{E / P} \times J^0_{E / P}\, ,\, U_{\overline F} / J_{\overline M , \overline F} )\,\rightarrowtail\, 
K^J_2 ( J_{E / P} \times P \times J^0_{E / P}\, ,\, U_{\overline F} / J_{\overline M , \overline F} )\> , $$
\par\medskip\noindent
(iii) any element in the trefoil intersection of the images of the groups
$\, K^J_2 ( J_{E / P}\, ,\, U_{\overline F} / J_{\overline M , \overline F} )\, $, 
$\, K^J_2 ( P\, ,\, U_{\overline F} / J_{\overline M , \overline F} )\, $, and
$\, K^J_2 ( J^0_{E / P}\, ,\, U_{\overline F} / J_{\overline M , \overline F} )\, $ in the enveloping group
$\, K^J_2 ( J_{E / P} \times P \times J^0_{E / P}\, ,\, U_{\overline F} / J_{\overline M , \overline F} )\, $ 
lifts to an element in the corresponding trefoil intersection with 
$K^J_2$-groups replaced by the extended $\widetilde K^J_2$-groups in each instance.
\par\medskip\noindent
We will say that $\, ( P , E )\, $ is {\it exact of second order} if 
$\, ( J_{E / P}\, ,\,\overline F )\, $ is an exact pair. $\, ( P , E )\, $ is said to be {\it fully exact} if it is exact of first and second order.
\par\bigskip\noindent  
The conditions of Definition 13 should not be easy to check in a general context. Condition (i) is the one we are mainly interested in because it is directly connected with certain problems of exactness. There is a counterexample to exactness of the following part of the Mapping Cone sequence
\smallskip
$$ K^J_2 ( C_{\varphi }\, ,\, U_F ) \longrightarrow K^J_2 ( N , F ) \longrightarrow K^J_2 ( M , F ) $$
\par\medskip\noindent
in the case of a free abelian group $\, C\, $ of rank greater than 2 sitting as a direct summand of some arbitrary group $\, B\, $ (say $\, B = C\, $, compare with Lemmas 5, 6 of section 1) on putting 
$\, F = U_C / J_{{\mathbb Z}_0\, ,\, C}\, ,\, N = J_{C / {\mathbb Z}_0}\, ,\, 
M = J_{C / {\mathbb Z}_0} \times {\mathbb Z}_0\, $ with $\, {\mathbb Z}_0 \simeq \mathbb Z\, $ a direct summand of $\, C\, $. The obstruction to exactness must come from conditions (i)--(iii) not being satisfied. We show the implications (i) $\Rightarrow$ (ii) and (ii) $\Rightarrow$ (iii) below, the implication (iii) $\Rightarrow$ (i) is given in the proof of the Mapping Cone Theorem. The following example gives a special case where the condition (iii) applies which is very important for the regular theory.
\par\bigskip\noindent
{\it Example.}\quad Assume that $\, P\, $ is a central subgroup of $\, E\, $ contained in $\, [ E , E ]\, $. We claim that the trefoil intersection of Definition 13 is empty in this case.  Such an element is in the kernel of the projection dividing by $\, P\, $ so it can be represented as an element 
$\, [\, U_P\, ,\, U_{J_{E / P} \times J^0_{E / P}}\, ]\, $ modulo trivial elements. Since $\, P\, $ is central it maps injectively to the quotient by the normal subgroup 
$\, [\, J_{E / P}\, ,\, U_{E / P}\, ] \times [\, J^0_{E / P}\, ,\, U_{E / P}\, ]\, $. For assume the element is in the kernel of the map induced by division of $\, [\, J_{E / P}\, ,\, U_{E / P}\, ]\, $. Then it is represented by an element of $\, [\, U_{[\, J_{E / P}\, ,\, U_{E / P}\, ]}\, ,\, 
U_{( J_{E / P} \times P ) \rtimes U_{E / P}}\, ]\, $ which becomes trivial dividing by $\, P\, $ so from the splitting $\, J_{E / P} \rtimes U_{E / P} \longrightarrow ( J_{E / P} \times P ) \rtimes U_{E / P}\, $ which is the identity on $\, J_{E / P}\, $ one finds that the element is represented by a commutator 
$\, [\, U_{[\, J_{E / P}\, ,\, U_{E / P}\, ]}\, ,\, U_P\, ]\, $ which gives a trivial element due to the fact that 
$\, P\, $ is central. Then by an analogous argument in the quotient 
$\, K^J_2 ( \overline J_{E / P} \times P \times J^0_{E / P}\, ,, 
( \overline J_{E / P} \times P ) \rtimes U_{E / P} )\, $ one finds that one can also divide by 
$\, [\, J^0_{E / P}\, ,\, U_{E / P}\, ]\, $  leaving us with an extension 
$\,\overline E\, $ of $\, E\, $ by 
\smallskip
$$ {\overline J_{E / P}} \times {\overline J^0_{E / P}}\> :=\>
 {J_{E / P} \over [\, J_{E / P}\, ,\, U_{E / P}\, ]} \times  {J_{E / P} \over [\, J_{E / P}\, ,\, U_{E / P}\, ]}\> . $$
\par\bigskip\noindent
The diagonal of the kernel is now a normal subgroup which may be divided out to obtain the quotient 
$\, {\overline J_{E / P}} \times E\, $.  Since $\, P \subset [ E , E ]\, $ one gets 
$\, [\, U_P\, ,\, U_{\overline J_{E / P}}\, ]\, \subseteq\, 
[\, [\, U_{\overline J_{E / P}}\, ,\, U_E\, ] , U_E \, ]\, $ modulo trivial elements in the quotient, which is trivial so that the element can be represented by an element of
$\, [\, U_P\, ,\, U_{\overline J^{\Delta }_{E / P}}\, ]\, $ where $\, \overline J^{\Delta }_{E / P}\, $ denotes the diagonal. Dividing by $\, \overline J_{E / P}\, $ the element becomes trivial. On the other hand there is a splitting $\, \overline E \simeq \overline J_{E / P} \times \widetilde E\, $ which is the identity on $\, P\, $ and sends the image of $\, \overline J^0_{E / P}\, $ in $\,\widetilde E\, $ to the diagonal 
$\, \overline J^{\Delta }_{E / P}\, $ showing that our element is trivial. Hence the trefoil intersection is empty settling condition (iii).    
\par\smallskip\noindent
Another instance where the trefoil intersection is empty is when $\, P\, $ has a weak core for $\, E\, $. In this case Lemma 12 gives a formula for the intersection of the images of 
$\, K^J_2 ( P\, ,\, ( J_{E / P} \times P ) \rtimes U_{E / P} )\, $ and of
$\, K^J_2 ( J_{E / P} \times J^0_{E / P}\, ,\, ( J_{E / P} \times P ) \rtimes U_{E / P} )\, $ which is 
$\, {\mathcal C}_P \otimes ( J^{ab}_{E / P} \times J^{0 , ab}_{E / P} )\, =\, 
( {\mathcal C}_P \otimes J^{ab}_{E / P} ) \times ( {\mathcal C}_P \otimes J^{0 , ab}_{E / P} )\, $ where 
$\, {\mathcal C}_P\, $ denotes the core of the abelianization of $\, P\, $. From this one induces that the intersection of the images of $\, K^J_2 ( P\, ,\, ( J_{E / P} \times P ) \rtimes U_{E / P} )\, $ and 
$\, K^J_2 ( J_{E / P}\, ,\, ( J_{E / P} \times P ) \rtimes U_{E / P} )\, $ resp. 
$\, K^J_2 ( J^0_{E / P}\, ,\, ( J_{E / P} \times P ) \rtimes U_{E / P} )\, $ in 
$\, K^J_2 ( J_{E / P} \times P \times J^0_{E / P}\, ,\, ( J_{E / P} \times P ) \rtimes U_{E / P} )\, $ is equal to 
$\, {\mathcal C}_P \otimes J^{ab}_{E / P}\, $ resp. $\, {\mathcal C}_P \otimes J^{0 , ab}_{E / P}\, $, so the trefoil intersection is empty and in particular $\, ( P , E )\, $ is exact.
\par\smallskip\noindent 
To become a little more familiar with the notion of an exact pair we include a brief discussion in a general setting. Let us look at condition (iii) in more detail. By the extension
\smallskip
$$ 1 \longrightarrow ( J_{\overline M , \overline F}\, ,\, U_{\overline F} ) \longrightarrow 
( U_{\overline M , \overline F}\, ,\, U_{\overline F} ) \longrightarrow 
( \overline M\, ,\, U_{\overline F} / J_{\overline M , \overline F} ) \longrightarrow 1 $$
\par\bigskip\noindent
and the associated exact sequence of $K^J_2$- and $K_1$-groups the map 
$\, K^J_2 ( J_{E / P} \times P\, ,\, U_{\overline F} / J_{\overline M , \overline F} )\, \rightarrowtail\, 
K_1 ( J_{\overline M , \overline F}\, ,\, U_{\overline F} )\, $ is injective, hence also the map 
$\, K^J_2 ( J_{E / P} \times P\, ,\, U_{\overline F} / J_{\overline M , \overline F} )\,\rightarrow\, 
K^J_2 ( J_{E / P} \times P \times J^0_{E / P}\, ,\, U_{\overline F} / J_{\overline M , \overline F} )\, $. There is a symmetry of $\, U_{\overline F} / J_{\overline M , \overline F}\, $ exchanging the copies $\, J_{E / P}\, $ and $\, J^0_{E / P}\, $ and leaving invariant the normal subgroup $\, P\, $ as well as the diagonal copy of $\  U_{E / P}\, $ so that also the map 
$\, K^J_2 ( P \times J^0_{E / P}\, ,\, U_{\overline F} / J_{\overline M , \overline F} )\,\rightarrowtail\, 
K^J_2 ( J_{E / P} \times P \times J^0_{E / P}\, ,\, U_{\overline F} / J_{\overline M , \overline F} )\, $ is injective. From the extension 
\medskip
$$ 1 \longrightarrow ( J_{\overline M , \overline F}\, ,\, U_{\overline F} ) \longrightarrow 
( U_{P , \overline F} + J_{\overline M , \overline F}\, ,\, U_{\overline F} ) \longrightarrow 
( P\, ,\, U_{\overline F} / J_{\overline M , \overline F} ) \longrightarrow 1 $$
\par\medskip\noindent
and the fact that $\, K^J_2 ( J_{\overline M , \overline F}\, ,\, U_{\overline F} )\,\twoheadrightarrow\, 
K^J_2 ( U_{P , \overline F} + J_{\overline M , \overline F}\, ,\, U_{\overline F} )\, $ is surjective (the group on the right can be identified with 
$\, K^J_2 ( J_{J_{E / P}\, ,\, U_{E / P}}\, ,\, U_{U_{E / P}} )\, $ so the map splits) one deduces that the maps from $\, K^J_2 ( P\, ,\, U_{\overline F} / J_{\overline M , \overline F} )\, $ to  
$\, K^J_2 ( J_{E / P} \times P\, ,\, U_{\overline F} / J_{\overline M , \overline F} )\, $ and 
$\, K^J_2 ( P \times J^0_{E / P}\, ,\, U_{\overline F} / J_{\overline M , \overline F} )\, $ are also injective. 
Considering the extension
\medskip
$$ 1 \longrightarrow ( J_{\overline M , \overline F}\, ,\, U_{\overline F} ) \longrightarrow 
( C_{\rho }\, ,\, U_{\overline F} ) \longrightarrow 
( J_{E / P}\, ,\, U_{\overline F} / J_{\overline M , \overline F} ) \longrightarrow 1 $$
\par\medskip\noindent
and surjectivity of $\, K^J_2 ( J_{\overline M , \overline F}\, ,\, U_{\overline F} )\, \twoheadrightarrow\,
K^J_2 ( C_{\rho }\, ,\, U_{\overline F} )\,\simeq\, K^J_2 ( J_{P , E}\, ,\, U_E )\, $ (construct a section 
$\, s : E \nearrow \overline F\, $ which is the identity on $\, P\, $, it induces a splitting
$\, ( J_{P , E}\, ,\, U_E ) \rightarrow  (J_{\overline M , \overline F}\, ,\, U_{\overline F} )\, $) 
one gets injectivity of the maps from  
$\, K^J_2 ( J_{E / P}\, ,\, U_{\overline F} / J_{\overline M , \overline F} )\, $ to  
$\, K^J_2 ( J_{E / P} \times P\, ,\, U_{\overline F} / J_{\overline M , \overline F} )\, $ and 
$\, K^J_2 ( J_{E / P} \times J^0_{E / P}\, ,\, U_{\overline F} / J_{\overline M , \overline F} )\, $,
and by symmetry injectivity of the maps $\, K^J_2 ( J^0_{E / P}\, ,\, U_{\overline F} / J_{\overline M , \overline F} )\,\rightarrowtail\, 
K^J_2 ( P \times J^0_{E / P} \, ,\, U_{\overline F} / J_{\overline M , \overline F} )\, $ as well as 
$\, K^J_2 ( J^0_{E / P}\, ,\, U_{\overline F} / J_{\overline M , \overline F} )\,\rightarrowtail\, 
K^J_2 ( J_{E / P} \times J^0_{E / P} \, ,\, U_{\overline F} / J_{\overline M , \overline F} )\, $. So the only map which is not necessarily injective is the one from 
$\, K^J_2 ( J_{E / P} \times J^0_{E / P} \, ,\, U_{\overline F} / J_{\overline M , \overline F} )\, $ to 
$\, K^J_2 ( J_{E / P} \times P \times J^0_{E / P} \, ,\, U_{\overline F} / J_{\overline M , \overline F} )\, $. If this map also were injective the lifting problem (iii) could be solved in the following manner. Take the unique preimages of an element in the trefoil intersection in the groups 
$\, K^J_2 ( J_{E / P} \times P\, ,\, U_{\overline F} / J_{\overline M , \overline F} )\, $ and 
$\, K^J_2 ( P \times J^0_{E / P}\, ,\, U_{\overline F} / J_{\overline M , \overline F} )\, $ respectively and lift them to elements of 
$\,\widetilde K^J_2 ( J_{E / P} \times P\, ,\, U_{\overline F} / J_{\overline M , \overline F} )\, $ in the intersection of the images of $\, \widetilde K^J_2 ( J_{E / P}\, ,\, U_{\overline F} / J_{\overline M , \overline F} )\, $ and $\, \widetilde K^J_2 ( P\, ,\, U_{\overline F} / J_{\overline M , \overline F} )\, $, and of 
$\, \widetilde K^J_2 ( P \times J^0_{E / P}\, ,\, U_{\overline F} / J_{\overline M , \overline F} )\, $ in the intersection of the images of $\, \widetilde K^J_2 ( P\, ,\, U_{\overline F} / J_{\overline M , \overline F} )\, $ and of $\, \widetilde K^J_2 ( J^0_{E / P}\, ,\, U_{\overline F} / J_{\overline M , \overline F} )\, $ respectively (this is always possible). Take the difference of some corresponding preimages in 
$\, \widetilde K^J_2 ( P\, ,\, U_{\overline F} / J_{\overline M , \overline F} )\, $. The resulting element is contained in the kernel of the map to $\, K^J_2 ( P\, ,\, U_{\overline F} / J_{\overline M , \overline F} )\, $, i.e. it is represented by an element of 
$ [\, U_P\, , J_{U_{\overline F} / J_{\overline M , \overline F}}\, ] $ modulo 
$\, [\, J_{P , U_{\overline F} / J_{\overline M , \overline F}}\, ,\, 
U_{U_{\overline F} / J_{\overline M , \overline F}}\, ]\, $ and becomes trivial dividing by the normal subgroup $\, J_{E / P} \times J^0_{E / P}\, $. Then its image in 
$\, \widetilde K^J_2 ( J_{E / P} \times P \times J^0_{E / P}\, ,\, 
U_{\overline F} / J_{\overline M , \overline F} )\, $ comes from an element in 
$\, \widetilde K^J_2 ( J_{E / P} \times J^0_{E / P}\, ,\, U_{\overline F} / J_{\overline M , \overline F} )\, $ in the kernel of the map to 
$\, K^J_2 ( J_{E / P} \times J^0_{E / P}\, ,\, U_{\overline F} / J_{\overline M , \overline F} )\, $ (it is here that one uses injectivity of 
$\, K^J_2 ( J_{E / P} \times J^0_{E / P}\, ,\, U_{\overline F} / J_{\overline M , \overline F} )\,\rightarrowtail\, 
K^J_2 ( J_{E / P} \times P \times J^0_{E / P}\, ,\, U_{\overline F} / J_{\overline M , \overline F} )\, $). This means it can be represented by an element of 
$\, [\, U_{J_{E / P}}\, ,\, J_{U_{\overline F} / J_{\overline M , \overline F}}\, ]\, +\, 
[\, U_{J^0_{E / P}}\, ,\, J_{U_{\overline F} / J_{\overline M , \overline F}}\, ]\, $ modulo 
$\, [\, J_{J_{E / P} \times J^0_{E / P}\, ,\, U_{\overline F} / J_{\overline M , \overline F}}\, ,\, 
U_{U_{\overline F} / J_{\overline M , \overline F}}\, ]\, $. Also this element becomes trivial dividing by the normal subgroup $\, P\, $. Consider the split exact sequence
\medskip
$$ 1 \rightarrow ( P , ( J_{E / P} \times P ) \rtimes U_{E / P} ) \rightarrow 
( J_{E / P} \times P \times J^0_{ E / P} , ( J_{E / P} \times P ) \rtimes U_{E / P} ) \rightarrow $$
$$\qquad\qquad\qquad\qquad\qquad\qquad\qquad \rightarrow 
( J_{E / P} \times J^0_{E / P}\, ,\,  J_{E / P} \rtimes U_{E / P} ) \rightarrow 1\> . $$
\par\medskip\noindent
The splitting can be chosen to map the image of $\, J_{E / P}\, $ in the quotient identically to itself and the image of $\, J^0_{E / P}\, $ to $\, P \times J^0_{E / P}\, $. Then  
a simple argument shows that it can be represented by an element of 
$\, [\, U_{J_{E / P}}\, ,\, J_{P\, ,\, U_{\overline F} / J_{\overline M , \overline F}}\, ]\, +\, 
[\, U_{J^0_{E / P}}\, ,\, J_{P\, ,\, U_{\overline F} / J_{\overline M , \overline F}}\, ]\, +\, 
( [\, U_P\, ,\, J_{U_{\overline F} / J_{\overline M , \overline F}}\, ]\,\cap\, 
U_{J^0_{E / P}\, ,\, U_{\overline F} / J_{\overline M , \overline F}} )\, $. Since the first two terms become trivial in 
$\, \widetilde K^J_2 ( J_{E / P} \times P \times J^0_{E / P}\, ,\, ( J_{E / P} \times P ) \rtimes U_{E / P} )\, $ they are negligible. We may then add/subtract the preimage of the last term in 
$\, \widetilde K^J_2 ( J^0_{E / P}\, ,\, U_{\overline F} / J_{\overline M , \overline F} )\, $ from the preimage of the original lift to  
$\, \widetilde K^J_2 ( P \times J^0_{E / P}\, ,\, U_{\overline F} / J_{\overline M , \overline F} )\, $ to obtain an element in the trefoil intersection of all three groups in 
$\, \widetilde K^J_2 ( J_{E / P} \times P \times J^0_{E / P}\, ,\, 
U_{\overline F} / J_{\overline M , \overline F} )\, $, solving (iii). The implication (i) $\Rightarrow$ (ii) also follows easily from the injectivity results obtained above.
\par\medskip\noindent
Inductively, one could say that a pair $\, ( P , E )\, $ is exact of order $\, n\, $ if its (first) derived pair
$\, ( J_{E / P}\, ,\, \overline F )\, $ is exact of order $\, n-1\, $. However one has the result that any pair 
$\, ( Q , H )\, $ is exact of third order. Of course this needs to be proved. Put $\, E = H / Q\, $. Then one is concerned with proving exactness of the second derived pair 
$\, ( J_H\, ,\, J_E \rtimes U_H )\, $. Considering the third derived pair 
$\, ( J_{Q \rtimes U_E}\, ,\, J_H \rtimes U_{Q \rtimes U_E} )\, $ one notes that 
the map 
\smallskip
$$ K^J_2 ( J_{Q\, ,\, H} \rtimes ( J_{Q \rtimes U_E} + U_{Q\, ,\, Q \rtimes U_E} )\, ,\, 
J_H \rtimes U_{Q \rtimes U_E} ) \rightarrowtail\qquad\qquad\qquad\qquad\qquad\qquad $$ 
$$\qquad\qquad\qquad\qquad\qquad\qquad 
K^J_2 ( J_H \rtimes ( J_{Q \rtimes U_E} + U_{Q\, ,\, Q  \rtimes U_E} )\, ,\, 
J_H \rtimes U_{Q \rtimes U_E} ) $$
\par\medskip\noindent
is injective by existence of a splitting.
This implies that the image of an element in the kernel of
\smallskip
$$ K^J_2 ( J_{Q \rtimes U_E}\, ,\, J_H \rtimes U_{Q \rtimes U_E} ) \longrightarrow 
K^J_2 ( J_{Q \rtimes U_E} \times J_H\, ,\, J_H \rtimes U_{Q \rtimes U_E} ) $$ 
\par\medskip\noindent
must be trivial in 
$\, K^J_2 ( J_{Q\, ,\, H} \rtimes ( J_{Q \rtimes U_E} + U_{Q\, ,\, Q \rtimes U_E} )\, ,\, 
J_H \rtimes U_{Q \rtimes U_E} )\, $. Put 
$\, L = J_{Q \rtimes U_E} + U_{J_E\, ,\, Q \rtimes U_E}\, ,\, 
M = J_{Q \rtimes U_E} + U_{Q\, ,\, Q \rtimes U_E}\, $ and consider the diagram
\smallskip
$$ \vbox{\halign{ #&#&#\cr 
\hfil $( J_{Q \rtimes U_E}\, ,\, J_H \rtimes U_{Q \rtimes U_E} )$\hfil & 
\hfil $\longrightarrow $\hfil & 
\hfil $( J_{Q\, ,\, H} \rtimes M\, ,\, 
J_H \rtimes U_{Q \rtimes U_E} )$\hfil \cr
\hfil $\Bigm\downarrow $\hfil &&\hfil $\Bigm\downarrow $\hfil \cr
\hfil $( L\, ,\, J_H \rtimes U_{Q \rtimes U_E} )$\hfil & 
\hfil $\longrightarrow $\hfil & 
\hfil $( J_{Q\, ,\, H} \rtimes ( L + M )\, ,\, J_H \rtimes U_{Q \rtimes U_E} )$\hfil \cr
\hfil $\Bigm\downarrow $\hfil &&\hfil $\Bigm\downarrow $\hfil \cr
\hfil $( J_E\, ,\, J_E \rtimes U_H )$\hfil &\hfil $\longrightarrow $\hfil & 
\hfil $( J_E\, ,\, J_E \rtimes U_E )$\hfil \cr }} \> . $$
\par\medskip\noindent
The middle horizontal map and the upper right vertical map are injective in $K^J_2$ by existence of a splitting and both quotients in the lower row are exact. Then from the Mapping Cone Theorem the element in 
$\, K^J_2 ( J_{Q \rtimes U_E} \, ,\, J_H \rtimes U_{Q \rtimes U_E} )\, $ lifts to 
$\, K^J_3 ( J_E\, ,\, J_E \rtimes U_H ) \simeq 
K^J_3 ( J_E\, ,\, J_E \rtimes U_E ) \simeq K^J_3 ( J_E\, ,\, U_E )\, $, these identities following from the fact that $\, J_E\, $ is a free group (see below).  The boundary map to 
$\, K^J_2 ( J_{Q\, ,\, H} \rtimes M\, ,\, J_H \rtimes U_{Q \rtimes U_E} )\, $ factors over 
$\, K^J_3 ( J_E \rtimes U_E )\, $ so the image in this group lifts to the group 
$\, K^J_3 ( J_H \rtimes U_{Q \rtimes U_E} )\, $ by Corollary 5.1 and its image in 
$\, K^J_3 ( J_E \rtimes U_H )\, $ lifts again to 
$\, K^J_3 ( J_E\, ,\, J_E \rtimes U_H )\, $, so that we may assume the original lift in this group to become trivial in $\, K^J_3 ( J_E \rtimes U_E )\, $ by taking the difference with the image of this element. But then considering the compatible splittings 
\smallskip
$$ \vbox{\halign{ #&#&#&#&#\cr
\hfil $( J_E\, ,\, J_E \rtimes U_H )$\hfil &\hfil $\longrightarrow $\hfil &
\hfil $J_E \rtimes U_H$\hfil &\hfil $\longrightarrow $ &
\hfil $U_H$\hfil \cr
\hfil $\Bigm\downarrow\Bigm\uparrow $\hfil &&\hfil $\Bigm\downarrow\Bigm\uparrow $\hfil && 
\hfil $\Bigm\downarrow\Bigm\uparrow $\hfil \cr
\hfil $( J_E\, ,\, J_E \rtimes U_E )$\hfil &\hfil $\longrightarrow $\hfil & 
\hfil $J_E \rtimes U_E$\hfil &\hfil $\longrightarrow $\hfil &\hfil $U_E$\hfil \cr }} $$
\par\medskip\noindent
one gets that already its image in $\, K^J_3 ( J_E \rtimes U_H )\, $ must be trivial, and hence the original element in 
$\, K^J_2 ( J_{Q \rtimes U_E}\, ,\, J_H \rtimes U_{Q \rtimes U_E} )\, $ so that 
$\, ( J_H\, ,\, J_E \rtimes U_H )\, $ is exact, and $\, ( Q , H )\, $ is exact of third order. 
\par\smallskip\noindent
To see that $\, K^J_3 ( J_E\, ,\, J_E \rtimes U_H ) \simeq K^J_3 ( J_E\, ,\, U_E )\, $ assume one is given a regular surjection $\, ( P , E ) \twoheadrightarrow ( P , E / Q )\, $ and that $\, P\, $ is free. Then the surjective map 
$\, K^J_2 ( J_{P , E} , U_E ) \rightarrow K^J_2 ( J_{P , E / Q} , U_{E / Q} )\, $ factors over the injective map to $\, K^J_2 ( J_{P , E / Q}\, ,\, U_{Q , E / P} \rtimes U_{E / Q} )\, $ since the kernel of 
$\, J_{P , E} \twoheadrightarrow J_{P , E / Q}\, $ has a weak core for $\, E\, $. We want to show that the kernel of $\, K^J_2 ( J_{P , E / Q}\, ,\, U_{Q , E / P} \rtimes U_{E / Q} ) \twoheadrightarrow 
K^J_2 ( J_{P , E / Q}\, ,\, U_{E / Q} )\, $ injects into   
$\, K^J_2 ( U_{P , E / Q}\, ,\, U_{Q , E / P} \rtimes U_{E / Q} )\, $. Since both $\, U_{P , E / Q}\, $ and
 $\, U_{Q , E / P}\, $ have a core for the quotient $\, U_{E / (P \times Q )}\, $, Lemmas 11 and 12 show that the group $\, K^J_2 ( U_{P , E / Q}\, ,\, U_{Q , E / P} \rtimes U_{E / Q} )\, $ is naturally isomorphic to the group $\, K^J_2 ( U_{Q , E / P}\, ,\, U_{Q , E / P} \rtimes U_{E / Q} )\, $ and equals the tensor product 
$\, {\mathcal C}_{U_{Q , E / P}} \otimes ( U_{P , E / Q} )^{ab}\, $ where the first factor denotes the abelianization of the core of $\, U_{Q , E / P}\, $. It is also clear that the kernel of 
\smallskip 
$$ K^J_2 ( J_{P , E / Q}\, ,\, U_{Q , E / P} \rtimes U_{E / Q} ) \twoheadrightarrow 
K^J_2 ( J_{P , E / Q}\, ,\, U_{E / Q} ) $$ 
\par\medskip\noindent
is a quotient of 
$\, {\mathcal C}_{U_{Q , E / P}} \otimes ( J_{P , E / Q} )^{ab}\, $ so that since $\, P\, $ is free, any element mapping to zero in $\, {\mathcal C}_{U_{Q , E / P}} \otimes ( U_{P , E / Q} )^{ab}\, $ must come from the image of $\, {\mathcal C}_{U_{Q , E / P}} \otimes 
( [ J_{P , E / Q}\, ,\, U_{P , E / Q} ] / [ J_{P , E / Q}\, ,\, J_{P , E / Q} ] )\, $. The image of this subgroup in  
$\, K^J_2 ( J_{P , E / Q}\, ,\, U_{E / Q} )\, $ is readily seen to be trivial whence the result follows. This implies that the image of $\, K^J_2 ( J_{P , E}\, ,\, U_E )\, $ cannot intersect with this kernel since the map to $\, K^J_2 ( U_{P , E / Q}\, ,\, U_{Q , E / P} \rtimes U_{E / Q} )\, $ factors over 
$\, K^J_2 ( U_{P , E}\, ,\, U_E ) = 0\, $. Note that the proof also applies only assuming that 
$\, K^J_2 ( P ) = 0\, $ (which implies that $\, J_{P , E / Q} \cap [\, U_{P , E / Q}\, ,\, U_{P , E / Q}\, ] = 
[\, J_{P , E / Q}\, ,\, U_{P , E / Q}\, ]\, $).
\par\smallskip\noindent
Another useful result is that $\, ( P\, ,\, U_E / J_{P , E} ) \simeq ( P\, ,\, P \rtimes U_{E / P} )\, $ is always exact, more generally $\, ( P\, ,\, P \rtimes U )\, $ is exact whenever $\, U\, $ is a free group. Namely, consider the split extension
\medskip
$$ 1 \longrightarrow ( J_U \, ,\, (  P \times J_U ) \rtimes U ) \rightarrow\qquad\qquad\qquad\qquad\qquad\qquad\qquad $$
$$ \rightarrow ( J_U \times P \, ,\, 
( P \times J_U ) \rtimes U ) \rightarrow $$
$$\qquad\qquad\qquad\qquad\qquad\qquad\qquad\qquad\qquad\qquad 
\rightarrow ( P\, ,\,  P  \rtimes U ) \longrightarrow 1 $$
\par\bigskip\noindent
from which follows that the map $\, K^J_2 ( J_U\, ,\, 
( P \times J_U ) \rtimes U )\,\rightarrow\, 
K^J_2 ( J_U \times P \, ,\, 
( P \times J_U ) \rtimes U )\, $ is (split) injective since one gets a well defined map from the kernel of 
$\, K^J_2 ( P\, ,\, ( P \times J_U ) \rtimes U )\,\rightarrow\, 
K^J_2 ( P\, ,\, P \rtimes U )\, $ to the kernel of
$\, K^J_2 ( J_U \, ,\, 
( P \times J_U ) \rtimes U )\,\rightarrow\, 
K^J_2 ( J_U \, ,\,  J_U \rtimes U )\, $. A similar argument shows that 
$\, ( R\, ,\, U )\, $ is an exact pair whenever $\, U\, $ is a free group, so in particular the suspension 
$\, ( J_{P , E}\, ,\, U_E )\, $ of any normal pair $\, ( P , E )\, $ is exact. In fact, the pair 
$\, ( R\, ,\, U )\, $ can be shown to be fully exact. Second order exactness is immediate since the pair 
$\, ( J_{U / R}\, ,\, R \rtimes U_{U / R} )\, =\, 
( J_{U / R}\, ,\, J_{U / R} \rtimes U )\, $ is obviously exact. 
If $\, P\, $ and $\, U\, $ are free then the pair $\, ( P\, ,\, P \rtimes U )\, $ is also exact of second order. Consider the diagram
\smallskip
$$ \vbox{\halign{ #&#&#\cr
\hfil $( J_{P \rtimes U}\, ,\, J_U \rtimes U_{P \rtimes U} )$\hfil &\hfil $\longrightarrow $\hfil & 
\hfil $( J_{P \rtimes U} \times J_U\, ,\, J_U \rtimes U_{P \rtimes U} )$\hfil \cr
\hfil $\Bigm\downarrow $\hfil &&\hfil $\Bigm\downarrow $\hfil \cr
\hfil $( J_{P \rtimes U} + U_{P ,\, P \rtimes U} ,\, J_U \rtimes U_{P \rtimes U} )$\hfil & 
\hfil $\longrightarrow $\hfil & 
\hfil $( ( J_{P \rtimes U} + U_{P ,\, P \rtimes U} ) \times J_U ,\, J_U \rtimes U_{P \rtimes U} )$\hfil \cr
\hfil $\Bigm\downarrow $\hfil &&\hfil $\Bigm\downarrow $\hfil \cr
\hfil $( P\, ,\, P \rtimes U_U )$\hfil &\hfil $\longrightarrow $\hfil &\hfil $( P\, ,\, P \rtimes U )$\hfil \cr }}\> . $$
\par\medskip\noindent
The middle horizontal map is injective in $K^J_2$ by existence of a splitting, and both pairs in the bottom row are exact. Since $\, P\, $ is free one gets an isomorphism 
$\, K^J_3 ( P\, ,\, P \rtimes U_U ) \simeq K^J_3 ( P\, ,\, P \rtimes U )\, $. Then by an argument as above 
(splitting of $\, P \rtimes U_U \twoheadrightarrow P \rtimes U\, $) one proves injectivity of the upper horizontal map in $K^J_2$.
\par\smallskip\noindent
Assume that one is given a split extension
\smallskip
$$ 1 \longrightarrow ( N , F ) \longrightarrow ( M , F ) \longleftrightarrow ( P , E ) \longrightarrow 1 $$
\par\medskip\noindent
such that $\, ( N , F )\, $ and $\, ( P , E )\, $ are exact. Then also $\, ( M , F )\, $ is exact. First note that from a simple diagram chase exactness of $\, ( N , F )\, $ entails injectivity of 
$\, K^J_2 ( J_{E / P}\, ,\, M \rtimes U_{E / P} ) \rightarrowtail 
K^J_2 ( J_{E / P} \times N\, ,\, M \rtimes U_{E / P} )\, $. Then any element in the kernel of 
$\, K^J_2 ( J_{E / P}\, ,\, M \rtimes U_{E / P} ) \rightarrow 
K^J_2 ( J_{E / P} \times M\, ,\, M \rtimes U_{E / P} )\, $ must become trivial in 
$\, K^J_2 ( J_{E / P}\, ,\, P \rtimes U_{E / P} )\, $ by exactness of $\, ( P , E )\, $ so its (injective) image in 
$\, K^J_2 ( J_{E / P} \times N\, ,\, M \rtimes U_{E / P} )\, $ lifts to an element of 
$\, K^J_2 ( N\, ,\, M \rtimes U_{E / P} )\, $. By splitting of the quotient map 
$\, ( J_{E / P} \times M\, ,\, M \rtimes U_{E / P} ) \rightarrow ( J_{E / P} \times P\, ,\, P \rtimes U_{E / P} )\, $ the map $\, K^J_2 ( N\, ,\, M \rtimes U_{E / P} ) \rightarrowtail 
K^J_2 ( J_{E / P} \times M\, ,\, M \rtimes U_{E / P} )\, $ is injective whence the result follows. 
\par\smallskip\noindent
Suppose that one is given a regular surjection 
$\, ( P , E )\,\twoheadrightarrow\, ( P , \overline E )\, $ such that the kernel $\, P'\, $ of 
$\, E\,\twoheadrightarrow\, \overline E\, $ has a weak core for $\, E / P\, $. Then, if $\, P\, $ is free and 
$\, ( P , \overline E )\, $ is exact, then also $\, ( P , E )\, $ is exact. Assume that $\, ( P , \overline E )\, $ is exact and $\, P'\, $ is free. Consider the commutative diagram
\smallskip
$$ \vbox{\halign{ #&#&#&#&#\cr
\hfil $\bigl( J_{E / P}\, ,\, P \rtimes U_{E / P} \bigr)$\hfil & \hfil $\rightarrow $\hfil & 
\hfil $\bigl( J_{E / P} + U_{P' , E / P}\, ,\, P \rtimes U_{E / P} \bigr)$\hfil &\hfil $\rightarrow $\hfil & 
\hfil $\bigl( P' , E \bigr)$\hfil \cr
\hfil $\Bigm\downarrow $\hfil && \hfil $\Bigm\downarrow $\hfil &&\cr
\hfil $\bigl( J_{E / P} \times P\, ,\, P \rtimes U_{E / P} \bigr)$\hfil &
\hfil $\rightarrow $\hfil & 
\hfil $\bigl( ( J_{E / P} + U_{P' , E / P} ) \times P\, ,\, P \rtimes U_{E / P} \bigr)$\hfil 
&& \cr
\hfil $\Bigm\downarrow $\hfil && \hfil $\Bigm\downarrow $\hfil &&\cr
\hfil $\bigl( P , E \bigr) $\hfil & \hfil $\rightarrow $\hfil & \hfil $\bigl( P , \overline E \bigr)$\hfil &&\cr }} \> . $$
\par\medskip\noindent
Since $\, ( P' , E / P )\, $ has a weak core for $\, E / P\, $ one gets a $E / P$-normal splitting of the map 
$\, U_{P' , E / P} \rightarrow P'\, $ modulo the subgroup $\, [\, U_{P' , E / P}\, ,\, J_{E / P}\, ]\, $ and it is easy to see that the map $\, K^J_2 ( J_{E / P}\, ,\, P \rtimes U_{E / P} ) \rightarrow 
K^J_2 ( J_{E / P} \times P\, ,\, P \rtimes U_{E / P} )\, $ is injective on the image of 
$\, K^J_2 ( [\, U_{P' , E / P}\, ,\, J_{E / P}\, ]\, ,\, P \rtimes U_{E / P} )\, $ so we may divide the diagram by this image. Then since $\, ( P' , E )\, $ is exact the resulting map on the quotients of 
the first upper horizontal map is injective in $K^J_2$ by existence of a splitting to the following projection in $K^J_3$. By exactness of $\, ( P , \overline E )\, $ the upper vertical map in the middle is injective on 
$K^J_2$-groups, so that also the map 
$\, K^J_2 ( J_{E / P}\, ,\, P \rtimes U_{E / P} ) \rightarrowtail 
K^J_2 ( J_{E / P} \times P\, ,\, P \rtimes U_{E / P} )\, $ must be injective and $\, ( P , E )\, $ is exact.  
\par\smallskip\noindent
It is clear that the pair $\, ( J_E , J_E \rtimes U_E )\, $ is fully exact for arbitrary $\, E\, $, and from this one can induce that also $\, ( J_{P , E}\, ,\, J_E \rtimes U_E )\, $ is exact by an application of the $P \times Q$-Lemma (Lemma 22 of section 6) which states that exactness of
$\, ( J_{P , E}\, ,\, J_E \rtimes U_E )\, $ follows from full exactness of 
$\, ( J^0_E , J_E \rtimes U_E ) \simeq 
( J_E , J_E \rtimes U_E )\, $, and full exactness of
$\, ( J_{P , E}\, ,\, ( J_E \rtimes U_E ) / J^0_E ) \simeq ( J_{P , E} , U_E )\, $ since $\, J^0_E\, $ is free. 
In the same way one proves that the pairs 
$\, ( J^n_{P , E}\, ,\, U^{n+1}_E / J^{n+1}_{E , E} ) \simeq 
( J^n_{P , E}\, ,\, J_{U^n_E / J^n_{E , E}} \rtimes U^n_E )\, $ are exact. For example by the $P \times Q$-Lemma exactness of $\, ( J^2_{P , E}\, ,\, J_{J_E \rtimes U_E} \rtimes U^2_E )\, $ follows from full exactness of the pairs $\, ( J^2_{P , E}\, ,\, U^2_E )\, $ and 
$\, ( J_{J_E \rtimes U_E}\, ,\, J_{J_E \rtimes U_E} \rtimes U^2_E )\, $.    
\par\smallskip\noindent
What is the motivation to consider these cases ? We are preparing for the proof of the Mapping Cone Theorem. Many of the examples considered above are useful in proving exactness of various parts of the Mapping Cone sequence. Consider for instance the extension
\medskip
$$ 1 \rightarrow ( N , U_F / J_{M , F} ) \rightarrow ( N \times J_{E / P}\, ,\, U_F / J_{M , F} ) \rightarrow 
( J_{E / P}\, ,\, U_E / J_{P , E} ) \rightarrow 1 \> . $$
\par\medskip\noindent
When the quotient pair is exact ( so that $\, ( P , E )\, $ is exact of second order one gets an injective map 
$\, K^J_2 ( N\, ,\, U_F / J_{M , F} )\,\rightarrowtail\, K^J_2 ( N \times J_{E / P}\, ,\, U_F / J_{M , F} )\, $ 
which is half of the proof of exactness of the Mapping Cone sequence at 
$\, K^J_2 ( C_{\varphi }\, ,\, U_F )\, $ associated to the inclusion $\, ( N , F ) \subseteq ( M , F )\, $. For the other half one needs to assume that $\, ( M , F )\, $ is exact.
\par\bigskip\noindent
The next result is the major outcome of this section. It seems that there is no analogue of the long exact sequence given by the Theorem for noninjective maps.   
\par\bigskip\noindent
{\bf Theorem 5.}\quad (Mapping Cone Theorem)\quad  Let
$\, ( N , F )\,\buildrel\phi\over\rightarrowtail\, ( M , G )\, $ be an injective morphism and 
$\, (\, C_\phi\, ,\, U_\phi\, )\, $ its mapping cone. Let 
$\, u^n_{\phi }\, $ be the induced map
$\, ( J^n_{N , F}\, ,\, U^n_F )\,\rightarrowtail\, ( J^n_{M , G}\, ,\, U^n_G )\, $. 
For $\, M_F\> =\> {\phi }^{-1} ( M \cap \phi (F) )\, $ 
let $\, \varphi  : ( N , F )\,\subseteq\, ( M_F , F )\, $ be the restriction of $\,\phi\, $ and put 
$\, ( P , E ) := ( M_F / N\, ,\, F / N )\, $. 
Assume that $\, ( M_F , F )\,\subseteq\, ( M , G )\, $ is preexcisive and that 
$\, ( P , E )\, $ is exact. Then the following sequence is halfexact at 
$\, K^J_2 ( N , F )\, $
\smallskip
$$ K^J_2 ( C_{\phi }\, ,\, U_{\phi } ) \longrightarrow  K^J_2 ( N , F ) \longrightarrow K^J_2 ( M , G ) \> . $$
\par\bigskip\noindent
If $\, ( M_F , F )\,\subseteq\, ( M , G )\, $ is preexcisive and the pairs 
$\, ( N , F )\, ,\, ( M_F , F )\, $ and  $\, ( P , E )\, $ are fully exact there is associated the following long exact sequence of $K^J_*$-groups
\smallskip
$$ \>\>\>\cdots\>\>\> \longrightarrow K^J_2 (\, J^{n+2}_{M , G}\, ,\, U_{u^{n+1}_{\phi }}\, )
\buildrel i_*\over\longrightarrow K^J_2 (\, C_{u^{n+1}_{\phi }}\, ,\, U_{u^{n+1}_{\phi }}\, )
\buildrel p_*\over\longrightarrow K^J_{n+3} (\, N , F\, )\, $$
$$ \buildrel {\varphi }_*\over\longrightarrow K^J_2 (\, J^{n+1}_{M , G}\, ,\, U_{u^n_{\phi }}\, )
\buildrel i_*\over\longrightarrow\quad\cdots\quad\>\cdots\quad\>\cdots\quad\>\cdots\quad\>\cdots\quad  
\buildrel p_*\over\longrightarrow K^J_3 (\, N , F \,) $$
$$ \buildrel {\varphi }_*\over\longrightarrow\,
K^J_2 (\, J_{M , G}\, ,\, U_{\phi }\, ) \buildrel i_*\over\longrightarrow 
K^J_2 (\, C_{\phi }\, ,\, U_{\phi } \, ) \buildrel p_*\over\longrightarrow K^J_2 (\, N , F \, )
\buildrel {\varphi }_*\over\longrightarrow K^J_2 (\, M , G \, ) $$
$$ \longrightarrow {\widetilde K}_1 (\, C_{\phi }\, ,\, U_{\phi }\, )\longrightarrow K_1 (\, N , F\, )\longrightarrow K_1 (\, M , G\, )\>\>  .\quad \qquad\qquad\qquad\qquad  $$
\par\bigskip\noindent 
{\it Proof.}\quad Before going on an observation borrowed from the proof of the Excision Theorem gives an injective map
\medskip 
$$ K^J_n ( J_{M_F , F}\, ,\, U_F )\,\rightarrowtail\, K^J_2 ( J^{n-1}_{M , G}\, ,\, U_{u^{n-2}_{\phi }} ) $$
\par\medskip\noindent 
for $\, n \geq 2\, $ since there is a projection 
$\, U_{u^{n-1}_{\phi }} \twoheadrightarrow U^n_F\, $ projecting $\, J^n_{M , G}\, $ onto 
$\, J^n_{M_F , F}\, $. One also gets injective maps 
$\, K^J_2 ( C_{u^n_{\varphi }}\, ,\, U^{n+1}_F )\,\rightarrowtail\, 
K^J_2 ( C_{u^n_{\phi }}\, ,\, U_{u^n_{\phi }} )\, $ by existence of a splitting.
\par\smallskip\noindent
We first show exactness of the last part of the Mapping Cone sequence to the right of 
$\, K^J_2 ( N , F )\, $ which is supposed to hold without any assumptions on the normal pairs involved.
Let us consider the case of $\, K^J_2 ( M , G )\, $. 
$\, {\widetilde C}_\phi\,\cap\, [\, {\overline M}^{ _{\sps G}} , G\, ]\, $ is an $F$-central extension of 
$\, N \cap [ M , G ]\, $, so there is associated a canonical map
\smallskip
$$ \bigl[\, {\overline N}^{ _{\sps F}} , F\,\bigr]\,\buildrel {\overline i}\over\largerightarrow\,
{\widetilde C}_\phi\,\cap\, \bigl[\, {\overline M}^{ _{\sps G}} , G\,\bigr]    $$
\par\medskip\noindent
and the composition
\medskip
$$ K^J_2 ( N , F )\,\subseteq\,
\bigl[\, {\overline N}^{ _{\sps F}} , F\,\bigr]\,\largerightarrow\,
{\widetilde C}_\phi\,\cap\, \bigl[\, {\overline M}^{ _{\sps G}} , G\,\bigr]
\,\largerightarrow\, {\widetilde K}_1 \bigl(\, C_\phi\, ,\, U_\phi\,\bigr)       $$
\par\medskip\noindent
is zero. On the other hand if an element of $\, K^J_2 ( M , G )\, $ is in 
$\, [\, {\widetilde C}_\phi\, ,\, F\, ] = {\overline i} (\, [\, {\overline N}^{ _{\sps F}}\, ,\, F\, ]\, )\, $ it must lie in 
$\, {\overline i} (\, K^J_2 ( N , F )\, )\, $. Next, consider the case of 
$\, {\widetilde K}_1 (\, C_\phi \, ,\, U_\phi\, )\, $. The quotient of this group by the image of 
$\, K^J_2 ( M , G )\, $ is isomorphic to 
$\, (\, N \cap [ M , G ]\, )\, /\, [ N , F ]\, $ which injects into $\, K_1 ( N , F )\, $. The case of $\, K_1 ( N , F )\, $ is now obvious.
\par\smallskip\noindent
We proceed by showing halfexactness at $\, K^J_2 ( N , F )\, $. It is immediate that the composition
\smallskip
$$ K^J_2 ( C_{\phi }\, ,\, U_{\phi } ) \longrightarrow K^J_2 ( N , F ) \longrightarrow K^J_2 ( M , G ) $$ 
\par\bigskip\noindent
is trivial because it factors over the cone $\, ( U_{M , G}\, ,\, U_G )\, $. If an element $\, x\, $ in 
$\, K^J_2 ( N , F )\, $ is in the kernel of the map to $\, K^J_2 ( M , G )\, $ it is zero already in 
$\, K^J_2 ( M_F , F )\, $ by injectivity of $\, K^J_2 ( M_F , F ) \rightarrowtail K^J_2 ( M , G )\, $. Then it suffices to consider the case $\, G = F\, $. Put $\,\overline F = U_F / J_{M_F , F}\, $. Then there is a compatible splitting
\medskip
$$ \bigl( U_{N , F}\, ,\, U_{M_F , F}\, ,\, U_F \bigr) \buildrel\sigma\over\longrightarrow 
\bigl( U_{N , \overline F}\, ,\, U_{M_F ; \overline F}\, ,\, U_{\overline F} \bigr) $$
\par\medskip\noindent
which sends $\, J_{N , F}\, $ to $\, J_{N , \overline F}\, $ and $\, J_{M_F , F}\, $ to 
$\, J_{M_F , \overline F}\, $ lifting $\, x\, $ to an element in $\, K^J_2 ( N , \overline F )\, $ such that its image in $\, K^J_2 ( M_F , \overline F )\, $ is in the kernel of the map to 
$\, K^J_2 ( M_F , F )\, $ and is represented by an element of the group
\medskip
$$ \bigl(\, \bigl[\, J_{M_F , \overline F}\, ,\, U_{\overline F}\,\bigr]\, +\, 
\bigl[\, U_{M_F , \overline F}\, ,\, J_{\overline F}\,\bigr]\, +\,
\bigl[\, U_{M_F , \overline F}\, ,\, U_{J_{F / M_F}\, ,\, \overline F}\,\bigr]\,\bigr)\>\cap\> 
J_{N , \overline F} \> . $$
\par\medskip\noindent
Then construct a compatible splitting
\smallskip
$$ \bigl( U_{N , \overline F}\, ,\, U_{M_F , \overline F}\, ,\, U_{\overline F} \bigr) \buildrel {\sigma }'\over\longrightarrow \bigl( U_{C_{\varphi }\, ,\, U_F}\, ,\, U_{U_{M_F , F}\, ,\, U_F}\, ,\, U^2_F \bigr)\> . $$
\par\medskip\noindent
It will map our representative to an element $\,\xi\, $ of the group
\medskip
$$ \bigl(\, \bigl[\, J_{U_{M_F , F}\, ,\, U_F}\, +\, U_{J_{M_F , F}\, ,\, U_F}\, ,\, U^2_F\,\bigr]\, +\, 
\bigl[\, U_{U_{M_F , F}\, ,\, U_F}\, ,\, J_{U_F}\, +\, U_{J_{M_F , F}\, ,\, U_F}\,\bigr]\, $$
$$\> +\, \bigl[\,  U_{U_{M_F , F}\, ,\, U_F}\, ,\, U_{J_F\, ,\, U_F}\,\bigr]\,\bigr)\>\cap\> 
\bigl[\, U_{C_{\varphi }\, ,\, U_F}\, ,\, U^2_F\,\bigr]  $$
\par\medskip\noindent
so that collapsing to $\, C_{\varphi }\, $ gives an element in
\medskip
$$ \bigl(\, \bigl[\, J_{M_F , F}\, ,\, U_F\,\bigr]\, +\, \bigl[\, U_{M_F}\, ,\, J_F\,\bigr]\, \bigr)\>\cap\> 
\bigl[\, C_{\varphi }\, ,\, U_F\, \bigr] \> . $$
\par\medskip\noindent
Take the difference of $\,\xi\, $ with a lift $\, p\, $ of its image in $\, J_{M_F , F}\, $ to 
$\, U_{J_{M_F , F}\, ,\, U_F}\, $ to get a preimage in $\, J_{C_{\varphi }\, ,\, U_F}\, $ of our representative in 
\medskip
$$ J_{N , \overline F}\>\cap\> \bigl[\, U_{N , \overline F}\, ,\, U_{\overline F}\, \bigr]\> . $$
\par\medskip\noindent
Now assume temporarily that the pair $\, ( J_{M_F , F}\, ,\, U_F )\, $ is partially  admissible on 
$\, R_{0 , M_F , F}\, $ which contains the image of $\, p\, $ (note that this is automatically the case when 
$\, M_F = F\, $). Let $\, D := C_{\varphi } + J_F\, $ and choose equivariant lifts
\medskip
$$ {\overline R}^{(J_{M_F , F}\, ,\, U_F )}_{0 , M_F , F} \buildrel\beta\over\longrightarrow 
{\overline C_{\varphi }}\qquad ,\qquad 
{\overline R}^{(J_{M_F , F}\, ,\, U_F )}_{0 , M_F , F}\buildrel {\beta }_0\over\longrightarrow 
\overline J_F $$
\par\medskip\noindent
where the group on the left is the unique "semiuniversal" $U_F$-central extension of 
$\, R_{0 , M_F , F}\, $ by 
$\, K^J_2 ( J_{M_F , F} , U_F )\, $ contained in the canonical $U_F$-central extension of $\, J_{M_F , F}\, $ and the target groups are some $U_F$-central extensions of $\, C_{\varphi }\, $ and $\, D\, $ by 
$\, K^J_2 ( C_{\varphi }\, ,\, U_F )\, $ and $\, K^J_2 ( D , U_F )\, $ respectively (images of the canonical 
$U_F$-central extensions under some almost canonical map). Also choose some almost canonical maps $\, {\alpha }_0\, ,\, {\alpha }_1\, $ as depicted in the diagram below
\medskip
$$ \vbox{\halign{ #&#&#&#&#\cr
$U_{J_{M_F , F}\, ,\, U_F}$ & $\buildrel {\alpha }_0\over\largerightarrow$ & 
$\overline R^{( J_{M_F , F}\, ,\, U_F )}_{0 , M_F , F}$ & $\buildrel {\beta }_0\over\largerightarrow$ &
$\overline J_F$ \cr
\hfil $\Bigm\downarrow$\hfil && \hfil $\Bigm\downarrow\beta$\hfil && 
\hfil $\Bigm\downarrow {\beta }_2$\hfil \cr
$U_{C_{\varphi }\, ,\, U_F}$ & $\buildrel {\alpha }_1\over\largerightarrow$ & 
\hfil $\overline C_{\varphi }\quad$\hfil & 
$\buildrel {\beta }_1\over\largerightarrow$ & \hfil $\overline D\quad$\hfil \cr }}   $$
\par\medskip\noindent
where the maps $\, {\beta }_1\, ,\, {\beta }_2\, $ are well defined and canonical on commutators with elements from $\, U_F\, $ (and also $\,\overline D\, $ is canonical only on the commutator subgroup 
$\, [\, \overline D\, ,\, U_F\, ]\, $). We claim that under the assumed conditions the image of our representative $\,\xi p^{-1} \in J_{C_{\varphi }\, ,\, U_F}\,  $ mapped to 
$\, K^J_2 ( C_{\varphi }\, ,\, U_F )\, $ under $\, {\alpha }_1\, $ is a preimage for the element 
$\, x \in K^J_2 ( N , F )\, $ up to a difference term which is again in the image of 
$\, K^J_2 ( C_{\varphi }\, ,\, U_F )\, $. Note that the map to $\, K^J_2 ( N , F )\, $ factors over 
$\, K^J_2 ( D , U_F )\, $. Put $\,\overline p = {\alpha }_0 ( p )\, $. We will show that the image of 
$\, \beta ( \overline p )\, $ in $\, K^J_2 ( N , F ) \subseteq [\, \overline N\, ,\, F\, ]\, $ is contained in the image of $\, K^J_2 ( C_{\varphi }\, ,\, U_F )\, $. From this it will follow that the element $\, x\, $ is in the image of $\, K^J_2 ( C_{\varphi }\, ,\, U_F )\, $ since the image of $\,\xi\, $ under $\, {\alpha }_1\, $ followed by the map $\, \overline C_{\varphi } \rightarrow \overline N\, $
is equal to the original element, both maps being canonical on commutators so that the composition 
\medskip
$$ \bigl[\, U_{C_{\varphi }\, ,\, U_F}\, ,\, U^2_F\, \bigr]\> \longrightarrow\> \overline C_{\varphi }\> \longrightarrow\> 
\overline N $$
\par\medskip\noindent
equals
\medskip
$$ \bigl[\, U_{C_{\varphi }\, ,\, U_F}\, ,\, U^2_F\, \bigr]\> \longrightarrow\> 
\bigl[\, U_{N , F}\, ,\, U_F\,\bigr]\> \longrightarrow\> \overline N \> . $$
\par\medskip\noindent
Also note that $\, \beta ( \overline p )\, $ differs from $\, {\alpha }_1 ( p )\, $ only modulo 
$\, K^J_2 ( C_{\varphi }\, ,\, U_F )\, $. We now claim that $\, {\beta }_1 \circ \beta ( \overline p )\, $ differs from $\, {\beta }_2 \circ {\beta }_0 ( \overline p )\, $ only by an element in the image of 
$\, K^J_2 ( C_{\varphi }\, ,\, U_F )\, $ plus the image of $\, K^J_2 ( J_F\, ,\, U_F )\, $ the latter of which maps to zero in $\, K^J_2 ( N , F )\, $. From this the result will follow since 
$\, {\beta }_0 ( \overline p ) \in [\, \overline J_F\, ,\, U_F\, ]\, $ maps to zero in $\, K^J_2 ( N , F )\, $. Let 
$\, \widetilde D\, $ be the subgroup generated by $\, J_F\, $ and 
$\, U_{M_F , F}\, $ so that $\, D \subseteq \widetilde D\, $. One gets that the element 
$\, {\beta }_1 \circ \beta ( \overline p )\, $ maps to zero in 
$\, K^J_2 ( D / C_{\varphi }\, ,\, U_F / C_{\varphi } )\, \subseteq  
\overline { D / C_{\varphi } }^{ _{\sps U_F / C_{\varphi }}}\, $ by halfexactness of the canonical composition
\medskip
$$ \bigl[\, \overline C_{\varphi }\, ,\, U_F\,\bigr]\>\longrightarrow\> 
\bigl[\, \overline D\, ,\, U_F\,\bigr]\>\longrightarrow\> 
\bigl[\, \overline { D / C_{\varphi }}\, ,\, U_F / C_{\varphi }\,\bigr] \> . $$
\par\medskip\noindent
We claim that the same is true for $\, {\beta }_2 \circ {\beta }_0 ( \overline p )\, $ up to an element in the image of $\, K^J_2 ( J_F\, ,\, U_F )\, $ from which follows by exactness that 
$\, {\beta }_2 \circ {\beta }_0 ( \overline p )\, $ can differ from $\, {\beta }_1 \circ \beta ( \overline p )\, $ only by an element in the images of $\, K^J_2 ( C_{\varphi }\, ,\, U_F )\, $ plus $\, K^J_2 ( J_F\, ,\, U_F )\, $ as suggested. The image of $\, p\, $ in 
$\, \overline J_F\, $ is equal to $\, \alpha ( q )\, $ up to an element in $\, K^J_2 ( J_F\, ,\, U_F )\, $ where 
$\,\alpha : U_{J_F\, ,\, U_F} \rightarrow \overline J_F\, $ is an almost canonical map and 
$\, q\, $ can be taken from the subgroup 
\smallskip
$$ \bigl[\, U_{J_{M_F , F}\, ,\, U_F}\, ,\, U^2_F\,\bigr]\, +\, 
\bigl[\, U_{J_F\, ,\, U_F}\, ,\, U_{U_{M_F , F}\, ,\, U_F}\,\bigr]   \> , $$
\par\bigskip\noindent
i.e. it is in the subgroup $\, [\, U_{U_{M_F , F}\, ,\, U_F}\, ,\, U^2_F\, ]\, $ which maps to zero in 
$\, K^J_2 ( D / C_{\varphi }\, ,\, U_{F / M_F} ) \subseteq 
[\, \overline { \widetilde D / U_{M_F , F}}\, ,\, U_{F / M_F}\, ]\, $ by exactness. Then the image of the difference of $\, {\beta }_1 \circ \beta ( \overline p )\, $ and $\,\alpha ( q )\, $ in 
$\, K^J_2 ( D / C_{\varphi }\, ,\, U_F / C_{\varphi } )\, $ is in the kernel of the regular surjection induced by division of a complementary copy of $\, P\, $. Since the map 
$\, K^J_2 ( P\, ,\, U_F / C_{\varphi } )\,\rightarrow K_1 ( C_{\varphi }\, ,\, U_F )\, $ is injective from Lemma 13 and Corollary 5.1 below one also gets injectivity of 
$\, K^J_2 ( P\, ,\, U_F / C_{\varphi } )\,\rightarrow\, 
K^J_2 ( ( D / C_{\varphi } ) \times P\, ,\, U_F / C_{\varphi } )\, $ which factors the former map. This in turn means that one has a well defined projection of the kernel of 
$\, K^J_2 ( D / C_{\varphi }\, ,\, U_F / C_{\varphi } )\,\twoheadrightarrow\, 
K^J_2 ( D / C_{\varphi }\, ,\, U_{F / M_F} )\, $ onto the kernel of 
$\, K^J_2 ( P\, ,\, U_F / C_{\varphi } )\,\twoheadrightarrow\, K^J_2 ( P , E )\, $ (compare with Lemma 11). Suppose this map was an isomorphism. Then also the kernel of the former map would inject into 
$\, K_1 ( C_{\varphi }\, ,\, U_F )\, $, but certainly the image of our difference element is zero in this group because it lifts to an element of $\, K^J_2 ( D\, ,\, U_F )\, $ whose image in 
$\, K_1 ( C_{\varphi }\, ,\, U_F )\, $ is trivial by exactness (see Corollary 5.1 below). Thus we must show injectivity of the map 
\medskip
$$ K^J_2 ( D / C_{\varphi }\, ,\, U_F / C_{\varphi } )\,\longrightarrow\, 
K^J_2 ( ( D / C_{\varphi } ) \times P\, ,\, U_F / C_{\varphi } )  $$
\par\medskip\noindent
which is again a problem of exactness since the mapping cone of this inclusion is $K^J_2$-equivalent to the $J$-suspension of the quotient $\, ( P , E )\, $ and the map 
$\, K^J_2 ( J_{(D / C_{\varphi }) \times P\, ,\, U_F / C_{\varphi }}\, ,\, U_{U_F / C_{\varphi }} )\,\twoheadrightarrow\,
K^J_2 ( J_{P , E}\, ,\, U_E )\, $ is surjective by existence of a splitting (so that it easily follows that
the map $\, K^J_2 ( J_{P , E}\, ,\, U_E ) \,\rightarrow\, K^J_2 ( D / C_{\varphi }\, ,\, U_F / C_{\varphi } )\, $ which is part of the Mapping Cone sequence corresponding to the inclusion as above is trivial). It seems that we havn't gained a lot by this argument, trading one exactness problem for another, but one has to see that the second extension is of a semisplit type (i.e. it splits in the first argument) which gives much more insight in the nature of the problem. In fact, if the extension were completely split, then the injectivity would follow immediately. We may now apply the same machinery used above to analyse the new problem replacing $\, N\, $ by 
$\,\overline N :=  J_F / J_{M_F , F} \simeq J_{E / P}\, $, $\, M\, $ by $\, \overline M := J_{E / P} \times P\, $ and $\, F\, $ by $\, \overline F := U_F / C_{\varphi } \simeq P \rtimes U_{E / P}\, $. However we want to avoid the technical condition that $\, ( J_{\overline M , \overline F}\, ,\, U_{\overline F} )\, $ is partially admissible on $\, R_{0 , \overline M , \overline F}\, $ which is very hard to check in many cases. Luckily it isn't really needed here. Consider the following diagram where $\, C_{\rho }\, $ denotes the mapping cone of the inclusion 
$\, ( \overline N , \overline F ) \subseteq ( \overline M , \overline F )\, $ and $\, J^0_{E / P}\, $ is the complementary copy of $\, J_{E / P}\, $ as in Definition 13
\medskip
$$\vbox{\halign{#&#&#&#&#&#&#&#&#\cr
&&&&&&\hfil ${\overline J}^{ _{\sps U_{\overline F}}}_{\overline F}$\hfil & \hfil $\longrightarrow$\hfil 
& \hfil ${\overline {J^0_{E / P}}}^{ _{\sps U_{\overline F} / J_{\overline M , \overline F}}}$\hfil \cr
&&&&&& \hfil $\Bigm\downarrow$\hfil && \hfil $\Bigm\downarrow\quad$\hfil \cr
&&&& \hfil ${\overline C_{\rho }}^{ _{\sps U_{\overline F}}}$\hfil & \hfil $\longrightarrow $\hfil & 
\hfil ${\overline {J_{\overline F} + C_{\rho }}}^{ _{\sps U_{\overline F}}}$\hfil & 
\hfil $\longrightarrow $ \hfil & \hfil ${\overline {J^0_{E / P}}}^{ _{\sps U_{\overline F} / C_{\rho }}}$\hfil \cr
&&&& \hfil $\Bigm\downarrow $\hfil && \hfil $\Bigm\downarrow $\hfil &&\cr
&&&& \hfil ${\overline {\overline N}}^{ _{\sps U_{\overline F} / J_{\overline M , \overline F}}}$\hfil & 
\hfil $\longrightarrow $\hfil & \hfil ${\overline {\overline N}}^{ _{\sps \overline F}}$\hfil &&\cr 
&&&& \hfil $\Bigm\downarrow $\hfil && \hfil $\Bigm\downarrow $\hfil &&\cr
\hfil  ${\overline J}^{ _{\sps U_{\overline F}}}_{\overline F}$\hfil & \hfil $\longrightarrow $\hfil & 
\hfil ${\overline {J^0_{E / P}}}^{ _{\sps U_{\overline F} / J_{\overline M , \overline F}}}$\hfil & 
\hfil $\longrightarrow $\hfil & 
\hfil ${J_{\overline M , \overline F}\over [\, J_{\overline M , \overline F}\, ,\, U_{\overline F}\, ]}$\hfil &
\hfil $\longrightarrow $\hfil & \hfil ${J_{\overline F}\over [\, J_{\overline F}\, ,\, U_{\overline F}\, ]}$\hfil &&\cr
\hfil $\Bigm\downarrow $\hfil && \hfil $\Bigm\downarrow $\hfil && \hfil $\Bigm\downarrow $\hfil &&&&\cr
\hfil ${\overline {J_{\overline F} + C_{\rho }}}^{ _{\sps U_{\overline F}}}$\hfil & \hfil $\longrightarrow $\hfil & 
\hfil ${\overline {J^0_{E / P}}}^{ _{\sps U_{\overline F} / C_{\rho }}}$\hfil & \hfil $\longrightarrow $\hfil & 
\hfil ${C_{\rho }\over [\, C_{\rho }\, ,\, U_{\overline F}\, ]}$\hfil &&&&\cr  }}\> . $$
\par\medskip\noindent
All rows and columns are exact (on the canonical commutator subgroups). We need to show that the difference $\,\delta\, $ between the (new) elements 
$\, {\beta }_1 \circ \beta ( \overline p)\, $ and $\, {\beta }_2 \circ {\beta }_0 ( \overline p )\, $ which is contained in 
$\, K^J_2 ( J_{\overline F} + C_{\rho }\, ,\, U_{\overline F} ) \subseteq 
[\, {\overline {J_{\overline F} + C_{\rho }}}^{ _{\sps U_{\overline F}}}\, ,\, U_{\overline F}\, ]\, $ maps to an element in the image of $\, K^J_2 ( C_{\rho }\, ,\, U_{\overline F} )\, $ in 
$\, K^J_2 ( \overline N\, ,\, \overline F ) \subseteq 
[\, {\overline {\overline N}}^{ _{\sps \overline F}}\, ,\, \overline F\, ]\, $. In case the equivariant lifts 
$\,\beta\, $ and $\, {\beta }_0\, $ do not exist we may divide $\, \overline D\, $ by the image of 
$\, K^J_2 ( C_{\rho }\, ,\, U_{\overline F} )\, $ and 
$\, K^J_2 ( J_{\overline F}\, ,\, U_{\overline F} )\, $ to get a (noncommutative) diagram as above and measure the difference of $\, {\beta }_1 \circ \beta ( \overline p )\, $ and 
$\, {\beta }_2 \circ {\beta }_0 ( \overline p )\, $ in the corresponding quotient of 
$\, K^J_2 ( D\, ,\, U_{\overline F} )\, $. One notes that this difference is zero in 
$\, K^J_2 ( \overline N\, ,\, \overline F / P ) \simeq K^J_2 ( J_{E / P}\, ,\, U_{E / P} )\, $ because the map
$\, K^J_2 ( \overline N\, ,\, \overline F )\,\rightarrow\, K^J_2 ( \overline N , \overline F / P )\, $ factors over the group $\, K^J_2 ( \overline M , \overline F )\, $ so that 
$\, [\, {\overline C_{\rho }}\, ,\, U_{\overline F}\, ] \longrightarrow 
[\, {\overline {\overline N}}^{U_{E / P}}\, ,\, U_{E / P}\, ]\, $ factors over 
$\, [\, U_{\overline M , \overline F}\, ,\, U_{\overline F}\, ]\, $ where the image of 
$\, {\beta }_1 \circ \beta ( \overline p )\, $ is represented by an element of 
$\, [\, U_{\overline M , \overline F}\, ,\, J_{\overline F}\, ]\, +\, 
[\, J_{\overline M , \overline F}\, ,\, U_{\overline F}\, ]\, $ which maps to zero. Also 
$\, {\beta }_2 \circ {\beta }_0 (\overline p )\, $, which is trivial already in 
$\, K^J_2 ( \overline N , \overline F )\, $ maps to zero. In any case the image of $\,\delta\, $ in 
(the quotient by the image of $\, K^J_2 ( C_{\varphi }\, ,\, U_{\overline F} )\, $ of) 
$\, K^J_2 ( \overline N , \overline F )\, $  can be lifted to an element (in the quotient by the image of 
$\, K^J_2 ( C_{\varphi }\, ,\, U_{\overline F} )\, $) of 
$\, K^J_2 ( \overline N , U_{\overline F} / J_{\overline M , \overline F} )\, $ whose image in 
$\, K_1 ( J_{\overline M , \overline F}\, ,\, U_{\overline F} )\, $ is in the kernel of the map to 
$\, K_1 ( J_{\overline F}\, ,\, U_{\overline F} )\, $. This means that as an element of the group 
$\, \overline {U_{\overline F}} := U_{\overline F} / [\, J_{\overline M , \overline F}\, ,\, U_{\overline F}\, ]\, $ it is represented as an element of 
\smallskip
$$  \bigl[\, \overline {J^0_{E / P}}\, ,\, \overline { ( J_{E / P} \times P ) \rtimes U_{E / P}}\, \bigr]\> \cap\> 
\bigl[\, \overline J_{E / P}\, ,\, \overline { ( J_{E / P} \times P ) \rtimes U_{E / P}}\, \bigr]\>  .   $$
\par\medskip\noindent
The image of (any) lift in 
$\, K^J_2 ( \overline N , U_{\overline F} / ( J_{\overline M , \overline F} + U_{P\, ,\, \overline F} ) )\, $ lies in the kernel of the map to $\, K^J_2 ( \overline N , U_{E / P} )\, $ so that its image in 
$\, K_1 ( J_{\overline N\, ,\, U_{E / P}}\, ,\, U_{ U_{E / P}} )\, $ is represented by an element of 
$\, [\, \overline {J^0_{E / P}}\, ,\, \overline {J_{E / P}}\, ]\, $ which by splitting of the map 
$\, \overline F\,\rightarrow\, U_{E / P}\, $ means the element in 
$\, K_1 ( J_{\overline M , \overline F}\, ,\, U_{\overline F} )\, $ is represented by an element of 
$\, [\, \overline {J^0_{E / P}}\, ,\, \overline { J_{E / P} \times P}\, ]\, =\, 
[\, \overline {J^0_{E / P}}\, ,\, \overline P\, ]\, +\, [\, \overline {J^0_{E / P}}\, ,\, \overline {J_{E / P}}\, ]\, $. The second part can be lifted to an element in the kernel of the map
$\, K^J_2 ( \overline N , U_{\overline F} / J_{\overline M , \overline F} )\,\twoheadrightarrow\, 
K^J_2 ( \overline N , \overline F )\, $ so without loss of generality we may assume that the image of the lift of $\,\delta\, $ in $\, K_1 ( J_{\overline M , \overline F}\, ,\, U_{\overline F} )\, $ by the procedure above is represented by an element of $\, [\, \overline {J^0_{E / P}}\, ,\, \overline P\, ]\, $ which can be lifted to an element of $\, K^J_2 ( P\, ,\, U_{\overline F} / J_{\overline M , \overline F} )\, $ in the kernel of the map to 
$\, K^J_2 ( P\, ,\, \overline F )\, $ and in the kernel of the map to
$\, K^J_2 ( P\, ,\, U_{\overline F} / C_{\rho } )\, $ since the latter group maps injectively to 
$\, K_1 ( C_{\rho }\, ,\, U_{\overline F} )\, $, i.e. it defines an element in the trefoil intersection for the pair 
$\, ( P , E )\, $ in the sense of Definition 13. By assumption this element lifts to an element in the trefoil intersection of the corresponding extended $K^J_2$-groups. Then one gets a well defined map to the kernel of   
$\, K^J_2 ( J^0_{E / P}\, ,\, U_{\overline F} / J_{\overline M , \overline F} )\,\twoheadrightarrow\, 
K^J_2 ( J^0_{E / P}\, ,\, U_{\overline F} / C_{\rho } )\, $ and the "image" of our element $\,\delta\, $ in 
$\, K_1 ( J_{\overline M , \overline F}\, ,\, U_{\overline F} )\, $ lifts to an element of 
$\, K^J_2 ( J^0_{E / P}\, ,\, U_{\overline F} / J_{\overline M , \overline F} )\, $ in this kernel. Let 
$\,{\delta }_0\, $ denote the "image" of $\,\delta\, $ in 
$\, K^J_2 ( {\overline N}\, ,\, U_{\overline F} / J_{\overline M , \overline F} )\, $ (i.e. any lift of the "image" of $\,\delta\, $ in $\, K_1 ( J_{\overline M , \overline F}\, ,\, U_{\overline F} )\, $ in this group will do) and 
$\, {\delta }_1\, $ the "image" of $\,\delta\, $ in 
$\, K^J_2 ( J^0_{E / P}\, ,\, U_{\overline F} / J_{\overline M , \overline F} )\, $ as above and consider their difference $\,\Delta\, $ in the group 
$\, K^J_2 ( J_{E / P} \times J^0_{E / P}\, ,\, U_{\overline F} / J_{\overline M , \overline F} )\, $. Then 
$\,\Delta\, $ maps to zero in $\, K_1 ( J_{\overline M , \overline F}\, ,\, U_{\overline F} )\, $ by the natural boundary map, it also maps to zero in $\, K^J_2 ( J^0_{E / P}\, ,\, U_{\overline F} / C_{\rho } )\, $ and it maps to the image of $\,\delta\, $ in $\, K^J_2 ( \overline N\, ,\, \overline F )\, $ modulo the image of 
$\, K^J_2 ( C_{\varphi }\, ,\, U_{\overline F} )\, $. By the first property it can be lifted to an element of 
$\, K^J_2 ( C_{\rho } + J_{\overline F}\, ,\, U_{\overline F} )\, $ which by the second property can be lifted to 
$\, K^J_2 ( C_{\rho }\, ,\, U_{\overline F} )\, $ showing that the image of 
$\,\delta\, $ in $\, K^J_2 ( \overline N , \overline F )\, $ is in the image of 
$\, K^J_2 ( C_{\rho }\, ,\, U_{\overline F} )\, $ as suggested. Replacing $\, M_F\, $ with $\, F\, $ and using the exactness results given above for the associated extension at $\, K^J_2 ( N , F )\, $ together with exactness at $\, K^J_3 ( E / P )\, $ for the sequence associated with the extension $\, E\, $ given below and a diagram chase argument one sees that the working assumption that 
$\, ( J_{M_F\, ,\, F}\, ,\, U_F )\, $ be partially admissible on $\, R_{0 , M_F , F}\, $ is not really needed.
\par\smallskip\noindent
To prove exactness of the whole Mapping Cone sequence we note that $\, ( J_F\, ,\, U_F )\, $ is partially admissible on $\, R_{0 , F , F}\, $ and that the pairs $\, ( N , F )\, ,\, ( M_F , F )\, $ and $\, ( P , E )\, $ are exact. We first do the special case $\, M = F = G\, $ and begin with exactness at 
$\, K^J_n ( F , F )\, $ for $\, n \geq 3\, $. Let $\, u^n_{\varphi }\, $ denote the inclusion 
$\, ( J^n_{N , F}\, ,\, U^n_F ) \subseteq ( J^n_{F , F}\, ,\, U^n_F )\, $. It is clear that the composition 
\smallskip
$$ K^J_2 ( J_{J^n_{N , F}\, ,\, U^n_F}\, ,\, U^{n+1}_F )\,\longrightarrow\, 
K^J_2 ( J_{J^n_{F , F}\, ,\, U^n_F}\, ,\, U^{n+1}_F )\,\longrightarrow\, 
K^J_2 ( C_{u^n_{\varphi }}\, ,\, U^{n+1}_F ) $$
\par\medskip\noindent
is zero because it factors over the cone $\, K^J_2 ( U_{J^n_{N , F}\, ,\, U^n_F}\, ,\, U^{n+1}_F )\> =\> 0\, $ by Lemma 13. Let 
$\, C^n_j\, $ be the mapping cone of 
$\, (\, J_{J^n_{F , F}\, ,\, U^n_F}\, ,\, U^{n+1}_F\, )\,\subseteq\, (\, C_{u^n_{\varphi }}\, ,\, U^{n+1}_F\, )\, $ and $\, D^n_j\, $ the subgroup of 
$\, C^n_j\, $ generated by $\, U_{J_{J^n_{N , F}\, ,\, U^n_F}\, ,\, U^{n+1}_F} \, $ and 
$\, J_{U_{J^n_{N , F}\, ,\, U^n_F}\, ,\, U^{n+1}_F}\, $. Since 
$\, ( J_{C_{u^n_{\varphi }}\, ,\, U^{n+1}_F}\, ,\, U^{n+2}_F )\, $ is admissible by Lemma 14 and 
$\, ( N , F )\, $ is exact by assumption (respectively for $\, n \geq 1\, $ the pair 
$\, ( C_{u^n_{\varphi }}\, /\, J^{n+1}_{F , F}\, ,\, U^{n+1}_F\, /\, J^{n+1}_{F , F} )\, $ is of the form 
$\, ( J^n_{N , F}\, ,\, U^{n+1}_F / J^{n+1}_{F , F} )\, $) an element in the kernel of the map 
$\, K^J_2 ( J^{n+1}_{F , F}\, ,\, U^{n+1}_F )\,\rightarrow\, 
K^J_2 ( C_{u^n_{\varphi }}\, ,\, U^{n+1}_F )\, $ can be lifted to $\, K^J_2 ( C^n_j\, ,\, U^{n+2}_F )\, $.
One gets a natural diagram
\smallskip
$$  \vbox{\halign{ #&#&#\cr
$ \bigl(\, D^n_j\, /\, U_{J^{n+1}_{N , F}\, ,\, U^{n+1}_F}\, ,\, U_{U^{n+1}_F\, /\, J^{n+1}_{N , F}}\,\bigr) $ & $ \longrightarrow $ & 
$ \bigl(\, C^n_j\, /\, U_{J^{n+1}_{F , F}\, ,\, U^{n+1}_F}\, ,\, U_{U^{n+1}_F\, /\, J^{n+1}_{F , F }}\,\bigr) $ \cr
\hfil $ \Bigm\Vert $\hfil && \hfil $ \Bigm\Vert $\hfil \cr
$ \bigl(\, J_{J^n_{N , F}\, ,\, U^{n+1}_F\, /\, J^{n+1}_{N , F}}\, ,\, U_{U^{n+1}_F\, /\, J^{n+1}_{N , F}}\,\bigr) $ & $ \longrightarrow $ & 
$ \bigl(\, J_{J^n_{N , F}\, ,\, U^{n+1}_F\, /\, J^{n+1}_{F , F}}\, ,\, U_{U^{n+1}_F\, /\, J^{n+1}_{F , F}}\,\bigr) $ \cr  }}  $$
\smallskip
so the horizontal map is surjective. But the group on the left side is naturally isomorphic with 
 $\, K^J_2 (\, D^n_j\, ,\, U^{n+2}_F\, )\, $ and the group on the right is isomorphic to 
 $\, K^J_2 (\, C^n_j\, ,\, U^{n+2}_F\, )\, $. Then the image of 
 $\, K^J_2 (\, C^n_j\, ,\, U^{n+2}_F\, )\, $ in $\, K^J_2 (\, J^{n+1}_{F , F}\, ,\, U^{n+1}_F\, )\, $ is the same as the image of 
 $\, K^J_2 (\, D^n_j\, ,\, U^{n+2}_F\, )\, $ , and the map 
 $\, (\, D^n_j\, ,\, U^{n+2}_F\, )\,\longrightarrow\, (\, J^{n+1}_{F , F}\, \, U^{n+1}_F\, )\, $ factors over 
 $\, (\, J^{n+1}_{N , F}\, ,\, U^{n+1}_F\, )\, $, whence the result. 
 \par\smallskip\noindent
Next we do exactness at $\, K^J_2 ( C_{u^n_{\varphi }}\, ,\, U^{n+1}_F )\, $.
By the exact sequence
\smallskip
$$ 1 \longrightarrow \bigl(\, J^{n+1}_{F , F}\, ,\, U^{n+1}_F\, \bigr)
\longrightarrow \bigl(\, C_{u^n_{\varphi }}\, ,\, U^{n+1}_F\, \bigr) \longrightarrow 
\bigl(\, J^n_{N , F}\, ,\, U^{n+1}_F / J^{n+1}_{F , F}\, \bigr) \longrightarrow 1   $$
\par\medskip\noindent
and Lemma 9 one must prove that the image of the quotient map in 
$\, K^J_2 (\, J^n_{N , F} , U^{n+1}_F\, /\, J^{n+1}_{F , F}\, )\, $ intersects trivially with the kernel of the regular surjection 
\smallskip
$$K^J_2 \bigl(\, J^n_{N , F}\, ,\, J_{U^n_F / J^n_{F , F}} \rtimes U^n_F\, \bigr)\,
\longrightarrow\, K^J_2 \bigl( J^n_{N , F}\, ,\, U^n_F \bigr)\>   .$$ 
\par\medskip\noindent
This is checked step by step starting with the lowest dimension $\, n = 1\, $ (the case $\, n = 0\, $ is of course trivial).  For $\, n = 1\, $ consider the diagram
\smallskip
$$ \vbox{\halign{ #&#&#&#&#\cr
\hfil $\bigl( J_{N , F}\, ,\, J^0_F \rtimes U_F \bigr)$\hfil &\hfil $\longrightarrow $\hfil &\hfil 
$\bigl( J_{N , F} \times J^0_F\, ,\, J^0_F \rtimes U_F \bigr)$\hfil &\hfil $\longrightarrow $\hfil & 
\hfil $\bigl( J^0_F\, ,\, J^0_F \rtimes U_F / J_{N , F} \bigr)$\hfil \cr
\hfil $\Bigm\downarrow $\hfil &&\hfil $\Bigm\downarrow $\hfil &&\hfil $\Bigm\downarrow $\hfil \cr
\hfil $\bigl( J_F\, ,\, J^0_F \rtimes U_F \bigr)$\hfil &\hfil $\longrightarrow $\hfil &
\hfil $\bigl( J_F \times J^0_F\, ,\, J^0_F \rtimes U_F \bigr)$\hfil &\hfil $\longrightarrow $\hfil & 
\hfil $\bigl( J_F\, ,\, U_F \bigr)$\hfil \cr
\hfil $\Bigm\downarrow $\hfil &&\hfil $\Bigm\downarrow $\hfil &&\cr
\hfil $\bigl( J_E\, ,\, J_E \rtimes U_F \bigr)$\hfil &\hfil $\longrightarrow $\hfil &
\hfil $\bigl( J_E\, ,\, U_F / J_{N , F} \bigr)$\hfil &&\cr }} \> . $$
\par\medskip\noindent
By existence of a splitting the first horizontal map in the middle row is injective, so any element in the kernel of the first upper horizontal map lifts to $\, K^J_3 ( J_E\, ,\, J_E \rtimes U_F )\, $ by exactness of this pair, which since $\, J_E\, $ is free is isomorphic to $\, K^J_3 ( J_E\, ,\, U_F / J_{N , F} ) \simeq 
K^J_2 ( C_{u_{\varphi }}\, ,\, U^2_F )\, $. By exactness of $\, ( J_F \times J^0_F\, ,\, J^0_F \rtimes U_F )\, $ and second order exactness of $\, ( J_E\, ,\, U_F / J_{N , F} )\, $ such an element lifts again to the 
$K^J_3$-group of the former pair and the image of this lift in 
$\, K^J_3 ( J_F\, ,\, U_F )\, $ can be lifted to $\, K^J_3 ( J_{N , F} \times J^0_F\, ,\, J^0_F \rtimes U_F )\, $ so that taking the difference with the image of this lift one can assume that the lift in 
$\, K^J_3 ( J_F \times J^0_F\, ,\, J^0_F \rtimes U_F )\, $ lifts to 
$\, K^J_3 ( J_F\, ,\, J^0_F \rtimes U_F )\, $ and by taking the difference with the image of this lift one gets a lift of the original element which is trivial, hence the map 
$\, K^J_2 ( J_{N , F}\, ,\, J^0_F \rtimes U_F ) \rightarrowtail 
K^J_2 ( J_{N , F} \times J^0_F\, ,\, J^0_F \rtimes U_F )\, $ is injective. 
This implies that the kernel of 
$\, K^J_2 ( J_{N , F}\, ,\, U^2_F / J^2_{F , F} )\,\twoheadrightarrow\, K^J_2 ( J_{N , F}\, ,\, U_F )\, $ is contained in the kernel of 
$\, K^J_2 ( J_F\, ,\, U^2_F / J^2_{F , F} )\,\twoheadrightarrow\, K^J_2 ( J_F\, ,\, U_F )\, $ since the latter is equal to the kernel of 
$\, K^J_2 ( J^0_F\, ,\, U^2_F / J^2_{F , F} )\,\twoheadrightarrow\, K^J_2 ( J_F\, ,\, U_F )\, $ by symmetry, and the kernel of 
$\, K^J_2 ( J^0_F\, ,\, U^2_F / J^2_{F , F} )\,\twoheadrightarrow K^J_2 ( J^0_F\, ,\, J_E \rtimes U_F )\, $ projects onto (in fact is equal to) the kernel of the first map. But the map from 
$\, K^J_2 ( J_{C_{\varphi }\, ,\, U_F}\, ,\, U^2_F )\, $ to 
$\, K^J_2 ( J_F\, ,\, U^2_F / J^2_{F , F} )\, $ factors over 
$\, K^J_2 ( J_{U_F}\, ,\, U^2_F ) = 0\, $, settling the case $\, n = 1\, $. The cases of higher 
$\, n\, $ are treated similar. For example, if $\, n = 2\, $ one reduces injectivity of 
$\, K^J_2 ( J^2_{N , F}\, ,\, J_{J_F \rtimes U_F} \rtimes U^2_F ) \rightarrow 
K^J_2 ( J^2_{N , F} \times J_{J_F \rtimes U_F}\, ,\, J_{J_F \rtimes U_F} \rtimes U^2_F )\, $ to injectivity of 
$\, K^J_2 ( J^2_{N , F}\, ,\, J_{U_F} \rtimes U^2_F ) \rightarrowtail 
K^J_2 ( J^2_{N , F} \times J_{U_F}\, ,\, J_{U_F} \rtimes U^2_F )\, $ shown above. This means that an element in the kernel of the first map lifts to an element in the kernel of 
\smallskip
$$ K^J_2 ( J_{J_F\, ,\, J_F \rtimes U_F}\, ,\, J_{J_F \rtimes U_F} \rtimes U^2_F ) \twoheadrightarrow 
K^J_2 ( J_{J_F\, ,\, J_F \rtimes U_F}\, ,\, J_{J_F \rtimes U_F} \rtimes U^2_F / J^2_{N , F} ) $$
\par\medskip\noindent
since the latter is contained in the kernel of 
\smallskip
$$ K^J_2 ( J_{J_F\, ,\, J_F \rtimes U_F}\, ,\, J_{J_F \rtimes U_F} \rtimes U^2_F ) \twoheadrightarrow 
K^J_2 ( J_{J_F\, ,\, J_F \rtimes U_F}\, ,\, J_{J_F \rtimes U_F} \rtimes U_F ) $$
\par\medskip\noindent
which is naturally isomorphic to the tensor product of the core of the abelianization of 
$\, J_{J_F\, ,\, J_F \rtimes U_F}\, $ for the quotient $\, U_F\, $ with the abelianization of $\, J_{U_F}\, $, and since the quotient $\, J_{U_F / J_{N , F}}\, $ is free one gets the result by a standard argument (see above). But the map 
\smallskip
$$ K^J_2 ( J_{J_F\, ,\, J_F \rtimes U_F}\, ,\, J_{J_F \rtimes U_F} \rtimes U^2_F ) \rightarrowtail 
K^J_2 ( J_{J_F\, ,\, J_F \rtimes U_F} \times J_{U_F}\, ,\, J_{J_F \rtimes U_F} \rtimes U^2_F ) $$
\par\medskip\noindent
is injective by existence of a splitting. Now one proceeds as above. The cases of higher $\, n\, $ are left to the reader.
\par\smallskip\noindent
Consider next the case of $\, K^J_n ( N , F )\, $. For $\, n = 2\, $ the proof above applies to give exactness at $\, K^J_2 ( N , F )\, $ by the fact that $\, ( J_F\, ,\, U_F )\, $ is partially admissible on 
$\, R_{0 , F , F}\, $. For $\, n = 3\, $ $\, ( J_F / J_{N , F}\, ,\, U_F / J_{N , F} ) \simeq 
( J_E\, ,\, N \rtimes U_E )\, $ is exact by assumption that $\, ( N , F )\, $ is exact of second order. In higher dimensions the quotients $\, ( J^n_{F , F} / J^n_{N , F}\, ,\, U^n_F / J^n_{N , F} )\, $ are all of the form 
$\, ( J_{J^{n-1}_{F , F} / J^{n-1}_{N , F}\, ,\, U^{n-1}_F / J^{n-1}_{N , F}}\, ,\, 
J_{U^{n-1}_F / J^{n-1}_{N , F}} \rtimes U^{n-1}_F )\, $. One gets a copy of $\, J^{n-1}_{N , F}\, $ complementary to the normal factor $\, J_{U^{n-1}_F / J^{n-1}_{N , F}}\, $ in each case, so the group 
$\, J_{U^{n-1}_F / J^{n-1}_{N , F}} \rtimes U^{n-1}_F\, $ can also be written as 
$\, J^{n-1}_{N , F} \rtimes U_{U^{n-1}_F / J^{n-1}_{N , F}}\, $. Then the $P \times Q$-Lemma yields exactness of these pairs by full exactness of the pairs 
$\, ( J^{n-1}_{N , F}\, ,\, J^{n-1}_{N , F} \rtimes U_{U^{n-1}_F / J^{n-1}_{N , F}} )\, $ and 
$\, ( J_{U^{n-1}_F / J^{n-1}_{N , F}}\, ,\, U_{U^{n-1}_F / J^{n-1}_{N , F}} )\, $. Then from Lemma 14 exactness of the special Mapping Cone sequence at $\, K^J_n ( N , F )\, $ follows. 
\par\smallskip\noindent
Next, consider the case $\, M = M_F\, $ and $\, F = G\, $. Let $\, u^n_{\varphi }\, $ denote the inclusion 
$\, ( J^n_{N , F}\, ,\, U^n_F )\,\subseteq\, ( J^n_{M , F}\, ,\, U^n_F )\, $. It follows from the proof of 
Corollary 5.1 below that the natural split surjective map 
\medskip
$$ K^J_2 ( C_{u^{n-1}_{\varphi }}\, ,\, U^n_F )\,\rightarrow\, K^J_2 ( J^n_{P , E}\, ,\, U^n_E )  $$
\par\medskip\noindent  
is an isomorphism. Then the exactness of the Mapping Cone sequence for the inclusion 
$\, ( N , F )\,\subseteq\, ( M , F )\, $ follows from exactness of the respective sequences for the inclusions 
$\, ( N , F )\,\subseteq\, ( F , F )\, ,\, ( M , F )\,\subseteq\, ( F , F )\, $ and $\, ( P , E )\,\subseteq\, ( E , E )\, $ respectively, by a diagram chase.
\par\smallskip\noindent
Finally consider the general case. The composition
\medskip
$$ K^J_2 \bigl(\, C_{u^n_{\phi }}\, ,\, U_{u^n_{\phi }}\, \bigr)\,\longrightarrow\, 
K^J_{n+2} \bigl(\, N\, ,\, F\, \bigr)\,\longrightarrow\, 
K^J_2 \bigl(\, J^n_{M , G}\, ,\, U_{u^{n-1}_{\phi }}\,\bigr)   $$
\par\medskip\noindent
is trivial because the map $\, K^J_2 ( C_{u^n_{\phi }}\, ,\, U_{u^n_{\phi }} ) \rightarrow 
K^J_{n+2} ( N , F )\, $ factors over the group $\, K^J_2 ( C_{u^n_{\varphi }}\, ,\, U^{n+1}_F )\, $ (there is a splitting of the natural inclusion $\, ( C_{u^n_{\varphi }}\, ,\, U^{n+1}_F ) \subseteq 
( C_{u^n_{\phi }}\, ,\, U_{u^n_{\phi }} )\, $ compatible with the evaluation map), and the image of this group in 
$\, K^J_{n+2} ( M_F , F )\, $ is trivial, while the map $\, K^J_{n+2} ( N , F ) \rightarrow 
K^J_2 ( J^n_{M , G}\, ,\, U_{u^{n-1}_{\phi }} )\, $ factors over the latter. 
Then if an element from $\, K^J_{n+2} ( N , F )\, $ maps to zero in 
$\, K^J_2 (\, J^n_{M , G}\, , U_{u^{n-1}_{\phi }} )\, $ it is already trivial in $\, K^J_{n+2} (\, M_F\, , F\, )\, $ by existence of a splitting for the inclusion 
$\, ( J^n_{M_F , F}\, ,\, U^n_F ) \subseteq ( J^n_{M , G}\, ,\, U_{u^{n-1}_{\phi }} )\, $.  So we have reduced exactness at $\, K^J_n ( N , F )\, $ to the case 
$\, M = M_F\, ,\, F = G\, $. It is clear that the composition 
\smallskip
$$ K^J_{n+2} \bigl(\, N\, ,\, F\, \bigr)\,\longrightarrow\, K^J_2 \bigl(\, J^n_{M , G}\, ,\, U_{u^{n-1}_{\phi }}\, \bigr)\,\longrightarrow\, K^J_n \bigl(\, C_{u^{n-1}_{\phi }}\, ,\, U_{u^{n-1}_{\phi  }}\,\bigr)     $$
\smallskip
is zero because it factors over $\, K^J_2 (\, U_{J^{n-1}_{N , F}\, ,\, U^{n-1}_F}\, ,\, U^n_F\, )\> =\> 0\, $. Since the map 
$\, (\, J^n_{N , F}\, ,\, U^n_F\, )\,\rightarrow\, (\, J^n_{M , G}\, ,\, U_{u^{n-1}_{\phi }}\, )\, $ factors over 
$\, (\, J^n_{M_F , F}\, ,\, U^n_F\, )\, $ it suffices to consider the case $\, M = M_F\, ,\, F = G\, $, and the case 
$\, N = M_F\, $ separately and then glue them together to get the general result. Assume that
$\, N = M_F\, $. Let $\,{\rho }^n\, $ denote the inclusion 
$\, (\, J^n_{M_F , F}\, ,\, U^n_F\, )\,\subseteq\, (\, J^n_{M , G}\, ,\, U^n_G\, )\, $ and 
$\, C_{j^n}\, $ be the mapping cone of the inclusion
$\, (\, J^{n+1}_{M , G}\, ,\, U_{{\rho }^n}\, )\,\subseteq\, (\, C_{{\rho }^n}\, ,\, U_{{\rho }^n}\, )\, $. 
By our argument above an element in the kernel of 
$\, K^J_2 (\, J^{n+1}_{M , G}\, ,\, U_{{\rho }^n}\, )\,\rightarrow\, K^J_2 ( C_{{\rho }^n}\, ,\, U_{{\rho }^n} )\, $ can be lifted to 
$\, K^J_2 (\, C_{j^n}\, ,\, U_{U_{{\rho }^n}}\, )\, $. The map from 
$\, K^J_2 (\, C_{j^n}\, ,\, U_{U_{{\rho }^n}} )\, $ to $\, K^J_2 ( J^{n+1}_{M , G}\, ,\, U_{{\rho }^n} )\, $ factors over $\, K^J_2 ( J^{n+1}_{M_F , F}\, ,\, U^{n+1}_F )\, $ because $\, (\, C_{j^n}\, ,\, U^2_{{\rho }^n}\, )\, $ is the extension of 
$\, (\, J_{J^n_{M_F\, ,\, F}\, ,\, U_{{\rho }^n}\, /\, J^{n+1}_{M , G}}\, ,\, 
U_{U_{{\rho }^n}\, /\, J^{n+1}_{M , G}}\, )\, $ by the cone 
$\, (\, U_{J^{n+1}_{M , G}\, ,\, U_{{\rho }^n} }\, ,\, U^2_{{\rho }^n} \, )\, $ so that the induced map on 
$K^J_2$-groups is injective (by Lemma 9) and surjective (by existence of a splitting). Then the projection 
\medskip
$$  K^J_2 (\, J_{J^n_{M_F\, ,\, F}\, ,\, U_{{\rho }^n}\, /\, J^{n+1}_{M , G}}\, ,\, 
U_{U_{{\rho }^n}\, /\, J^{n+1}_{M , G}}\, )\,\longrightarrow\,
K^J_2 (\, J^{n+1}_{M_F , F}\, ,\, U^{n+1}_F\, )   $$
\par\medskip\noindent
makes a commutative diagram with the projection  
\medskip
$$ K^J_2 ( J_{J^n_{M , G}\, ,\, U_{{\rho }^n} / J^{n+1}_{M , G}}\, ,\, U_{U_{{\rho }^n} / J^{n+1}_{M , G}} ) \longrightarrow K^J_2 ( J^{n+1}_{M , G}\, ,\, U_{{\rho }^n} )\> , $$
\par\medskip\noindent 
settling the case 
$\, N = M_F\, $. Then as
$\, ( C_{u^{n-1}_{\varphi }}\, ,\, U^n_F )\,\subseteq\, ( C_{u^{n-1}_{\phi }}\, ,\, U_{u^{n-1}_{\phi }} )\, $ admits a splitting the map 
$\, K^J_2 (\, C_{u^{n-1}_{\varphi }}\, ,\, U^n_F\, )\,\rightarrowtail\, 
K^J_2 (\, C_{u^{n-1}_{\phi }}\, ,\, U_{u^{n-1}_{\phi }}\, )\, $ is injective for all $\, n \geq 2\, $. Then one has reduced to the case $\, M = M_F\, ,\, F = G\, $. Exactness at 
$\, K^J_2 ( C_{u^n_{\phi }}\, ,\, U_{u^n_{\phi }} )\, $ follows since the kernel of
\medskip
$$ K^J_2 ( J^n_{N , F}\, ,\, U_{u^n_{\phi }} / J^{n+1}_{M , G} )\,\twoheadrightarrow\, 
K^J_2 ( J^n_{N , F}\, ,\, U^n_F ) $$
\par\medskip\noindent
is contained in the kernel of 
\medskip
$$ K^J_2 ( J^n_{M_F , F}\, ,\, U_{u^n_{\phi }} / J^{n+1}_{M , G} )\,\twoheadrightarrow\, 
K^J_2 ( J^n_{M_F , F}\, ,\, U^n_F ) $$
\par\medskip\noindent
by our argument above and the map to the latter group factors over 
$\, K^J_2 ( C_{{\rho }^n}\, ,\, U_{{\rho }^n } )\, $ which is trivial by the special case of the long exact sequence applied to the extension
\medskip
$$ 1 \longrightarrow ( C_{{\rho }^n}\, ,\, U_{{\rho }^n } ) \longrightarrow 
( U_{{\rho }^n}\, ,\, U_{{\rho }^n} ) \longrightarrow ( U_{{\overline\rho }^n}\, ,\, U_{{\overline\rho }^n} ) 
\longrightarrow 1 $$
\par\medskip\noindent
where $\,\overline {\rho }^n\, $ is the quotient inclusion $\, U^n_{F / M_F}\,\subseteq\, U^n_{G / M}\, $ 
\qed
\par\bigskip\noindent
An important special case of the Theorem is when $\, G = F\, $. Then it makes sense to compare the 
$\, K^J$-theory of $\, ( P , E ) = ( M/N , F/N )\, $ and of $\, (\, C_\varphi\, ,\, U_\varphi\, )\, $. One gets the following
\par\bigskip\noindent
{\bf Corollary 5.1.}\quad Assume given an exact sequence
\medskip
$$1\,\largerightarrow\, ( N , F )\,\largerightarrow\, ( M , F )\,\largerightarrow\, ( P , E )\,\largerightarrow\, 1$$ 
\par\medskip\noindent
of fully exact pairs. Then there is associated a long exact sequence of $\, K^J_*$-groups
\medskip
$$ \quad\cdots\>\>\>\buildrel p_*\over\longrightarrow K^J_{n+1} (P , E ) \buildrel {\delta }_*\over\longrightarrow 
K^J_n ( N , F )n \buildrel i_*\over\longrightarrow K^J_n ( M , F ) \buildrel p_*\over\longrightarrow 
K^J_n ( P , E )  $$ 
$$  \buildrel {\delta }_*\over\longrightarrow\>\quad\cdots \quad\cdots\quad\cdots\quad 
\buildrel {\delta }_*\over\longrightarrow K^J_2 ( N , F ) \buildrel i_*\over\longrightarrow K^J_2 ( M , F ) 
\buildrel p_*\over\longrightarrow K^J_2 ( P , E )    $$
$$  \buildrel {\delta }_*\over\longrightarrow K_1 ( N , F ) 
\buildrel i_*\over\longrightarrow K_1 ( M , F )\buildrel p_+\over\longrightarrow K_1 ( P , E ) 
\longrightarrow 0\>  .\qquad\qquad\qquad\qquad  $$
\par\bigskip\noindent
First order exactness of $\, ( P , E )\, $ is only needed and sufficient for exactness at $\, K^J_2 ( N , F )\, $, first order exactness of 
$\, ( M , F )\, $ and second order exactness of $\, ( P , E )\, $ are only needed and sufficient for exactness at $\, K^J_3 ( P , E )\, $, and first order exactness of $\, ( N , F )\, $ and second order exactness of 
$\, ( M , F )\, $ are only needed and sufficient for exactness at $\, K^J_3 ( M , F )\, $, second order exactness of $\, ( N , F )\, $ is only needed and sufficient for exactness at $\, K^J_3 ( N , F )\, $. At all other places the sequence is exact in any instance. If $\, ( N , F )\, $ is partially admissible on 
$\, N \cap [ M , F ]\, $ the sequence splits at $\, K^J_2 ( P , E )\, $.
\par\bigskip\noindent
{\it Proof.}\quad Let $\, (\, C_\varphi\, ,\, U_F\, )\, $ be the mapping cone of the inclusion 
$\, ( N , F ) \subseteq ( M , F )\, $. Then one gets an exact sequence
\smallskip
$$ 1\,\longrightarrow\, \bigl(\, U_{N,F}\, ,\, U_F\, \bigr)\,\longrightarrow\, \bigl(\, C_\varphi\, ,\, U_F\, \bigr)\,\longrightarrow\,\bigl(\, J_{P,E}\, ,\, U_E\,\bigr)\,\longrightarrow\, 1 $$
\smallskip
showing first of all that 
$\, K^J_2 (\, C_\varphi\, ,\, U_F\, )\,\simeq\, K^J_2 (\, J_{P,E}\, ,\, U_E\, )\,\simeq\, K^J_3 ( P , E )\, $. 
By a similar argument one gets an isomorphism
\medskip
$$ K^J_2 ( C_{u^n_{\varphi }}\, ,\, U^{n+1}_F )\buildrel\sim\over\longrightarrow 
K^J_2 ( J_{J^n_{M , F} / J^n_{N , F}\, ,\, U^n_F / J^n_{N , F}}\, ,\, U_{U^n_F / J^n_{N , F}} )\> . $$
\par\medskip\noindent
We claim that the split surjective map
\medskip
$$ K^J_2 ( J_{J^n_{M , F} / J^n_{N , F}\, ,\, U^n_F / J^N_{N , F}}\, ,\, U_{U^n_F / J^n _{N , F}} )\,\rightarrow\, K^J_2 ( J^{n+1}_{P , E}\, ,\, U^{n+1}_E ) $$
\par\medskip\noindent
is also injective. Consider first the case $\, n = 1\, $. Then the pair 
$\, ( J_{M , F} / J_{N , F}\, ,\, U_F / J_{N , F} )\, $ can be written as 
$\, ( J_{P , E}\, ,\, N \rtimes U_E )\, $ and since the subgroup 
$\, U_{J_{P , E}\, ,\, U_F / J_{N , F}} \cap U_{N\, ,\, U_F / J_{N , F}}\, $ has a weak core for 
$\, U_{U_F / J_{N , F}}\, $ one gets an injective map
\medskip
$$ K^J_2 ( J_{J_{P , E}\, ,\, N \rtimes U_E}\, ,\, U_{N \rtimes U_E} )\,\rightarrowtail\, 
K^J_2 ( J^2_{P , E}\, ,\, {U_{N \rtimes U_E}\over U_{J_{P , E}\, ,\, N \rtimes U_E} \cap 
U_{N\, ,\, N \rtimes U_E} } ) \> . $$
\par\medskip\noindent
The map 
\medskip
$$ K^J_2 ( J^2_{P , E}\, ,\, {U_{N \rtimes U_E}\over U_{J_{P , E}\, ,\, N \rtimes U_E} \cap 
U_{N\, ,\, N \rtimes U_E} } )\,\twoheadrightarrow\, K^J_2 ( J^2_{P , E}\, ,\, U^2_E ) $$
\par\medskip\noindent
is a regular surjection whose kernel is contained in the kernel of
\medskip
$$ K^J_2 ( U_{J_{P , E}\, ,\, N \rtimes U_E}\, ,\, {U_{N \rtimes U_E}\over U_{J_{P , E}\, ,\, N \rtimes U_E} \cap U_{N\, ,\, N \rtimes U_E} } )\,\twoheadrightarrow\, 
K^J_2 ( U_{J_{P , E}\, ,\, U_E}\, ,\, U^2_E ) $$
\par\medskip\noindent
since the quotient $\, ( J_{P , E}\, ,\, J_{P , E} \rtimes U_{M \rtimes U_{E / P}} )\, $ is exact and the image of the $K^J_2$-group of its suspension in 
$\, K^J_2 ( J^2_{P , E}\, ,\, U_{N\, ,\, U_F / J_{M , F}} \rtimes U^2_E )\, $ intersects trivially with the kernel of the map to $\, K^J_2 ( J^2_{P , E}\, ,\, U^2_E )\, $ (compare the analogous argument in the examples following Definition 13). The cases of higher $\, n\, $ are quite similar. For example when $\, n = 2\, $ one gets $\, J_{J^2_{M , F} / J^2_{N , F}\, ,\, U^2_F / J^2_{N , F}}\,\simeq\, 
J_{J_{J_{P , E}\, ,\, U_F / J_{N , F}}\, ,\, U^2_F / J^2_{N , F}}\, $. One can first reduce to 
$\, ( J_{J_{J_{P , E}\, ,\, U_F / J_{N , F}}\, ,\, U_{U_F / J_{N , F}}}\, ,\, U^2_{U_F / J_{N , F}} )\, $ by an argument as above. Since the kernel of 
$\,  J_{J_{P , E}\, ,\, U_F / J_{N , F}}\,\twoheadrightarrow\, J_{J_{P , E}\, ,\, U_E} \, $ satisfies the conditions of a strict splitting for the enveloping group 
$\, U_{U_F / J_{N , F}}\, $, the kernel of 
$\, J_{J_{J_{P , E}\, ,\, U_F / J_{N , F}}\, ,\, U_{U_F / J_{N , F}}}\,\twoheadrightarrow\, 
J_{J_{J_{P , E}\, ,\, U_E}\, ,\, U_{N\, ,\, U_F / J_{M , F}} \rtimes U_E}\, $ which is equal to the suspension of the former kernel also has a weak core for $\, U^2_{U_F / J_{N , F}}\, $ by Theorem 3, so one gets an injective map
\medskip
$$ K^J_2 ( J_{J_{J_{P , E}\, ,\, N \rtimes U_E}\, ,\, U_{N \rtimes U_E}}\, ,\, U^2_{N \rtimes U_E} )\,\rightarrowtail\, $$
$$ K^J_2 ( J_{J_{J_{P , E}\, ,\, U_E}\, ,\, U_{N , U_F / J_{M , F}} \rtimes U^2_E}\, ,\, 
( U_{J_{P , E} , N \rtimes U_E} \cap U_{N , N \rtimes U_E} ) \rtimes 
U_{U_{N , U_F / J_{M , F}} \rtimes U^2_E} )\> . $$
\par\medskip\noindent 
The image of this map intersects trivially with the kernel of the regular surjection
\medskip
$$ K^J_2 ( J_{J_{J_{P , E}\, ,\, U_E}\, ,\, U_{N , U_F / J_{M , F}} \rtimes U^2_E}\, ,\, 
( U_{J_{P , E} , N \rtimes U_E} \cap U_{N , N \rtimes U_E} ) \rtimes 
U_{U_{N , U_F / J_{M , F}} \rtimes U^2_E} ) $$
$$ \twoheadrightarrow\, 
K^J_2 ( J_{J_{J_{P , E}\, ,\, U_E}\, ,\, U_{N , U_F / J_{M , F}} \rtimes U^2_E}\, ,\, 
U_{U_{N , U_F / J_{M , F}} \rtimes U^2_E} ) $$
\par\medskip\noindent
by a standard argument. Then again the kernel of
\medskip 
$$ J_{J_{J_{P , E}\, ,\, U_E}\, ,\, U_{N\, ,\, U_F / J_{M , F}} \rtimes U^2_E}\,\twoheadrightarrow\, 
J_{J_{J_{P , E}\, ,\, U_E}\, ,\, U^2_E} $$
\par\medskip\noindent
has a weak core for 
$\, U_{U_{N\, ,\, U_F / J_{M , F}} \rtimes U^2_E}\, $, so one gets an injective map
\medskip
$$ K^J_2 ( J_{J_{J_{P , E}\, ,\, U_E}\, ,\, U_{N , U_F / J_{M , F}} \rtimes U^2_E}\, ,\, 
U_{U_{N , U_F / J_{M , F}} \rtimes U^2_E} )\, $$
$$ \rightarrowtail\, 
K^J_2 ( J^3_{P , E}\, ,\, U_{N\, ,\, U_F / J_{M , F}} \rtimes U^3_E ) $$
\par\medskip\noindent
the image of which intersects trivially with the kernel of the regular surjection
\medskip
$$ K^J_2 ( J^3_{P , E}\, ,\, U_{N\, ,\, U_F / J_{M , F}} \rtimes U^3_E )\,\twoheadrightarrow\, 
K^J_2 ( J^3_{P , E}\, ,\, U^3_E ) \> . $$
\par\medskip\noindent
It should now be clear how to generalize this argument to arbitrary values of $\, n\, $.
\par\smallskip\noindent
For $\, n = 1\, $ consider the surjective map 
\smallskip
$$ \bigl[\, {\overline M}^{ _{\sps F}}\, ,\, F\, \bigr]\,\longrightarrow\, 
\bigl[\, {\overline P}^{ _{\sps E}}\, ,\, E\, \bigr] \>  .  $$
\smallskip
The kernel is equal to the image of $\, [\, {\overline N}^{ _{\sps F}}\, ,\, F\, ]\, $ in 
$\, [\, {\overline M}^{ _{\sps F}}\, ,\, F\, ]\, $ which is $\, [\, {\widetilde C}_\varphi\, ,\, F\, ]\, $, giving an injective map $\, {\widetilde K}_1 (\, C_\varphi\, ,\, U_F\, )\,\longrightarrow\, K^J_2 ( P , E )\, $. But then an element of $\, [\, {\overline M}^{ _{\sps F}}\, ,\, F\, ]\, $ maps to $\, K^J_2 ( P , E )\, $ if and only if it is contained in $\, {\widetilde C}_\varphi\, $, so the map 
$\, {\widetilde K}_1 (\, C_\varphi\, ,\, U_F\, )\,\longrightarrow\, K^J_2 ( P , E )\, $ is an isomorphism. The boundary maps $\, {\delta }_*\, $ are defined with respect to these identities, and the exactness of the sequence up to $\, K_1 ( N , F )\, $ follows from Theorem 5. The rest of the sequence is immediate. The splitting at $\, K^J_2 ( P , E )\, $ can be seen by the following argument. Let $\, x \in K^J_2 ( P , E )\, $ be given. Choose a compatible splitting 
$\, (\, U_{P,E}\, ,\, U_E\, )\,\buildrel\sigma\over\longrightarrow\, (\, U_{M,F}\, ,\, U_F\, )\, $, it will send 
$\, J_{P,E}\, $ to $\, J_{M,F}\, +\, U_{N,F}\, $ (and $\, J_E\, $ to $\, J_F\, +\, U_{N,F}\, $). Let $\,\xi\, $  be a representative for $\, x\, $ in $\, J_{P,E} \cap [\, U_{P,E}\, ,\, U_E\, ]\, $ and $\, p \in U_{N,F}\, $ an element such that $\, \xi p^{-1} \in J_{M,F}\, $. Choose an "almost canonical map" 
\smallskip
$$ U_{N,F}\,\buildrel {\alpha }_N\over\longrightarrow\, {\overline N}^{ _{\sps F}}\,\supseteq\, 
{\overline {N \cap [ M , F ]}}^{ _{\sps ( N , F )}} \>  ,   $$
\smallskip
an equivariant lift 
$\, {\overline {N \cap [ M , F ]}}^{ _{\sps ( N , F )}} \,\longrightarrow\, [\, {\overline M}^{ _{\sps F}}\, ,\, F\, ]\, $ and 
$\, [\, {\overline M}^{ _{\sps F}}\, ,\, F\, ]\,\buildrel {\overline p}_*\over\longrightarrow\, 
 [\, {\overline P}^{ _{\sps E}}\, ,\, E\, ]\, $. Then choose an "almost canonical map" 
$\, U_{M,F}\,\buildrel {\alpha }_M\over\longrightarrow\, {\overline M}^{ _{\sps F}}\, $ induced from 
$\, {\alpha }_N\, $, whose existence is guaranteed by Lemma 10. Put 
$\, {\pi }_M ( x ) = {\overline p} ( {\alpha }_M ( \sigma ( \xi ) p^{-1} ) )\, $. Then $\, {\pi }_M ( x )\, $ is in the image of $\, K^J_2 ( M , F )\, $. Moreover the map $\, {\pi }_M\, $ does not depend on the choice of 
$\, p\, $. If $\, p'\, $ is another such element, then $\, p {p'}^{-1} \in J_{N,F}\, $ and is projected to the image of $\, K^J_2 ( N , F )\, $ in $\, K^J_2 ( M , F )\, $ by $\, {\alpha }_M\, $, because $\, {\alpha }_M\, $ is induced from $\, {\alpha }_N\, $. Hence the map $\, {\pi }_M\, $ is well defined on representatives $\,\xi\, $ in the first place. One checks that it maps trivial elements 
$\, \xi \in [\, J_{P,E}\, ,\, U_E\, ]\, +\, [\, U_{P,E}\, ,\, J_E\, ]\, $ to the identity, so it drops to a map from 
$\, K^J_2 ( P , E )\, $ to itself. It is easy to check that this map is the identity on the image of 
$\, K^J_2 ( M , F )\, $, thus it gives a projection onto this image, hence a splitting. Of course the splitting can also be seen in the following way: if $\, ( N , F )\, $ is partially admissible on 
$\, N \cap [ M , F ]\, $ it means that the kernel of the map 
$\, K_1 ( N , F )\,\longrightarrow\, K_1 ( M , F )\, $ is free abelian, so it lifts to $\, K^J_2 ( P , E )\, $ 
\qed
\par\bigskip\noindent
{\bf Remark.} \,  From Corollary 5.1 one deduces that the pair $\, (\, C_\rho\, ,\, U_\rho\, )\, $, the 
mapping cone of a preexcisive inclusion has trivial $\, K^J_*$-groups for all $\, n \geq 2\, $, so that Theorem 5 tells us that $\, K^J_n (\, J_{M_F , F}\, ,\, U_F\, )\,\simeq\, K^J_n (\, J_{M , G}\, ,\, U_\phi\, )\, $ and then also 
$\, K^J_n (\, C_{\varphi }\, ,\, U_F\, )\,\simeq\, K^J_n (\, C_\phi\, ,\, U_\phi\, )\, $ for all  
$\, n \geq 2\, $. This is a bit disappointing in a way, since it means that the sequence of Theorem 5 really reduces to that of Corollary 5.1 (at least to the left of $\, K^J_2 ( N , F )\, $). On the other hand the isomorphism above adds to the computability of the theory and may also prove useful in other respects
\par\smallskip\noindent
\bigskip\bigskip\bigskip
\par\noindent
\hfil\hfil \Large{\bf{ Section VI -- Regular $K$-theory.}} \hfil 
\bigskip\bigskip

\par\noindent
Building on the results of the last section we are now going to venture into regular $K$-theory. After a couple of preparatory results highlighted by the Special Deconstruction Theorem (Lemma 16), and completing the proof of Theorem 2, we show that the ad-hoc definition of the regular $K$-groups presented in Definition 14 leads to a properly defined theory (Proposition 1) which has long eaxct sequences analogous to those of (and embedding into) $K^J$-theory (Theorem 6). We begin with 
\par\bigskip\noindent
{\bf Definition 14.}\quad A chain $\, \{\, R^n_{N , F} \}^m_{n=0}\, $ of subgroups of $\, U^n_F\, $ will be called {\fndef regular } of length m for the pair $\, ( N , F )\, $, iff 
$\, ( R^0_{N , F} , U^0_F )\, =\, ( N , F )\, $ and for every 
$\, n = 0\, ,\ldots , m-1\, $ the pair 
$\, ( R^n_{N , F} , U^n_F )\, $ admits a universal $U^n_F$-central extension 
$\, {\overline R}^{ _{\sps n}}_{N , F}\, $ relative to 
$\, ( J_{R^{n-1}_{N , F}\, ,\, U^{n-1}_F}\, ,\, U^n_F )\, $ such that the central kernel of the projection onto 
$\, R^n_{N , F}\, $ is contained in the k-fold commutator subgroup 
$\, [\, [ \cdots [ {\overline R}^{ _{\sps n}}_{N , F}\, , U^n_F\, ],\cdots ], U^n_F\, ]\, $ modulo the canonical subgroup 
$\, [\, J_{R^{n-1}_{N , F}\, ,\, U^{n-1}_F}\, ,\, 
J_{R^{n-1}_{N , F}\, ,\, U^{n-1}_F}\, ] \subseteq 
{\overline R}^{ _{\sps n}}_{N , F}\, $ for any k and $\, R^{n+1}_{N , F}\, $ is the kernel of an almost canonical map
\smallskip
$$ U_{R^n_{N , F}, U^n_F} \longrightarrow {\overline R}^n_{N , F}\> .  $$
\par\smallskip\noindent
The pair $\, ( N , F )\, $ will be called {\it 0-regular} if it is admissible and the regular part of its semiuniversal $F$-central extension is equal to $\, K^J_2 ( N , F )\, $. It is called
{\it m-regular}, if it is (m-1)-regular and for any given regular chain of length m-1 as above, any kernel $\, R^m_{N , F}\, $ of some almost canonical map 
\smallskip
$$ U_{R^{m-1}_{N , F}\, ,\, U^{m-1}_F} \longrightarrow {\overline R}^{m-1}_{N , F} $$
\par\medskip\noindent
extends the regular chain $\,\{ R^n_{N , F} \}\, $ to a regular chain of length m. $\, ( N , F )\, $ will be called a {\it regular pair } if it is m-regular for all values of m. One has a notion of {\it regular equivalence} of two regular pairs generated by basic relations of the following form. Assume given a map $\, ( N , F ) \rightarrow ( M , G )\, $ of regular pairs. Then 
$\, ( N , F )\, $ is said to be equivalent to $\, ( M , G )\, $ if for every regular chain $\,\{ R^n_{N , F} \}\, $ of 
$\, ( N , F )\, $ there is an induced regular chain $\,\{ R^n_{M , G} \}\, $ of $\, ( M , G )\, $ and for every regular chain $\,\{ R^n_{M , G} \}\, $ there is a restricted regular chain $\,\{ R^n_{N , F} \}\, $ such that the induction (restriction) maps  
$\,{\varphi }_n :  ( R^n_{N , F} , U^n_F )\longrightarrow ( R^n_{M , G} , U^n_G )\, $ induce an isomorphism on the $K^J_2$-groups for all values of $\, n \geq 0\, $. 
If $\, ( N , F )\, $ is m-regular define its n-th {\it suspension} (or {\it $R$-suspension}) to be the pair $\, ( R^n_{N , F} , U^n_F )\, $ for $\, n \leq m\, $ and define its regular  $K$-groups 
$\, K_n (\, N , F\, )\, $ for $\, n = 2,\ldots , m+1\, $ by
$$ K_n ( N , F ) := K^J_2 ( R^{n-2}_{N , F}\, ,\, U^{n-2}_F\, )\> . $$
It will be shown in the next section that the groups $\, K_n ( N , F )\, $ are well defined up to natural isomorphism and functorial. These definitions can be extended to the {\it almost regular } case, i.e. the case where $\, ( N , F )\, $ contains a regular subpair $\, ( N_r , F )\, $ such that $\, N /  N_r\, $ is 
$(F /  N_r)$-nilpotent and any two such subpairs are regular equivalent by putting 
$\, K_n ( N , F ) = K_n ( N_r , F )\, $ in this case.
\par\bigskip\noindent
{\bf Lemma 15.}\quad Let $\, ( N , F )\, $ be a normal (resp. regular) pair. Then the pairs 
$\, ( J^n_{N , F}\, ,\, L^n_F )\, $ (resp. $\, ( J_{R^{n-1}_{N , F}\, ,\, U^{n-1}_F}\, ,\, L^n_F )\, $) are admissible for all $\, n \geq 1\, $ with $\, K^J_2 ( J^n_{N , F}\, ,\, L^n_F ) = 0\, $ (resp. 
$\, K^J_2 ( J_{R^{n-1}_{N , F}\, ,\, U^{n-1}_F}\, ,\, L^n_F ) = 0\, $).
\par\bigskip\noindent
{\it Proof.}\quad For $\, n = 1\, $, $\, L^1_F = J_F\, $ and $\, ( J_{N , F}\, ,\, J_F )\, $ is a cone, hence admissible with $\, K^J_2 ( J_{N , F}\, ,\, J_F ) = 0\, $. For $\, n = 2\, $ consider the inclusion 
$\, \varphi : ( J_{N , F}\, ,\, J_F ) \subset ( J_{N , F}\, ,\, U_F )\, $. The mapping cone 
$\, ( C_{\varphi }\, ,\, U_{\varphi } ) = ( U_{J_{N , F}\, ,\, U_F}\, ,\, L^2_F )\, $ is a cone, so 
$\, K^J_2 ( C_{\varphi }\, ,\, U_{\varphi } ) = 0\, $. Then consider the exact sequence
\smallskip
$$ 1 \longrightarrow ( J^2_{N , F}\, ,\, L^2_F ) \longrightarrow ( C_{\varphi }\, ,\, L^2_F ) \longrightarrow 
( J_{N , F}\, ,\, L^2_F / J^2_{N , F} ) \longrightarrow 1\> . $$
\par\medskip\noindent
The long exact sequence of Corollary 5.1 shows that $\, K_1 ( J^2_{N , F}\, ,\, L^2_F )\, $ is an extension of a free abelian group by $\, K^J_2 ( J_{N , F}\, ,\, L^2_F / J^2_{N , F} )\, $ which is free abelian since 
$\, L^2_F / U_{J_{N , F}\, ,\, U_F}\, $ is free (compare with the proof of Lemma 14). Then 
$\, ( J^2_{N , F}\, ,\, L^2_F )\, $ is admissible and Theorem 5 gives that 
$\, K^J_2 ( J^2_{N , F}\, ,\, L^2_F ) \rightarrowtail K^J_2 ( C_{\varphi }\, ,\, L^2_F )\, $ is injective, so 
$\, K^J_2 ( J^2_{N , F}\, ,\, L^2_F ) = 0\, $. One proceeds inductively, considering the inclusions 
$\, ( J^k_{N , F}\, ,\, L^k_F ) \subseteq ( J^k_{N , F}\, ,\, U^k_F )\, $. For $\, R_{N , F}\, $ consider the inclusion $\, {\varphi }_R : ( R_{N , F}\, ,\, J_F ) \subseteq (R_{N , F}\, ,\, U_F )\, $. Again 
$\, ( C_{{\varphi }_R}\, ,\, U_{{\varphi }_R} ) = ( U_{R_{N , F}\, ,\, U_F}\, ,\, L^2_F )\, $ is a cone and as above it follows that $\, ( J_{R_{N , F}\, ,\, U_F}\, ,\, L^2_F )\, $ is admissible. Theorem 5 gives halfexactness of 
\smallskip
$$ K^J_2 ( J_{R_{N , F}\, ,\, J_F}\, ,\, U_{J_F} ) \longrightarrow K^J_2 ( J_{R_{N , F}\, ,\, U_F}\, ,\, L^2_F ) \longrightarrow K^J_2 ( C_{{\varphi }_R}\, ,\, L^2_F ) $$
\par\medskip\noindent
and the Excision Theorem gives injectivity of the map
\smallskip 
$$ K^J_2 ( J_{R_{N , F}\, ,\, J_F}\, ,\, U_{J_F} ) \rightarrowtail K^J_2 ( J_{R_{N , F}\, ,\, U_F}\, ,\, U^2_F ) $$
\par\medskip\noindent 
which implies that $\, K^J_2 ( J_{R_{N , F}\, ,\, J_F}\, ,\, U_{J_F} ) = 0\, $, since as will be shown below the map $\, K^J_2 ( J_{R_{N , F}\, ,\, U_F}\, ,\, U^2_F ) \rightarrowtail K^J_2 ( J^2_{N , F}\, ,\, U^2_F )\, $ is injective, so that also 
$\, K^J_2 ( J_{R_{N , F}\, ,\, J_F}\, ,\, U_{J_F} ) \rightarrowtail 
K^J_2 ( J_{J_{N , F}\, ,\, J_F}\, ,\, U_{J_F} ) = 0\, $ must be injective. The argument that 
$\, ( J_{R^{n-1}_{N , F}\, ,\, U^{n-1}_F}\, ,\, L^n_F )\, $ is admissible with 
$\, K^J_2 ( J_{R^{n-1}_{N , F}\, ,\, U^{n-1}_F}\, ,\, L^n_F ) = 0\, $ is much the same\qed
\par\bigskip\noindent
{\bf Remark.} The pair $\, ( R^n_{N , F}\, ,\, L^n_F )\, $ need not be admissible in general, nor 
$\, K^J_2 ( R^n_{N , F}\, ,\, L^n_F ) = 0\, $ (this is an unsolved question so far), but in any case the image of $\, K^J_2 ( R^n_{N , F}\, ,\, L^n_F )\, $ in $\, K^J_2 ( R^n_{N , F}\, ,\, U^n_F )\, $ is trivial because the latter injects into $\, K^J_2 ( J_{R^{n-1}_{N , F}\, ,\, U^{n-1}_F}\, ,\, U^n_F )\, $ as will be shown below.
\par\medskip\noindent 
{\bf Lemma 16.}\quad ( Special Deconstruction Theorem )
Let $\, N\, $ be a full subgroup of a perfect group $\, F\, $ and $\, {\mathcal U} \subset U^m_F\, $ a normal subgroup such that $\, {\mathcal U} + L^m_F = U^m_F\, $ and $\, {\mathcal U}\, $ contains 
$\, J_{R^{m-1}_{N , F}\, ,\, U^{m-1}_F}\, $. For $\, m = 1\, $ one also assumes that $\, {\mathcal U}\, $ contains $\, U_{N , F}\, $. Assume that $\, ( R^m_{N , F}\, ,\, {\mathcal U} )\, $ admits a relative semiuniversal ${\mathcal U}$-central extension with respect to the inclusion into 
$\, ( J_{R^{m-1}_{N , F}\, ,\, U^{m-1}_F}\, ,\, {\mathcal U} )\, $. Put
 $\, {\mathcal U}' = {\mathcal U} + J^m_{N , F}\, $. The map 
$\, K^J_2 ( R^m_{N , F}\, ,\, {\mathcal U} ) \twoheadrightarrow K^J_2 ( R^m_{N , F}\, ,\, U^m_F )\, $ is surjective and its image in $\, K_2^J ( J^m_{N , F} , U_F^m )\, $ is equal to the image of $\, K^J_2 ( R^m_{N , F}\, ,\, {\mathcal U} )\, $ in 
$\, K^J_2 ( J^m_{N , F}\, ,\, {\mathcal U}' ) \subseteq K_2^J ( J^m_{N , F} , U^m_F )\, $.
\par\bigskip\noindent
{\it Proof.}\quad In the following put $\, R = R^m_{N , F}\> ,\> 
J = J_{R^{m-1}_{N , F}\, ,\, U^{m-1}_F}\> ,\> U = U^m_F\, $ and 
$\, L = L^m_F\, $. We first show injectivity of 
$\, K^J_2 ( J^m_{N , F}\, ,\, {\mathcal U}' ) \rightarrowtail K^J_2 ( J^m_{N , F}\, ,\, U )\, $ from which the last statement follows from injectivity of the map   
$\, K^J_2 ( R^m_{N , F}\, ,\, U^m_F ) \rightarrowtail K^J_2 ( J^m_{N , F}\, ,\, U^m_F )\, $ (see below). 
Let 
\smallskip
$$ I^{k , m}_F = U^k_{( J_{U^{m-k-1}_F}\, ,\, U^{m-k}_F )}\quad ,\quad k = 0 ,\cdots , m-1\> , $$ 
\par\medskip\noindent
be one of the  m  "fingers" of $\, L^m_F\, $, so that $\, L^m_F = {\bigcup }_{k = 0}^{m-1} I^{k , m}_F\, $. Inductively put $\, {\mathcal U}^{(0)} = {\mathcal U}'\, ,\, 
{\mathcal U}^{(k+1)} = {\mathcal U}^{(k)} + I^{k , m}_F\, $. Then as $\, J^{m-k}_{N , F}\, $ has a core for 
$\, J_{U^{m-k-1}_F}\, $, $\, J^m_{N , F}\, $ has a weak core for $\, I^{k , m}_F\, $ (compare with sections 3 and 4). From this one gets $\, K^J_2 ( J^m_{N , F}\, ,\, {\mathcal U}^{(k)} \cap I^{k , m}_F ) = 0\, $ and the Excision Theorem gives an injective map 
$\, K^J_2 ( J^m_{N , F}\, ,\, {\mathcal U}^{(k)} ) \rightarrowtail 
K^J_2 ( J^m_{N , F}\, ,\, {\mathcal U}^{(k+1)} )\, $ for each $\, k\, $. 
In order to show surjectivity of 
$\, K^J_2 ( R\, ,\, {\mathcal U} ) \rightarrow K^J_2 ( R\, ,\, U )\, $ one can assume the (relative) semiuniversal ${\mathcal U}$- (and $U$-)central extensions of $\, R^m_{N , F}\, $ both being embedded in the semiuniversal $U$-central extension $\, {\overline J}^{ _U}\, $ of 
$\, J = J_{R^{m-1}_{N , F}\, ,\, U^{m-1}_F}\, $. The strategy is now as follows. Dividing the extension 
$\, {\overline J}^{ _U}\, $ by the image of $\, K^J_2 ( R\, ,\, {\mathcal U} )\, $ one obtains a 
${\mathcal U}$-normal copy $\, \widehat R\, $ of $\, R\, $ and an $L^m_F$-normal copy 
$\, \Breve R\, $ of $\, R\, $ in the quotient since $\, ( J\, ,\, L )\, $ is admissible with 
$\, K^J_2 ( J\, ,\, L ) = 0\, $ by Lemma 15. One gets a canonical $U$-normal copy of 
$\, [\, R\, ,\, {\mathcal U}\, ]\, $ inside $\, {\overline J}^{ _U} / K^J_2 ( R , {\mathcal U} )\, $ contained in 
$\, \widehat R\, $. The $L$-normal copy $\, \Breve R\, $ sits inside an $L$-normal copy 
$\, \Breve J\, $ of $\, J\, $ so it contains the canonical $U$-normal copy of 
$\, [\, J\, ,\, L\, ] = [\, {\Breve J}\, ,\, L\, ]\, $. Both $\, [\, \widehat R\, ,\, {\mathcal U} ]\, $ and 
$\, [\, {\Breve J}\, ,\, L\, ]\, $ are normal for the action of $\, U\, $, because $\, {\mathcal U}\, $ and 
$\, L\, $ are normal subgroups of $\, U\, $ (by the Jacobi identity). The problem is to "glue together" the two subgroups over (some part of) their intersection. This intersection certainly contains 
$\, [\, R\, ,\, {\mathcal U} \cap L\, ]\, $ and hence also $\, [\, [\, R\, ,\, {\mathcal U}\, ] , L\, ]\, $ which is congruent to $\, [\, [\, R\, ,\, L\, ] , {\mathcal U}\, ]\, $ modulo the former subgroup. For 
$\, x \in U^m_F\, $ let $\,\xi\, $ be its image in $\, F\, $ modulo $\, L^m_F\, $. Consider the m-fold iterated section $\, \xi \mapsto u_{\xi } \mapsto u_{u_{\xi }} \mapsto\cdots \, $ and let $\, s^m ( \xi )\, $ denote its image in $\, U^m_F\, $. Modulo $\, L^m_F\, $ each representative $\, s^m ( \xi )\, $ is equivalent to another representative $\, u^m ( \xi )\, $ in $\, {\mathcal U}\, $. Since $\, F\, $ is perfect these elements can be chosen to lie in $\, [\, {\mathcal U}\, ,\, {\mathcal U}\, ]\, $. Denote this new set $\, {\mathcal U}^{ _F}\, $. Now consider the usual core of $\, J = J_{R^{m-1}_{N , F}\, ,\, U^{m-1}_F}\, $ for $\, J_{U^{m-1}_F}\, $ consisting of the set 
$\, \{\, u_z\, B_1\, u^{-1}_z\,\} ÷\cup \{\, u_z\, B_2\, u^{-1}_z\,\} \cup \{\, u_z\, B_3\, u^{-1}_z\,\}\, $ as in Lemma 3 of section 1. The set $\, B_3' = \{\, u_c\, u_y\, u^{-1}_{cy}\, \} \subset B_3\, $ can be completed to a basis of $\, U^m_F\, $, so that its normalization defines a cone with respect to both $\, {\mathcal U}\, $ and $\, U\, $. One can divide by this cone without changing $\, K^J_2 ( J\, ,\, U )\, $ and consider the reduced form of $\, U^m_F\, $ which is freely generated by the union of the two sets 
$\, \{\, u_c\,\vert\, c \in R^{m-1}_{N , F}\,\} \cup \{\, u_z\,\vert\, z = s ( \overline z )\,\}\, $, the latter of which is a semicanonical lift of the standard basis of $\, U_{U^{m-1}_F / R^{m-1}_{N , F}}\, $. First assume that 
$\, m \geq 2\, $. Modulo $\, [\, J\, ,\, J\, ]\, $ (the $B_3'$-reduced) $\, J\, $ is generated by the set 
$\, \{\, u^{\pm }_{x_1}\cdots u^{\pm }_{x_n}\, B_{1 , 2}\, u^{\mp }_{x_n}\cdots u^{\mp}_{x:_1}\,\}\, $ 
(Lemma 4) which can also be written in the form 
\smallskip
$$ \bigl\{\> u^{\pm }_{x_1}\cdots u^{\pm }_{x_n}\> s^m ( \zeta )\> B_{1 , 2}\> s^m ( \zeta )^{-1}\> 
u^{\mp }_{x_n}\cdots u^{\mp }_{x_1}\>\bigr\} $$
\par\medskip\noindent
with $\, u^{\pm }_{x_1}\cdots u^{\pm }_{x_n} \in {\overline L} = L^m_F \cap \langle\, \{\, u_z\,\}\,\rangle \simeq L^m_F / U_{R^{m-1}_{N , F}\, ,\, U^{m-1}_F}\, $. It is not clear that these generators are independent modulo $\, [\, J\, ,\, J\, ]\, $, but it is fairly easy to check that they are independent in $\, J\, $. Let $\, \{\, f_k\,\}\, $ denote a basis of $\, R^{m-1}_{N , F}\, $, put 
$\, B^0_1 = \{\, u_c\, u_{f_k}\, u^{-1}_{cf_k}\,\} \subseteq B_1\, $ and 
$\, B^0_2 = \{\, [\, u_x\, ,\, u_{f_k}\, ]\, u_{f_k}\, u^{-1}_{(f_k)_x}\,\} \subseteq B_2\, $. Then a simple calculation gives that $\, J\, /\, [\, J\, ,\, J\, ]\, $ is equivalently generated by the set 
\smallskip
$$ \bigl\{\> u^{\pm }_{x_1}\cdots u^{\pm }_{x_n}\> s^m ( \zeta )\> u_e\> B^0_{1 , 2}\> u^{-1}_e\> 
s^m ( \zeta )^{-1}\> 
u{\mp }_{x_n}\cdots u^{\mp }_{x_1}\> \bigr\} $$
\par\medskip\noindent
which is again linear independent in $\, J\, $ (this is a bit more complicated to check, one separates the subcores corresponding to $\,\{\, u_e\, B^0_1\, u^{-1}_e\,\}\, $ and $\,\{\, u_e\, B^0_2\, u^{-1}_e\,\}\, $, the first of which is obviously equivalent to the set $\, B_1\, $, so that modulo its normalization it suffices to check that the set $\,\{\, u_e\, B^0_2\, u{-1}_e\,\}\, $ is an equivalent basis for the subcore generated by 
$\, B_2\, $). We may now replace the coefficients $\, s^m ( \zeta )\, $ by 
$\, u^m ( \zeta ) \in [\, {\mathcal U}\, ,\, {\mathcal U}\, ]\, $. Then the resulting set 
\smallskip
$$ \bigl\{\> u^{\pm }_{x_1}\cdots u^{\pm }_{x_n}\> u^m ( \zeta )\> u_e\> B^0_{1 , 2}\> u^{-1}_e\> 
u^m ( \zeta )^{-1}\> u^{\mp }_{x_n}\cdots u^{\mp }_{x_1}\> \bigr\} $$
respectively
$$ \bigl\{\> u^{\pm }_{x_1}\cdots u^{\pm }_{x_n}\> u_e\> [\, u^m ( \zeta )\, ,\, B^0_{1 , 2}\, ]\> u^{-1}_e\> 
u^{\mp }_{x_n}\cdots u^{\mp }_{x_1}\> \bigr\} $$
$$ \cup\quad \bigl\{\> u^{\pm }_{x_1}\cdots u^{\pm }_{x_n}\> u_e\> B^0_{1 , 2}\> u^{-1}_e\> 
u^{\mp }_{x_n}\cdots u^{\mp }_{x_1}\>\bigr\} $$
\par\medskip\noindent
still generates $\, J\, $ modulo $\, [\, J\, ,\, J\, ]\, $.  Consider the two subgroups generated by 
\smallskip
$$ \bigl\{\> \bigl[\, v\, ,\, \bigl[\, u^m ( \zeta )\, ,\, B^0_{1 , 2}\,\bigr]\, \bigr]\> \bigr\}\quad \cup\quad 
\bigl\{\> v\, B^0_{1 , 2}\, v^{-1}\>\bigr\} \> , $$
\par\medskip\noindent
$\, v \in \{\, u^{\pm }_{x_1}\cdots u^{\pm }_{x_n}\, u_e\,\}\, $ on one side and 
\smallskip
$$ \bigl\{\> \bigl[\, v\, ,\, \bigl[\, u^m ( \zeta )\, ,\, B^0_{1 , 2}\, \bigr]\, \bigr]\> \bigr\} \quad \cup\quad 
\bigl\{\> \bigl[\, u^m ( \zeta )\, ,\, B^0_{1 , 2}\, \bigr]\>\bigr\} $$
\par\medskip\noindent
on the other side. The second subgroup is contained in $\, [\, R\, ,\, {\mathcal U}\, ]\, $ and lifts to 
$\, \widehat R\, $, whereas the intersection of the first subgroup with $\, [\, J\, ,\, L\, ]\, $ lifts to 
$\, [\, \Breve J\, ,\, L\, ]\, \subseteq\, \Breve R\, $. Given a relation in the generators
\smallskip 
$$ \bigl\{\> \bigl[\, v\, ,\, \bigl[\, u^m ( \zeta )\, ,\, B^0_{1 , 2}\, \bigr]\, \bigr]\> \bigr\}\> ,\> 
\bigl\{\> \bigl[\, u^m ( \zeta )\, ,\, B^0_{1 , 2}\, \bigr]\> \bigr\}\> ,\> 
\bigl\{\> v\> B^0_{1 , 2}\> v^{-1}\> \bigr\} $$ 
\par\medskip\noindent
write it in the form 
\smallskip
$$ x_1\> x_2\> x_3\> \equiv\> 0 $$
\par\medskip\noindent
modulo the subgroup $\, [\, R\, ,\, J\, ]\, $ with $\, x_1\, $ in the first part, $\, x_2\, $ in the second, and 
$\, x_3\, $ in the intersection of the range of the third part and $\, [\, J\, ,\, L\, ]\, $. Then 
$\, x_1\, x_2\, $ lifts to $\, \widehat R\, $, so also $\, x_3\, $ lifts to $\, \widehat R\, $, and since 
$\, x_1\, x_3\, $ lifts to $\,\Breve R\, $, so does $\, x_2\, $. Then $\, x_1\, $ as well as $\, x_2\, $ and 
$\, x_3\, $ are contained in the subgroup corresponding to $\, \widehat R \cap \Breve R\, $, and lifting the subgroup generated by elements $\, x_1\, x_2\, $ to $\, \widehat R\, $ and the subgroup generated by 
$\, x_1\, x_3\, $ to $\,\Breve R\, $ gives a compatible lift in the sense that any relation 
$\, x_1\, x_2\, x_3\, \equiv\, 0\, $ (mod $\, [\, R\, ,\, J\, ]\, $) lifts to the trivial element. Also one checks that both subgroups are normal for the action of $\, U\, $ modulo $\, \widehat R \cap \Breve R\, $, so the lift extends to a normal lift of the normalization of this subgroup of $\, J\, $. This normal lift can be glued in a compatible way to the canonical image of $\, [\, R\, ,\, J\, ] \subseteq \widehat R \cap \Breve R\, $  and the resulting subgroup is equal to a normal lift of $\, [\, J\, ,\, U\, ]\, $. To see this one first checks that it contains $\, [\, J\, ,\, J\, ]\, $. Since the subgroup of $\, J\, $ generated by the three parts as above generates $\, J\, $ modulo $\, [\, J\, ,\, J\, ]\, $ it follows that its commutator subgroup generates 
$\, [\, J\, ,\, J\, ]\, $ modulo $\, [\, J\, ,\, [\, J\, ,\, J\, ]\, ]\, \subseteq [\, R\, ,\, J\, ]\, $, but $\, [\, R\, ,\, J\, ]\, $ has been included. Next one checks that it contains $\, [\, J\, ,\, \widetilde U\, ]\, $ where $\, \widetilde U\, $ is generated by $\, J\, $ together with $\, \sigma ( U_{U^{m-1}_F / R^{m-1}_{N , F}} )\, $. Modulo 
$\, [\, J\, ,\, J\, ]\, $ we may replace the coefficients $\, u^m ( \zeta )\, $ by expressions of the form 
$\, v\, s^m ( \zeta )\, $ where 
$\, v \in L^m_F\, $ is as above so that modulo $\, [\, J\, ,\, J\, ]\, $, the first part and the intersection of the third part with $\, [\, J\, ,\, L\, ]\, $, we may replace the first set by 
$\, \{\> [\, \overline v\, ,\, [\, s^m ( \zeta )\, ,\, B_{1 , 2}\, ]\, ]\> \}\, $ with $\, \overline v \in \overline L^m_F\, $ and the second set can be transformed to the form $\, \{\> [\, s^m ( \zeta )\, ,\, B_{1 , 2}\, ]\,\}\, $ leaving invariant the third part. Then one easily finds from the fact that the subgroup $\, J_0\, $ of $\, J\, $ generated by 
$\,\{\> \overline v\> s^m ( \zeta )\> B_{1 , 2}\> s^m ( \zeta )^{-1}\> {\overline v}^{-1}\>\}\, $ has a weak core for $\, \overline U = \sigma ( U_{U^{m-1}_F / R^{m-1}_{N , F}} )\, $ that our subgroup contains 
$\, [\, J\, ,\, \widetilde U\, ]\, $. Modulo this group $\, [\, J\, ,\, U\, ]\, $ is generated by commutators 
$\, \{\> [\, u_e\, ,\, B^0_{1 , 2}\, ]\>\}\, $ which are contained in the intersection of the third part with 
$\, [\, J\, ,\, L\, ]\, $. Modulo $\, [\, J\, ,\, J\, ]\, $ any element in the image of 
$\, K^J_2 ( R , U ) \rightarrowtail K^J_2 ( J , U )\, $ is regular, whence the result 
$\, K^J_2 ( R , {\mathcal U} ) \twoheadrightarrow K^J_2 ( R , U )\, $ follows.
\par\noindent
Now suppose that $\, m = 1\, $ and that $\, {\mathcal U}\, $ contains $\, U_{N , F}\, $. Then one considers the subgroup $\, J_0\, $ of (the $B_3'$-reduced) $\, J\, $ generated by the set 
\smallskip
$$ \bigl\{\> v^{\pm }_1\cdots v^{\pm }_n\> u ( \zeta )\> B_{1 , 2}\> 
u ( \zeta )^{-1}\> v^{\mp }_n\cdots v^{\mp }_1\>\bigr\} $$
\par\medskip\noindent
where $\, v_1 ,\cdots , v_n\, $ are basis elements of $\, \lambda ( J_{F / N} )\, $ and $\, u ( \zeta )\, $ is a chosen representative of $\, \zeta \in F / N\, $ in 
$\, \sigma ( U_{F / N} ) \cap [\, {\mathcal U}\, ,\, {\mathcal U}\, ]\, $. It is readily checked that this set generates $\, J\, $ modulo $\, [\, R\, ,\, J\, ]\, $ and that $\, J_0\, $ is normal modulo 
$\, [\, J_0\, ,\, R\, ]\, $ (or $\, [\, J_0\, ,\, [\, R\, ,\, J\, ]\, ]\, $). The generators can be seen to be independent in $\, J\, $. Dividing $\, J_0\, $ by the subgroup generated by the set 
\smallskip
$$ \bigl\{\> \bigl[\, v\, ,\, \bigl[\, u ( \zeta )\, ,\, B_{1 , 2}\, \bigr]\, \bigr]\> \bigr\}    $$
\par\medskip\noindent
with $\, v = v^{\pm }_1 \cdots v^{\pm }_n\, $ the basis of $\, J_0\, $ reduces to the two parts 
\smallskip
$$ \bigl\{\> v^{\pm }_1\cdots v^{\pm }_n\> \bigl[\, v_0\, ,\, B_{1 , 2}\, \bigr]\> v^{\mp }_n\cdots v^{\mp }_1\>\bigr\}\quad \cup\quad 
\bigl\{\> u ( \zeta )\> B_{1 , 2}\> u ( \zeta )^{-1}\>\bigr\}\> . $$
\par\medskip\noindent
The first part is a subset of a basis of $\, J\, $ and its intersection with $\, [\, J\, ,\, J\, ]\, $ is equal to its commutator subgroup. Also since $\, N\, $ is full in $\, F\, $ each element $\, u_e\, $ is congruent to an element  
$\, {\Pi }_{\lambda }\> [\, u_{e_{\lambda }}\, ,\, u_{f_{\lambda }}\, u ( {\zeta }_{\lambda } )\, ]\, $ modulo 
$\, J_{N , F}\, $. Then one checks as before that the part of $\, J_0\, $ generated by 
$\, \{\> [\, v\, ,\, [\, u ( \zeta )\, ,\, B_{1 , 2}\, ]\, ]\>\}\>\cup\>\{\> u ( \zeta )\> B_{1 , 2}\> u ( \zeta )^{-1}\>\}\, $ (or rather its intersection with $\, [\, R\, ,\, {\mathcal U}\, ]\, $) is normal modulo a subgroup of $\, R\, $ lifting to 
$\, \widehat R \cap \Breve R\, $, so also the intersection of this normalization with the first part consists of elements lifting to $\, \widehat R \cap \Breve R\, $. The same is true for the normalization of the first part with respect to the second, so that altogether one gets a compatible normal lift of the normalization of 
$\, J_0\, $ to $\, {\overline J}^{ _U}\, /\, K^J_2 ( R , {\mathcal U} )\, $ which coincides on the intersection with $\, [\, R\, ,\, J\, ]\, $ with the canonical lift of this subgroup, so that the lift extends to a normal lift of some cocentral subgroup of $\, R\, $ into $\, {\overline J}^{ _U}\, /\, K^J_2 ( R , {\mathcal U} )\, $. Now modulo $\, [\, R\, ,\, J\, ]\, $ any element in the image of 
$\, K^J_2 ( R , U ) \rightarrowtail K^J_2 ( J , U )\, $ is regular (since $\, N\, $ is full), so that 
$\, K^J_2 ( R , {\mathcal U} ) \twoheadrightarrow K^J_2 ( R , U )\, $ surjective follows\qed
\par\bigskip\noindent 
If $\, {\mathcal U}\, $ does not contain $\, U_{N , F}\, $ for $\, m = 1\, $ the result may fail to be true but since modulo $\, J_{N , F}\, $ the group $\, U_{N , F}\, $ is congruent to $\, [\, U_{N , F}\, ,\, U_F\, ]\, $ this result suffices for our purposes.
\par\bigskip\noindent
{\bf Lemma 17.}\quad Let $\, C \subseteq B\, $ be a central subgroup. Then $\, K^J_2 ( C , B )\, $ is a quotient of $\, C \otimes B^{ab} = C \otimes ( B / [ B , B ] )\, $ and an extension of 
$\, C \otimes ( B / C )^{ab}\, $. The kernel of the map 
$\, C \otimes B^{ab}\,\twoheadrightarrow\, K^J_2 ( C , B )\, $ is generated by the image of elements 
$\, \{\, c \otimes c\,\vert\, c \in C\,\}\, $ in $\, C \otimes B^{ab}\, $.
\par\bigskip\noindent
{\it Proof.}\quad  Since $\, C\, $ is central one has $\, [\, U_{C , B}\, ,\, U_B\, ]\, =\, 
J_{C , B} \cap [\, U_{C , B}\, ,\, U_B\, ]\, $. We use the basis of $\, J_{C , B}\, $ determined by the core for a semicanonical lift $\, \lambda : J_{B / C} \rightarrow J_B\, $ as in section 1. Dividing by 
$\, [\, J_{C , B}\, ,\, J_{C , B}\, ]\, $ one gets a free abelian group generated by the elements of the core, i.e. 
$\, \{\> u_z\, B_1\, u^{-1}_z\>\} \cup \{\> u_z\, B_2\, u^{-1}_z\>\} \cup \{\> u_z\, B_3\, u^{-1}_z\>\}\, $, which can be reduced to $\, B_1 \cup B_2 \cup B_3\, $ modulo 
$\, [\, J_{C , B}\, ,\, \langle \{\> u_z\> \} \rangle\, ]\, $. Since $\, B_3\, $ is congruent to $\, B_1\, $ modulo the normalization of $\, B_3'\, $ which can be divided out, since $\, B_3' \cup \{\> u_c\>\} \cup \{\> u_z\>\}\, $ is a basis of $\, U_B\, $, so the intersection of the normalization of $\, B_3'\, $ with $\, [\, U_B\, ,\, U_B\, ]\, $ is contained in $\, [\, J_{C , B}\, ,\, U_B\, ]\, $, one can pass to the reduced form of $\, U_B\, $ and is left with a quotient of the free abelian group generated by $\, B_1\, $ and $\, B_2\, $, intersected with 
$\, [\, U_B\, ,\, U_B\, ]\, $. Note that since $\, C\, $ is central all elements in $\, B_2\, $ are of the form 
$\, \{\> [\, u_x\, ,\, u_d\, ]\>\} \subseteq [\, U_B\, ,\, U_B\, ]\, $. Also all commutators 
$\, [\, u_c\, ,\, u_d\, ] = ( u_c\, u_d\, u^{-1}_{cd} )\, ( u_{dc}\, u^{-1}_c\, u^{-1}_d )\, $ where $\, c , d \in C\, $ and say $\, c < d\, $ for some chosen order on $\, C\, $, give a free set of generators for 
$\, J_C \cap [\, U_C\, ,\, U_C\, ]\, $ modulo $\, [\, [\, U_C\, ,\, U_C\, ] , U_C\, ] \subseteq [\, J_C\, ,\, U_C\, ]\, $. By the formulas given in the proof of Theorem 4 there are only three more types of relations left which are the image of  $\, [\, u_e\, ,\, B_1\, ]\, $ or 
\smallskip
$$ \bigl[\, u_e\, ,\, u_{cd}\, \bigr]\> \equiv\> \bigl[\, u_e\, ,\, u_c\,\bigr]\> \bigl[\, u_e\, ,\, u_d\,\bigr]\> , $$
\par\medskip\noindent
the image of $\, [\, B_2\, ,\, u_c\, ]\, $ given by the relations 
\smallskip
$$ \bigl[\, u_x\, ,\, u_{cd}\,\bigr]\>\equiv\> \bigl[\, u_x\, ,\, u_c\,\bigr]\> \bigl[\, u_x\, ,\, u_d\,\bigr]\> , $$
\par\medskip\noindent
and the image of $\, [\, u_c\, ,\, u_x\, u_y\, u^{-1}_{xy}\, ]\, $ modulo the former relations which is 
\smallskip
$$ \bigl[\, u_{c_{x , y}}\, ,\, u_c\,\bigr]\> \bigl[\, u_{(xy)}\, ,\, u_c\,\bigr]\>\equiv\> 
\bigl[\, u_x\, ,\, u_c\,\bigr]\> \bigl[\, u_y\, ,\, u_c\,\bigr]\> , $$
\par\medskip\noindent
so that grouping together the elements $\,\{\> [\, u_c\, ,\, u_d\, ]\>\}\, $ and $\, \{\> [\, u_x\, ,\, u_d\, ]\>\}\, $ in the form $\, \{\> [\, u_{cx}\, ,\, u_d\, ]\>\}\, $ one gets that $\, K^J_2 ( C , B )\, $ is a quotient of the tensor product $\, C \otimes B^{ab}\, $ which is presented by generators 
$\, \{\> c \otimes dx\>\} = \{\> c \otimes b\>\}\, $ subject to the relations 
$\, ( c_1\, c_2 \otimes b ) = ( c_1 \otimes b )\> ( c_2 \otimes b )\, $ and 
$\, ( c \otimes b_1\, b_2 ) = ( c \otimes b_1 )\> ( c \otimes b_2 ) = ( c \otimes b_2 )\> ( c \otimes b_1 ) = 
( c \otimes b_2\, b_1 )\, $. The only additional relations come from the identities 
$\, c \otimes c \equiv 0\, $, whence the result follows\qed
\par\bigskip\noindent 
The preceding Lemma gives an extremely important clou about the relation of the regular theory with $K^J$-theory, namely one has the following result.
\par\bigskip\noindent
{\bf Lemma 18.}\quad Let $\, n \geq 4\, $. Then $\, K_n ( N , F ) = K^J_n ( N , F )\, $ whenever 
$\, K_{n-1} ( N , F )\, $ is a torsion group, and in general $\, K^J_n ( N , F ) / K_n ( N , F )\, $ is free abelian for any regular pair $\, ( N , F )\, $. Also in case of a perfect group $\, F\, $ one has 
$\, K^J_3 ( F ) = K_3 ( F )\, $.
\par\bigskip\noindent
{\it Proof.}\quad Consider the extension 
\smallskip
$$ 1 \longrightarrow \bigl( R^n_{N , F}\, ,\, U^n_F \bigr) \longrightarrow 
\bigl( J_{R^{n-1}_{N , F}\, ,\, U^{n-1}_F}\, ,\, U^n_F \bigr) \longrightarrow 
\bigl( C_n\, ,\, U^n_F / R^n_{N , F} \bigr) \longrightarrow 1 $$
\par\medskip\noindent
with $\, C_n = K_{n+1} ( N , F ) \simeq J_{R^{n-1}_{N , F}\, ,\, U^{n-1}_F}\, /\, R^n_{N , F}\, $. If $\, C_n\, $ is torsion so is $\, K^J_2 ( C_n\, ,\, U^n_F\, /\, R^n_{N , F} )\, $ by Lemma 17, so that by exactness the quotient $\, K^J_{n+2} ( N , F )\, /\, K_{n+2} ( N , F )\, =\, 
K^J_2 ( J_{R^{n-1}_{N , F}\, ,\, U^{n-1}_F}\, ,\, U^n_F )\, /\, K^J_2 ( R^n_{N , F}\, ,\, U^n_F )\, $ is a torsion group (here we assume that 
$\, K_{n+2} ( N , F ) \rightarrowtail K^J_{n+2} ( N , F )\, $ is injective which follows from the next two Lemmas below). Now let $\, n \geq 2\, $ and consider the extension 
\smallskip
$$ 1 \longrightarrow \bigl(\,\bigl[\, J_{R^{n-1}_{N , F}\, ,\, U^{n-1}_F}\, ,\, U^n_F\,\bigr]\, ,\, U^n_F\,\bigr) \longrightarrow \bigl(\, J_{R^{n-1}_{N , F}\, ,\, U^{n-1}_F}\, ,\, U^n_F\,\bigr) \rightarrow\qquad\qquad\quad $$
$$\qquad\qquad\qquad\qquad\qquad\qquad\qquad \rightarrow 
\bigl(\, A\, ,\, {U^n_F \over [\, J_{R^{n-1}_{N , F}\, ,\, U^{n-1}_F}\, ,\, U^n_F\,]}\,\bigr) \longrightarrow 1\> . $$
\par\medskip\noindent
$\, A\, $ is a free abelian central subgroup of 
$\, U^n_F\, /\, [\, J_{R^{n-1}_{N , F}\, ,\, U^{n-1}_F}\, ,\, U^n_F\, ]\, $ by Lemma 14, so Lemma 17 gives that 
$\, K^J_2 ( A\, ,\, U^n_F\, /\, [\, J_{R^{n-1}_{N , F}\, ,\, U^{n-1}_F}\, ,\, U^n_F\, ] )\, $ is a quotient of 
$\, A \otimes (U^n_F)^{ab}\, $. Write 
\smallskip
$$  A\> =\> A' \times A''\> =\> A' \times 
{J_{R^{n-1}_{N , F}\, ,\, U^{n-1}_F}\,\cap\, [\, U^n_F\, ,\, U^n_F\, ] \over 
[\, J_{R^{n-1}_{N , F}\, ,\, U^{n-1}_F}\, ,\, U^n_F\, ]} \> . $$
\par\medskip\noindent
Then it is clear that $\, K^J_2 ( A\, ,\, U^n_F\, /\, [\, J_{R^{n-1}_{N , F}\, ,\, U^{n-1}_F}\, ,\, U^n_F\, ] )\, $ is a direct product with one factor equal to $\, A'' \otimes (U^n_F)^{ab}\, $ and the other one a quotient of 
$\, A' \otimes (U^n_F)^{ab}\, $. The image of $\, K^J_2 ( J_{R^{n-1}_{N , F}\, ,\, U^{n-1}_F}\, ,\, U^n_F )\, $ modulo $\, K^J_2 ( [\, J_{R^{n-1}_{N , F}\, ,\, U^{n-1}_F}\, ,\, U^n_F\, ]\, ,\, U^n_F )\, $ is the same as modulo $\, K^J_2 ( R^n_{N , F}\, ,\, U^n_F )\, $, because the map 
$\, K^J_2 ( [\, J_{R^{n-1}_{N , F}\, ,\, U^{n-1}_F}\, ,\, U^n_F\, ]\, ,\, U^n_F ) \twoheadrightarrow 
K^J_2 ( R^n_{N , F}\, ,\, U^n_F )\, $ is easily seen to be surjective by the fact that 
$\, ( R^n_{N , F}\, ,\, U^n_F )\, $ is 0-regular (modulo 
$\, [\, J_{R^{n-1}_{N , F}\, ,\, U^{n-1}_F}\, ,\, J_{R^{n-1}_{N , F}\, ,\, U^{n-1}_F}\, ]\, $). $A'\, $ can be lifted to $\, R^n_{N , F}\, /\, [\, J_{R^{n-1}_{N , F}\, ,\, U^{n-1}_F}\, ,\, U^n_F\, ]\, $ when chosen correctly, so the image of $\, K^J_2 ( J_{R^{n-1}_{N , F}\, ,\, U^{n-1}_F}\, ,\, U^n_F )\, $ does not intersect with the second factor and is isomorphic to a subgroup of $\, A'' \otimes (U^n_F)^{ab}\, $ which is a free abelian group, from which follows that $\, K^J_n ( N , F )\, /\, K_n ( N , F )\, $ is free abelian for $\, n \geq 4\, $ and trivial whenever $\, K_{n-1} ( N , F )\, $ is a torsion group. 
\par\noindent
If $\, F\, $ is perfect then so is its universal central extension $\,\overline F\, $ by $\, K_2 ( F ) =: C\, $. Then  $\, C \otimes (\overline F )^{ab} = 0\, $ and $\, K^J_3 ( F ) = K_3 ( F )\, $\qed
\par\bigskip\noindent
Before going further we must digress on certain aspects of $K^J$-theory for central subgroups.
\par\bigskip\noindent
{\bf Lemma 19.}\quad If $\, \mathbb Z\, $ is a central subgroup of an arbitrary group $\, B\, $, then 
$\, K^J_n ( \mathbb Z\, ,\, B ) = 0\, $ for all $\, n \geq 3\, $. For arbitrary central subgroups $\, C\, $ one has 
$\, K^J_n ( C , B ) = 0\, $ for all $\, n \geq 4\, $.
\par\bigskip\noindent
{\it Proof.}\quad First let $\, C \simeq \mathbb Z\, $ and consider the core $\,\{\> u_z\, B\, u^{-1}_z\>\}\, $ of $\, J_{\mathbb Z , B}\, $ for $\, \lambda ( J_{B / \mathbb Z} )\, $ as given by Lemma 3 of section 1. Let 
$\, B_3' = \{\> u_c\, u_y\, u^{-1}_{cy}\>\} \subseteq B_3\, $. Then its normalization is a cone hence irrelevant for our purposes and without loss of generality we may pass to the reduced form of $\, U_B\, $ generated by $\, U_{\mathbb Z}\, $ and $\, U_{B / \mathbb Z}\, $ (i.e. a semicanonical lift of the latter corresponding to a given section $\, s : B / \mathbb Z \nearrow B\, $). Then consider the normalization of the $B_2$-subcore. Since $\,\mathbb Z\, $ is central it has a basis of the form 
\smallskip
$$ \bigl\{\> u^{\pm }_{e_1}\cdots u^{\pm }_{e_m}\> u^{\pm }_{z_1}\cdots u^{\pm }_{z_k}\> 
\bigl[\, u_x\, ,\, u_d\,\bigr]\> u^{\mp }_{z_k}\cdots u^{\mp }_{z_1}\> u^{\mp }_{e_m}\cdots u^{\mp }_{e_1}\>
\bigr\} $$
\par\medskip\noindent
so that one gets a strict splitting of the corresponding subgroup $\, J_0\, $ into $\, U_{J_0\, ,\, U_B}\, $ modulo $\, [\, U_{J_0\, ,\, U_B}\, ,\, J_{U_B}\, ]\, $ and $\, K^J_n ( J_0\, ,\, U_B ) = 0\, $ for all $\, n \geq 2\, $ follows. The quotient by $\, J_0\, $ is isomorphic to $\, U_{\mathbb Z} \times U_{B / \mathbb Z} \simeq 
( J_{\mathbb Z} \rtimes {\mathbb Z} ) \times U_{B / \mathbb Z}\, $. Since the quotient of the extension 
\smallskip
$$ 1 \longrightarrow \bigl( J_{\mathbb Z}\, ,\, U_{\mathbb Z} \times U_{B / \mathbb Z} \bigr) \longrightarrow \bigl( U_{\mathbb Z}\, ,\, U_{\mathbb Z} \times U_{B / \mathbb Z} \bigr) \longrightarrow 
\bigl( {\mathbb Z}\, ,\, {\mathbb Z} \times U_{B / \mathbb Z} \bigr) \longrightarrow 1 $$
\par\medskip\noindent
is fully exact and $\, ( J_{\mathbb Z}\, ,\, U_{\mathbb Z} \times U_{B / \mathbb Z} )\, $ is second order exact as $\, J_{\mathbb Z}\, $ is free (this follows from the $P \times Q$-Lemma, see below), the maps
\smallskip 
$$ K^J_n ( J_{\mathbb Z}\, ,\, U_{\mathbb Z} \times U_{B / \mathbb Z} ) \rightarrowtail 
K^J_n ( U_{\mathbb Z}\, ,\, U_{\mathbb Z} \times U_{B / \mathbb Z} ) $$ 
\par\medskip\noindent
are injective for $\, n \geq 2\, $. On the other hand the image of 
$\, K^J_n ( J_{\mathbb Z , B}\, ,\, U_B )\, $ in 
$\, K^J_n ( U_{\mathbb Z}\, ,\, U_{\mathbb Z} \times U_{B / \mathbb Z} )\, $ is trivial as the map factors over $\, K^J_n ( U_{\mathbb Z , B}\, ,\, U_B )\, $. One can now apply Corollary 5.1 to the extension 
\smallskip
$$ 1 \longrightarrow \bigl( J_0\, ,\, U_B \bigr) \longrightarrow \bigl( J_{\mathbb Z , B}\, ,\, U_B \bigr) \longrightarrow \bigl( J_{\mathbb Z}\, ,\, U_{\mathbb Z} \times U_{B / \mathbb Z} \bigr) \longrightarrow 1$$
\par\medskip\noindent
and conclude that the corresponding sequence of $K^J_*$-groups is exact at all places 
$\, K^J_n ( J_{\mathbb Z , B}\, ,\, U_B )\, $ (which only involves full exactness of the pairs 
$\, ( J_0\, ,\, U_B )\, $ and $\, ( J_{\mathbb Z , B}\, ,\, U_B )\, $). It follows that 
$\, K^J_n ( J_{\mathbb Z , B}\, ,\, U_B ) = 0\, $ for all $\, n \geq 3\, $. 
\par\noindent
Now let $\, C\, $ be a central direct summand of $\, B\, $. To begin with assume that 
$\, C \simeq {\mathbb Z} / n {\mathbb Z}\, $ is a cyclic group. Let
\smallskip
$$ 1 \longrightarrow n {\mathbb Z} \longrightarrow {\mathbb Z} \longrightarrow C \longrightarrow 1 $$
\par\medskip\noindent
be a resolution of $\, C\, $ and consider the associated extension of $\, B\, $ (so that 
$\, {\mathbb Z}\, $ sits as a direct summand of $\, \widetilde B = {\mathbb Z} \times B / C\, $). Then 
$\, K^J_n ( {\mathbb Z}\, ,\, \widetilde B ) = 0\, $ for $\, n \geq 3\, $, so that 
$\, K^J_n ( C , B ) \simeq K^J_{n-1} ( n {\mathbb Z}\, ,\, \widetilde B )\, $ for $\, n \geq 4\, $ which is trivial by the argument above. Then Corollary 5.1 gives that $\, K^J_n ( C , B ) = 0\, $ for any finitely generated abelian group $\, C\, $ which sits as a direct summand of $\, B\, $. By continuity of the $K^J_*$-functors one gets that the same is true for arbitrary central direct summands, For let $\,\{\> C_{\lambda }\>\}\, $ be a net of finitely generated abelian groups converging to $\, C\, $ and 
$\,\{\> B_{\lambda } = C_{\lambda } \times B / C\>\}\, $ the associated net converging to $\, B\, $. It follows that $\, K^J_n ( C , B ) = {\displaystyle\lim_{\to}}\> K^J_n ( C_{\lambda }\, ,\, B_{\lambda } )\, $.
\par\noindent
Finally , let $\, C\, $ be an arbitrary central subgroup of $\, B\, $. Let $\, J_0\, $ denote the normalization of $B_2$- and $B_3'$-elements in $\, J_{C , B}\, $ as above, so that $\, K^J_n ( J_0\, ,\, U_B )\ = 0\, $ for all 
$\, n \geq 2\, $. Again the quotient is isomorphic to $\, U_C \times U_{B / C}\, $. Consider the extension 
\smallskip
$$ 1 \longrightarrow \bigl( J_C\, ,\, U_C \times U_{B / C} \bigr) \longrightarrow 
\bigl( U_C\, ,\, U_C \times U_{B / C} \bigr) \longrightarrow \bigl( C\, ,\, C \times U_{B / C } \bigr) \longrightarrow 1 $$
\par\medskip\noindent
which does not in general split. However one gets $\, K^J_n ( C\, ,\, C \times U_{B / C} ) = 0\, $ for all 
$\, n \geq 4\, $ by the argument above, hence an injective map 
$\, K^J_n ( J_C\, ,\, U_C \times U_{B / C} ) \rightarrowtail K^J_n ( U_C\, ,\, U_C \times U_{B / C} )\, $ for all $\, n \geq 3\, $ (since $\, ( J_C\, ,\, U_C \times U_{B / C} )\, $ is second order exact !) which implies that the image of $\, K^J_n ( J_{C , B}\, ,\, U_B )\, $ in $\, K^J_n ( J_C\, ,\, U_C \times U_{B / C} )\, $ is trivial for 
$\, n \geq 3\, $. By exactness of 
\smallskip
$$ K^J_n ( J_0\, ,\, U_B ) \longrightarrow K^J_n ( J_{C , B}\, ,\, U_B ) \longrightarrow 
K^J_n ( J_C, ,\, U_C \times U_{B / C} ) $$
\par\medskip\noindent
it follows that $\, K^J_n ( C , B ) = 0\, $ for $\, n \geq 4\, $\qed
\par\bigskip\noindent
{\bf Lemma 20}\quad Let $\, C\, $ be a central subgroup of $\, B\, $ and assume that 
$\, C \subseteq [\, B\, ,\, B\, ]\, $. Then $\, ( C , B )\, $ is fully exact and $\, K^J_n ( C , B ) = 0\, $ for 
$\, n \geq 3\, $.
\par\bigskip\noindent
{\it Proof.}\quad We begin with the second statement. The cases $\, n \geq 4\, $ follow from the preceding Lemma. First assume that $\, C \simeq {\mathbb Z} / n {\mathbb Z}\, $ is a cyclic group of order $\, n\, $. Using the $B_3'$-reduced form of $\, U_B\, $, the core of $\, J_{C , B}\, $ for $\, \lambda ( J_{B / C} )\, $ is reduced to the four components $\, [\, u_z\, ,\, B_2\, ]\, ,\, [\, u_z\, ,\, B_1\, ]\, ,\, B_2\, $ and $\, B_1\, $. Let 
$\, B_{1 , 0} \subseteq B_1\, $ denote the subset $\,\{\> u_1\, u_k\, u^{-1}_{k+1}\>\}\, $ where 
$\, 1\, $ is a generator of $\, C\, $ and $\, k \in \{\> 1 , \cdots , n-1\>\}\, $. Consider the union of the following two sets 
\smallskip
$$\quad \bigl\{\> w^{\pm }_1\cdots w^{\pm }_m\> u_z\> \bigl[\, u_c\, ,\, B_2\,\bigr]\> u^{-1}_z\> w^{\mp }_m\cdots w^{\mp }_1\>\bigr\} $$
$$ \cup\quad \bigl\{\> w^{\pm }_1\cdots w^{\pm }_m\> u_z\> \bigl[\, u_w\, ,\, B_{1 , 0}\,\bigr]\> u^{-1}_z\> w^{\mp }_m\cdots w^{\mp }_1\>\bigr\}\quad $$
\par\medskip\noindent
where $\, w_1 , \cdots , w_m\, $ are basis elements of $\, \lambda ( J_{B / C} )\, $, $\, w , z \in s ( B / C )\, $, $\, c \in C\, $. From the congruence 
\smallskip
$$ \bigl[\, u_w\, ,\, u_c\, u_d\, u^{-1}_{cd}\,\bigr]\> \equiv\> \bigl[\, u_w\, ,\, u_c\,\bigr]\> 
\bigl[\, u_w\, ,\, u_d\,\bigr]\> \bigl[\, u_{cd}\, ,\, u_w\,\bigr]\> \bigl[\, u_c\, ,\, \bigl[\, u_w\, ,\, u_d\,\bigr]\,\bigr]$$
\par\medskip\noindent
modulo $\, [\, J_{C , B}\, ,\, J_{C , B}\, ]\, $ relating the elements $\, [\, u_c\, ,\, B_2\,]\, $ and 
$\, [\, u_w\, ,\, B_1\, ]\, $ modulo the $B_2$-subcore one finds that the set $\,\{\> [\, u_c\, ,\, B_2\, ]\>\} \cup 
\{\> [\, u_w\, ,\, B_{1 , 0}\, ]\>\}\, $ is linear independent and generates the $( B_2 \cup [\, u_z\, ,\, B_1\, ] )$-subcore modulo elements of $\, B_{2 , 0} = \{\> [\, u_x\, ,\, u_1\, ]\>\}\, $ and $\, B_1\, $ (and for the former one gets that $\, [\, u_x\, ,\, u_1\, ]^n\, $ is contained in the range of 
$\, \{\> [\, u_c\, ,\, B^x_2\, ]\>\} \cup \{\> [\, u_x\, ,\, B_{1 , 0}\, ]\>\}\, $ modulo 
$\, [\, B_1\, ,\, B^x_2\, ] \subseteq [\, J_{C , B}\, ,\, J_{C , B}\, ]\, $. A full basis of the 
$( B_2 \cup [\, u_z\, ,\, B_1\, ] )$-subcore is given (modulo $\, [\, J_{C , B}\, ,\, J_{C , B}\, ]\, $) by the set 
$\,\{\> u_e\, B_2\, u^{-1}_e\>\}\, $ which can be changed according to the assignment
\smallskip
$$ u_e\> B_2\> u^{-1}_e\quad \mapsto\quad \bigl[\, u_e\, ,\, B_2\,\bigr] ,\quad e \neq 1 , $$
$$ B_2\,\backslash\, B_{2 , 0}\quad \mapsto\quad \bigl[\, u_w\, ,\, B^+_{1 , 0}\,\bigr] , \qquad\quad$$
$$ B_{2 , 0}\qquad\>\>\> \mapsto\quad\> B_{2 , 0} \qquad\qquad\qquad$$
\par\medskip\noindent
where $\, B^+_{1 , 0} = \{\> u_1\, u_k\, u^{-1}_{k+1}\,\vert\, k = 2 ,\cdots , n-1\>\}\, $ and the second assignment is given explicitely by 
\smallskip
$$ \bigl[\, u_w\, ,\, u_k\,\bigr]\quad \mapsto\quad \bigl[\, u_w\, ,\, u_1\, u_k\, u^{-1}_{k+1}\,\bigr]\> . $$
\par\medskip\noindent
Hence our set will be linear independent modulo $\, [\, J_{C , B}\, ,\, J_{C , B}\, ]\, $ if and only if the corresponding set of conjugates of $\,\{\> u_e\, B_2\, u^{-1}_e\>\}\, $ is linear independent modulo 
$\, [\, J_{C , B}\, ,\, J_{C , B}\, ]\, $. Then modulo $\, [\, J_{C , B}\, ,\, J_{C , B}\, ]\, $ the coefficients 
$\, w^{\pm }_1\cdots w^{\pm }_m\, u_z\, $ can be replaced by coefficients of the form 
$\, u^{\pm }_{z_1}\cdots u^{\pm }_{z_k}\, $ in a one-to-one fashion and the set 
\smallskip
$$ \bigl\{\> u^{\pm }_{z_1}\cdots u^{\pm }_{z_k}\> u_e\> B_2\> u^{-1}_e\ 
u^{\mp }_{z_k}\cdots u^{\mp }_{z_1}\>\bigr\} $$
\par\medskip\noindent
is part of a basis for the free abelian group $\, [\, U_{C , B}\, ,\, U_B\, ]\, /\, [\, J_{C , B}\, ,\, J_{C , B}\, ]\, $, so that our original set is also linear independent modulo $\, [\, J_{C , B}\, ,\, J_{C , B}\, ]\, $ and this remains true if we replace the second set by 
\smallskip
$$ \bigl\{\> w^{\pm }_1\cdots w^{\pm }_m\> \bigl[\, u_z\, ,\,\bigl[\, u_w\, ,\, B_{1 , 0}\,\bigr]\,\bigr]\> 
w^{\mp }_m\cdots w^{\mp }_1\>\bigr\} \> . $$
\par\medskip\noindent
Let $\, J_0\, $ denote the subgroup of $\, J_{C , B}\, $ generated by $\, [\, J_{C , B}\, ,\, J_{C , B}\, ]\, $, the two sets above together with 
\smallskip
$$ \bigl\{\> w^{\pm }_1\cdots w^{\pm }_m\> \bigl[\, u_w\, ,\, B_{1 , 0}\,\bigr]\> w^{\mp }_m\cdots w^{\mp }_1\>\bigr\} $$
\par\medskip\noindent
and 
\smallskip
$$ \bigl\{\> w^{\pm }_1\cdots w^{\pm }_m\> B_{1 , 0}\> w^{\mp }_m\cdots w^{\mp }_1\>\bigr\} \> . $$
\par\medskip\noindent
Then $\, J_0\, $ generates $\, J_{C , B}\, $ up to torsion elements of order dividing $\, n\, $, which are contained in the image of the set 
\smallskip
$$ \bigl\{\> w^{\pm }_1\cdots w^{\pm }_m\> u_z\> B_{2 , 0}\; u^{-1}_z\> 
w^{\mp }_m\cdots w^{\mp }_1\>\bigr\}  $$
\par\medskip\noindent
(because $\, B_1\, $ is generated modulo $\, [\, J_{C , B}\, ,\, J_B\, ]\, $ by $\, B_{1 , 0}\, $ and elements 
$\, [\, [\, u_z\, ,\, u_w\, ]\, ,\, B_{1 , 0}\, ]\, $). Let $\, \widetilde J_0\, $ be the subgroup generated by 
$\, [\, J_{C , B}\, ,\, J_{C , B}\, ]\, $, the first (modified) two sets, together with the subbasis corresponding to the $\{\> u_z\, B_{2 , 0}\, u^{-1}_z\>\}$-subcore. It obviously contains the subbasis for the 
$[\, u_z\, ,\, B_2\, ]$-component and on adding the set 
\smallskip
$$ \bigl\{\> w^{\pm }_1\cdots w^{\pm }_m\> \bigl[\, u_w\, ,\, B_{1 , 0}\,\bigr]\> 
w^{\mp }_m\cdots w^{\mp }_1\>\bigr\} $$
\par\medskip\noindent
it would also contain the subbasis corresponding to the $B_2$-component. Moreover the quotient of $\, J_{C , B}\, $ by this larger set will be isomorphic to $\, J_C\, /\, [\, J_C\, ,\, J_C\, ]\, $, so that also the quotient of $\, J_{C , B}\, $ by $\, \widetilde J_0\, $ is free abelian. Then also the quotient of $\, J_0\, $ by the subgroup containing $\, [\, J_{C , B}\, ,\, J_{C , B}\, ]\, $ and the first two sets is free abelian. Modulo the subgroup generated by the $[\, u_z\, ,\, B_2\, ]$-component and the $[\, u_c\, ,\, B_2\, ]$-component of the basis of $\, J_{C , B}\, $ , any additional element in this subgroup is a $\lambda ( J_{B / C} )$-conjugate of an element $\, [\, u_z\, ,\, [\, u_w\, ,\, B_1\, ]\, ]\, $ leading to relations of the form 
\smallskip
$$ \bigl[\> u_z\, u_w\, u^{-1}_{zw}\, ,\, \bigl[\, u_{(zw)}\, ,\, u_{c_{z , w}}\, B_1\, u^{-1}_{c_{z , w}}\,\bigr]\, 
u_{c_{z , w}}\, B_1\, u^{-1}_{c_{z , w}}\,\bigr]\quad \equiv\qquad\qquad\qquad\qquad $$
$$\qquad\qquad\qquad \bigl[\, u_{c_{z , w}}\, B_1\, u^{-1}_{c_{z , w}}\, ,\, u_{(zw)}\,\bigr]\> \bigl[\, B_1\, ,\, u_{c_{z , w}}\,\bigr]\> 
\bigl[\, u_z\, ,\, B_1\,\bigr]\> \bigl[\, u_w\, ,\, B_1\,\bigr] $$
\par\medskip\noindent
(and $\lambda ( J_{B / C} )$-conjugates thereof), so the elements 
\smallskip
$$ \bigl\{\> w^{\pm }_1\cdots w^{\pm }_m\> B_{1 , 0}\> w^{\mp }_m\cdots w^{\mp }_1\> \bigr\} $$
\par\medskip\noindent
are still linear independent in the quotient and may be added to our existing set. Since each $\, u_c\, $ can be factored modulo $\, J_{C , B}\, $ into a product 
$\, w^{\pm }_{j_1}\cdots w^{\pm }_{j_k} \in \lambda ( J_{B / C} )\, $ and a product 
$\, \Pi_{\lambda , \mu }\> [\, u_{z_{\lambda }} , u_{z_{\mu }}\, ] \in 
[ \sigma ( U_{B / C} ) , \sigma ( U_{B / C} ) ] \, $, the enveloping subgroup also contains all elements
\smallskip 
$$ \bigl\{\> w^{\pm }_1\cdots w^{\pm }_m\> u_c\> B_{1 , 0}\> u^{-1}_c\> 
w^{\mp }_m\cdots w^{\mp }_1\>\bigr\} \> , $$
\par\medskip\noindent
and hence the $B_1$-component of the basis of $\, J_{C , B}\, $ (since $\, [\, J_C\, ,\, U_C\, ]\, $ is generated modulo $\, [\, J_C\, ,\, J_C\, ]\, $ by elements $\, [\, u_c\, ,\, B_{1 , 0}\, ]\, $ and 
$\, J_C\, /\, [\, J_C\, ,\, U_C\, ]\, $ can be identified with the subgroup of $\, J_C\, /\, [\, J_C\, ,\, J_C\, ]\, $ generated by $\, B_{1 , 0}\, $). Then modulo our extended subgroup the relations above restricted to the 
$[\, u_w\, ,\, B_{1 , 0}\, ]$-component are
\smallskip
$$ \bigl[\, u_z\, u_w\, u^{-1}_{zw}\, ,\, \bigl[\, u_{(zw)}\, ,\, B_{1 , 0}\,\bigr]\,\bigr]\quad \equiv\quad 
\bigl[\, B_{1 , 0}\, ,\, u_{(zw)}\,\bigr]\> \bigl[\, u_z\, ,\, B_{1 , 0}\,\bigr]\> 
\bigl[\, u_w\, ,\, B_{1 , 0}\,\bigr] $$
\par\medskip\noindent
which has to be read as a relation modulo the $( [\, u_z\, ,\, B_2\, ] \cup [\, u_c\, ,\, B_2\, ] \cup B_1 )$-component with arbitrary $\lambda ( J_{B / C} )$-coefficients. Clearly the quotient of the subgroup generated by the subbasis for the $[\, u_w\, ,\, B_{1 , 0}\, ]$-component by these relations is still free abelian, so we may add a corresponding basis of the form 
\smallskip
$$ \bigl\{\> w^{\pm }_1\cdots w^{\pm }_m\> \bigl[\, u_w\, ,\, B_{1 , 0}\, \bigr]\> w^{\mp }_m\cdots w^{\mp }_1\>\bigr\} $$
\par\medskip\noindent
with $\, w_m = u_x\, u_y\, u^{-1}_{xy}\, $ and $\, (xy) \neq w\, $or else negative exponent of $\, w_m\, $. We will now construct a $U_B$-normal lift of $\, J_0\, $ to $\, U^2_{C , B}\, $ modulo 
$\, [\, J^2_{C , B}\, ,\, U^2_B\, ]\, +\, [\, U^2_{C , B}\, ,\, J_{U_B}\, ]\, $ showing that the image of 
$\, K^J_2 ( J_0\, ,\, U_B )\, $ in $\, K^J_2 ( J_{C , B}\, ,\, U_B )\, $ is trivial. Recall that modulo 
$\, [\, J_{C , B}\, ,\, J_{C , B}\, ]\, $ we have constructed the following basis of $\, J_0\, $ consisting of the four sets
\smallskip
$$ \bigl\{\> w^{\pm }_1 \cdots w^{\pm }_m\> u_z\> \bigl[\, u_c\, ,\, B_2\,\bigr]\> u^{-1}_z\> 
w^{\mp }_m\cdots w^{\mp }_1\>\bigr\} $$
$$ \cup\quad \bigl\{\> w^{\pm }_1\cdots w^{\pm }_m\> \bigl[\, u_z\, ,\, \bigl[\, u_w\, ,\, B_{1 , 0}\,\bigr]\,\bigr]\> w^{\mp }_m\cdots w^{\mp }_1\> \bigr\} $$
$$ \cup\quad \bigl\{\> w^{\pm }_1\cdots w^{\pm }_m\> B_{1 , 0}\> w^{\mp }_m\cdots w^{\mp }_1\>\bigr\} $$
$$ \cup\quad \bigl\{\> w^{\pm }_1\cdots w^{\pm }_m\> \bigl[\, u_w\, ,\, B_{1 , 0}\, \bigr]\> 
w^{\mp }_m\cdots w^{\mp }_1\> \bigr\} $$
\par\medskip\noindent
where in the last set $\, w^{\pm }_m\, $ is subject to the restrictions described above. For 
$\, [\, J_{C , B}\, ,\, J_{C , B}\, ]\, $ we choose the canonical lift obtained from an arbitrary lift of 
$\, J_{C , B}\, $. Then the corresponding quotient of $\, U^2_{C , B}\, $ is abelian, and we may extend by a lift of the free abelian group given by the basis above, lifting elements $\, B_2\, $ and $\, B_{1 , 0}\, $ by the rule 
\smallskip
$$ \bigl[\, u_x\, ,\, u_d\,\bigr]\quad \mapsto\quad u_{[\, u_x\, ,\, u_d\,]}\qquad ,\qquad 
u_1\, u_k\, u^{-1}_{k+1}\quad \mapsto\quad u_{u_1\, u_k\, u^{-1}_{k+1}} $$
\par\medskip\noindent
and completing normally for the canonical lift $\, U_B \rightarrow U^2_B\, $. One easily checks that this gives a $U_C$- (hence also $U_B$-)normal lift of the first two sets. Also writing $\, u_c\, $ as a product in 
$\, \lambda ( J_{B / C} )\, $ and $\, [\, \sigma ( U_{B / C} )\, ,\, \sigma ( U_{B / C} )\, ]\, $ as above one checks that the lift is normal also on the third set (modulo the second and fourth) and on the fourth set (modulo the first and second set). 
\par\noindent
The quotient of $\, U_{C , B}\, $ by $\, J_0\, $ is a torsion group generated by elements of order dividing $\, n\, $ (certain conjugates of $\, B_{2 , 0}\, $ and $\, u_c\, $) such that the image of $\, J_{C , B}\, $ is central in the quotient $\, U_B / J_0\, $ and contained in $\, [\, U_B / J_0\, ,\, U_B / J_0\, ]\, $. One gets an extension 
\smallskip
$$ 1 \rightarrow \bigl( J_{C , B} / J_0\, ,\, U_B / J_0 \bigr) \rightarrow 
\bigl( U_{C , B} / J_0\, ,\, U_B / J_0 \bigr) \rightarrow \bigl( C , C \times U_{B / C} \bigr) \rightarrow 1\> . $$
\par\medskip\noindent
Lemma 17 gives $\, K^J_2 ( J_{C , B} / J_0\, ,\, U_B / J_0 )\,\simeq\, {J_{C , B} \over J_0} \otimes 
( C \times U^{ab}_{B / C} ), $ and a similar argument (separating the variables coming from 
$\, J_{C , B} / J_0\, $ and $\, C\, $ respectively) for 
$\, K^J_2 ( ( J_{C , B} / J_0 ) \rtimes C\, ,\, U_B / J_0 )\, $ yields that the map 
\smallskip
$$ K^J_2 ( J_{C , B} / J_0\, ,\, U_B / J_0 ) \rightarrowtail K^J_2 ( U_{C , B} / J_0\, ,\, U_B / J_0 )$$
\par\medskip\noindent
is injective. In the following let $\, \Gamma = J_{C , B} / J_0\, $, so that 
$\, U_{C , B} / J_0 \simeq \Gamma \rtimes C\, $ and put $\, \widetilde B = U_B / J_0\, $. We denote elements of $\, U_{\Gamma }\, $ by letters $\, u_{\gamma }\, ,\, u_{\delta }\, $ etc., elements of $\, U_C\, $ by letters $\, u_c\, ,\, u_d\, $ etc. and elements of $\, \sigma ( U^2_{B / C} )\, $ by letters 
$\, u_x\, ,\, u_y\, ,\, u_z\, $ etc.. Then, computing $\, K^J_2 ( \Gamma \rtimes C\, ,\, \widetilde B )\, $ one notes the following facts: first one can pass to the reduced form of $\, U_{\widetilde B}\, $ generated by the three components $\, U_{\Gamma }\, ,\, U_C\, $ and $\, \sigma ( U^2_{B / C} )\, $, also since there is a splitting of $\, U_{B / C}\, $ into $\, \widetilde B\, $ one gets that the elements 
$\, ( \gamma c )_{x , y}\, $ are all trivial and $\, \lambda ( J_{U_{B / C}} ) = \sigma ( J_{U_{B / C}} )\, $. Then modulo $\, [\, J_{\Gamma \times C\, ,\, \widetilde B}\, ,\, J_{\Gamma \times C\, ,\, \widetilde B}\,]\, $ and the $[\, u_z\, ,\, B_2\, ]$- and $[\, u_z\, ,\, B_1\, ]$-subcores the (abelian) basis of 
$\, J_{\Gamma \times C\, ,\, \widetilde B}\, $ is reduced to the $B_2$- and $B_1$-subcores subject only to the following remaining relations, which are the images of $\, [\, B_1\, ,\, u_{\gamma }\, ]\> ,\> 
[\, B_1\, ,\, u_c\, ]\> ,\> [\, B^{\delta }_2\, ,\, u_{\gamma }\, ]\> ,\> [\, B^d_2\, ,\, u_{\gamma }\, ]\> ,\> 
[\, B^{\delta }_2\, ,\, u_1\, ]\> ,\> [\, B^d_2\, ,\, u_1\, ]\, $, $\, [\, u_1\, ,\, u_x\, u_y\, u^{-1}_{xy}\, ]\, $ and 
$\, [\, u_{\gamma }\, ,\, u_x\, u_y\, u^{-1}_{xy}\, ]\, $. Elements of $\, B_1\, $ are of the three types 
$\, ( u_{\gamma }\, u_{\delta }\, u^{-1}_{\gamma \delta } )\, ,\, ( u_c\, u_d\, u^{-1}_{cd} )\, $ and 
$\, [\, u_c\, ,\, u_{\gamma }\, ]\, $ in the reduced form of $\, U_{\widetilde B}\, $. The relations 
$\, [\, u_{\gamma }\, ,\, u_x\, u_y\, u^{-1}_{xy}\, ]\, ,\, [\, u_{\gamma }\, ,\, 
u_{\alpha }\, u_{\beta }\, u^{-1}_{\alpha \beta }\, ]\, ,\, [\, u_{\gamma }\, ,\, [\, u_x\, ,\, u_{\delta }\, ]\, ]\, $ are already relations for $\, K^J_2 ( \Gamma \, ,\, \widetilde B )\, $ so they may be dropped from consideration. Then in the reduced $\, U_{\widetilde B}\, $ the relations 
$\, [\, u_1\, ,\, [\, u_x\, ,\, u_{\delta }\, ]\, ]\, $ are contained in $\, [\, u_{\gamma }\, ,\, B^d_2\, ]\, $ modulo some extra terms from $\, [\, U_{\Gamma }\, ,\, U_{\Gamma }\, ]\, $, the image of which is trivial already in $\, K^J_2 ( \Gamma\, ,\, \widetilde B )\, $, so that relations $\, [\, u_1\, ,\, B^{\delta }_2\, ]\, $ may be dropped from consideration. The relation $\, [\, u_1\, ,\, u_c\, u_d\, u^{-1}_{cd}\, ]\, $ reduces the subbasis 
$\, \{\> u_c\, u_d\, u^{-1}_{cd}\>\}\, $ to the form $\,\{\> u_1\, u_k\, u^{-1}_{k+1}\>\}\, $, otherwise does not interfere with the representatives $\, [\, u_{\gamma }\, ,\, u_c\, ]\, $ and $\, [\, u_{\gamma }\, ,\, u_x\, ]\, $ for the image of $\, K^J_2 ( \Gamma\, ,\, \widetilde B )\, $. Then relation 
$\, [\, u_{\gamma }\, ,\, u_c\, u_d\, u^{-1}_{cd}\, ]\, $ reduces the subbasis 
$\,\{\> [\, u_c\, ,\, u_{\delta }\, ]\>\}\, $ to $\,\{\> [\, u_1\, ,\, u_{\delta }\, ]\>\}\, $ subject to the (independent) relations $\, [\, u_1\, ,\, u_{\delta }\, ]^n = 0\, $. Adding relations 
$\, [\, u_{\gamma }\, ,\, [\, u_1\, ,\, u_{\delta }\, ]\, ]\, ,\, [\, u_1\, ,\, [\, u_c\, ,\, u_{\delta }\, ]\, ]\, $ and 
$\, [\, u_1\, ,\, u_{\gamma }\, u_{\delta }\, u^{-1}_{\gamma \delta }\, ]\, $ reduces the subgroup generated by 
$\,\{\> [\, u_1\, ,\, u_{\delta }\, ]\>\}\, $ to $\, \Gamma \otimes C \subseteq \Gamma \otimes 
( \widetilde B / \Gamma )^{ab}\, $. The relation $\, [\, u_1\, ,\, B^d_2\, ]\, $ can be written modulo the former relations as 
\smallskip
$$ 0\quad \equiv\quad \bigl[\, u_1\, ,\, u_x\,\bigr]\> \bigl[\, u_k\, ,\, u_x\,\bigr]\> \bigl[\, u_x\, ,\, u_{k+1}\,\bigr]\> \bigl[\, u_{[\, x\, ,\, k\, ]}\, ,\, u_1\,\bigr] $$
$$\qquad \equiv\quad \bigl[\, u_1\, ,\, u_x\,\bigr]\> \bigl[\, u_k\, ,\, u_x\,\bigr]\> 
\bigl[\, u_x\, ,\, u_{k+1}\,\bigr]\> \bigl[\, u_{[\, x\, ,\, 1\, ]}\, ,\, u_1\,\bigr]^k \> . $$
\par\medskip\noindent
Since the elements $\,\{\> [\, u_1\, ,\, u_x\, ]\> [\, u_k\, ,\, u_x\, ]\> [\, u_x\, ,\, u_{k+1}\, ]\>\}\, $ are linear independent for $\, k = 1 , \cdots , n-1\, $, these relations cannot interfere with the representatives 
$\,\{\> [\, u_1\, ,\, u_{\delta }\, ]\>\}\, $ or $\,\{\> [\, u_x\, ,\, u_{\delta }\, ]\>\}\, $. Finally, relation 
$\, [\, u_1\, ,\, u_x\, u_y\, u^{-1}_{xy}\, ]\, $ reads 
\smallskip
$$ 0 \equiv \bigl( \bigl[\, u_{xy}\, ,\, u_1\,\bigr]\, u^{-1}_{[\, xy\, ,\, 1\, ]} \bigr)\, 
\bigl( u_{[\, x\, ,\, 1\, ]}\, \bigl[\, u_1\, ,\, u_x\,\bigr] \bigr)\, 
\bigl( u_{[\, y\, ,\, 1\, ]}\, \bigl[\, u_1\, ,\, u_y\,\bigr] \bigr)\, 
\bigl[\, u_{[\, y\, ,\, 1\, ]}\, ,\, u_x\,\bigr] $$
\par\medskip\noindent
modulo the former relations. Consider the subgroup of the $B^1_2$-component generated by triples
\smallskip 
$$ \bigl\{\> \bigl( \bigl[\, u_{xy}\, ,\, u_1\,\bigr]\, u^{-1}_{[\, xy\, ,\, 1\, ]} \bigr)\, 
\bigl( u_{[\, y\, ,\, 1\, ]}\, \bigl[\, u_1\, ,\, u_y\,\bigr] \bigr)\, 
\bigl( u_{[\, x\, ,\, 1\, ]}\, \bigl[\, u_1\, ,\, u_x\,\bigr] \bigr)\>\bigr\} \> . $$
\par\medskip\noindent
Since $\, U_{B / C}\, $ is a free group it has a basis consisting of special triples of the form 
\smallskip
$$ \bigl\{\> \bigl( \bigl[\, u_{x f^{\pm }_{\lambda }}\, ,\, u_1\,\bigr]\, u^{-1}_{[\, x f^{\pm }_{\lambda }\, ,\, 1\, ]} \bigr)\, 
\bigl( u_{[\, f_{\lambda }\, ,\, 1\, ]}\, \bigl[\, u_1\, ,\, u_{f_{\lambda }}\,\bigr] \bigr)^{\pm }\, 
\bigl( u_{[\, x\, ,\, 1\, ]}\, \bigl[\, u_1\, ,\, u_x\,\bigr] \bigr)\>\} $$
\par\medskip\noindent
where $\,\{\> f_{\lambda }\>\}\, $ is a basis of $\, U_{B / C}\, $ and the degree of the element 
$\, x f^{\pm }_{\lambda }\, $ with respect to this basis is larger than the degree of $\, x\, $. From this presentation one checks that if a combination of relations 
$\,\{\> [\, u_1\, ,\, u_x\, u_y\, u^{-1}_{xy}\, ]\>\}\, $ should lead to a trivial expression modulo the extra 
$[\, u_{[\, y\, ,\, 1\, ]}\, ,\, u_x\, ]$-terms, then the corresponding combination of these extra terms is trivial when mapped to $\, \Gamma \otimes U^{ab}_{B / C}\, $. Since we have exhausted thereby all our relations one concludes that the map 
$\, K^J_2 ( \Gamma\, ,\, \widetilde B ) \rightarrowtail K^J_2 ( (\Gamma \times C)\, ,\, \widetilde B )\, $ is injective and $\, K^J_2 ( J_{C , B}\, ,\, U_B ) = 0\, $ follows.
\par\noindent
We have already seen that $\, ( C , B )\, $ is exact of first order. Assume that 
$\, C \simeq {\mathbb Z} / n {\mathbb Z}\, $ is a (possibly infinite) cyclic group. We first show that 
$\, ( C\, ,\, C \times U )\, $ is second order (hence fully) exact whenever $\, U\, $ is a free group (this is clear if $\, C \simeq {\mathbb Z}\, $). Consider the diagram 
\smallskip
$$ \vbox{\halign{ #&#&#\cr
\hfil $\bigl( J_{C \times U}\, ,\, J_U \rtimes U_{C \times U} \bigr)$\hfil & \hfil $\rightarrow $\hfil & 
\hfil $\bigl( J_{C \times U} \times J_U\, ,\, J_U \rtimes U_{C \times U} \bigr)$\hfil \cr
\hfil $\Bigm\downarrow $\hfil && \hfil $\Bigm\downarrow $\hfil \cr
\hfil $\bigl( J_{C \times U} + U_{C , C \times U} ,\ J_U \rtimes U_{C \times U} \bigr)$\hfil & 
\hfil $\rightarrow $\hfil & 
\hfil $\bigl( ( J_{C \times U} + U_{C , C \times U} ) \times J_U , J_U \rtimes U_{C \times U} \bigr)$ \cr
\hfil $\Bigm\downarrow $\hfil && \hfil $\Bigm\downarrow $\hfil \cr
\hfil $\bigl( C\, ,\, C \times U_U \bigr)$\hfil & \hfil $\rightarrow $\hfil & 
\hfil $\bigl( C\, ,\, C \times U \bigr)$\hfil \cr  }} \> . $$
\par\medskip\noindent
Both pairs in the bottom row are exact, and the middle horizontal map induces an injective map in $K^J_2$ by existence of a splitting. Then an element in the kernel of 
$\, K^J_2 ( J_{C \times U}\, ,\, J_U \rtimes U_{C \times U} ) \rightarrow 
K^J_2 ( J_{C \times U} \times J_U\, ,\, J_U \rtimes U_{C \times U} )\, $ lifts to 
$\, K^J_3 ( C\, ,\, C \times U_U )\, $ which we claim to be isomorphic to $\, K^J_3 ( C\, ,\, C \times U )\, $ and to $\, K^J_3 ( C ) \simeq {\mathbb Z} / n {\mathbb Z}\, $. To see this note that an analogous description of a basis of $\, J_{C\, ,\, C \times U_U}\, $, respectively of $\, J_{C\, ,\, C \times U}\, $ as in Lemma 5 for $\, C\, $ free abelian yields that in case $\, C\, $ is generated by a single element one gets no relative terms in $\, K^J_3 ( C\, ,\, C \times U_U )\, $, i.e. the split injection 
$\, K^J_3 ( C ) \rightarrowtail K^J_3 ( C\, ,\, C \times U_U )\, $ is an isomorphism. It is easy to compute $\, K^J_3 ( C ) \simeq K^J_2 ( n {\mathbb Z}\, ,\, {\mathbb Z} ) \simeq {\mathbb Z} / n {\mathbb Z}\, $ from Lemma 17. Also since $\, U\, $ is free one gets a compatible splitting 
\smallskip
$$ \vbox{\halign{ #&#&#&#&#&#&#&#&#\cr
\hfil $1$\hfil & \hfil $\>\longrightarrow\> $\hfil & \hfil $\bigl( C\, ,\, C \times U_U \bigr)$\hfil & 
\hfil $\>\longrightarrow\> $\hfil & \hfil $C \times U_U$\hfil & \hfil $\>\longrightarrow\> $\hfil &
\hfil $U_U$\hfil & \hfil $\>\longrightarrow\> $\hfil & \hfil $1$\hfil \cr
&& \hfil $\Bigm\downarrow\Bigm\uparrow $\hfil && \hfil $\Bigm\downarrow\Bigm\uparrow $\hfil && 
\hfil $\Bigm\downarrow\Bigm\uparrow $\hfil &&\cr
\hfil $1$\hfil & \hfil $\>\longrightarrow\> $\hfil &\hfil $\bigl( C\, ,\, C \times U \bigr)$\hfil & 
\hfil $\>\longrightarrow\> $\hfil & \hfil $C \times U$\hfil & \hfil $\>\longrightarrow\> $\hfil & 
\hfil $U$\hfil & \hfil $\>\longrightarrow\> $\hfil & $1$\cr }}  $$
\par\medskip\noindent
so that by a standard argument (compare the examples in section 5) one gets that
$\, K^J_2 ( J_{C \times U}\, ,\, J_U \rtimes U_{C \times U} ) \rightarrowtail 
K^J_2 ( J_{C \times U} \times J_U\, ,\, J_U \rtimes U_{C \times U} )\, $ is injective and 
$\, ( C\, ,\, C \times U )\, $ is fully exact. In particular $\, ( C\, ,\, C \times U_{B / C} )\, $ is fully exact, so the $P \times Q$-Lemma gives that $\, ( J_{B / C} \times C\, ,\, C \times U_{B / C} ) = 
( J^{\Delta }_{B / C} \times C\, ,\, C \times U_{B / C} )\, $ is fully exact (where $\, J^{\Delta }_{B / C}\, $ denotes the normal copy contained in $\, U_{B / C}\, $) since also $\, ( J_{B / C}\, ,\, U_{B / C} )\, $ is fully exact. This in turn implies, again by the $P \times Q$-Lemma, that 
$\, ( J_{B / C}\, ,\, C \times U_{B / C} )\, $ is exact, so that $\, ( C , B )\, $ is second order (hence fully) exact 
(if $\, C \subseteq [\, B\, ,\, B\, ]\, $). One now uses induction to prove that the same is true for a finitely generated central subgroup $\, C\, $ with $\, C \subseteq [\, B\, ,\, B\, ]\, $. Assume the statement is true for all $\, C'\, $ generated by less or equal than k elements and let $\, C = C_0 \times C'\, $ be generated by k+1  elements such that $\, C_0 \simeq {\mathbb Z}\, $ or $\, C_0 \simeq {\mathbb Z} / n {\mathbb Z}\, $. Consider the extension
\smallskip
$$ 1 \longrightarrow \bigl( C_0\, ,\, B \bigr) \longrightarrow \bigl( C\, ,\, B \bigr) \longrightarrow 
\bigl( C'\, ,\, B' \bigr) \longrightarrow 1 \> . $$
\par\medskip\noindent
The $P \times Q$-Lemma yields that $\, ( C , B )\, $ is fully exact and Corollary 5.1 then yields that 
$\, K^J_3 ( C , B ) = 0\, $. By continuity of the $K^J_3$-functor one gets $\, K^J_3 ( C , B ) = 0\, $ for arbitrary central subgroups with $\, C \subseteq [\, B\, ,\, B\, ]\, $, and much in the same way one finds that 
$\, ( C , B )\, $ is second order exact by considering the system of inclusions
\smallskip
$$ \vbox{\halign{ #&#&#\cr
\hfil $K^J_2 ( J_B\, ,\, J_{B / C_{\lambda }} \rtimes U_B )$\hfil & \hfil $\rightarrowtail $\hfil & 
\hfil $K^J_2 ( J_B \times J_{B / C_{\lambda }}\, ,\, J_{B / C_{\lambda }} \rtimes U_B )$ \hfil \cr
\hfil $\Bigm\downarrow $\hfil && \hfil $\Bigm\downarrow $\hfil \cr
\hfil $K^J_2 ( J_B\, ,\, J_{B / C} \rtimes U_B )$\hfil & \hfil $\rightarrow $ \hfil & 
\hfil $K^J_2 ( J_B \times J_{B / C}\, ,\, J_{B / C} \rtimes U_B )$\hfil \cr }} $$
\par\medskip\noindent
converging to the lower horizontal map\qed
\par\bigskip\noindent  
With the two preceding Lemmas at hand showing (by use of Corollary 5.1) that the map 
$\, K^J_2 ( R^n_{N , F}\, ,\, U^n_F ) \rightarrowtail K^J_2 ( J^n_{N , F}\, ,\, U^n_F )\, $ is injective together with the Special Deconstruction Theorem the proof of Theorem 2 is now complete. The next result shows that the regular theory really "exists" and does not depend on the choices employed in its construction.
\par\bigskip\noindent
{\bf Proposition 1.}\quad Let $\, N\, $ be a full subgroup of a perfect group $\, F\, $. Define inductively 
$\, R^n_{0 , N , F} := R_{0 , R^{n-1}_{0 , N , F}\, ,\, U^{n-1}_F}\, $ for $\, n \geq 1\, $. Then the natural maps 
$\, K^J_2 ( R^n_{0 , N , F}\, ,\, U^n_F ) \twoheadrightarrow K^J_2 ( R^n_{N , F}\, ,\, U^n_F )\, $ are surjective, and the maps 
$\, K^J_2 ( R^n_{N , F}\, ,\, U^n_F ) \rightarrowtail K^J_2 ( J^n_{N , F}\, ,\, U^n_F )\, $ are injective. In particular the groups $\,\{ K_n ( N , F ) \}\, $ are well defined for any normal subgroup 
$\, N \subseteq F\, $ and behave functorial under morphisms $\, ( N , F ) \rightarrow ( N' , F' )\, $. They can be identified with the images of $\, K^J_2 ( R^{n-2}_{0 , N , F}\, ,\, U^{n-2}_F )\, $ in 
$\, K^J_2 ( J^{n-2}_{N , F}\, ,\, U^{n-2}_F )\, $ and also with the abelian regular part of 
$\, K^J_2 ( J^{n-2}_{N , F}\, ,\, U^{n-2}_F )\, $.
\par\bigskip\noindent
{\it Proof.}\quad Let $\, N\, $ be a normal subgroup of $\, F\, $ and $\, N_r = [ N , F ]\, $. By definition 
$\, K_n ( N , F ) = K_n ( N_r , F )\, $ so it suffices to assume $\, N = N_r\, $. Then $\, K_n ( N , F )\, $ is defined in terms of an arbitrary regular chain of $\, ( N , F )\, $ and we must show that the different choices are related by some natural isomorphisms of the $K_*$-groups. First of all since 
$\, K^J_2 ( R_{N , F}\, ,\, U_F ) = K_3 ( N , F )\, $ determines the abelian regular part of 
$\, K^J_2 ( J_{N , F}\, ,\, U_F )\, $ and injectivity of the map 
$\, K^J_2 ( R_{N , F}\, ,\, U_F ) \rightarrowtail K^J_2 ( J_{N , F}\, ,\, U_F )\, $ which follows from Lemma 20 and Corollary 5.1, the former group is well defined independent of the choice of regular suspension 
$\, ( R_{N , F}\, ,\, U_F )\, $. Put $\, D = R_{N , F} / R_{0 , N , F} \simeq J_{N , F} / J_{0 , N , F}\, $ so that 
$\, D\, $ is a free abelian central direct summand of 
$\, U_F / R_{0 , N , F} \simeq D \times ( {\overline N}^{ _F} \rtimes U_{F / N} )\, $. The image of 
$\, K^J_2 ( R_{N , F}\, ,\, U_F )\, $ in $\, K^J_2 ( D\, ,\, U_F / R_{0 , N , F} )\, $ must be trivial because the latter does not contain any nontrivial regular elements. It follows that 
$\, K^J_2 ( R_{0 , N , F}\, ,\, U_F ) \twoheadrightarrow K^J_2 ( R_{N , F}\, ,\, U_F )\, $ is surjective (and by a similar argument $\, K^J_2 ( R_{0 , R^{n-1}_{N , F}\, ,\, U^{n-1}_F}\, ,\, U^n_F ) \twoheadrightarrow 
K^J_2 ( R^n_{N , F}\, ,\, U^n_F )\, $ is surjective). Then any suspension $\, ( R^2_{N , F}\, ,\, U^2_F )\, $ may be restricted to $\, ( R_{R_{0 , N , F}\, ,\, U_F}\, ,\, U^2_F )\, $ which is the kernel of an almost canonical map 
\smallskip
$$ U_{R_{0 , N , F}\, ,\, U_F} \longrightarrow {\overline R_{0 , N , F}} \subseteq {\overline R}^{ _F}_{N , F} 
$$
\par\medskip\noindent
and $\, \overline R_{0 , N , F}\, $ denotes the preimage of $\, R_{0 , N , F}\, $ in the (relative) universal $U_F$-central extension $\, {\overline R}^{ _F}_{N , F}\, $. One gets injective maps
\smallskip 
$$ K^J_2 ( R_{R_{0 , N , F}\, ,\, U_F}\, ,\, U^2_F ) \rightarrowtail 
K^J_2 ( J_{R_{0 , N , F}\, ,\, U_F}\, ,\, U^2_F ) \buildrel\sim\over\rightarrow 
K^J_2 ( J_{J_{0 , N , F}\, ,\, U_F}\, ,\, U^2_F ) $$ 
and 
\smallskip
$$ K^J_2 ( R^2_{N , F}\, ,\, U^2_F ) \rightarrowtail 
K^J_2 ( J_{R_{N , F}\, ,\, U_F}\, ,\, U^2_F ) \buildrel\sim\over\rightarrow K^J_2 ( J^2_{N , F}\, ,\, U^2_F ) $$ \par\medskip\noindent
from Corollary 5.1. Consider the commutative diagram 
\smallskip
$$ \vbox{\halign{ #&#&#&#&#\cr
\hfil $K^J_2 ( R_{R_{0 , N , F}\, ,\, U_F}\, ,\, U^2_F )$\hfil & \hfil $\rightarrowtail $\hfil & 
\hfil $K^J_2 ( R^2_{N , F}\, ,\, U^2_F )$\hfil & \hfil $\rightarrow $\hfil & 
\hfil $K^J_2 ( J_{D , B}\, ,\, {\overline R_{0 , N , F}} \rtimes U_B )$\hfil \cr
\hfil $\Bigm\downarrow $\hfil && \hfil $\Bigm\downarrow $\hfil && 
\hfil $\Bigm\downarrow $\hfil \cr
\hfil $K^J_2 ( J_{R_{0 , N , F}\, ,\, U_F}\, ,\, U^2_F )$ \hfil & \hfil $\rightarrowtail $\hfil & 
\hfil $K^J_2 ( J_{R_{N , F}\, ,\, U_F}\, ,\, U^2_F )$ \hfil & \hfil $\rightarrow $\hfil & 
\hfil $K^J_2 ( J_{D , B}\, ,\, R_{0 , N , F} \rtimes U_B )$ \hfil \cr
\hfil $\Bigm\downarrow $\hfil && \hfil $\Bigm\downarrow $\hfil &&\cr
\hfil $K^J_2 ( C , {\overline R_{0 , N , F}} \rtimes U_B )$\hfil & \hfil $\twoheadrightarrow $ \hfil & 
\hfil $K^J_2 ( C , {\overline R_{N , F}} \rtimes U_{B / D} )$\hfil &&\cr }} $$
\par\medskip\noindent
with $\, B = U_F / R_{0 , N , F}\, $. We claim that the image of $\, K^J_2 ( R^2_{N , F}\, ,\, U^2_F )\, $ in 
$\, K^J_2 ( J_{D , B}\, ,\, R_{0 , N , F} \rtimes U_B )\, $ is trivial. By full exactness of 
$\, ( R_{N , F}\, ,\, U_F )\, $ and $\, ( R_{0 , N , F}\, ,\, U_F )\, $ the image of 
$\, K^J_2 ( J_{R_{N , F}\, ,\, U_F}\, ,\, U^2_F )\, $ does not intersect the kernel of 
$\, K^J_2 ( J_{D , B}\, ,\, R_{0 , N , F} \rtimes U_B ) \twoheadrightarrow K^J_2 ( J_{D , B}\, ,\, U_B )\, $ from Corollary 5.1, so we may replace 
$\, K^J_2 ( J_{D , B}\, ,\, R_{0 , N , F} \rtimes U_B )\, $ by $\, K^J_3 ( D , B ) = 
K^J_2 ( J_{D , B}\, ,\, U_B )\, $. This group consists of the split injective image of $\, K^J_3 ( D )\, $ and a purely relative part both of which are free abelian groups (compare with section 1, Lemmas 5 and 6). From the extension 
\smallskip
$$ 1 \longrightarrow \bigl( J_{0 , N , F}\, ,\, J_{N , F} \bigr) \longrightarrow \bigl( J_{N , F}\, ,\, J_{N , F} \bigr) \longrightarrow \bigl( D\, ,\, D \bigr) \longrightarrow 1 $$
\par\medskip\noindent
one finds that $\, K^J_3 ( D ) \simeq K^J_2 ( J_{0 , N , F}\, ,\, J_{N , F} )\, $ which injects into the group 
$\, K^J_2 ( J_{0 , N , F}\, ,\, U_F ) \simeq K^J_3 ( U_F / J_{0 , N , F} )\, $ by splitting of the corresponding map $\, D \subseteq U_F / J_{0 , N , F}\, $, so that the image of 
$\, K^J_2 ( J_{R_{N , F}\, ,\, U_F}\, ,\, U^2_F )\, $ does not intersect with 
$\, K^J_3 ( D )\, $. Let $\, J_0\, $ be the normalization of the $B_1$-(and $B_3'$)-subcore of 
$\, J_{D , B}\, $. Then the relative part of $\, K^J_2  J_{D , B}\, ,\, U_B )\, $ passes to 
$\, K^J_2 ( J_{D , B} / J_0\, ,\, U_B / J_0 )\, $ and in fact the proof of Lemma 6 gives that any such element lifts to $\, K^J_2 ( [\, U_{D , B}\, ,\, U_B\, ] / J_0\, ,\, U_B / J_0 )\, $. Consider the extension
\smallskip
$$ 1 \longrightarrow \bigl( {J_0 + [\, U_{D , B}\, ,\, U_B\, ] \over J_0}\, ,\, {U_B \over J_0} \bigr) \longrightarrow 
\bigl( {J_0 + [\, U_B\, ,\, U_B\, ] \over J_0}\, ,\, {U_B \over J_0} \bigr) \rightarrow\qquad\qquad $$ 
$$ \qquad\qquad\qquad\qquad\qquad\qquad\qquad 
\rightarrow \bigl( [\, U_{B / D}\, ,\, U_{B / D}\, ]\, ,\, A \times U_{B / D} \bigr) \longrightarrow 1 $$
\par\medskip\noindent
with $\, A = U_{D , B} / ( J_0 + [\, U_{D , B}\, ,\, U_B\, ] )\, $ a central direct summand in the quotient. Restricting to the preimage $\, \widetilde U_B\, $ of the factor $\, U_{B / D}\, $ one gets a split extension and hence an injective map 
\smallskip
$$ K^J_2 ( {[\, U_{D , B}\, ,\, U_B\, ] \over J_0}\, ,\, \widetilde U_B ) \rightarrowtail 
K^J_2 ( {[\, U_B\, ,\, U_B\, ] \over J_0}\, ,\, \widetilde U_B ) \> . $$
\par\medskip\noindent
It is clear that the image of the relative part lifts to 
$\, K^J_2 ( {[\, U_{D , B}\, ,\, U_B\, ] \over J_0}\, ,\, \widetilde U_B )\, $ so it injects into 
$\, K^J_2 ( [\, U_B\, ,\, U_B\, ] / J_0\, ,\, \widetilde U_B )\, $, respectively into the group  
$\, K^J_2 ( ( J_{D , B} + [\, U_B\, ,\, U_B\, ] ) / J_0\, ,\, \widetilde {\widetilde U_B} )\, $ where 
$\, \widetilde {\widetilde U_B} = \widetilde U_B + ( J_{D , B} / J_0 )\, $ by the same argument. The Special Deconstruction Theorem gives a surjection 
$\, K^J_2 ( R^2_{N , F}\, ,\, R^2_{N , F} + [\, U^2_F\, ,\, U^2_F\, ] ) \twoheadrightarrow 
K^J_2 ( R^2_{N , F}\, ,\, U^2_F )\, $. Now the map 
$\, K^J_2 ( R^2_{N , F}\, ,\, R^2_{N , F} + [\, U^2_F\, ,\, U^2_F\, ] ) \rightarrow 
K^J_2 ( ( J_{D , B} + [\, U_B\, ,\, U_B\, ] ) / J_0\, ,\, \widetilde {\widetilde U_B} )\, $ factors over 
$\, K^J_2 ( R^2_{N , F} + [\, U^2_F\, ,\, U^2_F\, ] ) = 0\, $ leading to the conclusion that the image of 
$\, K^J_2 ( R^2_{N , F}\, ,\, U^2_F )\, $ in $\, K^J_2 ( J_{D , B}\, ,\, R_{0 , N , F} \rtimes U_B )\, $ is trivial. 
Then any element in the cokernel of 
$\, K^J_2 ( R_{R_{0 , N , F}\, ,\, U_F}\, ,\, U^2_F ) \rightarrow K^J_2 ( R^2_{N , F}\, ,\, U^2_F )\, $ maps to the kernel of $\, K^J_2 ( J_{D , B}\, ,\, {\overline R_{0 , N , F}} \rtimes U_B ) \twoheadrightarrow 
K^J_2 ( J_{D , B}\, ,\, R_{0 , N , F} \rtimes U_B )\, $. Also by surjectivity of 
$\, K^J_2 ( R_{0 , R_{N , F}\, ,\, U_F}\, ,\, U^2_F ) \twoheadrightarrow K^J_2 ( R^2_{N , F}\, ,\, U^2_F )\, $ the map from $\, K^J_2 ( R^2_{N , F}\, ,\, U^2_F )\, $ to $\, K^J_2 ( J_{D , B}\, ,\, {\overline R_{0 , N , F}} \rtimes U_B )\, $ factors over the image of 
$\, K^J_2 ( R_{0 , D , B}\, ,\, {\overline R_{0 , N , F}} \rtimes U_B )\, $. Suppose that
\smallskip 
$$ K^J_2 ( R_{R_{0 , N , F}\, ,\, U_F}\, ,\, U^2_F ) \longrightarrow K^J_2 ( R^2_{N , F}\, ,\, U^2_F ) $$ 
\par\medskip\noindent
is not surjective. Then one finds an element of $\, K^J_2 ( R^2_{N , F}\, ,\, U^2_F )\, $ which considered as an element of $\, K^J_2 ( J_{R_{0 , N , F}\, ,\, U_F}\, ,\, U^2_F )\, $ maps to the kernel of 
$\, K^J_2 ( C\, ,\, {\overline R_{0 , N , F}\, ,\, U_F}\, ,\, U_B ) \twoheadrightarrow 
K^J_2 ( C\, ,\, ( {\overline R_{0 , N , F}} \rtimes U_B ) / R_{0 , D , B} )\, $. But 
$\, R_{0 , D , B} \subseteq [\, U_B\, ,\, U_B\, ]\, $, so this kernel is trivial. Then 
\smallskip
$$ K^J_2 ( R_{R_{0 , N , F}\, ,\, U_F}\, ,\, U^2_F ) \buildrel\sim\over\longrightarrow 
K^J_2 ( R^2_{N , F}\, ,\, U^2_F ) $$
\par\medskip\noindent
is an isomorphism and one may continue by restricting $\, ( R^2_{R_{0 , N , F} , U_F} , U^3_F )\, $ from 
$\, ( R^3_{N , F}\, ,\, U^3_F )\, $ respectively $\, ( R^n_{R_{0 , N , F}\, ,\, U_F}\, ,\, U^{n+1}_F )\, $ from 
$\, ( R^{n+1}_{N , F}\, ,\, U^{n+1}_F )\, $ by induction, showing that the induced map on $K^J_2$-groups is an isomorphism for each $\, n\, $. If $\, R'_{N , F}\, $ is any other suspension of $\, ( N , F )\, $ the 
chain $\,\{ ( R^n_{R_{0 , N , F}\, ,\, U_F}\, ,\, U^{n+1}_F ) \}\, $ can be induced to a regular chain 
$\,\{ ( R^n_{R'_{N , F} , U_F}\, ,\, U^{n+1}_F ) \}\, $ beginning with $\, ( R'_{N , F} , U_F )\, $ so that 
$\, K^J_2 ( R^n_{N , F}\, ,\, U^n_F )\, $ is isomorphic to 
$\, K^J_2 ( R^{n-1}_{R'_{N , F}\, ,\, U_F}\, ,\, U^n_F )\, $ by symmetry of the above argument. Then one has reduced the problem by one dimension to the question whether all possible continuations of a given suspension $\, ( R_{N , F}\, ,\, U_F )\, $ to a regular chain $\,\{ ( R^n_{N , F}\, ,\, U^n_F ) \}\, $ of 
$\, ( N , F )\, $ have (naturally) isomorphic $K^J_2$-groups. Inductively, assume this is the case for all 
$\, k \leq n-1\, $. Put $\, D_n = R^n_{N , F} / R_{0 , R^{n-1}_{N , F} , U^{n-1}_F}
\simeq J_{R^{n-1}_{N , F}\, ,\, U^{n-1}_F} / J_{0 , R^{n-1}_{N , F}\, ,\, U^{n-1}_F}\, $, 
$\, C_n = K^J_2 ( R^{n-1}_{N , F}\, ,\, U^{n-1}_F )\, $ and 
$\, B_n = U^n_F / R_{0 , R^{n-1}_{N , F}\, ,\, U^{n-1}_F}\, $. Then $\, D_n\, $ is still a free abelian central direct summand of $\, B_n\, $. To get a projection onto $\, D_n\, $ divide $\, U^n_F\, $ by 
$\, J_{0 , J_{R^{n-2}_{N , F}\, ,\, U^{n-2}_F}\, ,\, U^{n-1}_F}\, $ giving a semidirect product of an extension of $\, U_{C_n\, ,\, B_n}\, /\, [\, U_{C_n\, ,\, B_n}\, ,\, U_{B_n}\, ]\, $ by $\, D_n \times R^{n-1}_{N , F}\, $ with 
$\, U_{U^{n-1}_F / J_{R^{n-2}_{N , F}\, ,\, U^{n-2}_F}}\, $. Dividing out the normal copy of 
$\, R^{n-1}_{N , F}\, $ one obtains a central extension of 
$\, U_{U^{n-1}_F / J_{R^{n-2}_{N , F}\, ,\, U^{n-2}_F}}\, $ by 
$\, D_n \times ( U_{C_n\, ,\, B_n} / [\, U_{C_n\, ,\, B_n}\, ,\, U_{B_n}\, ] )\, $. So all the above arguments apply in this situation to show that any two suspensions of $\, ( R^{n-1}_{N , F}\, ,\, U^{n-1}_F )\, $ have continuations with isomorphic $K^J_2$-groups, hence the $K^J_2$-group of any given n-fold suspension is isomorphic to the $K^J_2$-group of any other n-fold suspension. But this is not enough, since there may be different possibilities to identify them in the manner above and it is not a priori clear that these yield the same result. But we need the identification to be natural. This will be the case if they can all be identified with the image of $\, K^J_2 ( R^n_{0 , N , F}\, ,\, U^n_F )\, $ in 
$\, K^J_2 ( J^n_{N , F}\, ,\, U^n_F )\, $, i.e. if the map 
$\, K^J_2 ( R^n_{0 , N , F}\, ,\, U^n_F ) \rightarrow K^J_2 ( R^n_{N , F}\, ,\, U^n_F )\, $ is surjective for each $\, n\, $. By induction we may assume this to be the case for all values $\, k < n\, $. Then one can apply the constructions in the proof of Theorem 2 to show that 
$\, ( R^k_{0 , N , F}\, ,\, U^k_F )\, $ admits a (nonfaithful) relative semiuniversal $U^k_F$-central extension with respect to the inclusion into $\, ( J_{R^{k-1}_{0 , N , F}\, ,\, U^{k-1}_F}\, ,\, U^k_F )\, $, having the (relative) splitting properties analogous to the case of $\, ( R^k_{N , F}\, ,\, U^k_F )\, $, so that the image of $\, K^J_2 ( R^{n-k}_{R^k_{0 , N , F}\, ,\, U^k_F}\, ,\, U^n_F )\, $ in 
$\, K^J_2 ( J_{R^{n-k-1}_{R^k_{0 , N , F}\, ,\, U^k_F}\, ,\, U^{n-1}_F}\, ,\, U^n_F )\, $ consists of abelian regular elements, and by an argument as above the map 
\smallskip
$$ K^J_2 ( R^n_{0 , N , F}\, ,\, U^n_F ) \twoheadrightarrow 
K^J_2 ( R_{R^{n-1}_{0 , N , F}\, ,\, U^{n-1}_F}\, ,\, U^n_F ) \simeq K^J_2 ( R^n_{N , F}\, ,\, U^n_F ) $$
\par\medskip\noindent
is surjective and the latter isomorphism follows as above by induction. Since the assignments 
\smallskip
$$ \bigl( N , F \bigr) \rightsquigarrow \bigl( R^n_{0 , N , F}\, ,\, U^n_F \bigr)\quad ,\quad 
\bigl( N , F \bigr) \rightsquigarrow \bigl( J^n_{N , F}\, ,\, U^n_F \bigr) $$
\par\medskip\noindent
are functorial, so is the assignment $\, ( N , F ) \rightsquigarrow K_n ( N , F )\, $ for every $\, n \geq 2\, $.
\par\noindent
There is yet another way to describe the regular $K_*$-groups. One knows that the image of 
$\, K^J_2 ( R^n_{N , F}\, ,\, U^n_F )\, $ in $\, K^J_2 ( J^n_{N , F}\, ,\, U^n_F )\, $ consists of abelian regular elements. We claim that it also exhausts the subset of abelian regular elements. Let 
$\, R_{J^{n-1}_{N , F}\, ,\, U^{n-1}_F}\, $ denote the kernel of an almost canonical map 
\smallskip
$$ U_{J^{n-1}_{N , F}\, ,\, U^{n-1}_F} \longrightarrow {\overline J^{n-1}_{N , F}}^{ _{U^{n-1}_F}} $$
\par\medskip\noindent
induced from $\, R^n_{N , F}\, $. Consider the diagram 
\smallskip
$$ \vbox{\halign{ #&#&#&#&#&#&#&#&#\cr
\hfil $1$\hfil & \hfil $\longrightarrow $\hfil & \hfil $\bigl( R^n_{N , F}\, ,\, U^n_F \bigr)$\hfil & 
\hfil $\longrightarrow $\hfil & \hfil $\bigl( J_{R^{n-1}_{N , F}\, ,\, U^{n-1}_F}\, ,\, U^n_F \bigr)$\hfil & 
\hfil $\longrightarrow $\hfil & \hfil $\bigl( C\, ,\, {U^n_F \over R^n_{N , F}} \bigr)$\hfil & \hfil $\longrightarrow $\hfil & 
\hfil $1$\hfil \cr
&& \hfil $\Bigm\downarrow $\hfil && \hfil $\Bigm\downarrow $\hfil && \hfil $\Bigm\downarrow $\hfil &&\cr
\hfil $1$\hfil & \hfil $\longrightarrow $\hfil & \hfil $\bigl( R_{J^{n-1}_{N , F}\, ,\, U^{n-1}_F}\, ,\, U^n_F \bigr)$\hfil & \hfil $\longrightarrow $\hfil & \hfil $\bigl( J^n_{N , F}\, ,\, U^n_F \bigr)$\hfil & 
\hfil $\longrightarrow $\hfil & \hfil $\bigl( C'\, ,\, {U^n_F \over R_{J^{n-1}_{N , F}\, ,\, U^{n-1}_F}} \bigr)$\hfil & 
\hfil $\longrightarrow $\hfil & \hfil $1$\hfil \cr
&& \hfil $\Bigm\downarrow $\hfil && \hfil $\Bigm\downarrow $\hfil &&&&\cr
&& \hfil $\bigl( \widetilde R\, ,\, {U^n_F \over R^n_{N , F}} \bigr)$\hfil & \hfil $\longrightarrow $\hfil & 
\hfil $\bigl( \widetilde J\, ,\, {U^n_F \over J_{R^{n-1}_{N , F}\, ,\, U^{n-1}_F}} \bigr)$\hfil &&&&\cr }} \> . $$
\par\medskip\noindent
with $\, \widetilde R = R_{J^{n-1}_{N , F}\, ,\, U^{n-1}_F} / R^n_{N , F}\, $ and    
$\, \widetilde J = J_{J^{n-1}_{N , F} / R^{n-1}_{N , F}\, ,\, U^{n-1}_F / R^{n-1}_{N , F}}\, $. If the map 
$\, K^J_2 ( R^n_{N , F} \rightarrowtail K^J_2 ( R_{J^{n-1}_{N , F}\, ,\, U^{n-1}_F}\, ,\, U^n_F )\, $ wouldn't be an isomorphism, there would exist an element of 
$\, K^J_2 ( R_{J^{n-1}_{N , F}\, ,\, U^{n-1}_F}\, ,\, U^n_F )\, $ mapping to a nontrivial element of 
$\, K^J_2 ( J_{D , B}\, ,\, U^n_F / R^n_{N , F} )\, $ in the kernel of the map to 
$\, K^J_2 ( J_{D , B}\, ,\, U^n_F / J_{R^{n-1}_{N , F}\, ,\, U^{n-1}_F} )\, $ where 
$\, D = J^{n-1}_{N , F} / R^{n-1}_{N , F}\, $ and $\, B = U^{n-1}_F / R^{n-1}_{N , F}\, $. This element is still nontrivial in $\, K^J_2 ( J_{D , B} \times C\, ,\, U^n_F / R^n_{N , F} )\, $ as 
$\, J_{D , B} \times C\, $ is the image of $\, J^n_{N , F}\, $ modulo $\, R^n_{N , F}\, $, and is represented by a commutator in $\, [\, U_C\, ,\, U_{J_{D , B}}\, ]\, $, hence it can be lifted to an element of the canonical subgroup $\, [\, J_{R^{n-1}_{N , F}\, ,\, U^{n-1}_F}\, ,\, J^n_{N , F}\, ] \subseteq 
[\, J^n_{N , F}\, ,\, J^n_{N , F}\, ] \subseteq {\overline J^n_{N , F}}^{ _{U^n_F}}\, $ which (the first bracket) by the proof of Theorem 2 is contained in a normal copy of $\, {\overline R^n_{N , F}}^{ _{U^n_F}}\, $ inside 
$\, {\overline J^n_{N , F}}^{ _{U^n_F}}\, $. Since its image in $\, K_1 ( R^n_{N , F}\, ,\, U^n_F )\, $ is trivial, it is even contained in $\, [\, {\overline R^n_{N , F}}^{ _{U^n_F}}\, ,\, U^n_F\, ]\, $, from which follows by exactness that the image of the corresponding element of 
$\, K^j_2 ( J_{R^{n-1}_{N , F}\, ,\, U^{n-1}_F}\, ,\, U^n_F )\, $ in $\, K^J_2 ( C\, ,\, U^n_F / R^n_{N , F} )\, $, hence in $\, K^J_2 ( J_{D , B} \times C\, ,\, U^n_F / R^n_{N , F} )\, $ is trivial, giving a contradiction. Then the map 
\smallskip
$$ K^J_2 ( R^n_{N , F}\, ,\, U^n_F ) \buildrel\sim\over\longrightarrow 
K^J_2 ( R_{J^{n-1}_{N , F}\, ,\, U^{n-1}_F}\, ,\, U^n_F )  $$
\par\medskip\noindent
is an isomorphism. On the other hand any abelian regular element of 
$\, K^J_2 ( J^n_{N , F}\, ,\, U^n_F )\, $ lifts to $\, K^J_2 ( R_{J^{n-1}_{N , F}\, ,\, U^{n-1}_F}\, ,\, U^n_F )\, $\qed
\par\bigskip\noindent
The next two Lemmas are concerned with exactness. The result below shows that a regular pair is always fully exact, so that Corollary 5.1 applies to arbitrary extensions of regular pairs.
\par\bigskip\noindent  
{\bf Lemma 21.}\quad Assume given an extension 
\smallskip
$$ 1 \longrightarrow ( N , F ) \longrightarrow ( M , F ) \longrightarrow ( P , E ) \longrightarrow 1 $$
\par\medskip\noindent
such that $\, N \subseteq [ M , F ]\, $ and $\, ( P , E )\, $ is exact. Then $\, ( M , F )\, $ is exact. In particular this is the case if $\, M \subseteq [ F , F ]\, $. If $\, N\, $ is a full subgroup of a perfect group 
$\, F\, $ then $\, ( N , F )\, $ is fully exact.
\par\bigskip\noindent
{\it Proof.}\quad Put $\, \overline F = U_F / J_{M , F} \simeq M \rtimes U_{E / P}\, $. By exactness of 
$\, ( P , E )\, $ any element in the kernel of 
\smallskip
$$ K^J_2 ( J_{E / P}\, ,\, M \rtimes U_{E / P} ) \longrightarrow 
K^J_2 ( J_{E / P} \times M\, ,\, M \rtimes U_{E / P} ) $$
\par\medskip\noindent
is in the kernel of the regular surjection 
\smallskip
$$ K^J_2 ( J_{E / P}\, ,\, M \rtimes U_{E / P} ) \twoheadrightarrow K^J_2 ( J_{E / P}\, ,\, P \rtimes U_{E / P} ) 
$$
\par\medskip\noindent
hence representable by a commutator expression in $\, [\, U_{J_{E / P}}\, ,\, U_N\, ]\, $. Such an element is nontrivial if and only if the corresponding element in 
$\, K^J_2 ( J_{E / P}\, ,\, M^{ab} \rtimes U_{E / P} )\, $ is nontrivial, dividing $\, \overline F\, $ by 
$\, [ M , M ]\, $. So without loss of generality we may assume that $\, M\, $ is abelian, so in particular there exists a well defined action of $\, E / P\, $ on $\, M\, $ inducing the action of $\, U_{E / P}\, $ in 
$\, \overline F = M \rtimes U_{E / P}\, $ and the diagonal copy $\, J^{\Delta }_{E / P} \subseteq U_{E / P} \subseteq \overline F\, $ is normal with $\, \overline F / J^{\Delta }_{E / P} \simeq M \rtimes E / P\, $. We first show that $\, ( M , F )\, $ is exact if $\, ( M , M \rtimes E / P )\, $ is exact (the converse does not hold in general !). More precisely we show that an element in the trefoil intersection of the groups 
$\, K^J_2 ( J_{E / P}\, ,\, ( J_{E / P} \times M ) \rtimes U_{E / P} )\, ,\, 
K^J_2 ( M\, ,\, ( J_{E / P} \times M ) \rtimes U_{E / P} )\, $ and 
$\, K^J_2 ( J^0_{E / P}\, ,\, ( J_{E / P} \times M ) \rtimes U_{E / P} )\, $ in 
$\, K^J_2 ( J_{E / P} \times M \times J^0_{E / P}\, ,\, ( J_{E / P} \times M ) \rtimes U_{E / P} )\, $ defines in a natural way an element in the corresponding trefoil intersection with $\, J_{E / P}\, $ replaced by 
$\, J^{\Delta }_{E / P}\, $, where the latter is a normal copy of $\, J_{E / P}\, $ in 
$\, ( J_{E / P} \times M ) \rtimes U_{E / P}\, $ which maps to the diagonal copy 
$\, J^{\Delta }_{E / P} \subseteq \overline F\, $ modulo $\, J^0_{E / P}\, $. Consider an element in the former trefoil intersection which we would like to lift to the corresponding trefoil intersection of extended 
$K^J_2$-groups. We may assume this element to be represented by a commutator expression in 
$\, v \in [\, U_{J_{E / P}}\, ,\, U_M\, ]\, $. Consider a bijection of 
$\, \widetilde F = U_{\overline F} / J_{\overline M , \overline F} \simeq 
( J_{E / P} \times M ) \rtimes U_{E / P}\, $ to itself which sends the subgroup 
$\, J_{E / P}\, $ (equivarianty) to $\, J^{\Delta }_{E / P}\, $, also sending $\, M\, $ and $\, J^0_{E / P}\, $ identically to themselves and leaving invariant some given section $\, s ( E / P ) \subseteq U_{E / P}\, $. This map induces compatible homomorphisms 
\smallskip
$$ \bigl(\> U_{J_{E / P}\, ,\, \widetilde F}\> ,\> 
U_{M , \widetilde F}\> ,\> 
U_{J^0_{E / P}\, ,\, \widetilde F}\> ,\> 
U_{J_{E / P} \times M \times J^0_{E / P}\, ,\, \widetilde F}\> ,\> 
U_{\widetilde F}\> \bigr)\quad \buildrel\sigma\over\largerightarrow\qquad\qquad $$
$$ \qquad\qquad\qquad \bigl(\> U_{J^{\Delta }_{E / P}\, ,\, \widetilde F}\> ,\> 
U_{M , \widetilde F}\> ,\> 
U_{J^0_{E / P}\, ,\, \widetilde F}\> ,\> 
U_{J_{E / P} \times M \times J^0_{E / P}\, ,\, \widetilde F}\> ,\> 
U_{\widetilde F}\> \bigr) $$
\par\medskip\noindent
inducing
\smallskip
$$ \bigl(\> J_{J_{E / P}\, ,\, \widetilde F}\> ,\> 
J_{M , \widetilde F}\> ,\> 
J_{J^0_{E / P}\, ,\, \widetilde F}\> ,\> 
J_{J_{E / P} \times M \times J^0_{E / P}\, ,\, \widetilde F}\> \bigr) 
\quad\buildrel\sigma\over\largerightarrow\qquad\qquad\qquad $$ 
$$\qquad\qquad\qquad\qquad \bigl(\> J_{J^{\Delta }_{E / P}\, ,\, \widetilde F}\> ,\> 
J_{M , \widetilde F}\> ,\> 
J_{J^0_{E / P}\, ,\, \widetilde F}\> ,\> 
J_{J_{E / P} \times M \times J^0_{E / P}\, ,\, \widetilde F}\> \bigr) $$
\par\medskip\noindent
and $\, J_{\widetilde F}\, \longrightarrow\, 
J_{\widetilde F} + U_{M , \widetilde F}\, $ and in general differing from the identity of 
$\, U_{\widetilde F}\, $ only modulo 
$\, U_{M , \widetilde F}\, $. It sends the representative $\, v\, $ to 
$\, v^{\Delta } \in [\, U_{J^{\Delta }_{E / P}}\, ,\, U_M\, ]\, $ which is therefore in the common intersection of 
the groups $\, K^J_2 ( J^{\Delta }_{E / P}\, ,\, ( J_{E / P} \times M ) \rtimes U_{E / P} )\, $ and 
$\, K^J_2 ( M\, ,\, ( J_{E / P} \times M ) \rtimes U_{E / P} )\, $ (in 
$\, K^J_2 ( J_{E / P} \times M \times J^0_{E / P}\, ,\, ( J_{E / P} \times M ) \rtimes U_{E / P} )\, $) and differs from the original element only by an element represented by a commutator in 
$\, [\, U_M\, ,\, U_M\, ]\, $, thus in the image of $\, K^J_2 ( M )\, $. This element must become trivial in 
$\, K^J_2 ( M , M \rtimes E / P )\, $ hence in $\, K^J_2 ( M \rtimes E / P )\, $ (since the map 
$\, K^J_2 ( M ) \rightarrow K^J_2 ( M , M \rtimes E / P )\, $ factors over 
$\, K^J_2 ( M , ( J_{E / P} \times M ) \rtimes U_{E / P} )\, $ and $\, v\, $ is in the image of 
$\, K^J_2 ( J^0_{E / P}\, ,\, ( J_{E / P} \times M ) \rtimes U_{E / P} )\, $ while $\, v^{\Delta }\, $ is in the image of $\, K^J_2 ( J^{\Delta }_{E / P}\, ,\, ( J_{E / P} \times M ) \rtimes U_{E / P} )\, $). But the map 
$\, K^J_2 ( M ) \rightarrowtail K^J_2 ( M \rtimes E / P )\, $ is injective from the Excision Theorem. Then 
$\, v \equiv v^{\Delta }\, $ in 
$\, K^J_2 ( J_{E / P} \times M \times J^0_{E / P}\, ,\, ( J_{E / P} \times M ) \rtimes U_{E / P} )\, $, so that any element in the trefoil intersection corresponding to $\, J_{E / P}\, ,\, M\, $ and $\, J^0_{E / P}\, $ is automatically  in the image of $\, K^J_2 ( J^{\Delta }_{E / P}\, ,\, ( J_{E / P} \times M ) \rtimes U_{E / P} )\, $. Assuming that $\, v \equiv v^{\Delta }\, $ has a lift to the trefoil intersection of the extended $K^J_2$-groups for $\, J^{\Delta }_{E / P}\, ,\, M\, $ and $\, J^0_{E / P}\, $ the inverse of the equivariant bijection as above defines compatible maps 
\smallskip
$$ \widetilde K^J_ ( J^{\Delta }_{E / P}\, ,\, \widetilde F ) \longrightarrow 
\widetilde K^J_2 ( J_{E / P}\, ,\, \widetilde F )\> , $$
$$ \widetilde K^J_2 ( M\, ,\, \widetilde F ) \longrightarrow \widetilde K^J_2 ( M\, ,\, \widetilde F )\> , $$
$$ \widetilde K^J_2 ( J^0_{E / P}\, ,\, \widetilde F ) \longrightarrow 
\widetilde K^J_2 ( J^0_{E / P}\, ,\, \widetilde F )\> , $$
$$ \widetilde K^J_2 ( J_{E / P} \times M \times J^0_{E / P}\, ,\, \widetilde F ) \longrightarrow 
\widetilde K^J_2 ( J_{E / P} \times M \times J^0_{E / P}\, ,\, \widetilde F ) $$
\par\medskip\noindent
giving a lift of $\, [ v ]\, $ into the trefoil intersection of the extended $K^J_2$-groups on the right hand side. Returning to the setting given by the Lemma suppose given an element in the kernel of the map 
\smallskip
$$ K^J_2 ( J_{E / P} \times J^0_{E / P}\, ,\, \widetilde F ) \longrightarrow 
K^J_2 ( J_{E / P} \times M \times J^0_{E / P}\, ,\, \widetilde F ) $$
\par\medskip\noindent
represented by a commutator expression in
\smallskip 
$$ \bigl[\, U_{J_{E / P} \times J^0_{E / P}}\, ,\, U_N\, \bigr]\, \equiv\, 
\bigl[\, U_{J_{E / P}}\, ,\, U_N\, \bigr]\> \bigl[\, U_{J^0_{E / P}}\, ,\, U_N\, \bigr] $$ 
\par\medskip\noindent
and write it in the form 
$\, z = v \cdot w\, $ with $\, v \in [\, U_{J_{E / P}}\, ,\, U_N\, ]\, ,\, w \in [\, U_{J^0_{E / P}}\, ,\, U_N\, ]\, $. Then $\, v\, $ defines an element in the trefoil intersection of $\, J_{E / P}\, ,\, N\, $ and $\, J^0_{E / P}\, $ in  
$\, K^J_2 ( J_{E / P} \times M \times J^0_{E / P}\, ,\, ( J_{E / P} \times M ) \rtimes U_{E / P} )\, $ which we would like to lift to the corresponding trefoil intersection replacing $K^J_2$ by the extended $K^J_2$-functor. Without loss of generality we may assume that $\, F = M \rtimes E / P\, $. Let 
\smallskip
$$ 1 \longrightarrow R \longrightarrow U \longrightarrow E / P \longrightarrow 1 $$
\par\medskip\noindent
be a free resolution of $\, E / P\, $ and consider the extension $\, H = M \rtimes U\, $ of $\, F\, $ (where 
$\, R\, $ is supposed to commute with $\, M\, $). Then $\, ( M , H )\, $ is exact, so the lifting problem can be reduced to the problem of finding a lift of $\, [ v ]\, $ in the trefoil intersection of the subgroups 
$\, J_U\, ,\, M\, $ and 
$\, J^0_U\, $ in $\, K^J_2 ( J_U \times M \times J^0_U\, ,\, ( J_U \times M ) \rtimes U_U )\, $ which is trivial from Lemma 12 by the fact that $\, J_U\, $ (or $\, J^0_U\, $) has a core for the quotient $\, U\, $. In fact by the same reason already the corresponding trefoil intersection in 
$\, K^J_2 ( J_{E / P} \times M \times J^0_U\, ,\, ( J_{E / P} \times M ) \rtimes U_U )\, $ is trivial, so it suffices to find a lift into this group. Consider the image of $\, v\, $ in 
$\, K^J_2 ( J_{E / P} \times M \times J^0_U\, ,\, ( J_{E / P} \times M ) \rtimes U_U )\, $ and the diagram 
\smallskip
$$ \vbox{ \halign{ #&#&#\cr
\hfil $K^J_2 ( J_{E / P} \times M \times J^0_U\, ,\, ( J_{E / P} \times M ) \rtimes U_U )$\hfil & 
\hfil $\longrightarrow $\hfil & 
\hfil $K^J_2 ( J_{E / P} \times M \times J^0_{E / P}\, ,\, \widetilde F )$\hfil \cr
\hfil $\Bigm\downarrow $\hfil && \hfil $\Bigm\downarrow $\hfil \cr
\hfil $K^J_2 ( J_{E / P} \times M\, ,\, ( J_{E / P} \times M ) \rtimes U )$\hfil & \hfil $\longrightarrow $\hfil & 
\hfil $K^J_2 ( J_{E / P} \times M\, ,\, \overline F )$\hfil \cr }}\> . $$
\par\medskip\noindent
Then the image of $\, v\, $ in $\, K^J_2 ( J_{E / P} \times M\, ,\, ( J_{E / P} \times M ) \rtimes U )\, $ is represented by a commutator expression in $\, [\, U_{[\, J_{E / P}\, ,\, U_{E / P}\, ]}\, ,\, U_M\, ]\, $ (using the Jacobi identity). The right vertical map splits sending $\, J_{E / P}\, $ into the diagonal of 
$\, J_{E / P} \times J^0_{E / P}\, $ and $\, M\, $ to $\, M\, $ (note that by the factorization 
$\, F = M \rtimes E / P\, $ no perturbation of the diagonal into $\, M\, $ occurs). Divide 
$\, J^0_U \subseteq ( J_{E / P} \times M ) \rtimes U_U\, $ by the commutator subgroup 
$\, [\, J^0_U\, ,\, U_{R , U}\, ]\, $. The basic idea is to find a normal splitting of 
$\, [\, J_{E / P}\, ,\, U_{E / P}\, ]\, $ to $\, [\, J^0_U\, ,\, U_U\, ]\, $ modulo $\, [\, J^0_U\, ,\, U_{R , U}\, ]\, $ which in particular would give a "diagonal" splitting of the left vertical map restricted to 
$\, ( [\, J_{E / P}\, ,\, U_{E / P}\, ] \times M ) \rtimes U\, $ compatible with the splitting on the right side. However this can only occur when the image of 
$\, K^J_2 ( [\, J_{E / P}\, ,\, U_{E / P}\, ]\, ,\, U_{E / P} )\, $ in 
$\, K^J_2 ( J_{E / P}\, ,\, U_{E / P} ) = K^J_3 ( E / P )\, $ is trivial which of course is a too restrictive condition to be of much interest (and even this being the case it isn't clear that such a normal splitting modulo $\, [\, J^0_U\, ,\, U_{R , U}\, ]\, $ exists, thinking of $\, U\, $ as $\, U_{E / P}\, $ and $\, R\, $ as 
$\, J_{E / P}\, $). The strategy is to divide $\, [\, J_{E / P}\, ,\, U_{E / P}\, ]\, $ into smaller portions which do admit a normal lift modulo $\, [\, J^0_U\, ,\, U_{R , U}\, ]\, $ and dividing our representative in 
$\, [\, U_{[\, J_{E / P}\, ,\, U_{E / P}\, ]}\, ,\, U_M\, ]\, $ accordingly into a finite sum of elements each contained in one such portion, to lift them separately and summing up the lifts. To do so one must exhibit a preferred set of generators of $\, [\, J_{E / P}\, ,\, U_{E / P}\, ]\, $ modulo $\, [\, J_{E / P}\, ,\, J_{E / P}\, ]\, $. Modulo $\, [\, J_{E / P}\, ,\, J_{E / P}\, ]\, $ the group $\, [\, J_{E / P}\, ,\, U_{E / P}\, ]\, $ is generated by commutators $\, [\, u_z\, ,\, u_x\, u_y\, u^{-1}_{xy}\, ]\, $. The first task is to find a linear independent subset of these generators which still generates $\, [\, J_{E / P}\, ,\, U_{E / P}\, ]\, $ modulo 
$\, [\, J_{E / P}\, ,\, J_{E / P}\, ]\, $ and is suited for our purposes. Let a generating set 
$\,\{ v_{\lambda } \}\, $ of $\, E / P\, $ be given. We may assume that any identity 
$\, v_{\lambda } = v^{\pm }_{\mu }\, $ implies $\, \lambda = \mu\, $ and 
$\, v_{\lambda } = v^{\pm }_{\mu }\, v^{\pm }_{\nu }\, $ implies $\, \lambda = \mu = \nu\, $ and 
$\, v^3_{\lambda } = 1\, $. Such a set will be called minimal (although it may not be). Let 
$\, l : E / P \rightarrow {\mathbb N}\, $ denote the length function with respect to $\,\{ v_{\lambda } \}\, $. It is clear that modulo the generators $\, \{\> [\, u_z\, ,\, u_x\, u_{x^{-1}}\, ]\> \}\, $ and 
$\, [\, J_{E / P}\, ,\, J_{E / P}\, ]\, $ one can reduce to generaors 
$\, \{\> [\, u_z\, ,\, u_x\, u_y\, u^{-1}_{xy}\, ]\> \}\, $ satisfying $\, l ( x ) \leq l ( xy )\, $. We want to further reduce to such generators where $\, y \in \{ v_{\lambda } \}\, $ is a generator. To this end let 
$\, v^{\pm }_0\, $ be the generator appearing on the utmost right of a minimal expression for $\, y\, $, and put $\, y' = y v^{\mp }_0\, $ so that $\, l ( y') = l ( y ) - 1\, $. Writing 
\smallskip
$$ \bigl[\, u_z\, ,\, u_x\, u_y\, u^{-1}_{xy}\,\bigr]\, =\, 
\bigl[\, u_z\, ,\, u_x\, ( u_y\, u^{-1}_{v^{\pm }_0}\, u^{-1}_{y'} )\, u^{-1}_x\> ( u_x\, u_{y'}\, u^{-1}_{xy'} )\> 
( u_{xy'}\, u_{v^{\pm }_0}\, u^{-1}_{xy} )\,\bigr] $$
$$\equiv \bigl[ u_{y'}\, u_{v^{\pm }_0}\, u^{-1}_y ,\, u_{zx} \bigr]\, 
\bigl[ u_z\, , u_{y'}\, u_{v^{\pm }_0}\, u^{-1}_y \bigr]\, 
\bigl[ u_z\, , u_x\, u_{y'}\, u^{-1}_{xy'} \bigr]\, 
\bigl[ u_z\, , u_{xy'}\, u_{v^{\pm }_0}\, u^{-1}_{xy} \bigr] $$
\par\medskip\noindent
modulo $\, [\, J_{E / P}\, ,\, J_{E / P}\, ]\, $ one has replaced $\, [\, u_z\, ,\, u_x\, u_y\, u^{-1}_{xy}\, ]\, $ by generators in $\, \{\> [\, u_{z'}\, ,\, u_{x'}\, u_{y'}\, u^{-1}_{x'y'}\, ]\>\}\, $ with $\, l ( y' ) \leq l ( y )\, $, so inductively one gets generators of the form 
$\,\{\> [\, u_z\, ,\, u_x\, u_{v^{\pm }_0}\, u^{-1}_{xv^{\pm }_0}\, ]\>\}\, $. Also since 
$\, [\, u_{x^{-1}}\, ,\, u_x\, u_{x^{-1}}\, ] = ( u_{x^{-1}}\, u_x )\> ( u_x\, u_{x^{-1}} )^{-1}\, $ we may replace all generators $\, [\, u_z\, ,\, u_{v_0}\, u_{v^{-1}_0}\, ]\, $ by 
$\, [\, u_z\, ,\, u_{v^{-1}_0}\, u_{v_0}\, ]\, $ where $\, v_0 \in \{ v_{\lambda } \}\, $. If moreover 
$\, l ( x ) = l ( x v^{\pm }_0 )\, $ one may assume that the element $\, x\, $ is in some sense smaller than 
$\, x v^{\pm }_0\, $, for example by choosing an order on the subset of all elements of a given length. We will use the shortcut notation $\, (\, x\, ,\, y\, )\, $ for $\, u_x\, u_y\, u^{-1}_{xy}\, $ and 
$\, (\, z\, ,\, x\, ,\, y\, )\, $ for the element $\, [\, u_z\, ,\, u_x\, u_y\, u^{-1}_{xy}\, ]\, $. Then one has the general equality 
\smallskip
$$ \bigl(\, z\, ,\, x\, ,\, y\, \bigr)\> =\> \bigl(\, z\, ,\, x\, \bigr)\> \bigl(\, zx\, ,\, y\, \bigr)\> \bigl(\, z\, ,\, xy\, \bigr)^{-1}\> \bigl(\, x\, ,\, y\, \bigr)^{-1} $$
\par\medskip\noindent
which degenerates to 
\smallskip
$$ \bigl(\, z\, ,\, x\, ,\, x^{-1}\, \bigr)\> =\> \bigl(\, z\, ,\, x\, \bigr)\> \bigl(\, zx\, ,\, x^{-1}\, \bigr)\> 
\bigl(\, x\, ,\, x^{-1}\, \bigr)^{-1} $$
\par\medskip\noindent
if $\, y = x^{-1}\, $ and to 
\smallskip
$$ \bigl(\, x^{-1}\, ,\, x\, ,\, y\, \bigr)\> =\> \bigl(\, x^{-1}\, ,\, x\, \bigr)\> \bigl(\, x^{-1}\, ,\, xy\, \bigr)^{-1}\> 
\bigl(\, x\, ,\, y\, \bigr)^{-1}$$
\par\medskip\noindent
for $\, z = x^{-1}\, $, respectively to 
\smallskip
$$ \bigl(\, x^{-1}\, ,\, x\, ,\, x^{-1}\, \bigr)\> =\> \bigl(\, x^{-1}\, ,\, x\, \bigr)\> \bigl(\, x\, ,\, x^{-1}\,\bigr)^{-1}\> . $$
\par\medskip\noindent
We are done if we can establish the following Nielsen type property: given a relation in the generators 
$\,\{ (\, z\, ,\, x\, ,\, v^{\pm }_0\, )\}\, $ if any two basis elements of $\, J_{E / P}\, $ in the expression of 
$\, (\, z\, ,\, x\, ,\, v^{\pm }_0\, )\, $ as above are cancelled by according two basis elements of 
$\, (\, z'\, ,\, x'\, ,\, v^{\pm }_1\, )\, $ adjacent to $\, (\, z\, ,\, x\, ,\, v^{\pm }_0\, )\, $ on the left side or on the right in the relation (where $\, (\, z\, ,\, x\, ,\, v^{\pm }_0\, )\, $ is supposed to be of the first three types), then $\, (\, z'\, ,\, x'\, ,\, v^{\pm }_1\, )\, $ $\, =\, (\, z\, ,\, x\, ,\, v^{\pm }_0\, )\, $ and the former appears inverse to the latter. Let us first check the case where $\, (\, z'\, ,\, x'\, ,\, v^{\pm }_1\, )\, =\, 
(\, {x'}^{-1}\, ,\, x'\, ,\, v^{\pm }_1\, )\, $. If $\, (\, z\, ,\, x\, ,\, v^{\pm }_0\, )\, $ is also of this type the condition is readily checked. If $\, (\, z\, ,\, x\, ,\, v^{\pm }_0\, )\, =\, (\, z\, ,\, v^{-1}_0\, ,\, v_0\, )\, $ the cancelling of the two basis elements implies the equations $\, x' = z\, ,\, v^{\pm }_1 = v^{-1}_0\, ,\, {x'}^{-1} = z v^{-1}_0\, ,\, 
x'  v^{\pm }_1 = v_0\, $, i.e. $\, z v^{-1}_0 = z^{-1} = v_0 = v^{-2}_0\, $ which leads to the form 
$\, (\, v_0\, ,\, v^{-1}_0\, ,\, v^{-1}_0\, )\, $ for $\, (\, {x'}^{-1}\, ,\, x'\, ,\, v^{\pm }_1\, )\, $ which is not included in our set (assuming that $\, v_0 < v^{-1}_0 $). If $(\, z\, ,\, x\, ,\, v^{\pm }_0\, )\, $ is nondegenerate one gets $\, x' = z\, ,\, v^{\pm }_1 = x\, ,\, {x'}^{-1} = zx\, ,\, x' v^{\pm }_1 = v^{\pm }_0\, $, i.e. 
$\, zx = z^{-1} = v^{\pm }_0\, $. Then $\, x' v^{\pm }_1\, $, hence $\, x'\, $, is a generator and 
$\, x' = v^{\pm }_1 = v^{\mp }_0\, ,\, v^3_0 = 1\, $ follows. Then $\, (\, {x'}^{-1}\, ,\, x'\, ,\, v^{\pm }_1\, )\, =\, 
(\, v^{\pm }_0\, ,\, v^{\pm }_0\, ,\, v^{\mp }_0\, )\, $ and $\, (\, z\, ,\, x\, ,\, v^{\pm }_0\, )\, =\, 
(\, v^{\mp }_0\, ,\, v^{\mp }_0\, ,\, v^{\pm }_0\, )\, $ contradicting the assumption that 
$\, (\, z\, ,\, x\, ,\, v^{\pm }_0\, )\, $ is nondegenerate. Next assume that $\, (\, z\, ,\, x\, ,\, v^{\pm }_0\, )\, $ is of the first two types and  $\, (\, z'\, ,\, x'\, ,\, v^{\pm }_1\, )\, =\, (\, z'\, ,\, v^{-1}_1\, ,\, v_1\, )\, $. The case 
$\, (\, z\, ,\, v^{-1}_0\, ,\, v_0\, )\, $ is readily checked. Otherwise one gets 
$\, z' = x\, ,\, v^{-1}_1 = v^{-1}_0\, ,\, z' v^{-1}_1 = z\, ,\, v_1 = x v^{-1}_0\, $. Then $\, x v^{-1}_0\, $, hence 
$\, x\, $ is a generator and $\, x = v^{-1}_0\, ,\, v^3_0 = 1\, $ follows, but 
$\, (\, z\, ,\, v^{-1}_0\, ,\, v^{-1}_0\, )\, $ is not contained in our set. Finally suppose that both 
$\, (\, z\, ,\, x\, ,\, v^{\pm }_0\, )\, $ and $\, (\, z'\, ,\, x'\, ,\, v^{\pm }_1\, )\, $ are nondegenerate. The equations 
$\, x' = z\, ,\, v^{\pm }_1 = x\, ,\, z' = zx\, ,\, x' v^{\pm }_1 = v^{\pm }_0\, $ imply that $\, x'\, $ is a generator and $\, x' = v^{\pm }_1 = v^{\mp }_0\, $ giving a contradiction. 
\par\noindent
For each fixed $\, ( x , v^{\pm }_0 )\, $ consider the subset 
$\, A_{( x , v^{\pm }_0 )} = \{ (\, z\, ,\, x\, ,\, v^{\pm }_0\, )\}\, $ of our preferred set of generators 
$\, A = {\bigcup }_{( x , v^{\pm }_0 )}\, A_{( x , v^{\pm }_0 )}\, $. Denote the corresponding subgroup of 
$\, [\, J_{E / P}\, ,\, U_{E / P}\, ]\, $ generated by $\, A_{( x , v^{\pm }_o )}\, $ together with 
$\, [\, J_{E / P}\, ,\, J_{E / P}\, ]\, $ by $\, [\, J_{E / P}\, ,\, U_{E / P}\, ]_{( x , v^{\pm }_0 )}\, $. Then we may factor the element $\, v \in [\, U_{[J_{E / P}\, ,\, U_{E / P}\, ]}\, ,\, U_M\, ]\, $ into a finite sum of elements 
$\, v_{( x , v^{\pm }_0 )} \in [\, U_{[\, J_{E / P}\, ,\, U_{E / P}\, ]_{( x , v^{\pm }_0 )}}\, ,\, U_M\, ]\, $. We claim that there exists for each $\, ( x , v^{\pm }_0 )\, $ a normal lift of 
$\, [\, J_{E / P}\, ,\, U_{E / P}\, ]_{( x , v^{\pm }_0 )}\, $ to 
$\, [\, J^0_U\, ,\, U_U\, ] / [\, J^0_U\, ,\, U_{R , U}\, ]\, $. To do so one notes that by minimality of the generating set $\,\{ v_{\lambda } \}\, $ the element $\, x v^{\pm }_0\, $ is never contained in the set 
$\,\{ v^{\pm }_{\lambda } \}\, $ unless $\, x = v_0\, ,\, x v_0 = v^{-1}_0\, $ (this was already used above). Since $\, J_{E / P}\, $ is free any lift $\, J_{E / P} \rightarrow \widetilde J^0_U\, $ will be normal and canonical on the subgroup $\, [\, J_{E / P}\, ,\, J_{E / P}\, ]\, $ and equal to 
$\, [\, \lambda ( J_{E / P} )\, ,\, \lambda ( J_{E / P} )\, ] \equiv 
[\, \sigma ( J_{E / P} )\, ,\, \lambda ( J_{E / P} )\, ]\, $ where $\, \lambda : J_{E / P} \rightarrow J^0_U\, ,\, 
\sigma : U_{E / P} \rightarrow U_U\, $ are the semicanonical lifts induced from a given section 
$\, s : E / P \nearrow U\, $. We will lift $\, A_{( x , v^{\pm }_0 )}\, $ by the rule 
\smallskip
$$ \bigl[\, u_z\, ,\, u_x\, u_{v^{\pm }_0}\, u^{-1}_{x v^{\pm }_0}\,\bigr]\quad \mapsto\quad 
\bigl[\, \sigma ( u_z )\, ,\, \lambda ( u_x\, u_{v^{\pm }_0}\, u^{-1}_{ x v^{\pm }_0} )\,\bigr]\> . $$
\par\medskip\noindent
Writing out this lift on the basis elements gives 
\smallskip
$$ \bigl(\, z\, ,\, x\, ,\, v^{\pm }_0\, \bigr)\quad\mapsto\quad 
\sigma\, \bigl(\, z\, ,\, x\,\bigr)\> \lambda\, \bigl(\, zx\, ,\, v^{\pm }_0 \bigr)\> 
\widetilde\sigma\,\bigl(\, z\, ,\, x v^{\pm }_0\,\bigr)^{-1}\> \lambda\,\bigl(\, x\, ,\, v^{\pm }_0\,\bigr)^{-1} $$
\par\medskip\noindent
where $\,\widetilde\sigma\, (\, z\, ,\, x v^{\pm }_0\, )\, $ differs from $\, \lambda\, (\, z\, ,\, x v^{\pm }_0\, )\, $ modulo $\, U_{R , U} / [\, J^0_U\, ,\, U_{R , U}\, ]\, $. Obviously this lift extends to a normal lift of $\, [\, J_{E / P}\, ,\, U_{E / P}\, ]_{( x , v^{\pm }_0 )}\, $ if it is compatible with the canonical lift 
$\, [\, \sigma ( J_{E / P} )\, ,\, \lambda ( J_{E / P} )\, ]\, $ (modulo $\, [\, J^0_U\, ,\, U_{R , U}\, ]\, $). Assume given an element in the intersection of the range of $\, A_{( x , v^{\pm }_0 )}\, $ and 
$\, [\, J_{E / P}\, ,\, J_{E / P}\, ]\, $. The basis element $\, (\, z\, ,\, x v^{\pm }_0\, )\, $ can only appear as a basis element of a unique element $\, (\, z\, ,\, x\, ,\, v^{\pm }_0\, )\, $ if $\, x v^{\pm }_0\, $ is not a generator. This means that if such an element appears in the expression of an element of 
$\, \langle A_{( x , v^{\pm }_0 )} \rangle \cap [\, J_{E / P}\, ,\, J_{E / P}\, ]\, $, then the whole expression 
$\, (\, z\, ,\, x\, ,\, v^{\pm }_0\, )\, $ appears twice, cancelling modulo 
$\, [\, \langle A_{( x , v^{\pm }_0 )} \rangle\, ,\, \langle A_{( x , v^{\pm }_0 )} \rangle\, ]\, $ on which subgroup our lift identifies with the canonical lift of $\, [\, J_{E / P}\, ,\, J_{E / P}\, ]\, $. Then only the cases where 
$\, x = v_0\, ,\, v^3_0 = 1\, $ and the cases $\, x v^{\pm }_0 = 1\, $ remain. In the first case one can choose the section $\, s : E / P \nearrow U\, $ such that $\, s ( v^{-1}_0 ) = s ( v_0 ) s ( v_0 )\, $. Then 
$\,\widetilde\sigma (\, z\, ,\, x v^{\pm }_0\, )\, $ can be replaced by 
$\,\sigma (\, z\, ,\, x v^{\pm }_0\, )\, $ (also replacing $\, \lambda (\, zx\, ,\, v^{\pm }_0\, )\, $ by 
$\, \sigma (\, zx\, ,\, v^{\pm }_0\, )\, $ and 
$\, \lambda (\, x\, ,\, v^{\pm }_0\, ) = \sigma (\, x\, ,\, v^{\pm }_0\, )\,  $), so that one gets an expression consisting solely of elements $\, \sigma (\, z\, ,\, x\, )\, $ (and $\, \lambda (\, z\, ,\, x\, )\, $ for 
$\, x v^{\pm }_0 = 1\, $). Since the expression is contained in 
$\, [\, \widetilde J^0_U\, ,\, \widetilde U_U\, ]\, $ it follows that the difference to the corresponding expression in the canonical subgroup $\, [\, \lambda ( J_{E / P} )\, ,\, \lambda ( J_{E / P} )\, ]\, $ is contained in $\, [\, U_{R , U}\, ,\, J^0_U + U_{R , U}\, ] \cap J_{R , U} = [\, J^0_U\, ,\, U_{R , U}\, ]\, $ since 
$\, R\, $ is free. Thus one has constructed a normal splitting of 
$\, [\, J_{E / P}\, ,\, U_{E / P}\, ]_{( x , v^{\pm }_0 )}\, $ to $\, [\, \widetilde J^0_U\, ,\, \widetilde U_U\, ]\, $ also giving a diagonal splitting 
\smallskip
$$ K^J_2 ( [\, J_{E / P}\, ,\, U_{E / P}\, ]_{( x , v^{\pm }_0 )} \times M\, ,\, 
( [\, J_{E / P}\, ,\, U_{E / P}\, ]_{( x , v^{\pm }_0 )} \times M ) \rtimes U ) \longrightarrow $$
$$ K^J_2 ( [\, J_{E / P}\, ,\, U_{E / P}\, ]_{( x , v^{\pm }_0 )} \times M \times 
[\, \widetilde J^0_U\, ,\, \widetilde U_U\, ]\, ,\, ( [\, J_{E / P}\, ,\, U_{E / P}\, ]_{( x , v^{\pm }_0 )} \times M ) \rtimes \widetilde U_U ) $$
\par\medskip\noindent
compatible with the diagonal splitting of the right vertical map of the diagram above.  Extending 
$\,\langle A_{( x , v^{\pm }_0 )} \rangle\, $ to the normalization of the basis element $\, ( x , v^{\pm }_0 )\, $ of $\, J_{E / P}\, $ denote the  corresponding subgroup by $\, \overline {\langle A_{( x , v^{\pm }_0 )} \rangle }\, $, the $\, U_{E / P} / \overline {\langle A_{( x , v^{\pm }_0 )}\rangle }\, $ is freely generated by 
$\, \{\, u_z\,\vert\, z \neq x v^{\pm }_0\,\}\, $ such that $\, u_z\, $ maps to $\, u_z\, $ if 
$\, z \neq x v^{\pm }_0\, $ and $\, u_{x v^{\pm }_0 }\, $ maps to $\, u_x\, u_{v^{\pm }_0}\, $. A similar result holds in the cases $\, ( x , v^{\pm }_0 ) = ( v^{-1}_0 , v_0 )\, $ and $\, ( x , v^{\pm }_0 ) = ( v_0 , v_0 )\, ,\, 
v^3_0 = 1\, $ in that $\, U_{E / P} / \overline {\langle A_{( x , v^{\pm }_0 )} \rangle }\, $ is freely generated by $\,\{\, u_z\,\vert\, z \neq v^{-1}_0\,\}\, $ in these cases. In particular modulo 
$\, [\, J_{E / P}\, ,\, J_{E / P}\, ]\, $ there exists an equivariant bijection of $\, [\, J_{E / P}\, ,\, U_{E / P}\, ]\, $ to itself which is the identity for $\,\langle A_{( x , v^{\pm }_0 )} \rangle\, $ and sends the complement for 
$\, A_{( x , v^{\pm }_0 )}\, $ in $\, A\, $ into the range of the splitting 
$\, U_{E / P} / \overline {\langle A_{( x , v^{\pm }_0 )} \rangle} \rightarrow U_{E / P}\, $ (which is normal modulo $\, [\, J_{E / P}\, ,\, J_{E / P}\, ]\, $). Now let $\, v = \sum v_{(x , v^{\pm }_0 )}\, $ be a factorization of $\, v\, $ into a finite sum with each $\, v_{( x , v^{\pm }_0 )}\, $ contained in 
$\, [\, U_{\langle A_{( x , v^{\pm }_0 )} \rangle}\, ,\, U_M\, ]\, $ perturbing $\, v\, $ modulo 
$\, [\, U_{[\, J_{E / P}\, ,\, J_{E / P}\, ]}\, ,\, U_M\, ]\, $ if necessary. It is clear that the image of $\, v\, $ in the lower right group $\, K^J_2 ( J_{E / P} \times M\, ,\, M \rtimes U_{E / P} )\, $ of the diagram above is trivial. We claim that this is trus for each single $\, v_{( x , v^{\pm }_0 )}\, $. Suppose this is not the case. Dividing by the normalization of the complement $\, \langle A^c \rangle\, $ of $\, A\, $ to 
$\, A_{( x , v^{\pm }_0 )}\, $ the element $\, v_{( x , v^{\pm }_0 )}\, $ becomes trivial, so it defines an element in the trefoil intersection of $\, M\, ,\, \overline {\langle A_{( x , v^{\pm }_0 )} \rangle }\, $ and 
$\, \langle A^c \rangle\, $ (modulo $\, [\, \overline {\langle A_{( x , v^{\pm }_0 )}  \rangle }\, ,\, J_{E / P}\, ]\, $ say) in $\, K^J_2 ( \overline {J_{E / P}} \times M\, ,\, M \times \overline {U_{E / P}} )\, $. This element is equivalently represented by the inverse of the image of $\, v^c = v v^{-1}_{( x , v^{\pm }_0 )} \in 
[\, U_{\langle A^c \rangle }\, ,\, U_M\, ]\, $. Using the equivariant bijection sending $\, A^c\, $ to the complement $\, J^c_{E / P}\, $ for $\, A_{( x , v^{\pm }_0 )}\, $ in $\, J_{E / P}\, $ defined by the splitting above one gets a representative in $\, [\, U_{J^c_{E / P}}\, ,\, U_M\, ]\, $ which differs from the original element only modulo $\, [\, U_{\overline {\langle A_{( x , v^{\pm }_0 )} \rangle }}\, ,\, U_M\, ]\, $ hence it defines an element in the trefoil intersection of $\, M\, ,\, \overline {\langle A_{( x , v^{\pm }_0 )}\rangle }\, $ and $\, \overline {J^c_{E / P}}\, $ in the group 
$\, K^J_2 ( \overline { J_{E / P}} \times M\, ,\, M \rtimes \overline {U_{E / P}} )\, $ 
(modulo $\, [\, \overline {\langle A_{( x , v^{\pm }_0 )} \rangle }\, ,\, J_{E / P}\, ]\, $). Since 
$\, \overline  {\langle A_{( x , v^{\pm }_0 )} \rangle }\, $ has a weak core for the quotient $\, E / P\, $, Lemma 12 yields that the kernel of 
\smallskip
$$ K^J_2 ( \overline {\langle A_{( x , v^{\pm }_0 )}\rangle }\, ,\, M \rtimes \overline {U_{E / P}} ) \twoheadrightarrow 
K^J_2 ( \overline {\langle A_{( x , v^{\pm }_0 )}\rangle }\, ,\,\overline {U_{E / P}} / \overline {J^c_{E / P}} )$$
\par\medskip\noindent
is given by the formula $\, {\mathcal C}_{\overline {\langle A_{( x , v^{\pm }_0 )}\rangle }\, ,\, E / P} \otimes 
( M \times {\overline {J^c_{E / P}}}^{ab} )\, =\, 
( {\mathcal C}_{\overline {\langle A_{( x , v^{\pm }_0 )}\rangle }\, ,\, E / P} \otimes M ) \times 
( {\mathcal C}_{\overline {\langle A_{( x , v^{\pm }_0 )}\rangle }\, ,\, E / P} \otimes 
{\overline {J^c_{E / P}}}^{ab} )\, $ where 
$\, {\mathcal C}_{\overline {\langle A_{( x , v^{\pm }_0 )}\rangle }\, ,\, E / P}\, $ denotes the core of 
$\, \overline {\langle A_{( x , v^{\pm }_0 )}\rangle }\, $ (which is isomorphic to $\,\mathbb Z\, $). On the other hand one gets $\, \overline {U_{E / P}} / \overline {J^c_{E / P}} \simeq 
\overline {\langle A_{( x , v^{\pm }_0 )}\rangle } \rtimes E / P\, $. This implies that the following maps are injective (by existence of a splitting)
\smallskip
$$ K^J_2 ( \overline {\langle A_{( x , v^{\pm }_0 )}\rangle }\, ,\, M \rtimes \overline {U_{E / P}} ) 
\rightarrowtail K^J_2 ( \overline {\langle A_{( x , v^{pm }_0 )}\rangle } \times M \times 
\overline {J^c_{E / P}}\, ,\, M \rtimes \overline {U_{E / P}} )\> , $$
$$ K^J_2 ( \overline {\langle A_{( x , v^{\pm }_0 )}\rangle }\, ,\, 
M \rtimes \overline {U_{E / P}} / \overline {J^c_{E / P}} ) \rightarrowtail 
K^J_2 ( \overline {\langle A_{( x , v^{\pm }_0 )}\rangle } \times M\, ,\, 
M \rtimes \overline {U_{E / P}} / \overline {J^c_{E / P}} )\>  ,  $$
$$ K^J_2 ( \overline {\langle A_{( x , v^{\pm }_0 )}\rangle }\, ,\, \overline {U_{E / P}} ) \rightarrowtail 
K^J_2 ( \overline {\langle A_{( x , v^{\pm }_0 )}\rangle } \times \overline {J^c_{E / P}}\, ,\, 
\overline {U_{E / P}} ) \> . $$
\par\medskip\noindent
Then any element in the trefoil intersection as above is in the intersection of 
$\, {\mathcal C}_{\overline {\langle A_{( x , v^{\pm }_0 )}\rangle }\, ,\, E / P} \otimes 
{\overline {J^c_{E / P}}}^{ab}\, $ and 
$\, {\mathcal C}_{\overline {\langle A_{( x , v^{\pm }_0 )}\rangle }\, ,\, E / P} \otimes M\, $ hence trivial. Then if $\, \overline {[\, J_{E / P}\, ,\, U_{E / P}\, ]^c}\, $ denotes the intersection of 
$\, [\, \overline {J_{E / P}}\, ,\, \overline {U_{E / P}}\, ]\, $ with $\, \overline {J^c_{E / P}}\, $ our element is trivial already in the trefoil intersection of $\, M\, ,\, \overline {\langle A_{( x , v^{\pm }_0 )}\rangle }\, $ and 
$\, \overline {[\, J_{E / P}\, ,\, U_{E / P}\, ]^c}\, $ in 
$\, K^J_2 ( ( \overline {\langle A_{( x , v^{\pm }_0 )}\rangle } + 
[\,\overline {J_{E / P}}\, ,\,\overline {U_{E / P}}\, ] ) \times M\, ,\, M \rtimes \overline {U_{E / P}} )\, $ since it is trivial in 
$\, K^J_2 ( \overline {\langle A_{( x , v^{\pm }_0 )}\rangle }\, ,\, M \rtimes \overline {U_{E / P}} )\, $. By surjectivity of the map sending the kernel of 
\smallskip
$$ K^J_2 ( \overline {\langle A_{( x , v^{\pm }_0 )}\rangle }\, ,\, M \rtimes \overline {U_{E / P}} ) \twoheadrightarrow
K^J_2 ( \overline {\langle A_{( x , v^{\pm }_0 )}\rangle }\, ,\, 
M \rtimes ( \overline {U_{E / P}} / \overline {[\, J_{E / P}\, ,\, U_{E / P}\, ]^c} ) $$
\par\medskip\noindent
to the kernel of 
\smallskip
$$ K^J_2 ( \overline {[\, J_{E / P} , U_{E / P}\, ]^c} , M \rtimes \overline {U_{E / P}} ) \twoheadrightarrow 
K^J_2 ( \overline {[\, J_{E / P} , U_{E / P}\, ]^c} , M \rtimes 
\overline {U_{E / P}} / \overline {\langle A_{( x , v^{\pm }_0 )}\rangle } ) $$
\par\medskip\noindent
one gets that our representative also defines a trivial element in the group 
$\, K^J_2 ( \overline {[\, J_{E / P}\, ,\, U_{E / P}\, ]^c}\, ,\, M \rtimes \overline {U_{E / P}} )\, $ which is therefore represented as a commutator expression in 
\smallskip
$$ \bigl[\, J_{\overline {[\, J_{E / P}\, ,\, U_{E / P}\, ]^c}\, ,\, M \rtimes \overline {U_{E / P}}}\, ,\, 
U_{M \rtimes \overline {U_{E / P}}}\, \bigr]\, +\, 
\bigl[\, U_{\overline {[\, J_{E / P}\, ,\, U_{E / P}\, ]^c}}\, ,\, J_{M \rtimes \overline {U_{E / P}}}\,\bigr]\> . $$
\par\medskip\noindent
Using the inverse of the equivariant bijection as above this element maps into the trefoil intersection of 
$\, M\, ,\, \overline {\langle A_{( x , v^{\pm }_0 )}\rangle }\, $ and $\, \overline {\langle A^c \rangle }\, $ 
(modulo $\, [\, J_{E / P}\, ,\, J_{E / P}\, ]\, $) and is represented modulo a trivial element by a commutator in $\, [\, U_{\overline {\langle A^c \rangle }}\, ,\, U_{\overline {\langle A_{( x , v^{\pm }_0 )}\rangle }}\, ]\, $ (and also of course by the element $\, v_{( x , v^{\pm }_0 )}\, $). Due to the splitting of the map corresponding to division by $\, M\, $ this element must be trivial. Then all of the elements 
$\, [ v_{( x , v^{\pm }_0 )} ]\, $ are trivial in $\, K^J_2 ( J_{E / P} \times M\, ,\, M \rtimes U_{E / P} )\, $. Returning to our diagram above one gets that taking the difference of the original representative 
$\, v_{( x , v^{\pm }_0 )}\, $ and the diagonal lift $\, \Delta ( v_{( x , v^{\pm }_0 )} )\, $ as above defines a preimage in 
$\, K^J_2 ( J_{E / P} \times M \times \widetilde J^0_U\, ,\, ( J_{E / P} \times M ) \rtimes \widetilde U_U )\, $ for the element $\, [ v_{( x , v^{\pm }_0 )} ] \in 
K^J_2 ( J_{E / P} \times M \times J^0_{E / P}\, ,\, ( J_{E / P} \times M ) \rtimes U_{E / P} )\, $. It is clear by construction that this element lies in the common image of 
$\, K^J_2 ( M\, ,\, ( J_{E / P} \times M ) \rtimes \widetilde U_U )\,$ and  
$\, K^J_2 ( \widetilde {\langle A_{( x , v^{\pm }_0 )}\rangle }\, ,\, 
 J_{E / P} \times M ) \rtimes \widetilde U_U )\, $ where 
 $\, \widetilde {\langle A_{( x , v^{\pm }_0 )}\rangle }\, $ is a normal lift of 
 $\, \overline {\langle A_{( x , v^{\pm }_0 )}\rangle }\, $ in $\,\widetilde J^0_U\, $. Dividing by 
 $\, J_{E / P}\, $ gives an element of 
 $\, K^J_2 ( M \times \widetilde J^0_U\, ,\, M \rtimes \widetilde U_U )\, $ which becomes trivial dividing by 
 $\, U_{R , U} / [\, J^0_U\, ,\, U_{R , U}\, ]\, $, so it defines an element in the trefoil intersection of 
 $\, M\, ,\, \widetilde {\langle A_{( x , v^{\pm }_0 )}\rangle }\, $ and $\, \widetilde U_{R , U}\, $ in 
 $\, K^J_2 ( M \times \widetilde {\langle A_{( x , v^{\pm }_0 )}\rangle } \times \widetilde U_{R , U}\, ,\, 
 M \rtimes \widetilde U_U )\, $. Then $\, \widetilde {\langle A_{( x , v^{\pm }_0 )}\rangle }\, $ has a weak core for the quotient $\, U_{E / P} / \overline {\langle A_{( x , v^{\pm }_0 )}\rangle }\, $ so Lemma 12 gives the formula $\, {\mathcal C}_{\overline {\langle A_{( x , v^{\pm }_0 )}\rangle }\, ,\, U_{E / P}} \otimes 
( M \times \widetilde U_{R , U}^{ab} )\, $ for the kernel of the map
\smallskip
$$ K^J_2 ( \widetilde {\langle A_{( x , v^{\pm }_0 )}\rangle }\, ,\, M \rtimes \widetilde U_U ) \twoheadrightarrow 
K^J_2 ( \overline {\langle A_{( x , v^{\pm }_0 )}\rangle }\, ,\, U_{E / P} ) $$
\par\medskip\noindent
and the following maps are injective by existence of a splitting
\smallskip
$$ K^J_2 ( \widetilde {\langle A_{( x , v^{\pm }_0 )}\rangle }\, ,\, M \rtimes \widetilde U_U ) \rightarrowtail 
K^J_2 ( M \times \widetilde {\langle A_{( x , v^{\pm }_0 )}\rangle } \times 
\widetilde U_{R , U}\, ,\, M \rtimes \widetilde U_U ) \> , $$
$$ K^J_2 ( \widetilde {\langle A_{( x , v^{\pm }_0 )}\rangle }\, ,\, \widetilde U_U ) \rightarrowtail 
K^J_2 ( \widetilde {\langle A_{( x , v^{\pm }_0 )}\rangle } \times \widetilde U_{R , U}\, ,\,
\widetilde U_U )\> , $$ 
$$ K^J_2 ( \overline {\langle A_{( x , v^{\pm }_0 )}\rangle }\, ,\, 
M \rtimes U_{E / P} ) \rightarrowtail 
K^J_2 ( M \times \overline {\langle A_{( x , v^{\pm }_0 )}\rangle }\, ,\, M \rtimes U_{E / P} )\> . $$
\par\medskip\noindent
 showing that the trefoil intersection as above is trivial. Thus the element
 $\, [\, v_{( x , v^{\pm }_0 )}\cdot \Delta ( v_{( x , v^{\pm }_0 )} )^{-1}\, ]\, $ is a lift of 
 $\, [\, v_{( x , v^{\pm }_0 )} \,] \in 
 K^J_2 ( J_{E / P} \times M \times J^0_{E / P}\, ,\, ( J_{E / P} \times M ) \rtimes U_{E / P} )\, $ to the trefoil intersection of $\, M\, ,\, J_{E / P}\, $ and $\, \widetilde J^0_U\, $  in 
 $\, K^J_2 ( J_{E / P} \times M \times \widetilde J^0_U\, ,\, ( J_{E / P} \times M ) \rtimes \widetilde U_U )\, $. It is straightforward to find a lift of this element to the trefoil intersection of $\, J_{E / P}\, ,\, M\, $ and 
 $\, J^0_U\, $ in $\, K^J_2 ( J_{E / P} \times M \times J^0_U\, ,\, ( J_{E / P} \times M ) \rtimes U_U )\, $ which is trivial, so that $\, [\, v_{( x , v^{\pm }_0 )} \, ]\, $ and hence also 
 $\, [\, v \, ] \in K^J_2 ( J_{E / P} \times M \times J^0_{E / P}\, ,\, \widetilde F )\, $ is trivial. This means that 
 $\, [\, v\, ]\, $ is trivial already as an element of 
 $\, K^J_2 ( J_{E / P} \times J^0_{E / P}\, ,\, \widetilde F )\, $  since the map 
 \smallskip
 $$ K^J_2 ( J_{E / P}\, ,\, \widetilde F ) \rightarrowtail 
 K^J_2 ( J_{E / P} \times M \times J^0_{E / P}\, ,\, \widetilde F ) $$
 \par\medskip\noindent
 is injective. Then also $\, [\, w\, ]\, $ and hence $\, [\, z\, ]\, $ is trivial as an element of 
 $\, K^J_2 ( J_{E / P} \times J^0_{E / P}\, ,\, \widetilde F )\, $ showing that 
 \smallskip
 $$ K^J_2 ( J_{E / P} \times J^0_{E / P}\, ,\, \widetilde F ) \rightarrowtail 
 K^J_2 ( J_{E / P} \times M \times J^0_{E / P}\, ,\, \widetilde F ) $$
 \par\medskip\noindent
 is injective and $\, ( M , F )\, $ is exact.
 \par\noindent 
Assume now that $\, F\, $ is perfect and $\, N \subseteq F\, $ is a full subgroup. Let $\, E = F / N\, $ and consider the universal central extensions $\,\overline F\, $ and $\,\overline E\, $ of $\, F\, $ and $\, E\, $ respectively. One gets an extension 
\smallskip
$$ 1 \longrightarrow ( \overline N , \overline F ) \longrightarrow ( \overline F , \overline F ) \longrightarrow 
( \overline E , \overline E ) \longrightarrow 1 $$
\par \medskip\noindent
where $\,\overline N\, $ is the image of the universal $F$-central extension $\, {\overline N}^{ _F}\, $ of 
$\, N\, $ in $\,\overline F\, $. $\,\overline N\, $ is again a full subgroup of $\, \overline F\, $, so we begin by considering the case $\, ( N , F ) = ( \overline N , \overline F )\, $. Exactness follows from the first part of the Lemma. We have to show exactness of the pair 
$\, ( J_{\overline E}\, ,\, \overline N \rtimes U_{\overline E} )\, $. Dividing $\, J_{\overline E}\, $ by 
$\, [\, J_{\overline E}\, ,\, U_{\overline E}\, ]\, $this can be reduced using the first part to proving exactness of the pair $\, ( D , D \times \overline F )\, $, where 
$\, D = J_{\overline E} / [\, J_{\overline E}\, ,\, U_{\overline E}\, ]\, $ is a free abelian central direct summand of the quotient $\, \overline N \rtimes U_{\overline E} / [\, J_{\overline E}\, ,\, U_{\overline E}\, ] 
\simeq D \times \overline F\, $. Since $\, ( \overline F , \overline F )\, $ and $\, ( D , D )\, $ are fully exact, exactness of the pair $\, ( D , D \times \overline F )\, $ follows from full exactness of 
$\, ( \overline F , D \times \overline F )\, $ by the $P \times Q$-Lemma below, and exactness of 
$\, ( \overline F , D \times \overline F )\, $ was shown above, so that it remains to show that 
$\, ( \overline F , D \times \overline F )\, $ is second order exact. Consider the following diagram
\smallskip
$$ \vbox{ \halign {#&#&#\cr
\hfil $\bigl( J_{D \times \overline F}\, ,\, J_D \rtimes U_{D \times \overline F} \bigr)$\hfil & 
\hfil $\rightarrow $\hfil & 
\hfil $\bigl( J_{D \times \overline F} \times J_D\, ,\, J_D \rtimes U_{D \times \overline F} \bigr)$\hfil \cr
\hfil $\Bigm\downarrow $\hfil && \hfil $\Bigm\downarrow $\hfil \cr 
\hfil $\bigl( J_{D \times \overline F}  + U_{\overline F , D \times \overline F} ,
J_D \rtimes U_{D \rtimes \overline F} \bigr)$\hfil & \hfil $\rightarrow $\hfil & 
\hfil $\bigl( ( J_{D \times \overline F} + U_{\overline F , D \times \overline F} ) \times J_D , 
J_D \rtimes U_{D \times \overline F} \bigr)$\hfil \cr
\hfil $\Bigm\downarrow $\hfil && \hfil $\Bigm\downarrow $\hfil \cr
\hfil $\bigl( \overline F , \overline F \times U_D \bigr)$\hfil & \hfil $\rightarrow $\hfil & 
\hfil $\bigl( \overline F , \overline F \times D \bigr)$\hfil \cr }}\> . $$ 
\par\medskip\noindent
We need to sow injectivity of the upper horizontal map in $K^J_2$. Any element in its kernel will be trivial in $\, K^J_2 ( J_{D \times \overline F} + U_{\overline F , D \times \overline F}\, ,\, 
J_D \rtimes U_{D \times \overline F} )\, $, hence lifts to $\, K^J_3 ( \overline F , \overline F \times U_D )\, $ by exactness of $\, ( \overline F , \overline F \times U_D )\, $. One has 
$\, K^J_3 ( \overline F , \overline F \times U_D ) \simeq K^J_3 ( \overline F , \overline F \times D )\, $ by the natural map. This is seen in the following way. The map 
$\, K^J_2 ( J_{\overline F , \overline F \times U_D}\, ,\, U_{\overline F \times U_D} ) \rightarrow 
K^J_2 ( J_{\overline F , \overline F \times D}\, ,\, U_{\overline F \times D} )\, $ factors into the injective map 
\smallskip
$$ K^J_2 ( J_{\overline F , \overline F \times U_D}\, ,\, U_{\overline F \times U_D} ) \rightarrowtail 
K^J_2 ( J_{\overline F , \overline F \times D}\, ,\, U_{J_D , U_D} \rtimes U_{\overline F \times D} ) $$
\par\medskip\noindent
and the regular surjection 
\smallskip
$$ K^J_2 ( J_{\overline F , \overline F \times D}\, ,\, U_{J_D , U_D} \rtimes U_{\overline F \times D} ) 
\twoheadrightarrow K^J_2 ( J_{\overline F , \overline F \times D}\, ,\, U_{\overline F \times D} )\> . $$
\par\medskip\noindent
Since $\, U_{\overline F , \overline F \times D}\, $ and $\, U_{J_D , U_D}\, $ both have a core for the quotient $\, U_D\, $, Lemma 12 gives the formula 
$\, {\mathcal C}_{J_D , U_D} \otimes U^{ab}_{\overline F , \overline F \times D}\, $ for the kernel of the regular surjection 
\smallskip
$$ K^J_2 ( U_{\overline F , \overline F \times D}\, ,\, U_{J_D , U_D} \rtimes U_{\overline F \times D} ) \twoheadrightarrow 
K^J_2 ( U_{\overline F , \overline F \times D}\, ,\, U_{\overline F \times D} ) $$
\par\medskip\noindent
so that an element in the kernel of 
\smallskip 
$$ K^J_2 ( J_{\overline F , \overline F \times D}\, ,\, U_{J_D , U_D} \rtimes U_{\overline F \times D} ) \twoheadrightarrow 
K^J_2 ( J_{\overline F , \overline F \times D}\, ,\, U_{\overline F \times D} ) $$
\par\medskip\noindent
which is trivial in 
$\, K^J_2 ( U_{\overline F , \overline F \times D}\, ,\, U_{J_D , U_D} \rtimes U_{\overline F \times D} )\, $ 
is in the image of 
$\, {\mathcal C}_{J_D , U_D} \otimes ( J_{\overline F , \overline F \times D} \cap 
[\, U_{\overline F , \overline F \times D}\, ,\, U_{\overline F , \overline F \times D}\, ] )^{ab}\, =\, 
{\mathcal C}_{J_D , U_D} \otimes [\, J_{\overline F , \overline F \times D}\, ,\, U_{\overline F , \overline F \times D}\, ]^{ab}\, $ 
which subgroup is easily seen to lie in the kernel of the natural map to 
$\, K^J_2 ( J_{\overline F , \overline F \times D}\, ,\, U_{J_D , U_D} \rtimes U_{\overline F \times D} )\, $. But the map to  
$\, K^J_2 ( U_{\overline F , \overline F \times D}\, ,\, U_{J_D , U_D} \rtimes U_{\overline F \times D} )\, $ factors over $\, K^J_2 ( U_{\overline F , \overline F \times U_D}\, ,\, U_{\overline F \times U_D} ) = 0\, $, so it follows that the map $\, K^J_3 ( \overline F\, ,\, \overline F \times U_D ) \buildrel\sim\over\longrightarrow 
K^J_3 ( \overline F , \overline F \times D )\, $ is an isomorphism. Also the map 
$\, K^J_3 ( \overline F \times U_D ) \rightarrowtail K^J_3 ( \overline F \times D )\, $ is injective, since there is a projection $\, K^J_3 ( \overline F \times U_D ) \rightarrow K^J_3 ( \overline F ) \simeq 
K^J_3 ( \overline F , \overline F \times U_D )\, $ which shows that 
$\, K^J_3 ( \overline F , \overline F \times U_D ) \buildrel\sim\over\longrightarrow 
K^J_3 ( \overline F \times U_D )\, $ is injective hence an isomorphism, and by the same argument 
$\, K^J_3 ( \overline F , \overline F \times D ) \rightarrowtail K^J_3 ( \overline F \times D )\, $ is injective. Now consider the image of our lifted element in $\, K^J_3 ( \overline F \times D )\, $ which lifts to 
$\, K^J_3 ( J_D \rtimes U_{D \times \overline F} )\, $, and taking the difference with the original lift in 
$\, K^J_3 ( \overline F , \overline F \times U_D )\, $ one gets a lift which is trivial in 
$\, K^J_3 ( \overline F \times D )\, $, hence trivial, so the map 
\smallskip
$$ K^J_2 ( J_{D \times \overline F}\, ,\, J_D \rtimes U_{D \times \overline F} ) \rightarrowtail 
K^J_2 ( J_{D \times \overline F} \times J_D\, ,\, J_D \rtimes U_{D \times \overline F} ) $$
\par\medskip\noindent
is injective, and $\, ( J_D\, ,\, \overline F \times U_D )\, $ is exact. Next we want to induce injectivity of 
\smallskip
$$ K^J_2 ( J_F\, ,\, J_E \rtimes U_F ) \rightarrow K^J_2 ( J_F \times J_E\, ,\, J_E \rtimes U_F ) $$
\par\medskip\noindent
from injectivity of 
\smallskip
$$ K^J_2 ( J_{\overline F}\, ,\, J_{\overline E} \rtimes U_{\overline F} ) \rightarrowtail 
K^J_2 ( J_{\overline F} \times J_{\overline E}\, ,\, J_{\overline E} \rtimes U_{\overline F} )\> . $$
\par\medskip\noindent
Let $\, C_N\, $ be the kernel of $\, \overline N \twoheadrightarrow N\, $ and put 
$\, {\overline F}^{ _E} = \overline F / C_N\, $. Then injectivity of the last map implies injectivity of
\smallskip
$$ K^J_2 ( J_{{\overline F}^{ _E}}\, ,\, ( J_{\overline E} \times C_N ) \rtimes U_{{\overline F}^{ _E}} ) \rightarrowtail 
K^J_2 ( J_{{\overline F}^{ _E}} \times J_{\overline E}\, ,\, ( J_{\overline E} \times C_N ) \rtimes 
U_{{\overline F}^{ _E}} ) $$
\par\medskip\noindent
by Corollary 5.1 and a simple diagram chase. Since one gets a compatible splitting 
\smallskip
$$ \bigl( J_{{\overline F}^{ _E}}\, ,\, J_{\overline E}\, ,\, J_{\overline E} \rtimes U_{{\overline F}^{ _E}} \bigr) 
\rightarrow \bigl( J_{{\overline F}^{ _E}}\, ,\, J_{\overline E}\, ,\, 
( J_{\overline E} \times C_N ) \rtimes U_{{\overline F}^{ _E}} \bigr) $$
\par\medskip\noindent
this implies injectivity of 
\smallskip
$$ K^J_2 ( J_{{\overline F}^{ _E}}\, ,\, J_{\overline E} \rtimes U_{{\overline F}^{ _E}} ) \rightarrowtail 
K^J_2 ( J_{{\overline F}^{ _E}} \times J_{\overline E}\, ,\, 
J_{\overline E} \rtimes U_{{\overline F}^{ _E}} )\> . $$
\par\medskip\noindent
Next we want to show injectivity of 
\smallskip
$$ K^J_2 ( J_{{\overline F}^{ _E}} \times J_{C , \overline E}\, ,\, 
J_{\overline E} \rtimes U_{{\overline F}^{ _E}} ) \rightarrow 
K^J_2 ( J_{{\overline F}^{ _E}} \times J_{\overline E}\, ,\, 
J_{\overline E} \rtimes U_{{\overline F}^{ _E}} ) \> . $$
\par\medskip\noindent
Suppose this not being the case. Any element in the kernel gives a nontrivial element in the kernel of 
\smallskip
$$ K^J_2 ( J_{C , \overline E}\, ,\, N \rtimes U_{\overline E} ) \rightarrow 
K^J_2 ( J_{\overline E}\, ,\, N \rtimes U_{\overline E} ) $$
\par\medskip\noindent
which by injectivity of $\, K^J_2 ( N\, ,\, N \rtimes U_{\overline E} ) \rightarrowtail 
K^J_2 ( N \rtimes U_{\overline E} )\, $ must also be in the kernel of the map
\smallskip
$$ K^J_2 ( J_{C , \overline E}\, ,\, N \rtimes U_{\overline E} ) \rightarrow 
K^J_2 ( J_{C , \overline E} \times N\, ,\, N \rtimes U_{\overline E} ) $$
\par\medskip\noindent
(since $\, K^J_2 ( J_{C , \overline E}\, ,\, U_{\overline E} ) = 0\, $). From the argument above we know that the map 
\smallskip
$$ K^J_2 ( J_{\overline E}\, ,\, N \rtimes U_{\overline E} ) \rightarrowtail 
K^J_2 ( J_{\overline E} \times N\, ,\, N \rtimes U_{\overline E} ) $$
\par\medskip\noindent
is injective since $\, N\, $ is full. A similar argument (restricting to the normal subgroup 
$\, J_{C , \overline E} \subseteq J_{\overline E}\, $) then yields injectivity of 
\smallskip
$$ K^J_2 ( J_{C , \overline E}\, ,\, N \rtimes U_{\overline E} ) \rightarrowtail 
K^J_2 ( J_{C , \overline E} \times N\, ,\, N \rtimes U_{\overline E} )\> , $$
\par\medskip\noindent
implying that
\smallskip
$$ K^J_2 ( J_{{\overline F}^{ _E}} \times J_{C , \overline E}\, ,\, 
J_{\overline E} \rtimes U_{{\overline F}^{ _E}} ) \rightarrowtail 
K^J_2 ( J_{{\overline F}^{ _E}} \times J_{\overline E}\, ,\, J_{\overline E} \rtimes U_{{\overline F}^{ _E}} )$$
\par\medskip\noindent
and then 
\smallskip
$$ K^J_2 ( J_{{\overline F}^{ _E}}\, ,\, J_E \rtimes U_{{\overline F}^{ _E}} ) \rightarrowtail 
K^J_2 ( J_{{\overline F}^{ _E}} \times J_E\, ,\, J_E \rtimes U_{{\overline F}^{ _E}} ) $$
\par\medskip\noindent
is injective. Then also the map
\smallskip
$$ K^J_2 ( J_{{\overline F}^{ _E}} + U_{C , {\overline F}^{ _E}}\, ,\, J_E \rtimes U_{{\overline F}^{ _E}} ) \rightarrowtail 
K^J_2 ( ( J_{{\overline F}^{ _E}} + U_{C , {\overline F}^{ _E}} ) \times J_E\, ,\, 
J_E \rtimes U_{{\overline F}^{ _E}} ) $$
\par\medskip\noindent
is injective, since any nontrivial element in its kernel would give a nontrivial element in 
$\, K^J_2 ( C\, ,\, ( N \times C ) \rtimes U_E )\, $ (because the map 
$\, K^J_2 ( J_{{\overline F}^{ _E}} \times J_E\, ,\, J_E \rtimes U_{{\overline F}^{ _E}} ) \rightarrowtail 
K^J_2 ( ( J_{{\overline F}^{ _E}} + U_{C , {\overline F}^{ _E}} ) \times J_E\, ,\, 
J_E \rtimes U_{{\overline F}^{ _E}} )\, $ is injective due to the fact that 
$\, K^J_3 ( C\, ,\, {\overline F}^{ _E} ) = 0\, $). But the map 
\smallskip
$$ K^J_2 ( C\, ,\, ( N \times C ) \rtimes U_E ) \rightarrowtail 
K^J_2 ( C \times J_E\, ,\, ( N \times C ) \rtimes U_E ) $$
\par\medskip\noindent
is injective by splitting of the quotient, and the map to 
$\, K^J_2 ( C \times J_E\, ,\, ( N \times C ) \rtimes U_E )\, $ factors over 
$\, K^J_2 ( ( J_{{\overline F}^{ _E}} + U_{C , {\overline F}^{ _E}} ) \times J_E\, ,\, 
J_E \rtimes U_{{\overline F}^{ _E}} )\, $. Considering the diagram
\smallskip
$$ \vbox{\halign{ #&#&#\cr
\hfil $\bigl( U_{C , {\overline F}^{ _E}}\, ,\, J_E \rtimes U_{{\overline F}^{ _E}} \bigr)$\hfil & 
\hfil $=$\hfil & \hfil  $\bigl( U_{C , {\overline F}^{ _E}}\, ,\, J_E \rtimes U_{{\overline F}^{ _E}} \bigr)$\hfil \cr
\hfil $\Bigm\downarrow $\hfil && \hfil $\Bigm\downarrow $\hfil \cr
\hfil $\bigl( J_{{\overline F}^{ _E}} + U_{C , {\overline F}^{ _E}}\, ,\, 
J_E \rtimes U_{{\overline F}^{ _E}} \bigr)$ \hfil & \hfil $\longrightarrow $\hfil & 
\hfil $\bigl( ( J_{{\overline F}^{ _E}} + U_{C , {\overline F}^{ _E}} ) \times J_E\, ,\, 
J_E \rtimes U_{{\overline F}^{ _E}} \bigr)$\hfil \cr
\hfil $\Bigm\downarrow $\hfil && \hfil $\Bigm\downarrow $\hfil \cr
\hfil $\bigl( J_F\, ,\, J_E \rtimes U_E \bigr)$\hfil & \hfil $\longrightarrow $\hfil &
\hfil $\bigl( J_F \times J_E\, ,\, J_E \rtimes U_F \bigr)$\hfil \cr }}  $$
\par\medskip\noindent
it now follows from a simple diagram chase that the lower horizontal map is also injective in 
$K^J_2$, and $\, ( N , F )\, $ is fully exact\qed
\par\bigskip\noindent
Another useful result concerning exactness is the following which has been used in the proof of the long exact sequence of $K^J_*$-theory for the higher $K^J_*$-groups ($\, n \geq 4\, $) and many of the previous results. Recall that a semisplit extension is an extension of normal pairs that splits normally in the first variable.
\par\bigskip\noindent
{\bf Lemma 22.}\quad ($P \times Q$-Lemma)\quad Assume given a semisplit extension 
\smallskip
$$ 1 \longrightarrow \bigl( Q , H \bigr) \longrightarrow \bigl( P \times Q , H \bigr) \longrightarrow 
\bigl( P , E \bigr) \longrightarrow 1 $$
\par\medskip\noindent
such that $\, ( Q , H )\, $ and $\, ( P , E )\, $ are exact. Then also $\, ( P \times Q\,  , H )\, $ is exact. If moreover $\, ( Q , H )\, $ is fully exact and $\, ( P , E )\, $ is exact of second order, then also 
$\, ( P \times Q\,  , H )\, $ is exact of second order (this result passes to arbitrary -- not necessarily semisplit-- extensions, see the proof below), and if $\, ( P \times Q\, , H )\, $ and $\, ( Q , H )\, $ are fully exact with $\, ( Q , H / P )\, $ exact of second order and $\, Q\, $ is free (or $\, K^J_2 ( Q ) = 0\, $), then also 
$\, ( P , H )\, $ is exact. If $\, P\, $ is free (resp. $\, K^J_2 ( P ) = 0\, $) and $\, ( P \times Q\,  , H )\, $ is exact of second order then also 
$\, ( P , H )\, $ is exact of second order.
\par\bigskip\noindent
{\it Proof.}\quad Assume that $\, ( P , E )\, $ and $\, ( Q , H )\, $ are exact. Consider the diagram 
\smallskip
$$ \vbox{\halign{ #&#&#\cr
\hfil $K^J_2\, ( J_{P , E}\, ,\, Q \rtimes U_E )$\hfil &\hfil $=$ \hfil & 
\hfil $K^J_2\, ( J_{P , E}\, ,\, Q \rtimes U_E )$\hfil \cr
\hfil $\Bigm\downarrow $\hfil && \hfil $\Bigm\downarrow $\hfil \cr
\hfil $K^J_2\, ( J_E\, ,\, Q \rtimes U_E )$\hfil & \hfil $\rightarrowtail $\hfil & 
\hfil $K^J_2\, ( J_E \times Q\, ,\, Q \rtimes U_E )$\hfil \cr
\hfil $\Bigm\downarrow $\hfil && \hfil $\Bigm\downarrow $\hfil \cr
\hfil $K^J_2\, ( J_{E / P}\, ,\, ( P \times Q ) \rtimes U_{E / P} )$\hfil & \hfil $\rightarrow $\hfil & 
\hfil $K^J_2\, ( J_{E / P} \times Q\, ,\, ( P \times Q ) \rtimes U_{E / P} )$\hfil \cr }}\> . $$
\par\medskip\noindent
Any element in the kernel of the lower horizontal map is represented by a commutator in 
$\, [\, U_{J_{E / P}}\, ,\, U_Q\, ]\, $, hence can be lifted to 
$\, K^J_2 ( J_E\, ,\, Q \rtimes U_E ) \rightarrowtail K^J_2 ( J_E \times Q\, ,\, Q \rtimes U_E )\, $. By halfexactness of the $K^J_2$-functor this lift must come from 
$\, K^J_2 ( J_{P , E}\, ,\, Q \rtimes U_E )\, $, then the element in 
$\, K^J_2 ( J_{E / P}\, ,\, ( P \times Q ) \rtimes U_{E / P} )\, $ is trivial, so the lower horizontal map is injective. Now consider the diagram 
\smallskip
$$ \vbox{\halign{ #&#&#\cr
\hfil $K^J_2\, ( Q\, ,\, ( P \times Q ) \rtimes U_{E / P} )$\hfil & \hfil $=$\hfil & 
\hfil $K^J_2\, ( Q\, ,\, ( P \times Q ) \rtimes U_{E / P} )$\hfil \cr
\hfil $\Bigm\downarrow $\hfil && \hfil $\Bigm\downarrow $\hfil \cr
\hfil $K^J_2\, ( J_{E / P} \times Q , ( P \times Q ) \rtimes U_{E / P} )$\hfil & \hfil $\rightarrow $\hfil & 
\hfil $K^J_2\, ( J_{E / P} \times P \times Q , ( P \times Q ) \rtimes U_{E / P} )$\hfil \cr
\hfil $\Bigm\downarrow $\hfil && \hfil $\Bigm\downarrow $\hfil \cr
\hfil $K^J_2\, ( J_{E / P}\, ,\, P \rtimes U_{E / P} )$\hfil & \hfil $\rightarrowtail $\hfil &
\hfil $K^J_2\, ( J_{E / P} \times P\, ,\, P \rtimes U_{E / P} )$\hfil \cr }}\> . $$
\par\medskip\noindent
An element in the kernel of the middle horizontal map lifts to 
$\, K^J_2 ( Q\, ,\, ( P \times Q ) \rtimes U_{E / P} )\, $. But the map 
\smallskip
$$ K^J_2\, ( Q\, ,\, ( P \times Q ) \rtimes U_{E / P} ) \rightarrowtail 
K^J_2\, ( J_{E / P} \times P \times Q\, ,\, ( P \times Q ) \rtimes U_{E / P} ) $$
\par\medskip\noindent
is injective as the corresponding extension splits. Thus 
\smallskip
$$ K^J_2\, ( J_{E / P}\, ,\, ( P \times Q ) \rtimes U_{E / P} ) \rightarrowtail 
K^J_2\, ( J_{E / P} \times P \times Q\, ,\, ( P \times Q ) \rtimes U_{E / P} ) $$
\par\medskip\noindent
is injective and $\, ( P \times Q\, , H )\, $ is exact.
\par\noindent
Next suppose that $\, ( Q , H )\, $ is fully exact and that $\, ( P , E )\, $ is exact of second order. We first show that $\, (\, Q \,, \,( P \times Q ) \rtimes U_{E / P} )\, $ is exact. Put 
$\, M = J_{P , E} \rtimes U_{J_{E / P}\, ,\, P \rtimes U_{E / P}}\, $ and
$\, L =J_{P , E} \rtimes ( J_{P \rtimes U_{E / P}} + U_{J_{E / P}\, ,\, P \rtimes U_{E / P}} )\, $ which are to denote the obvious normal subgroups of the quotient
$\, U_{Q \rtimes U_E} / ( J^2_{P , E} + J_{Q , Q \rtimes U_E} )\, $.  From the diagram 
\smallskip
$$ \vbox{\halign{ #&#&#\cr
\hfil $K^J_2\, ( M\, ,\, ( J_{P , E} \times Q ) \rtimes U_{P \rtimes U_{E / P}} )$\hfil &\hfil $=$\hfil & 
\hfil $K^J_2\, ( M\, ,\, ( J_{P , E} \times Q ) \rtimes U_{P \rtimes U_{E / P}} )$\hfil \cr
\hfil $\Bigm\downarrow $\hfil && \hfil $\Bigm\downarrow $\hfil \cr
\hfil $K^J_2\, ( L\, ,\, ( J_{P , E} \times Q ) \rtimes U_{P \rtimes U_{E / P}} )$\hfil & \hfil $\rightarrow $\hfil & 
\hfil $K^J_2\, ( Q  \times L\, ,\, ( J_{P , E} \times Q ) \rtimes U_{P \rtimes U_{E / P}} )$\hfil \cr
\hfil  $\Bigm\downarrow\Bigm\uparrow $\hfil && \hfil $\Bigm\downarrow\Bigm\uparrow $\hfil \cr
\hfil $K^J_2\, ( J_E\, ,\, Q \rtimes U_E )$\hfil & \hfil $\rightarrowtail $\hfil &
\hfil $K^J_2\, ( J_E \times Q\, ,\, Q \rtimes U_E )$\hfil \cr }} $$
\par\medskip\noindent
one induces that the middle horizontal map is injective. From the diagram 
\smallskip
$$\vbox{\halign{ #&#&#\cr
\hfil $( J_{P \rtimes U_{E / P}}\, ,\, ( J_{P , E} \times Q ) \rtimes U_{P \rtimes U_{E / P}} )$\hfil & 
\hfil $=$\hfil & 
\hfil $( J_{P \rtimes U_{E / P}}\, ,\, ( J_{P , E} \times Q ) \rtimes U_{P \rtimes U_{E / P}} )$\hfil \cr 
\hfil $\Bigm\downarrow $\hfil && \hfil $\Bigm\downarrow $\hfil \cr
\hfil $( L\, ,\, ( J_{P , E} \times Q ) \rtimes U_{P \rtimes U_{E / P}} )$\hfil & \hfil $\rightarrow $\hfil & 
\hfil $( Q \times L\, ,\, ( J_{P , E} \times Q ) \rtimes U_{P \rtimes U_{E / P}} )$\hfil \cr
\hfil $\Bigm\downarrow\Bigm\uparrow $\hfil && \hfil $\Bigm\downarrow\Bigm\uparrow $\hfil \cr 
\hfil $( J_E\, ,\, Q \rtimes U_E )$\hfil & \hfil $\rightarrow $\hfil & 
\hfil $( J_E \times Q\, ,\, Q \rtimes U_E )$\hfil \cr }} $$
\par\medskip\noindent
one gets injectivity of the map 
\smallskip
$$K^J_2 ( J_{P \rtimes U_{E / P}}\, ,\, ( J_{P , E} \times Q ) \rtimes U_{P \rtimes U_{E / P}} ) \rightarrowtail \qquad\qquad\qquad\qquad $$ 
$$\qquad\qquad\qquad\qquad\qquad 
K^J_2 ( J_{P \rtimes U_{E / P}} \times Q\, ,\, ( J_{P , E} \times Q ) \rtimes U_{P \rtimes U_{E / P}} )\> . $$
\par\medskip\noindent
Then there is a compatible splitting
\smallskip
$$ \vbox{\halign{ #&#&#\cr
\hfil $( J_{P \rtimes U_{E / P}}\, ,\, Q \rtimes U_{P \rtimes U_{E / P}} )$\hfil & \hfil $\rightarrow $\hfil &
\hfil $( J_{P \rtimes U_{E / P}} \times Q\, ,\, Q \rtimes U_{P \rtimes U_{E / P}} )$\hfil \cr
\hfil $\Bigm\downarrow\Bigm\uparrow $\hfil && \hfil $\Bigm\downarrow\Bigm\uparrow $\hfil \cr
\hfil $( J_{P \rtimes U_{E / P}} , ( J_{P , E}\negthinspace\times\negthinspace Q )\negthinspace\rtimes\negthinspace 
U_{P \rtimes U_{E / P}} )$\hfil & 
\hfil $\rightarrow $\hfil & 
\hfil $( J_{P \rtimes U_{E / P}} \times Q , ( J_{P , E}\negthinspace\times\negthinspace Q ) \negthinspace\rtimes\negthinspace 
U_{P \rtimes U_{E / P}} )$\hfil \cr }} $$
\par\medskip\noindent
showing that the upper horizontal map is injective in $K^J_2$. This implies that 
$\, (\, Q\, , ( P \times Q ) \rtimes U_{E / P} )\, $ is exact. Then putting 
$\, \widetilde U =  ( J_{P , E} \times Q ) \rtimes U_{P \rtimes U_{E / P}}\, $ consider the diagram
\smallskip
$$ \vbox{\halign{ #&#&#\cr
\hfil $K^J_2\, ( J_{P , E}\, ,\, \widetilde U\, )$\hfil & \hfil $=$\hfil & 
\hfil $K^J_2\, ( J_{P , E}\, ,\, \widetilde U\, )$\hfil \cr
\hfil $\Bigm\downarrow $\hfil && \hfil $\Bigm\downarrow $\hfil \cr
\hfil $K^J_2\, ( J_{P \rtimes U_{E / P}} \times J_{P , E}\, ,\, \widetilde U\, )$\hfil & \hfil $\rightarrow $\hfil & 
\hfil $K^J_2\, ( J_{P \rtimes U_{E / P}} \times J_{P , E} \times Q\, ,\, \widetilde U\, )$\hfil \cr
\hfil $\Bigm\downarrow\Bigm\uparrow $\hfil && \hfil $\Bigm\downarrow\Bigm\uparrow $\hfil \cr
\hfil $K^J_2\, ( J_{P \rtimes U_{E / P}}\, ,\, Q \rtimes U_{P \rtimes U_{E / P}}\, )$\hfil & \hfil $\rightarrowtail $\hfil & \hfil $K^J_2\, ( J_{P \rtimes U_{E / P}} \times Q\, ,\, Q \rtimes U_{P \rtimes U_{E / P}}\, )$\hfil \cr }} $$
\par\medskip\noindent
showing that the middle horizontal map is also injective. Then the diagram 
\smallskip
$$ \vbox{\halign{ #&#&#\cr
\hfil $K^J_2\, ( J_{P \rtimes U_{E / P}}\, ,\, \widetilde U\, )$\hfil & \hfil $=$\hfil & 
\hfil $K^J_2\, ( J_{P \rtimes U_{E / P}}\, ,\, \widetilde U\, )$\hfil \cr
\hfil $\Bigm\downarrow $\hfil && \hfil $\Bigm\downarrow $\hfil \cr
\hfil $K^J_2\, ( J_{P \rtimes U_{E / P}} \times J_{P , E}\, ,\, \widetilde U\, )$\hfil & \hfil $\rightarrowtail $\hfil & 
\hfil $K^J_2\, ( J_{P \rtimes U_{E / P}} \times J_{P , E} \times Q\, ,\, \widetilde U\, )$\hfil \cr
\hfil $\Bigm\downarrow $\hfil && \hfil $\Bigm\downarrow $\hfil \cr
\hfil $K^J_2\, ( J_{P , E}\, ,\, Q \rtimes U_E\, )$\hfil & \hfil $\rightarrow $\hfil & 
\hfil $K^J_2\, ( J_{P , E} \times Q\, ,\, Q \rtimes U_E\, )$\hfil \cr }} $$
\par\medskip\noindent
shows that 
\smallskip
$$ K^J_2\, ( J_{P , E}\, ,\, Q \rtimes U_E\, ) \rightarrowtail 
K^J_2\, ( J_{P , E} \times Q\, ,\, Q \rtimes U_E\, ) $$
\par\medskip\noindent
is injective. Now suppose given an element in the kernel of the map
\smallskip
$$ K^J_2\, ( J_H\, ,\, J_{E / P} \rtimes U_H ) \rightarrow 
K^J_2\, ( J_H \times J_{E / P}\, ,\, J_{E / P} \rtimes U_H )\> . $$
\par\medskip\noindent
One notes that injectivity of 
\smallskip
$$ K^J_2\, ( J_E\, ,\, J_E \rtimes U_E ) \rightarrowtail K^J_2\, ( J_E \times J^0_E\, ,\, J_E \rtimes U_E ) $$
\par\medskip\noindent
and 
\smallskip
$$ K^J_2\, ( J_E\, ,\, J_{E / P} \rtimes U_E ) \rightarrowtail 
K^J_2\, ( J_E \times J_{E / P}\, ,\, J_{E / P} \rtimes U_E ) $$
\par\medskip\noindent
entails injectivity of 
\smallskip
$$ K^J_2\, ( J_H + U_{Q , H}\, ,\, J_E \rtimes U_H ) \rightarrowtail 
K^J_2\, ( ( J_H + U_{Q , H} ) \times J_E\, ,\, J_E \rtimes U_H ) $$
\par\medskip\noindent
and of 
\smallskip
$$ K^J_2\, ( J_H + U_{Q , H}\, ,\, J_{E / P} \rtimes U_H ) \rightarrowtail 
K^J_2\, ( ( J_H + U_{Q , H} ) \times J_{E / P}\, ,\, J_{E / P} \rtimes U_H ) \> . $$
\par\medskip\noindent
Also the maps 
\smallskip
$$ K^J_2\, ( U_{Q , H}\, ,\, J_E \rtimes U_H ) \rightarrowtail 
K^J_2\, ( ( J_H + U_{Q , H} ) \times J_E\, ,\, J_E \rtimes U_H )\> ,  $$
$$ K^J_2\, ( U_{Q , H}\, ,\, J_{E / P} \rtimes U_H ) \rightarrowtail 
K^J_2\, ( ( J_H + U_{Q , H} ) \times J_{E / P}\, ,\, J_{E / P} \rtimes U_H ) $$
\par\medskip\noindent
are injective by existence of a splitting. Then our element lifts to an element of 
$\, K^J_2 ( J_H\, ,\, J_E \rtimes U_H )\, $ mapping to the kernel of 
\smallskip
$$ K^J_2\, ( J_H + U_{Q , H}\, ,\, J_E \rtimes U_H ) \twoheadrightarrow 
K^J_2\, ( J_H + U_{Q , H}\, ,\, J_{E / P} \rtimes U_H ) $$
\par\medskip\noindent
which is represented by an element of $\, [\, U_{J_{P , E}}\, , U_{J_H + U_{Q , H}}\, ] \equiv 
[\, U_{J_{P , E}}\, , U_{J_H}\, ]$ $[\, U_{J_{P , E}}\, ,\, U_{U_{Q , H}}\, ]\, $. Modifying our lift modulo the image of the first bracket in $\, K^J_2 ( J_H\, ,\, J_E \rtimes U_H )\, $ we may assume that the lift is in the image of $\, K^J_2 ( U_{Q , H}\, ,\, J_E \rtimes U_H )\, $ in 
$\, K^J_2 ( J_H + U_{Q , H}\, ,\, J_E \rtimes U_H )\, $ and more precisely comes from the kernel of 
\smallskip
$$ K^J_2\, ( U_{Q , H}\, ,\, J_E \rtimes U_H ) \twoheadrightarrow 
K^J_2\, ( U_{Q , H}\, ,\, J_{E / P} \rtimes U_H ) $$
\par\medskip\noindent
which can be lifted to the kernel of 
\smallskip
$$ K^J_2\, ( J_{P , E}\, ,\, J_E \rtimes U_H ) \twoheadrightarrow 
K^J_2\, ( J_{P , E}\, ,\, J_E \rtimes U_E ) \> . $$ 
\par\medskip\noindent
Also in $\, K^J_2 ( ( J_H + U_{Q , H} ) \times J_{P , E}\, ,\, J_E \rtimes U_H )\, $ the element becomes trivial upon dividing by $\, J_H\, $, so its image in 
$\, K^J_2 ( J_H \times J_{P , E}\, ,\, J_E \rtimes U_H )\, $ is already in the kernel of the map 
\smallskip
$$ K^J_2\, ( J_H \times J_{P , E}\, ,\, J_E \rtimes U_H ) \rightarrow 
K^J_2\, ( J_{P , E}\, ,\, Q \rtimes U_E ) $$
\par\medskip\noindent
(by injectivity of $\, K^J_2 ( J_{P , E}\, ,\, Q \rtimes U_E ) \rightarrowtail 
K^J_2 ( J_{P , E} \times Q\, ,\, Q \rtimes U_E )\, $) and can be represented by an element of 
$\, [\, U_{J_{P , E}}\, ,\, U_{J_H}\, ]\, $. Again modifying the lift modulo this bracket one gets a lift which becomes trivial in $\, K^J_2 ( J_H + U_{Q , H}\, ,\, J_E \rtimes U_H )\, $, hence lifts to 
$\, K^J_3 ( Q\, ,\, Q \rtimes U_E )\, $. Then since the image of our element in $\, K^J_2 ( J_H\, ,\, U_H ) = 
K^J_3 ( H )\, $ is trivial and considering the extension 
\smallskip
$$ 1 \longrightarrow ( J_H \times J_E\, ,\, J_E \rtimes U_H ) \longrightarrow 
( J_E \rtimes U_H\, ,\, J_E \rtimes U_H ) \longrightarrow ( H , H ) \longrightarrow 1 $$
\par\medskip\noindent
one finds that the image of our lift in $\, K^J_2 ( J_H \times J_E\, ,\, J_E \rtimes U_H )\, $ is trivial. But the map 
\smallskip
$$ K^J_2 ( J_H\, ,\, J_E \rtimes U_H ) \rightarrowtail K^J_2 ( J_H \times J_E\, ,\, J_E \rtimes U_H ) $$
\par\medskip\noindent
is injective by assumption that $\, ( Q , H )\, $ is exact of second order. Thus 
\smallskip
$$ K^J_2\, ( J_H\, ,\, J_{E / P} \rtimes U_H ) \rightarrowtail 
K^J_2\, ( J_H \times J_{E / P}\, ,\, J_{E / P} \rtimes U_H )  $$
\par\medskip\noindent
is injective and $\, (\, P \times Q\, , H\, )\, $ is exact of second order. One notes that we did not use the assumption that the extension is semisplit, but we have retained the familiar $P \times Q$-notation to comply with the rest of the proof.
\par\noindent
Now assume that $\, (\, P \times Q\, , H )\, $ and $\, ( Q , H )\, $ are fully exact and that $\, ( Q , H / P )\, $ is exact of second order and $\, K^J_2 ( Q ) = 0\, $. An element in the kernel of 
\smallskip
$$ K^J_2\, J_{H / P}\, ,\, P \rtimes U_{H / P} ) \rightarrow K^J_2 ( J_{H / P} \times P\, ,\, P \rtimes U_{H / P} )$$
\par\medskip\noindent
is also in the kernel of the map to 
$\, K^J_2 ( ( J_{H / P} + U_{Q , H / P} ) \times P\, ,\, P \rtimes U_{H / P} )\, $ so considering the diagram
\smallskip
$$ \vbox{\halign{ #&#&#&#&#&#&#&#&#\cr
\hfil $1$\hfil &\hfil $\rightarrow $\hfil & \hfil $( J_{P , H} , U_H )$\hfil & \hfil $\rightarrow $\hfil & 
\hfil $(\, J_H\, ,\, U_H\, )$\hfil & \hfil $\rightarrow $\hfil & \hfil $(\, J_{H / P}\, ,\, P \rtimes U_{H / P}\, )$\hfil & 
\hfil $\rightarrow $\hfil & \hfil $1$\hfil \cr 
&& \hfil $\Bigm\Vert $\hfil && \hfil $\Bigm\downarrow $\hfil && \hfil $\Bigm\downarrow $\hfil &&\cr
\hfil $1$\hfil & \hfil $\rightarrow $\hfil & \hfil $( J_{P , H} , U_H )$\hfil & \hfil $\rightarrow $\hfil &
\hfil $( J_H\negthinspace +\negthinspace U_{P\negthinspace\times\negthinspace Q , H} , U_H )$\hfil & \hfil $\rightarrow $\hfil & 
\hfil $(\negthinspace (J_{H / P}\negthinspace +\negthinspace U_{Q , H / P})\negthinspace\times\negthinspace\negthinspace P , 
P \negthinspace\rtimes\negthinspace U_{H / P} )$\hfil & 
\hfil $\rightarrow $\hfil & \hfil $1$\hfil \cr  }} $$
\par\medskip\noindent
its image in $\, K_1 ( J_{P , H}\, ,\, U_H )\, $ must vanish, so it lifts to $\, K^J_2 ( J_H\, ,\, U_H )\, $. If the image of this lift in $\, K^J_2 ( J_H + U_{Q , H}\, ,\, U_H )\, $ would be nontrivial it would lift to 
$\, K^J_2 ( J_{P , H}\, ,\, U_H )\, $ and be zero in $\, K^J_2 ( J_{H / P}\, ,\, P \rtimes U_{H / P} )\, $. Thus considering the extension
\smallskip
$$ 1 \longrightarrow ( J_H\, ,\, U_H ) \longrightarrow ( J_H + U_{P\times Q , H}\, ,\, U_H ) \longrightarrow 
( P\times Q\, ,\, H ) \longrightarrow 1 $$
\par\medskip\noindent
with $\, ( P\times Q\, ,\, H )\, $ exact, the element in $\, K^J_2\, ( J_H\, ,\, U_H )\, $ lifts to the group
$\, K^J_3 ( P\times Q\, ,\, H )\, $. One can assume that its image in $\, K^J_3 ( Q\, ,\, H / P )\, $ is trivial by comparing the extension above with
\smallskip
$$ 1 \longrightarrow ( J_{H / P}\, ,\, U_{H / P} ) \longrightarrow ( J_{H / P} + U_{Q , H / P}\, ,\, U_{H / P} ) 
\longrightarrow ( Q\, ,\, H / P ) \longrightarrow 1\> . $$
\par\medskip\noindent
Certainly its image in $\, K^J_2 ( J_{H / P}\, ,\, U_{H / P} )\, $ is trivial, so it lifts to the group
$\, K^J_3 ( J_{H / P} + U_{Q , H / P}\, ,\, U_{H / P} ) \simeq K^J_3 ( J_H + U_{P\times Q , H}\, ,\, U_H )\, $ as $\, ( Q\, ,\, H / P )\, $ is exact of second order, and the image of the latter in $\, K^J_2 ( J_H\, ,\, U_H )\, $ is trivial. Now there is a split surjection $\, K^J_3 ( P\times Q\, ,\, H ) \twoheadrightarrow K^J_3 ( P , E )\, $ and the boundary map 
\smallskip
$$ K^J_3 ( P , H ) \longrightarrow K^J_3 ( P \times Q\, ,\, H ) \buildrel\delta\over\longrightarrow
K^J_2 ( J_H\, ,\, U_H ) $$
\par\medskip\noindent
is easily seen to map into the image of the natural map 
\smallskip
$$ K^J_3 ( P , H ) \simeq K^J_2 ( J_{P , H}\, ,\, U_H ) \longrightarrow K^J_2 ( J_H\, ,\, U_H ) $$
\par\medskip\noindent
so the image of any element from $\, K^J_3 ( P , H )\, $ in 
$\, K^J_2 ( J_{H / P}\, ,\, P \rtimes U_{H / P} )\, $ will be trivial, then by taking the difference of the lift in 
$\, K^J_3 ( P\times Q\, ,\, H )\, $ with a lift of its image in $\, K^J_3 ( P , E )\, $ to $\, K^J_3 ( P , H )\, $ we may assume that the image of the lift in $\, K^J_3 ( P , E )\, $ is trivial and lifts to $\, K^J_3 ( Q , H )\, $. This element may be further lifted to $\, K^J_3 ( J_E \times Q\, ,\, Q \rtimes U_E )\, $. Its image in 
$\, K^J_3 ( J_{E / P} \times Q\, ,\, Q \rtimes U_{E / P} )\, $ lifts to 
$\, K^J_3 ( J_{E / P}\, ,\, Q \rtimes U_{E / P} )\, $ (check that 
$\, ( J_{E / P} \times Q\, ,\, Q \rtimes U_{E / P} )\, $ is exact of second order and 
$\, ( J_{E / P}\, ,\, Q \rtimes U_{E / P} )\, $ is fully exact by the results above). Then 
\smallskip
$$ K^J_3 ( J_{E / P}\, ,\, Q \rtimes U_{E / P} ) \simeq K^J_3 ( J_{E / P}\, ,\, U_{E / P} ) \simeq 
K^J_3 ( J_{E / P}\, ,\, P \rtimes U_{E / P} ) $$
\par\medskip\noindent
and since the image of this element in $\, K^J_2 ( J_{P , E}\, ,\, U_E ) \simeq K^J_3 ( P , E )\, $ is trivial it lifts to $\, K^J_3 ( J_E\, ,\, U_E ) \simeq K^J_3 ( J_E\, ,\, Q \rtimes U_E )\, $ and by taking the difference of the original lift with this element we can assume that the lift in 
$\, K^J_3 ( J_E \times Q\, ,\, Q \rtimes U_E )\, $ is zero in 
$\, K^J_3 ( J_{E / P} \times Q\, ,\, Q \rtimes U_{E / P} )\, $. Then 
$\, K^J_3 ( J_{E / P} \times Q\, ,\, Q \rtimes U_{E / P} )\, $ is isomorphic to 
$\, K^J_3 ( (J_E + U_{P , E}) \times Q\, ,\, Q \rtimes U_E )\, $ by the extension 
\smallskip
$$ 1\negthinspace\negthinspace\rightarrow\negthinspace\negthinspace ( U_{P , E} , Q \rtimes U_E ) \negthinspace\rightarrow\negthinspace ((J_E + U_{P , E}) \times Q , Q \rtimes U_E ) \negthinspace\rightarrow\negthinspace 
( J_{E / P}\negthinspace\times\negthinspace Q , Q \rtimes U_{E / P} ) \negthinspace\rightarrow\negthinspace\negthinspace 1 $$
\par\medskip\noindent
and $\, K^J_3 ( U_{P , E}\, ,\, Q \rtimes U_E ) \simeq K^J_3 ( U_{P , E}\, ,\, U_E ) = 0\, $ provided that we can show that $\, ( U_{P , E}\, ,\, Q \rtimes U_E )\, $ is exact since 
$\, ( ( J_E + U_{P , E} ) \times Q\, ,\, Q \rtimes U_E )\, $ is second order exact from the fact that 
$\, K^J_2 ( Q ) = 0\, $ together with the results above. 
Putting 
$\, L = J_{Q \rtimes U_{E / P}} + U_{Q\, ,\, Q \rtimes U_{E / P}}\, $ consider the diagram
\smallskip
$$ \vbox{\halign{ #&#&#\cr
\hfil $\bigl(\, J_{Q \rtimes U_{E / P}}\, ,\, U_{P , E} \rtimes U_{Q \rtimes U_{E / P}} \,\bigr)$\hfil & \hfil $\longrightarrow $\hfil & 
\hfil $\bigl(\, J_{Q \rtimes U_{E / P}} \times U_{P , E}\, ,\, U_{P , E} \rtimes U_{Q \rtimes U_{E / P}} \,\bigr)$\hfil \cr
\hfil $\Bigm\downarrow $\hfil && \hfil $\Bigm\downarrow $\hfil \cr
\hfil 
$\bigl(\, L\, ,\, U_{P , E} \rtimes U_{Q \rtimes U_{E / P}} \,\bigr)$\hfil & \hfil $\longrightarrow $\hfil & \hfil 
$\bigl(\, L \times U_{P , E}\, ,\, U_{P , E} \rtimes U_{Q \rtimes U_{E / P}} \,\bigr)$\hfil \cr
\hfil $\Bigm\downarrow $ \hfil && \hfil $\Bigm\downarrow $\hfil \cr
\hfil $\bigl(\, Q\, ,\, Q \rtimes U_E \,\bigr)$\hfil & \hfil $\longrightarrow $\hfil & 
\hfil $\bigl(\, Q\, ,\, Q \rtimes U_{E / P} \,\bigr)$\hfil \cr }}\> . $$
\par\medskip\noindent
Both pairs in the bottom row are fully exact as well as the pair 
$\, ( L \times U_{P , E}\, ,\, U_{P , E} \rtimes U_{Q \rtimes U_{E / P}} )\, $ and since 
$\, U_{P , E} \rtimes U_{Q \rtimes U_{E / P}} = L \rtimes U_E\, $ the same is true for the pair 
$\, ( L\, ,\, U_{P , E} \rtimes U_{Q \rtimes U_{E / P}} )\, $. The map 
\smallskip
$$ K^J_2 ( L\, ,\, U_{P , E} \rtimes U_{Q \rtimes U_{E / P}} ) \rightarrowtail 
K^J_2 ( L \times U_{P , E}\, ,\, U_{P , E} \rtimes U_{Q \rtimes U_{E / P}} ) $$
\par\medskip\noindent
is obviously injective. Then a diagram chase gives injectivity of the upper horizontal map in $K^J_2$ by the isomorphism $\, K^J_3 ( Q\, ,\, Q \rtimes U_E ) \simeq K^J_3 ( Q\, ,\, Q \rtimes U_{E / P} )\, $. Then the image of our lift is trivial in the group $\, K^J_3 ( (J_E + U_{P , E}) \times Q\, ,\, Q \rtimes U_E )\, $ and lifts to $\, K^J_4 (P , E )\, $. Obviously the image of $\, K^J_4 ( P , E )\, $ in $\, K^J_3 ( P\times Q\, ,\, H )\, $ is trivial proving exactness of $\, ( P , H )\, $.
\par\noindent
Now let $\, K^J_2 ( P ) = 0\, $ and $\, ( P \times Q\, ,\, H )\, $ exact of second order. Consider the diagram 
\smallskip
$$ \vbox{\halign{ #&#&#\cr
\hfil $\bigl(\, J_H\, ,\, J_{H / P} \rtimes U_H\,\bigr)$\hfil & \hfil $\longrightarrow $\hfil & 
\hfil $\bigl(\, J_{Q\, ,\, H / P} \times (J_H  + U_{Q , H})\, ,\, J_{H / P} \rtimes U_H\,\bigr)$\hfil \cr
\hfil $\Bigm\downarrow $\hfil && \hfil $\Bigm\downarrow $\hfil \cr
\hfil $\bigl(\, J_H + U_{P , H }\, ,\, J_{H / P} \rtimes U_H\,\bigr)$\hfil & \hfil $\longrightarrow $\hfil & 
\hfil $\bigl(\, J_{Q , H / P} \rtimes (J_H + U_{P\times Q , H})\, ,\, J_{H / P} \rtimes U_H\,\bigr)$\hfil \cr
\hfil $\Bigm\downarrow $\hfil && \hfil $\Bigm\downarrow $\hfil \cr
\hfil $\bigl(\, P\, ,\, P \rtimes U_{H / P}\,\bigr)$\hfil &\hfil $\longrightarrow $\hfil &
\hfil $\bigl(\, P\, ,\, P \rtimes U_{E / P}\,\bigr)$\hfil \cr }} . $$
\par\medskip\noindent
Then the middle horizontal map is injective for $K^J_2$ by existence of a splitting. An element in the kernel of
\smallskip
$$ K^J_2 ( J_H\, ,\, J_{H / P} \rtimes U_H ) \rightarrow 
K^J_2 ( J_{Q\, ,\, H / P} \times (J_H + U_{Q , H})\, ,\, J_{H / P} \rtimes U_H )   $$
\par\medskip\noindent
lifts to $\, K^J_3 ( P\, ,\, P \rtimes U_{H / P} ) \simeq K^J_3 ( P\, ,\, P \rtimes U_{E / P} ) \simeq 
K^J_3 ( P \rtimes U_{E / P} )\, $. The first isomorphism follows since $\, K^J_2 ( P ) = 0\, $ and the second by the fact that $\, ( P\, ,\, P \rtimes U_{E / P} )\, $ is fully exact. A simple diagram chase gives a lift which is trivial so that the map 
\smallskip
$$ K^J_2 ( J_H\, ,\, J_{H / P} \rtimes U_H ) \rightarrowtail 
K^J_2 ( J_H \times J_{Q\, ,\, H / P}\, ,\, J_{H / P} \rtimes U_H ) $$
\par\medskip\noindent
which factors the upper horizontal map is injective. Now consider the diagram 
\smallskip
$$ \vbox{\halign{ #&#&#\cr
\hfil $\bigl(\, J_{Q\, ,\, H / P}\, ,\, J_{H / P} \rtimes U_H \,\bigr)$\hfil & \hfil $=$\hfil & 
\hfil $\bigl(\, J_{Q\, ,\, H / P}\, ,\, J_{H / P} \rtimes U_H \,\bigr)$\hfil \cr
\hfil $\Bigm\downarrow $\hfil && \hfil $\Bigm\downarrow $\hfil \cr
\hfil $\bigl(\, J_H \times J_{Q\, ,\, H / P}\, ,\, J_{H / P} \rtimes U_H\,\bigr)$\hfil & \hfil $\longrightarrow $\hfil & 
\hfil $\bigl(\, J_H \times J_{H / P}\, ,\, J_{H / P} \rtimes U_H \,\bigr)$\hfil \cr
\hfil $\Bigm\downarrow $\hfil  && \hfil $\Bigm\downarrow $\hfil \cr
\hfil $\bigl(\, J_H\, ,\, J_{E / P} \rtimes U_H\,\bigr)$\hfil & \hfil $\longrightarrow $\hfil & 
\hfil $\bigl(\, J_H \times J_{E / P}\, ,\, J_{E / P} \rtimes U_H\,\bigr)$\hfil \cr }} . $$
\par\medskip\noindent
Since the map 
\smallskip
$$ K^J_2 ( J_H\, ,\, J_{E / P} \rtimes U_H ) \rightarrowtail 
K^J_2 ( J_H \times J_{E / P}\, ,\, J_{E / P} \rtimes U_H ) $$
\par\medskip\noindent
is injective and $\, ( J_H \times J_{E / P}\, ,\, J_{E / P} \rtimes U_H )\, $ is fully exact an element in the kernel of the map 
\smallskip
$$ K^J_2 ( J_H \times J_{Q\, ,\, H / P}\, ,\, J_{H / P} \rtimes U_H ) \rightarrow 
K^J_2 ( J_H \times J_{H / P}\, ,\, J_{H / P} \rtimes U_H ) $$
\par\medskip\noindent
lifts to $\, K^J_2 ( J_{Q\, ,\, H / P}\, ,\, J_{H / P} \rtimes U_H )\, $ then to 
$\, K^J_3 ( J_H \times J_{E / P}\, ,\, J_{E / P} \rtimes U_H )\, $. Suppose that the original element in 
$\, K^J_2 ( J_H \times J_{Q\, ,\, H / P}\, ,\, J_{H / P} \rtimes U_H )\, $ came from 
$\, K^J_2 ( J_H\, ,\, J_{H / P} \rtimes U_H )\, $. Then the lifted element in the group
$\, K^J_2 ( J_{Q\, ,\, H / P}\, ,\, J_{H / P} \rtimes U_H )\, $ is trivial in 
$\, K^J_2 ( J_{Q\, ,\, H / P}\, ,\, P \rtimes U_{H / P} )\, $ so the image of the lift in 
$\, K^J_3 ( J_H \times J_{E / P}\, ,\, J_{E / P} \rtimes U_H )\, $ by the composition
\smallskip
$$ K^J_3 ( J_H\negthinspace\times\negthinspace J_{E / P} , J_{E / P} \negthinspace\rtimes\negthinspace U_H )\negthinspace\twoheadrightarrow\negthinspace 
K^J_3 ( J_{E / P} , (P\times Q)\negthinspace\rtimes\negthinspace U_{E / P} )\negthinspace\rightarrow\negthinspace K^J_3 ( (P\times Q)\negthinspace\rtimes\negthinspace U_{E / P} ) $$
\par\medskip\noindent
lifts to $\, K^J_3 ( P \rtimes U_{H / P} ) \simeq K^J_3 ( P\, ,\, P \rtimes U_{H / P} ) \simeq 
K^J_3 ( P\, ,\, (P\times Q) \rtimes U_{E / P} )\, $. But the intersection of the images of 
$\, K^J_3 ( P\, ,\, (P\times Q) \rtimes U_{E / P} )\, $ and 
$\, K^J_3 ( J_{E / P}\, ,\, (P\times Q) \rtimes U_{E / P} )\, $ in $\, K^J_3 ( (P\times Q) \rtimes U_{E / P} )\, $ is obviously trivial, so that the image of the lift in $\, K^J_3 ( (P\times Q) \rtimes U_{E / P} )\, $ by the composition as above must be trivial. Then the image in 
$\, K^J_3 ( J_{E / P}\, ,\, (P\times Q) \rtimes U_{E / P} )\, $ lifts to 
$\, K^J_4 ( H ) = K^J_3 ( J_H\, ,\, U_H ) \simeq K^J_3 ( J_H\, ,\, J_{E / P} \rtimes U_H ) \, $ the image of which is trivial in $\, K^J_2 ( J_H \times J_{Q\, ,\, H / P}\, ,\, J_{H / P} \rtimes U_H )\, $ so by taking the difference of the original lift with this element one obtains a lift in 
$\, K^J_3 ( J_H \times J_{E / P}\, ,\, J_{E / P} \rtimes U_H )\, $ trivial in 
$\, K^J_3 ( J_{E / P}\, ,\, (P\times Q) \rtimes U_{E / P} )\, $. This means that it lifts to 
$\, K^J_3 ( J_H\, ,\, J_{E / P} \rtimes U_H )\, $ showing that the original element in  
$\, K^J_3 ( J_H \times J_{Q\, ,\, H / P}\, ,\, J_{H / P} \rtimes U_H )\, $ was trivial and the map 
\smallskip
$$ K^J_2 ( J_H\, ,\, J_{H / P} \rtimes U_H ) \rightarrowtail 
K^J_2 ( J_H \times J_{H / P}\, ,\, J_{H / P} \rtimes U_H ) $$
\par\medskip\noindent
must be injective\qed  
\par\bigskip\noindent
Now we are ready to prove the main theorem of this section. We restrict ourselves to the case of 
(full subgroups of) perfect groups which provide the principal examples for the regular theory. Lemma 21 shows that a regular pair in this narrow sense is automatically fully exact, so the only problem is to show that the exactness results of Corollary 5.1 pass to the regular $K_*$-groups. It is conceivable that the result generalizes to abstract regular pairs (which should match the requirements of Proposition 1 of course).
\par\bigskip\noindent
{\bf Theorem 6.}\quad (Long Exact Sequence of Regular $K$-Theory)\quad Assume that $\, F\, $ is perfect and given an extension 
$$ 1 \longrightarrow \bigl( N , F \bigr) \longrightarrow \bigl( M , F \bigr) \longrightarrow 
\bigl( P , E \bigr) \longrightarrow 1\> . $$
\par\medskip\noindent
There is associated the following long exact sequence of $K_*$-groups
\medskip
$$ \qquad \cdots \buildrel {p_*}\over\longrightarrow K_{n+1}\, \bigl( P , E \bigr) 
\buildrel {{\delta }_*}\over\longrightarrow K_n\, \bigl( N , F \bigr) \buildrel {i_*}\over\longrightarrow 
K_n\, \bigl( M , F \bigr) \buildrel {p_*}\over\longrightarrow K_n\, \bigl( P , E \bigr) $$
$$ \buildrel {{\delta }_*}\over\longrightarrow \cdots \>\>\>\qquad\qquad\qquad \cdots 
\buildrel {{\delta }_*}\over\longrightarrow 
K_2\, \bigl( N , F \bigr) \buildrel {i_*}\over\longrightarrow  K_2\, \bigl( M , F \bigr) 
\buildrel {p_*}\over\longrightarrow K_2\, \bigl( P , E \bigr) $$
$$\buildrel {{\delta }_*}\over\longrightarrow K_1\, \bigl( N , F \bigr) \buildrel {i_*}\over\longrightarrow 
K_1\, \bigl( M , F \bigr) \buildrel {p_*}\over\longrightarrow K_1\,\bigl( P , E \bigr)\> .\qquad\qquad\qquad\qquad\qquad\quad $$
\par\bigskip\noindent
{\it Proof.}\quad We first reduce to the case where $\, N\, ,\, M\, $ and $\, P\, $ are full subgroups of $\, F\, $ and $\, E\, $ respectively. By definition of the regular $K_*$-groups one can assume that $\, M\, $ and 
$\, P\, $ are full subgroups, replacing $\, N\, $ with $\, N \cap [ M , F ]\, $, $\, M\, $ with $\, [ M , F ]\, $ and $\, P\, $ with $\, [ P , E ]\, $ if necessary. Then consider $\, \overline P = M / [\, N , F\, ] \subseteq \overline E = 
F / [\, N , F\, ]\, $. $\,\overline P\, $ is a full subgroup of $\,\overline E\, $ as a quotient of $\, M\, $ and we must show that $\, K_n ( \overline P , \overline E ) = K_n ( P . E )\, $ for $\, n \geq 3\, $. It is clear that 
$\, K^J_n ( \overline P , \overline E ) = K^J_n ( P , E )\, $ for all $\, n \geq 3\, $, so the map 
$\, K_n ( \overline P , \overline E ) \rightarrowtail K_n ( P , E )\, $ is necessarily injective. We begin with 
$\, n = 3\, $. Choose suspensions $\, ( R_{\overline P , \overline E}\, ,\, U_{\overline E} )\, $ and 
$\, ( R_{P , E}\, ,\, U_E )\, $ respectively and put 
$\, D = J_{\overline P , \overline E} / R_{\overline P , \overline E}\, ,\, D' = J_{P , E} / R_{P , E}\, $. Also let 
$\, C = N / [ N , F ] = ker \{ \overline P \twoheadrightarrow P \}\, $. Suppose the map 
$\, K_3 ( \overline P , \overline E ) \rightarrow K_3 ( P , E )\, $ is not surjective. Then a corresponding element of $\, K^J_2 ( J_{\overline P , \overline E}\, ,\, U_{\overline E} )\, $ maps to zero in 
$\, K^J_2 ( D'\, ,\, U_E / R_{P , E} )\, $ but not in 
$\, K^J_2 ( D\, ,\, U_{\overline E} / R_{\overline P , \overline E} )\, $. Since 
$\, K^J_2 ( C ,\overline E ) = 0\, $ one gets 
$\, J_{C , \overline E} \subseteq R_{\overline P , \overline E}\, $, so that 
$\, U_{\overline E} / R_{\overline P , \overline E}\, $ contains a central copy of $\, C\, $. There is no direct map 
\smallskip
$$ K^J_2 ( D\, ,\, U_{\overline E} / R_{\overline P , \overline E} ) \longrightarrow 
K^J_2 ( D'\, ,\, U_E / R_{P , E} ) $$
\par\medskip\noindent
but replacing $\, R_{\overline P , \overline E}\, $ with $\, R_{0 , \overline P , \overline E}\, $ and 
$\, R_{P , E}\, $ with $\, R_{0 , P , E}\, $ one gets a map
\smallskip
$$ K^J_2 ( D\, ,\, U_{\overline E} / R_{0 , \overline P , \overline E} ) \longrightarrow 
K^J_2 ( D'\, ,\, U_E / R_{0 , P , E} ) $$
\par\medskip\noindent
and surjective maps 
\smallskip
$$ K^J_2 ( D\, ,\, U_{\overline E} / R_{0 , \overline P , \overline E} ) \twoheadrightarrow 
K^J_2 ( D\, ,\, U_{\overline E} / R_{\overline P , \overline E} )\> , $$ 
$$ K^J_2 ( D'\, ,\, U_E / R_{0 , P , E} ) \twoheadrightarrow K^J_2 ( D'\, ,\, U_E / R_{P , E} ) $$
\par\medskip\noindent
showing that our element lifts to an element of 
$\, K^J_2 ( D\, ,\, U_{\overline E} / R_{0 , \overline P , \overline E} )\, $ mapping to zero in 
$\, K^J_2 ( D\, ,\, U_E / R_{0 , P , E} )\, $ which injects into the group 
$\, K^J_2 ( D'\, ,\, U_E / R_{0 , P , E} )\, $ since $\, D\, $ and $\, D'\, $ are in the image of 
$\, [\, U_E\, ,\, U_E\, ]\, $ and $\, ( U_E / J_{0 , P , E} )^{ab} = 
U^{ab}_E\, $ is free abelian. Such an element maps to an element in the kernel of 
$\, K^J_2 ( D\, ,\, U_{\overline E} / R_{\overline P , \overline E} ) \rightarrow
K^J_2 ( D\, ,\, U_E / \widetilde R_{P , E} )\, $ where $\,\widetilde R_{P , E}\, $ denotes the image of 
$\, R_{\overline P , \overline E}\, $ in $\, U_E\, $, i.e. it is represented by an element of 
$\, [\, U_D\, ,\, U_C\, ]\, $. Since $\, \overline E\, $ is perfect this element has a preimage in 
$\, [\, U_{J_{\overline P , \overline E}\, ,\, U_{\overline E}}\, ,\, U_{J_{\overline E}\, ,\, U_{\overline E}}\, ]\, +\, 
[\, [\, U_{J_{\overline P , \overline E}\, ,\, U_{\overline E}}\, ,\, U^2_{\overline E}\, ] , U^2_{\overline E}\, ]\, $ modulo $\, [\, U_{J_{\overline P , \overline E}\, ,\, U_{\overline E}}\, ,\, J_{U_{\overline E}}\, ]\, $ defining an element of $\, [\, \overline J_{\overline P , \overline E}\, ,\, U_{\overline E}\, ]\, $ with 
$\, \overline J_{\overline P , \overline E}\, $ some minimal (almost canonical) $U_{\overline E}$-central extension of $\, J_{\overline P , \overline E}\, $ by 
$\, K^J_2 ( J_{\overline P , \overline E}\, ,\, U_{\overline E} )\, $. It is clear that 
$\, [\, [\, \overline J_{\overline P , \overline E}\, ,\, U_{\overline E}\, ] , U_{\overline E}\, ]\, $ maps to zero in 
$\, K^J_2 ( D\, ,\, U_{\overline E} / R_{\overline P , \overline E} )\, $ since $\, D\, $ is central. Then one gets a preimage for our element in $\, K^J_2 ( D\, ,\, J_{\overline E} / R_{\overline P , \overline E} )\, $
$\, = K^J_2 ( D\, ,\, J_E / \widetilde R_{P , E} )\, $ which is an extension of $\, D \otimes J^{ab}_{E / P}\, $ by the image of $\, D \otimes D\, $ which can be neglected being in the kernel of the map to 
$\, K^J_2 ( D\, ,\, U_{\overline E} / R_{\overline P , \overline E} )\, $ (because $\, D\, $ is in the image of 
$\, [\, U_{\overline E}\, ,\, U_{\overline E}\, ]\, $). The remaining part of this kernel is an extension of 
$\, K_2 ( E / P )\, $ by $\, [\, J_{E / P}\, ,\, U_{E / P}\, ] / [\, J_{E / P}\, ,\, J_{E / P}\, ]\, $ tensored with 
$\, D\, $ and coincides with the kernel of the map to $\, K^J_2 ( D\, ,\, U_E / \widetilde R_{P , E} )\, $. Thus if our element is to be trivial in the latter group it is already zero in 
$\, K^J_2 ( D\, ,\, U_{\overline E} / R_{\overline P , \overline E} )\, $, proving that 
$\, K_3 ( \overline P , \overline E ) = K_3 ( P , E )\, $. One proceeds by restricting the pair 
$\, ( R_{R_{0 , \overline P , \overline E}\, ,\, U_{\overline E}}\, ,\, U^2_{\overline E} )\, $ from the pair 
$\, ( R_{R_{0 , P , E}\, ,\, U_E}\, ,\, U^2_E )\, $ which are supposed to be the kernels of some almost canonical maps of $\, U_{R_{0 , \overline P , \overline E}\, ,\, U_{\overline E}}\, $ and 
$\, U_{R_{0 , P , E}\, ,\, U_E}\, $ to the preimages of $\, R_{0 , \overline P , \overline E}\, $ and 
$\, R_{0 , P , E}\, $ in $\, {\overline R}^{ _{\sps U_{\overline E}}}_{\overline P , \overline E}\, $ and 
$\, {\overline R}^{ _{\sps U_E}}_{P , E}\, $ respectively (the latter is obtained from $\, R^2_{P , E}\, $ by restriction, say). Define 
$\, J := J_{R_{0 , P , E}\, ,\, U_E}\, ,\, 
\overline J  :=  J_{R_{0 , \overline P , \overline E}\, ,\, U_{\overline E}}\, $. Then
\smallskip
$$ K^J_2\, ( \overline J\, ,\, U^2_{\overline E} )\buildrel\sim\over\longrightarrow 
K^J_2 ( J\, ,\, U^2_E ) $$
\par\medskip\noindent
is an isomorphism (the map is surjective by existence of a splitting, to see this consider a section 
$\, ( P , E ) \nearrow ( \overline P , \overline E )\, $ which induces a splitting 
$\, ( R_{0 , P , E}\, ,\, U_E ) \rightarrow ( R_{0 , \overline P , \overline E} + 
[\, U_{C , \overline E}\, ,\, U_{\overline E}\, ]\, ,\, U_{\overline E} ) = 
( R_{0 , \overline P , \overline E}\, ,\, U_{\overline E} )\, $ since 
$\, [\, U_{C , \overline E}\, ,\, U_{\overline E}\, ] = 
J_{C , \overline E} \cap [\, U_{C , \overline E}\, ,\, U_{\overline E}\, ] = R_{0 , C , \overline E} \subseteq 
R_{0 , \overline P , \overline E}\, $). Injectivity follows from the inclusions
\smallskip 
$$ K^J_2 ( \overline J\, ,\, U^2_{\overline E} ) \rightarrowtail 
K^J_2 ( J_{R_{\overline P , \overline E}\, ,\, U_{\overline E}}\, ,\, U^2_{\overline E} ) \simeq 
K^J_2 ( J_{R_{P , E}\, ,\, U_E}\, ,\, U^2_E )\> , $$ 
$$ K^J_2 ( J\, ,\, U^2_E ) \rightarrowtail K^J_2 ( J_{R_{P , E}\, ,\, U_E}\, ,\, U^2_E )\> . $$
\par\medskip\noindent
Since $\, \overline R = R_{R_{0 , \overline P , \overline E}\, ,\, U_{\overline E}}\, $ is just the preimage of 
$\, R = R_{R_{0 , P , E}\, ,\, U_E}\, $ in $\,\overline J\, $ one also gets an isomorphism 
$\, K^J_2 ( \overline R\, ,\, U^2_{\overline E} ) \buildrel\sim\over\rightarrow K^J_2 ( R\, ,\, U^2_E )\, $ proving that $\, K_4 ( \overline P , \overline E ) \simeq K_4 ( P , E )\, $. The argument in higher dimensions is completely analogous restricting the pair 
$\, ( R^n_{R_{0 , \overline P , \overline E}\, ,\, U_{\overline E}}\, ,\, U^{n+1}_{\overline E} )\, $ from 
$\, ( R^n_{R_{0 , P , E}\, ,\, U_E}\, ,\, U^{n+1}_E )\, $ at each step by induction.
\par\noindent
We may now restrict to extensions of the form 
\smallskip
$$ 1 \longrightarrow \bigl( N , F \bigr) \longrightarrow \bigl( F , F \bigr) \longrightarrow 
\bigl( E , E \bigr) \longrightarrow 1 $$
\par\medskip\noindent
with $\, N\, $ a full subgroup of $\, F\, $. The general case will follow from a diagram chase argument. Up to the case $\, K_3 ( E ) = K^J_3 ( E )\, $ exactness of the associated long sequence of $K_*$-groups follows from Corollary 5.1. We begin with exactness at $\, K_{n+3} ( E )\, ,\, n \geq 1\, $. Recall from Proposition 1 that $\, K_{n+2} ( N , F )\, $ can be identified with the subgroup of abelian regular elements of $\, K^J_2 ( J^n_{ N , F}\, ,\, U^n_F )\, $ and with 
$\, K^J_2 ( R_{J^{n-1}_{N , F}\, ,\, U^{n-1}_F}\, ,\, U^n_F )\, $ where 
$\, R_{J^{n-1}_{N , F}\, ,\, U^{n-1}_F}\, $ is the kernel of an almost canonical map 
\smallskip
$$ U_{J^{n-1}_{N , F}\, ,\, U^{n-1}_F}\> \longrightarrow\> 
{\overline {J^{n-1}_{N , F}\, }}^{ _{ U^{n-1}_F}} $$
\par\medskip\noindent
induced from some $\, R^n_{N , F}\, $. To save notation we wil write $\, R^n_{N , F}\, $ for 
$\, R_{J^{n-1}_{N , F}\, ,\, U^{n-1}_F}\, $ etc. in the following argument. Let 
$\, ( C_{{\varphi }_n}\, ,\, U^{n+1}_F )\, $ be the mapping cone of the inclusion 
$\, ( R^n_{N , F}\, ,\, U^n_F ) \subseteq ( R^n_{F , F}\, ,\, U^n_F )\, $ and $\, C^r_{{\varphi }_n}\, $ the subgroup generated by $\, U_{R^n_{N , F}\, ,\, U^n_F}\, $ and $\, R_{R^n_{F , F}\, ,\, U^n_F} = 
R^2_{J^{n-1}_{F , F}\, ,\, U^{n-1}_F}\, $ which is induced in the usual way from some ordinary $(n+1)$-fold suspension of $\, F\, $. Also consider the mapping cone $\,\widetilde C_{{\varphi }_n}\, $ of the inclusion $\, ( J^n_{N , F}\, ,\, U^n_F ) \subseteq ( J^n_{N , F} + R^n_{F , F}\, ,\, U^n_F )\, $. Let 
$\,\widetilde R^{n+1}_{E , E}\, $ denote the image of $\, R^2_{J^{n-1}_{F , F}\, ,\, U^{n-1}_F}\, $ in 
$\, J^{n+1}_{E , E}\, $ and $\,\widetilde R^n_{E , E}\, $ the image of $\, R^n_{F , F}\, $ in $\, J^n_{E , E}\, $. By existence of a splitting it is easy to see that 
$\, K^J_2 ( \widetilde C_{{\varphi }_n}\, ,\, U^{n+1}_F )\ $ is regular isomorphic to 
$\, K^J_2 ( J_{\widetilde R^n_{E , E}\, ,\, U^n_E}\, ,\, U^{n+1}_E )\, $ so that 
$\, ( \widetilde C^r_{{\varphi }_n}\, ,\, U^{n+1}_F )\, $ is regular isomorphic to 
$\, K^J_2 ( \widetilde R^{n+1}_{E , E}\, ,\, U^{n+1}_E )\, $ which contains the image of 
$\, K^J_2 ( R_{0 , R^n_{E , E}\, ,\, U^n_E}\, ,\, U^{n+1}_E )\, $ and injects into 
$\, K^J_2 ( J^{n+1}_{E , E}\, ,\, U^{n+1}_E )\, $, so that it contains all abelian regular elements of the latter group, i.e. it contains $\, K_{n+3} ( E )\, $. Obviously, the map 
$\, K^J_2 ( C_{{\varphi }_n}\, ,\, U^{n+1}_F ) \rightarrowtail 
K^J_2 ( \widetilde C_{{\varphi }_n}\, ,\, U^{n+1}_F )\, $ is injective as the latter is equal to 
$\, K^J_{n+3} ( E )\, $. We will show that any abelian regular element of 
$\, K^J_2 ( \widetilde C_{{\varphi }_n}\, ,\, U^{n+1}_F )\, $ lifts to 
$\, K^J_2 ( C^r_{{\varphi }_n}\, ,\, U^{n+1}_F )\, $ which itself consists solely of abelian regular elements of $\, K^J_2 ( C_{{\varphi }_n}\, ,\, U^{n+1}_F )\, $. Then $\, K^J_2 ( C^r_{{\varphi }_n}\, ,\, U^{n+1}_F )\, $ must be isomorphic to $\, K_{n+3} ( E )\, $. For the first statement apply Corollary 5.1 to the extension 
\smallskip
$$ 1 \rightarrow \bigl( C , {U^n_F \over R^n_{N , F}} \bigr) \rightarrow 
\bigl( {R^n_{F , F} \over R^n_{N , F}} , {U^n_F \over R^n_{N , F}} \bigr) \rightarrow 
\bigl( {R^n_{F , F} \over R^n_{F , F} \cap J^n_{N , F}} , {U^n_F \over R^n_{F , F} \cap J^n_{N , F}} \bigr) \rightarrow 1 \> . $$
\par\medskip\noindent 
By repeated application of Lemma 22 the quotient pair is exact, since 
$\, ( R^n_{F , F} / ( R^n_{F , F} \cap J^n_{N , F} )\, ,\, U_{U^{n-1}_F / J^{n-1}_{N , F}} )\, $ is fully exact as well as the pair $\, ( \overline J^{n-1}_{N , F}\, ,\, U^n_F / ( R^n_{F , F} \cap J^{n-1}_{N , F} )\, $ and also 
$\, ( \overline J^{n-1}_{N , F}\, ,\, U^n_F / R^n_{F , F} )\, $ is exact of second order 
(with $\,\overline J^{n-1}_{N , F} = U_{J^{n-1}_{N , F}\, ,\, U^{n-1}_F} / ( R^n_{F , F} \cap J^n_{N , F} )\, $). It is clear that 
$\, ( \overline J^{n-1}_{N , F}\, ,\, U^n_F / ( R^n_{F , F} \cap J^n_{N , F} ) ) \simeq 
( \overline J^{n-1}_{N , F}\, ,\, \overline J^{n-1}_{N , F} \rtimes U_{U^{n-1}_F / J^{n-1}_{N , F}} )\, $ is exact and that $\, ( J^{n-1}_{N , F}\, ,\, U^n_F / J^n_{N , F} )\, $ is fully exact. Then assume given an extension 
\smallskip
$$ 1\, \longrightarrow\, \bigl( N , F \bigr)\, \longrightarrow\, \bigl( M , F \bigr)\, \longrightarrow\, 
\bigl( P , E \bigr)\, \longrightarrow\, 1 $$
\par\medskip\noindent
such that $\, N \subseteq [ M , F ]\, $ and $\, ( P , E )\, $ is exact of second order. Consider the diagram 
\smallskip
$$ \vbox{\halign{ #&#&#\cr
\hfil $\bigl( J_F\, ,\, J_{E / P} \rtimes U_F \bigr)$\hfil & \hfil $\longrightarrow $\hfil & 
\hfil $\bigl( J_F \times J_{E / P}\, ,\, J_{E / P} \rtimes U_F \bigr)$\hfil \cr
\hfil $\Bigm\downarrow $\hfil && \hfil $\Bigm\downarrow $\hfil \cr
\hfil $\bigl( J_F + U_{N , F}\, ,\, J_{E / P} \rtimes U_F \bigr)$\hfil & \hfil $\longrightarrow $\hfil & 
\hfil $\bigl( (J_F + U_{N , F} \times J_{E / P}\, ,\, J_{E / P} \rtimes U_F \bigr)$\hfil \cr
\hfil $\Bigm\downarrow $\hfil && \hfil $\Bigm\downarrow $\hfil \cr
\hfil $\bigl( N\, ,\, M \rtimes U_{E / P} \bigr)$\hfil & \hfil $\longrightarrow $\hfil & 
\hfil $\bigl( N\, ,\, F \bigr)$\hfil \cr }} \> . $$
\par\medskip\noindent
Then $\, ( N\, ,\, M \rtimes U_{E / P} )\, $ is exact and supposing that 
$\, K^J_3 ( N\, ,\, M \rtimes U_{E / P} ) = 0\, $ (for example if $\, N\, $ is central) one gets that the map 
\smallskip
$$ K^J_2 ( J_F\, ,\, J_{E / P} \rtimes U_F ) \rightarrowtail 
K^J_2 ( J_F + U_{N , F}\, ,\, J_{E / P} \rtimes U_F ) $$
\par\medskip\noindent
is injective. On the other hand it follows from injectivity of the map
$\, K^J_2 ( U_{N , F}\, ,\, J_{E / P} \rtimes U_F ) \rightarrowtail K^J_2 ( J_{E / P} \rtimes U_F )\, $ that the middle horizontal map 
\smallskip
$$ K^J_2 ( J_F + U_{N , F}\, ,\, J_{E / P} \rtimes U_F ) \rightarrowtail 
K^J_2 ( (J_F + U_{N , F}) \times J_{E / P}\, ,\, J_{E / P} \rtimes U_F ) $$
\par\medskip\noindent
is injective since $\, ( P , E )\, $ is exact of second order. Thus also $\, ( M , F )\, $ is exact of second order, which leads to the statement that the pair
$\, ( \overline J^{n-1}_{N , F}\, ,\, U^n_F / ( R^n_{F , F} \cap J^n_{N , F} ) )\, $ is fully exact. Also one gets that the pair
$\, ( J^{n-1}_{N , F}\, ,\, U^n_F / ( R^n_{F , F} + J^n_{N , F} ) )\, $ is exact of second order because the pair
$\, ( J^{n-1}_{N , F} \times J_{U^{n-1}_F} / ( R^n_{F , F} + J^n_{N , F} )\, ,\, 
U^n_F / ( R^n_{F , F} + J^n_{N , F} ) )\, $ is exact of second order from the $P\times Q$-Lemma (and the argument above) and $\, J^{n-1}_{N , F}\, $ is free. The argument above shows that also 
$\, ( \overline J^{n-1}_{N , F}\, ,\, U^n_F / R^n_{F , F} )\, $ is exact of second order. Then the pair 
$\, ( R^n_{F , F} / ( R^n_{F , F} \cap J^n_{N , F} )\, ,\, U^n_F / ( R^n_{F , F} \cap J^n_{N , F} ))\, $ is exact, and even fully exact since $\, R^n_{F , F} / ( R^n_{F , F} \cap J^n_{N , F} )\, $ is free and 
$\, ( R^n_{F , F} / ( R^n_{F , F} \cap J^n_{N , F} ) \times \overline J^{n-1}_{N , F}\, ,\, 
U^n_F / ( R^n_{F , F} \cap J^n_{N , F} ) )\, $ is exact of second order. Also the kernel 
$\, ( C\, ,\, U^n_F / R^n_{F , F} )\, $ is fully exact. Assuming by induction that the regular sequence is exact at all places to the right of $\, K_{n+3} ( E )\, $ (in particular at $\, K_{n+1} ( N , F )\, $) one gets that $\, C\, $ is contained in $\, [\, R^n_{F , F} / R^n_{N , F}\, ,\, U^n_F / R^n_{N , F}\, ]\, $ because it consists of elements in the kernel of the map $\, K_{n+1} ( N , F ) \rightarrow K_{n+1} ( F )\, $ which can be lifted to 
$\, K_{n+2} ( E )\, $ by induction assumption (and one knows that 
$\, K^J_2 ( R^n_{F , F} / R^n_{N , F}\, ,\, U^n_F / R^n_{N , F} )\, $ contains preimages for every element in 
$\, K_{n+2} ( E )\, $). Then Lemma 21 gives that the pair 
$\, ( R^n_{F , F} / R^n_{N , F}\, ,\, U^n_F / R^n_{N , F} )\, $ is (fully, by the argument above) exact. An application of Corollary 5.1 to the extension above yields that any abelian regular element of
\smallskip 
$$ K^J_2 ( J_{R^n_{F , F} / (R^n_{F , F} \cap J_{N , F})\, ,\, U^n_F / (R^n_{F , F} \cap J^n_{N , F})}\, ,\, 
U_{U^n_F / ( R^n_{F , F} \cap J^n_{N , F})} ) {\simeq }_{reg} 
K^J_2 ( \widetilde C_{{\varphi }_n}\, ,\, U^{n+1}_F ) $$ 
\par\medskip\noindent
which must be trivial in 
$\, K^J_2 ( C\, ,\, U^n_F / R^n_{N , F} )\, $ lifts to
\smallskip 
$$ K^J_2 ( J_{R^n_{F , F} / R^n_{N , F}\, ,\, U^n_F / R^n_{N , F}}\, ,\, U_{U^n_F / R^n_{N ,F}} ) 
{\simeq }_{reg} K^J_2 ( C_{{\varphi }_n}\, ,\, U^{n+1}_F )\> . $$ 
\par\medskip\noindent
To see that it lifts to 
$\, K^J_2 ( C^r_{{\varphi }_n}\, ,\, U^{n+1}_F )\, $ one considers the diagram
\smallskip
$$ \vbox{\halign{ #&#&#&#&#\cr
\hfil $\bigl(\, \widetilde R_{C\, ,\, {U^n_F / R^n_{N , F}}}\, ,\, 
U_{U^n_F / R^n_{N , F}}\, \bigr)$\hfil & 
\hfil $\rightarrow $\hfil & \hfil $\bigl(\, \widetilde J_{C\, ,\, {U^n_F / R^n_{N , F}}}\, ,\, 
U_{U^n_F / R^n_{N , F}} \,\bigr)$\hfil &&\cr
\hfil $\Bigm\downarrow $\hfil && \hfil $\Bigm\downarrow $\hfil &&\cr
\hfil $\bigl( {R_{R^n_{F , F}\, ,\, U^n_F} \over R_{R^n_{F , F} , U^n_F} \cap J_{R^n_{N , F} , U^n_F}}\, ,\, 
U_{\negthinspace {U^n_F / R^n_{N , F}}} \bigr)$\hfil & \hfil $\rightarrow $\hfil & 
\hfil $\bigl(\, {J_{R^n_{F , F}\, ,\, U^n_F} \over J_{R^n_{N , F}\, ,\, U^n_F}}\, ,\, 
U_{\negthinspace {U^n_F / R^n_{N , F}}} \,\bigr)$\hfil & \hfil $\rightarrow $\hfil & 
\hfil $\bigl(\, D\, , {U^{n+1}_F \over C^r_{{\varphi }_n}} \,\bigr)$\hfil \cr
\hfil $\Bigm\downarrow $\hfil && \hfil $\Bigm\downarrow $\hfil && \hfil $\Bigm\downarrow $\hfil \cr
\hfil $\bigl( {R_{R^n_{F , F} , U^n_F} \over R_{R^n_{F , F} , U^n_F} \cap J^{n+1}_{N , F}}\, ,\, 
U_{\negthinspace {U^n_F / J^n_{N , F}}} \bigr)$\hfil & \hfil $\rightarrow $\hfil & 
\hfil $\bigl( {J_{R^n_{F , F} , U^n_F} \over J_{R^n_{F , F} , U^n_F} \cap J^{n+1}_{N , F}}\, ,\, 
U_{\negthinspace {U^n_F / J^n_{N , F}}} \bigr)$\hfil & 
\hfil $\rightarrow $\hfil & 
\hfil $\bigl(\, D'\, , {U^{n+1}_F \over \widetilde C^r_{{\varphi }_n}} \,\bigr)$\hfil \cr }} $$
\par\medskip\noindent
where 
\smallskip
$$ \widetilde J_{C\, ,\, U^n_F / R^n_{N , F}} = 
( J_{R^n_{F , F}\, ,\, U^n_F} \cap J^{n+1}_{N , F} ) / J_{R^n_{N , F}\, ,\, U^n_F}\> , $$ 
$$ \widetilde R_{C\, ,\, U^n_F / R^n_{N , F}} = \widetilde J_{C\, ,\, U^n_F / R^n_{N , F}} \cap 
R_{R^n_{F , F} , U^n_F} / (R_{R^n_{F , F} , U^n_F} \cap J_{R^n_{N , F} , U^n_F})\> , $$
$$ D = C_{{\varphi }_n} / C^r_{{\varphi }_n}\> ,\> 
D' = \widetilde C_{{\varphi }_n} / \widetilde C^r_{{\varphi }_n}\> . $$ 
\par\medskip\noindent
Put
\smallskip
$$ D_0 = \widetilde J_{C\, ,\, U^n_F / R^n_{N , F}} / \widetilde R_{C\, ,\, U^n_F / R^n_{N , F}} = 
{J_{R^n_{F , F}\, ,\, U^n_F} \cap J_{J^n_{N , F}\, ,\, U^n_F} \over
J_{R^n_{N , F}\, ,\, U^n_F} + ( R_{R^n_{F , F}\, ,\, U^n_F} \cap J_{J^n_{N , F}\, ,\, U^n_F} )} \> . $$
\par\medskip\noindent 
Thus $\, D_0\, $ can be described as the intersection of the images of the groups
$\, K^J_2 ( J^n_{N , F}\, ,\, U^n_F )\, $ and $\, K^J_2 ( R^n_{F , F}\, ,\, U^n_F )\, $ in 
$\, K^J_2 ( J^n_{F , F}\, ,\, U^n_F )\, $ modulo the image of $\, K^J_2 ( R^n_{N , F}\, ,\, U^n_F )\, $. By our induction assumption on exactness of the regular sequence at all places to the right of the one in question one gets $\, D_0 = 0\, $ and 
$\, \widetilde R_{C\, ,\, U^n_F / R^n_{N , F}} = \widetilde J_{C\, ,\, U^n_F / R^n_{N , F}}\, $. Then 
$\, D = D'\, $. The abelian regular element in 
\smallskip
$$ K^J_2 ( J_{R^n_{F , F}\, ,\, U^n_F} / ( J_{R^n_{F , F}\, ,\, U^n_F} \cap J^{n+1}_{N , F})\, ,\, 
U_{U^n_F / J^n_{N , F}} ) $$
\par\medskip\noindent
lifts to 
\smallskip
$$ K^J_2 ( R_{R^n_{F , F}\, ,\, U^n_F} / (R_{R^n_{F , F}\, ,\, U^n_F} \cap J^{n+1}_{N , F})\, ,\, 
U_{U^n_F / R^n_{N , F}} ) $$
\par\medskip\noindent
and then to 
\smallskip
$$ K^J_2 ( R_{R^n_{F , F}\, ,\, U^n_F} / (R_{R^n_{F , F}\, ,\, U^n_F} \cap J_{R^n_{N , F}\, ,\, U^n_F})\, ,\, 
U_{U^n_F / R^n_{N , F}} ) \simeq K^J_2 ( C^r_{{\varphi }_n}\, ,\, U^{n-1}_F ) $$
\par\medskip\noindent
as suggested since its image in $\, K_1 ( \widetilde J_{C\, ,\, U^n_F / R^n_{N , F}} )\, $ is trivial. It remains to prove that the (injective) image of $\, K^J_2 ( C^r_{{\varphi }_n}\, ,\, U^{n+1}_F )\, $ in 
$\, K^J_2 ( J^{n+1}_{E , E}\, ,\, U^{n+1}_E )\, $ is equal to $\, K^J_2 ( R^{n+1}_{E , E}\, ,\, U^{n+1}_E )\, $, which we will combine with proving halfexactness of the sequence 
\smallskip
$$ K^J_2\, \bigl( R^{n+1}_{F , F}\, ,\, U^{n+1}_F \bigr) \longrightarrow 
K^J_2\, \bigl( C^r_{{\varphi }_n}\, ,\, U^{n+1}_F \bigr) \longrightarrow 
K^J_2\, \bigl( R^n_{N , F}\, ,\, U^n_F \bigr) $$
\par\medskip\noindent
which implies exactness of the regular sequence at $\, K_{n+3} ( E )\, $. In fact using the induction assumption that the regular sequence is exact to the right of $\, K_{n+3} ( E )\, $, in particular at 
$\, K_{n+2} ( N , F )\, $, a positive answer to the second problem implies the solution of the first one. Namely, any element $\, x \in K^J_2 ( C^r_{{\varphi }_n}\, ,\, U^{n+1}_F )\, $ maps to an element in the kernel of $\, K^J_2 ( R^n_{N , F}\, ,\, U^n_F ) \rightarrow K^J_2 ( R^n_{F , F}\, ,\, U^n_F )\, $ which can be lifted to $\, K^J_2 ( R^{n+1}_{E , E}\, ,\, U^n_E )\, $ and taking the difference with the corresponding element of $\, K^J_2 ( C^r_{{\varphi }_n}\, ,\, U^{n+1}_F )\, $ reduces the problem of whether 
$\, K^J_2 ( C^r_{{\varphi }_n}\, ,\, U^{n+1}_F )\, $ is equal to 
$\, K^J_2 ( R^{n+1}_{E , E}\, ,\, U^{n+1}_E )\, $ to the proof of halfexactness of the above sequence, since any element in the image of $\, K^{n+1}_{F , F}\, ,\, U^{n+1}_F )\, $ maps to 
$\, K^J_2 ( R^{n+1}_{E , E}\, ,\, U^{n+1}_E )\, $. Assume given an element $\, x\, $ in the kernel of 
\smallskip
$$ K^J_2 ( C^r_{{\varphi }_n}\, ,\, U^{n+1}_F ) \longrightarrow K^J_2 ( R^n_{N , F}\, ,\, U^n_F )\> . $$
\par\medskip\noindent
Then $\, x\, $ is zero already in $\, K^J_2 ( R^n_{N , F}\, ,\, U^{n+1}_F / J_{R^n_{F , F}\, ,\, U^n_F} )\, $. Let 
$\, \overline {\varphi }\, ( \overline R^n_{N , F} )\, $ denote the image of 
$\, {\overline R^n_{N , F}}^{ _{\sps U^n_F}}\, $ in $\, {\overline R^n_{F , F}}^{ _{\sps U^n_F}}\, $ (by some equivariant lift of $\, R^n_{N , F} \subseteq R^n_{F , F}\, $). Dividing by 
$\, R_{R^n_{F , F}\, ,\, U^n_F}\, $ induces a map 
$\, ( C^r_{{\varphi }_n}\, ,\, U^{n+1}_F ) \rightarrow 
( \overline {\varphi }\, ( \overline R^n_{N , F} )\, ,\, \overline U^n_F )\, $ with 
$\, \overline U^n_F = U^{n+1}_F / R_{R^n_{F , F}\, ,\, U^n_F}\, $. One must show that the image of $\, x\, $ in $\, K^J_2 ( \overline {\varphi }\, ( \overline R^n_{N , F} )\, ,\, \overline U^n_F )\, $ is zero. One notes that
\smallskip 
$$ K^J_2 ( J_{\overline R^n_{F , F}\, ,\, \overline U^n_F}\, ,\, U_{\overline U^n_F} ) \simeq 
K^J_2 ( R_{R^n_{F , F}\, ,\, U^n_F}\, ,\, U^{n+1}_F ) $$ 
\par\medskip\noindent
by comparison of the exact sequences 
\smallskip
$$ \vbox{\halign{ #&#&#&#&#&#&#&#&#\cr
\hfil $1$\hfil &\hfil $\longrightarrow $\hfil & 
\hfil $\bigl( J_{\overline R^n_{F , F} , \overline U^n_F}\, ,\, U_{\overline U^n_F} \bigr)$\hfil & 
\hfil $\longrightarrow $\hfil & 
\hfil $\bigl( U_{\overline R^n_{F , F}\, ,\, \overline U^n_F}\, ,\, U_{\overline U^n_F} \bigr)$\hfil & 
\hfil $\buildrel can\over\longrightarrow $\hfil & 
\hfil $\bigl( \overline R^n_{F , F}\, , 
{U_{\overline U^n_F}\over J_{\overline R^n_{F , F} , \overline U^n_F}} \bigr)$\hfil & 
\hfil $\longrightarrow $\hfil & \hfil $1$\hfil \cr 
&& \hfil $\Bigm\downarrow\Bigm\uparrow $\hfil && 
\hfil $ u^{-1}_{\alpha }\Bigm\downarrow\Bigm\uparrow u_{\alpha } $\hfil && 
\hfil $ can \Bigm\downarrow\Bigm\uparrow\quad $\hfil &&\cr
\hfil $1$\hfil & \hfil $\longrightarrow $\hfil & 
\hfil $\bigl( R_{R^n_{F , F} , U^n_F}\, ,\, U^{n+1}_F \bigr)$\hfil & \hfil $\longrightarrow $\hfil & 
\hfil $\bigl( U_{R^n_{F , F}\, ,\, U^n_F}\, ,\, U^{n+1}_F \bigr)$\hfil & 
\hfil $\buildrel \alpha \over\longrightarrow $\hfil & 
\hfil $\bigl( \overline R^n_{F , F}\, ,\, \overline U^n_F \bigr)$ \hfil & \hfil $\longrightarrow $\hfil 
\hfil $1$\hfil \cr }} $$
\par\medskip\noindent
and the fact that $\, K^J_3 ( \overline R^n_{F , F}\, ,\, \overline U^n_F ) \simeq 
K^J_3 ( \overline R^n_{F , F}\, ,\, U_{\overline U^n_F} / J_{\overline R^n_{F , F}\, ,\, \overline U^n_F} )\, $. To see this one notes that both pairs $\, ( \overline R^n_{F , F}\, ,\, \overline U^n_F )\, $ and 
$\, ( \overline R^n_{F , F}\, ,\, U_{\overline U^n_F} / J_{\overline R^n_{F , F}\, ,\, \overline U^n_F} )\, $ are fully exact being central extensions of $\, R^n_{F , F}\, $ with kernel contained in 
$\, [\, \overline R^n_{F , F}\, ,\, \overline U^n_F\, ]\, $. Then the maps 
\smallskip
$$ K^J_3 ( \overline R^n_{F , F}\, ,\, \overline U^n_F ) \rightarrowtail 
K^J_3 ( R^n_{F , F}\, ,\, U^{n+1}_F / J_{R^n_{F , F}\, ,\, U^n_F} ) \simeq 
K^J_3 ( R^n_{F , F}\, ,\, U_{\overline U^n_F / \sim } ) $$
\par\medskip\noindent
and 
\smallskip
$$ K^J_3 ( \overline R^n_{F , F}\, ,\, U_{\overline U^n_F} / J_{\overline R^n_{F , F}\, ,\, \overline U^n_F} ) \rightarrowtail K^J_3 ( R^n_{F , F}\, ,\, U_{\overline U^n_F / \sim } ) $$
\par\medskip\noindent
are injective with $\,\sim\, $ corresponding to division by the preimage of the kernel of 
$\, \overline R^n_{F , F} \rightarrow R^n_{F , F}\, $ in $\, U_{\overline R^n_{F , F}\, ,\, \overline U^n_F}\, $, so the map 
\smallskip
$$ K^J_3 ( \overline R^n_{F , F}\, ,\, U_{\overline U^n_F} / J_{\overline R^n_{F , F}\, ,\,\overline U^n_F} )
\rightarrowtail
K^J_3 ( \overline R^n_{F , F}\, ,\, \overline U^n_F ) $$
\par\medskip\noindent
is injective, and surjective by existence of a splitting. Given the almost canonical map $\,\alpha\, $ of the diagram one constructs the lift $\, u_{\alpha }\, $ and its left inverse $\, u^{-1}_{\alpha }\, $ in the following way. Choose sections $\, U^n_F / R^n_{F , F} \nearrow \overline U^n_F / \overline R^n_{F , F} \nearrow 
\overline U^n_F\, $. Then drop $\, \overline U^n_F / \overline R^n_{F , F} \nearrow 
\overline U^n_F\, $ to $\, U^n_F / R^n_{F , F} \nearrow U^n_F\, $. Define the corresponding bases 
$\, C_{R^n_{F , F}\, ,\, U^n_F} \cup C_{U^n_F / R^n_{F , F}}\, $ of $\, U^{n+1}_F\, $ and 
$\, C_{\overline R^n_{F , F}\, ,\, \overline U^n_F} \cup C_{\overline U^n_F / \overline R^n_{F , F}}\, $ of 
$\, U_{\overline U^n_F}\, $ with respect to these sections. For every 
$\, u_x\, u^{-1}_{\overline x} \in C_{R^n_{F , F}\, ,\, U^n_F}\, $ with 
$\, \alpha\, ( u_x\, u^{-1}_{\overline x} ) = {\alpha }_{x , \overline x}\, $ there is an element 
$\, u_{\widetilde x}\, u^{-1}_{\overline {\widetilde x}} \in C_{\overline R^n_{F , F}\, ,\, \overline U^n_F}\, $ such that 
$\, \overline {\widetilde x}\, $ is the image of $\, \overline x\, $ under the composition 
$\, U^n_F / R^n_{F , F} \nearrow \overline U^n_F / \overline R^n_{F , F} \nearrow \overline U^n_F\, $ and 
$\, \widetilde x\, \overline {\widetilde x}^{-1} = {\alpha }_{x , \overline x}\, $. The so defined section 
$\, U^n_F \nearrow \overline U^n_F\, $ is unique making a commutative diagram with the other sections. Define $\, u_{\alpha } ( u_x\, u^{-1}_{\overline x} ) = u_{\widetilde x}\, u^{-1}_{\overline {\widetilde x}}\, $ and $\, u_{\alpha } ( u_{\overline x} ) = u_{\overline {\widetilde x}}\, $ respectively. For the reverse map let 
$\, u^{-1}_{\alpha } ( u_{\overline {\widetilde x}} ) = u_{\overline x }\, $, (i.e. no matter whether 
$\, \overline {\widetilde x}\, $ is in the image of the section $\, U^n_F / R^n_{F , F} \nearrow 
\overline U^n_F\, $ or not), and representing each element of the kernel 
$\, \overline R^n_{F , F} \twoheadrightarrow R^n_{F , F}\, $ by some element of 
$\, J_{R^n_{F , F}\, ,\, U^n_F}\, $, say $\, k ( x )\, $, define 
$\, u^{-1}_{\alpha } ( u_{\widetilde x}\, u^{-1}_{\overline {\widetilde x}} ) = u_x\, u^{-1}_{\overline x}\cdot 
k ( \alpha ( u_x\, u^{-1}_{\overline x} ) \overline {\widetilde x}\, {\widetilde x}^{-1} )\, $ if $\, \widetilde x\, $ is not in the image of $\, U^n_F \nearrow \overline U^n_F\, $ and by 
$\, u^{-1}_{\alpha } ( u_{\widetilde x}\, u^{-1}_{\overline {\widetilde x}} ) = u_x\, u^{-1}_{\overline x}\, $ if it is. The map $\, u_{\alpha }\, $ induces 
$\, ( R_{R^n_{F , F}\, ,\, U^n_F}\, ,\, U^{n+1}_F ) \rightarrow 
( J_{\overline R^n_{F , F}\, ,\, \overline U^n_F}\, ,\, U_{\overline U^n_F} )\, $ and hence a splitting 
$\, ( \overline R^n_{F , F}\, ,\, \overline U^n_F ) \rightarrow 
( \overline R^n_{F , F}\, ,\, U_{\overline U^n_F} / J_{\overline R^n_{F , F}\, ,\, \overline U^n_F} )\, $ of the canonical map in the reverse direction. Let $\, C^r_{\overline {\varphi }_n}\, $ be the preimage of 
$\, C^r_{{\varphi }_n}\, $ under $\, u^{-1}_{\alpha }\, $. It follows that the kernel of $\, u^{-1}_{\alpha }\, $ restricted to $\, C^r_{\overline {\varphi }_n}\, $ is the same as restricted to 
$\, J_{\overline R^n_{F , F}\, ,\, \overline U^n_F}\, $ whence the groups 
$\, K^J_2 ( C^r_{\overline {\varphi }_n}\, ,\, U_{\overline U^n_F} )\, $ and 
$\, K^J_2 ( C^r_{{\varphi }_n}\, ,\, U^{n+1}_F )\, $ are identified under $\, u_{\alpha }\, ,\, 
u^{-1}_{\alpha }\, $. Also the composition 
\smallskip
$$ U_{R^n_{F , F}\, ,\, U^n_F} \buildrel u_{\alpha }\over\longrightarrow 
U_{\overline R^n_{F , F}\, ,\, \overline U^n_F} \buildrel can \over\longrightarrow 
U_{R^n_{F , F}\, ,\, U^n_F} $$
\par\medskip\noindent
is the identity, so that the second canonical map restricted to $\, C_{\overline {\varphi }_n} \supseteq 
C^r_{\overline {\varphi }_n}\, $ induces the same isomorphism 
$\, K^J_2 ( C_{\overline {\varphi }_n}\, ,\, U_{\overline U^n_F} ) \buildrel\sim\over\rightarrow 
K^J_2 ( C_{{\varphi }_n}\, ,\, U^{n+1}_F )\, $ as $\, u^{-1}_{\alpha }\, $. Let 
$\, \overline {R_{R^n_{F , F}\, ,\, U^n_F}}\, $ denote the pull back of the universal $U^{n+1}_F$-central extension of some $(n+1)$-fold suspension of $\, F\, $ to a $U^{n+1}_F$-central extension of 
$\, R_{R^n_{F , F}\, ,\, U^n_F}\, $ by $\, K_{n+3} ( F )\, $, and $\, \overline {C^r_{{\varphi }_n}}\, $ some minimal $U^{n+1}_F$-central extension by $\, K^J_2 ( C^r_{{\varphi }_n}\, ,\, U^{n+1}_F )\, $. Taking the pull backs of $\,\overline {R_{R^n_{F , F}\, ,\, U^n_F}}\, $ and $\,\overline {C^r_{{\varphi }_n}}\, $ by 
$\, u^{-1}_{\alpha }\, $ one gets $U_{\overline U^n_F}$-central extensions 
$\, \overline {J_{\overline R^n_{F , F}\, ,\, \overline U^n_F}}\, $ and 
$\,\overline {C^r_{\overline {\varphi }_n}}\, $ respectively making a commutative diagram of the following sort 
\smallskip
$$ \vbox{\halign{ #&#&#&#&#\cr
\hfil $\overline {J_{\overline R^n_{F , F}\, ,\, \overline U^n_F}}\,$\hfil & \hfil $\largerightarrow $\hfil & 
\hfil $\overline {C^r_{\overline {\varphi }_n}}$\hfil &&\cr
\hfil $\overline u^{-1}_{\alpha }\negthinspace\negthinspace\Bigm\downarrow\quad $\hfil && 
\hfil $\overline u^{-1}_{\alpha }\negthinspace\negthinspace\Bigm\downarrow\quad $\hfil & 
\hfil $\searrow $\hfil &\cr
\hfil $\overline {R_{R^n_{F , F}\, ,\, U^n_F}}\,$\hfil & \hfil $\largerightarrow $\hfil & 
\hfil $\overline {C^r_{{\varphi }_n}}$\hfil & \hfil $\largerightarrow $ \hfil & 
\hfil $\overline {\overline {\varphi } ( \overline R^n_{N , F} )}$\hfil \cr }} $$
\par\medskip\noindent
where the group $\, \overline {\overline {\varphi } ( \overline R^n_{N , F} )}\, $ denotes a minimal $\overline U^n_F$-central extension of $\, \overline {\varphi } ( \overline R^n_{N , F} )\, $ by 
$\, K^J_2 ( \overline {\varphi } ( \overline R^n_{N , F} )\, ,\, \overline U^n_F )\, $ and the horizontal maps are only expected to be meaningful on commutators. Put 
$\, J' = C^r_{{\varphi }_n} \cap J_{R^n_{F , F}\, ,\, U^n_F} = 
R_{R^n_{F , F}\, ,\, U^n_F} + J_{R^n_{N , F}\, ,\, U^n_F}\, $. Assume given an element 
$\, \eta \in J_{C^r_{{\varphi }_n}\, ,\, U^{n+1}_F} \cap 
[\, U_{C^r_{{\varphi }_n}\, ,\, U^{n+1}_F}\, ,\, U^{n+2}_F\, ]\, $ representing an element in the kernel of 
\smallskip
$$ K^J_2\, \bigl( C^r_{{\varphi }_n}\, ,\, U^{n+1}_F \bigr) \longrightarrow 
K^J_2\, \bigl( R^n_{N , F}\, ,\, U^{n+1}_F / J_{R^n_{F , F}\, ,\, U^n_F} \bigr)\> . $$
\par\medskip\noindent 
Then modulo 
$\, [\, U_{J'\, ,\, U^{n+1}_F}\, ,\, U^{n+2}_F\, ]\, $ it is represented up to a trivial element which is of no importance since it can be lifted to a trivial element, by a commutator in 
$\, [\, U_{R^n_{N , F}}\, ,\, U_C\, ]\, $ with $\, C = J_{R^n_{F , F}\, ,\, U^n_F} / J' \subseteq 
[\, U^{n+1}_F / J'\, ,\, U^{n+1}_F / J'\, ]\, $. This element lifts to a double commutator in 
$\, [\, [\, U_{C^r_{{\varphi }_n}\, ,\, U^{n+1}_F}\, ,\, U^{n+2}_F\, ] , U^{n+2}_F\, ]\, $ modulo a trivial element, so that 
$\, \eta\, $ can be taken from  
$\, [\, U_{J'\, ,\, U^{n+1}_F}\, ,\, U^{n+2}_F\, ]\, +\, 
[\, [\, U_{C^r_{{\varphi }_n}\, ,\, U^{n+1}_F}\, ,\, U^{n+2}_F\, ] , U^{n+2}_F\, ]\, $. Then 
$\, [\, U_{J'\, ,\, U^{n+1}_F}\, ,\, U^{n+2}_F\, ]\, $ again decomposes modulo trivial elements as an expression in  
$\, [\, U_{R_{R^n_{F , F}\, ,\, U^n_F}\, ,\, U^{n+1}_F}\, ,\, U^{n+2}_F\, ]\, +\, 
[\, [\, U_{C^r_{{\varphi }_n}\, ,\, U^{n+1}_F}\, ,\, U^{n+2}_F\, ] , U^{n+2}_F\, ]\, $, so that altogether one gets an expression in 
\smallskip
$$ J_{C^r_{{\varphi }_n}\, ,\, U^{n+1}_F}\, \cap\, 
\bigl(\, \bigl[\, U_{R_{R^n_{F , F}\, ,\, U^n_F}\, ,\, U^{n+1}_F}\, ,\, U^{n+2}_F\, \bigr]\, +\, 
\bigl[\,\bigl[\, U_{C^r_{{\varphi }_n}\, ,\, U^{n+1}_F}\, ,\, U^{n+2}_F\,\bigr] , U^{n+2}_F\, \bigr]\,\bigr) \> . $$
\par\medskip\noindent
Choose a compatible section 
\smallskip
$$ \bigl(\, J_{R^n_{F , F}\, ,\, U^n_F}\, ,\, C_{{\varphi }_n}\, ,\, U^{n+1}_F\,\bigr)\,\nearrow\, 
\bigl(\, J_{\overline R^n_{F , F}\, ,\, \overline U^n_F}\, ,\, C_{\overline {\varphi }_n}\, ,\, U_{\overline U^n_F}\,\bigr) $$
\par\medskip\noindent
for the canonical map $\, U_{\overline U^n_F} \twoheadrightarrow U^{n+1}_F\, $ induced by 
$\, \overline U^n_F \twoheadrightarrow U^n_F\, $. Since the element 
$\, [ \eta ] \in K^J_2 ( C^r_{{\varphi }_n}\, ,\, U^{n+1}_F )\, $ is in the image of 
$\, K^J_2 ( J_{R^n_{F , F}\, ,\, U^n_F}\, ,\, U^{n+1}_F )\, $, the image of $\, \eta\, $ under the induced splitting 
\smallskip
$$ \bigl(\, U_{J_{R^n_{F , F} , U^n_F} , U^{n+1}_F}\, , U_{C_{{\varphi }_n} , U^{n+1}_F}\, , 
U^{n+2}_F\,\bigr) \longrightarrow
\bigl( U_{J_{\overline R^n_{F , F} , \overline U^n_F} , U_{\overline U^n_F}}\, , 
U_{C_{\overline {\varphi }_n} , U_{\overline U^n_F}}\, , U^2_{\overline U^n_F}\,\bigr) $$
\par\medskip\noindent
gives an element $\, \overline \eta\, $ of 
$\, [\, \overline {C_{\overline {\varphi }_n}}\, ,\, U_{\overline U^n_F}\, ]\, $ in the image of 
$\, [\, \overline {J_{\overline R^n_{F , F}\, ,\, \overline U^n_F}}\, ,\, U_{\overline U^n_F}\, ]\, $. One must show that the image of $\,\overline\eta\, $ by the composition 
\smallskip
$$ \overline {J_{\overline R^n_{F , F}\, ,\, \overline U^n_F}} \buildrel can \over\longrightarrow 
\overline {J_{R^n_{F , F}\, ,\, U^n_F}} \longrightarrow \overline {C_{{\varphi }_n}} \longrightarrow 
\overline {\overline {\varphi } ( \overline R^n_{N , F} )} $$
\par\medskip\noindent
is zero. It is clear that the composition 
\smallskip
$$ \overline {J_{\overline R^n_{F , F}\, ,\, \overline U^n_F}} \buildrel u^{-1}_{\alpha }\over\longrightarrow 
\overline {R_{R^n_{F , F}\, ,\, U^n_F}} \longrightarrow \overline {C^r_{{\varphi }_n}} \longrightarrow 
\overline {\overline {\varphi } ( \overline R^n_{N , F} )} $$
\par\medskip\noindent
is zero by halfexactness at $\, [\, \overline {C^r_{{\varphi }_n}}\, ,\, U^{n+1}_F\, ]\, $. On the other hand for any element in $\, \overline u_{\alpha } ( \overline {R_{R^n_{F , F}\, ,\, U^n_F}} )\, $ the map 
$\, \overline {can}\, $ coincides with $\, \overline u^{-1}_{\alpha }\, $, so that only elements in the kernel of 
$\, \overline u^{-1}_{\alpha }\, $ are of interest. We may then decompose $\,\overline\eta\, $ into a part coming from $\, \overline {R_{R^n_{F , F}\, ,\, U^n_F}}\, $ and a part in the kernel of 
$\, \overline u^{-1}_{\alpha }\, $. Modifying $\,\eta\, $ modulo the terms in 
$\, [\, U_{R_{R^n_{F , F}\, ,\, U^n_F}\, ,\, U^{n+1}_F}\, ,\, U^{n+2}_F\, ]\, $ one can write 
\smallskip
$$ \overline\eta \, =\, {\prod }_{\lambda , \mu }\> 
\bigl[\, \bigl[\, z_{\lambda , \mu }\, ,\, g^{\lambda }_{\mu }\, \bigr] , f_{\lambda } \bigr] $$
\par\medskip\noindent
where $\, z_{\lambda , \mu } \in C^r_{{\varphi }_n} \subseteq \overline {C^r_{{\varphi }_n}}\, $ by some chosen splitting, $\, f_{\lambda }\, ,\, g^{\lambda }_{\mu } \in U_{\overline U^n_F}\, $. Then 
\smallskip
$$ ( \overline u_{\alpha } \circ \overline u^{-1}_{\alpha } ) ( \overline\eta )\, =\, 
{\prod }_{\lambda , \mu }\> 
\bigl[\, \bigl[\, w_{\lambda , \mu }\, ,\, g^{\lambda }_{\mu }\, \bigr] ,  f_{\lambda }\, \bigr] $$
\par\medskip\noindent
with 
$\, w_{\lambda , \mu } \in ( \overline u_{\alpha } \circ \overline u^{-1}_{\alpha } )\, ( z_{\lambda , \mu } )\, $ and the difference 
$\, \overline\eta\, ((\overline u_{\alpha } \circ \overline u^{-1}_{\alpha } )\, ( \overline\eta ))^{-1}\, $ is congruent to 
\smallskip
$$ {\prod }_{\lambda , \mu }\> \bigl[\, \bigl[\, z_{\lambda , \mu }\, ,\, g^{\lambda }_{\mu }\,\bigr]\, 
\bigl[\, g^{\lambda }_{\mu }\, ,\, w_{\lambda , \mu }\,\bigr]\, ,\, f_{\lambda }\,\bigr] $$
\par\medskip\noindent
modulo the subgroup 
$\, [\, [\, \overline {J_{\overline R^n_{F , F}\, ,\, \overline U^n_F}}\, ,\, U_{\overline U^n_F}\, ] , 
\overline {C^r_{\overline {\varphi }_n}}\, ]\, $ which maps to zero by the composition
\smallskip
$$ \bigl[\, \overline {J_{\overline R^n_{F , F}\, ,\, \overline U^n_F}}\, ,\, U_{\overline U^n_F}\, \bigr] 
\buildrel \overline {can}\over\longrightarrow 
\bigl[\, \overline {J_{R^n_{F , F}\, ,\, U^n_F}}\, ,\, U^{n+1}_F\,\bigr] \longrightarrow 
\bigl[\, \overline {\overline {\varphi } ( \overline R^n_{N , F} )}\, ,\, \overline U^n_F\,\bigr]\> . $$
\par\medskip\noindent
Since $\, [\, z_{\lambda , \mu }\, ,\, g^{\lambda }_{\mu }\, ] \, 
[\, g^{\lambda }_{\mu }\, ,\, w_{\lambda , \mu }\, ]\, 
[\, g^{\lambda }_{\mu }\, ,\, z_{\lambda , \mu }\, w^{-1}_{\lambda , \mu }\, ]\,\in\, 
ker\, u^{-1}_{\alpha }\, \cap\, [\, C^r_{\overline {\varphi }_n}\, ,\, C^r_{\overline {\varphi }_n}\, ]\, =\, 
[\, ker\, u^{-1}_{\alpha }\, ,\, C^r_{\overline {\varphi }_n}\, ]\, $ 
as $\, u^{-1}_{\alpha } ( C^r_{\overline {\varphi }_n} )\, $ is free 
$\, {\prod }_{\lambda , \mu }\> [\, [\, z_{\lambda , \mu }\, ,\, g^{\lambda }_{\mu }\, ]\, 
[\, g^{\lambda }_{\mu }\, ,\, w_{\lambda , \mu }\, ]\, ,\, f_{\lambda }\, ]\, $ is congruent to 
\smallskip
$$ {\prod }_{\lambda , \mu }\> \bigl[\, \bigl[\, z_{\lambda , \mu }\, w^{-1}_{\lambda , \mu }\, ,\, 
g^{\lambda }_{\mu }\, \bigr] , 
f_{\lambda }\, \bigr] $$ 
\par\medskip\noindent
modulo 
$\, [\, [\, ker\, u^{-1}_{\alpha }\, ,\, \overline {C^r_{\overline {\varphi }_n}}\, ] , U_{\overline U^n_F}\, ]\, $ which maps to zero by the composition above. This completes the proof of halfexactness at 
$\, K_{n+3} ( E )\, $.
\par\noindent 
Next we consider the case of $\, K_{n+2} ( N , F)\, $. Let $\, ( C_{\psi }\, ,\, U^{n+1}_F )\, $ be the mapping cone of the inclusion $\, ( R_{R^{n-1}_{F , F}\, ,\, U^{n-1}_F}\, ,\, U^n_F ) \subseteq 
( C^r_{{\varphi }_{n-1}}\, ,\, U^n_F )\, $ if $\, n > 1\, $, respectively of $\, ( J_F\, ,\, U_F ) \subseteq 
( C_{\varphi }\, ,\, U_F )\, $ with $\, C_{\varphi }\, $ the mapping cone of $\, ( N , F ) \subseteq ( F , F )\, $ if $\, n = 1\, $. Also let $\,\widetilde C_{\psi }\, $ be the mapping cone of 
$\, ( J_{R^{n-1}_{F , F}\, ,\, U^{n-1}_F}\, ,\, U^n_F ) \subseteq ( C_{{\varphi }_{n-1}}\, ,\, U^n_F )\, $ for 
$\, n > 1\, $ and define the subgroups $\, C^r_{{\psi }_n}\, $ and $\, \widetilde C^r_{{\psi }_n}\, $ accordingly, so that $\, C^r_{{\psi }_n}\, $ is generated by 
$\, U_{R_{R^{n-1}_{F , F}\, ,\, U^{n-1}_F}\, ,\, U^n_F}\, $ and $\, R_{C^r_{{\varphi }_{n-1}}\, ,\, U^n_F}\, $, the latter being the kernel of an almost canonical map 
$\, U_{C^r_{{\varphi }_{n-1}}\, ,\, U^n_F} \rightarrow 
{\overline {C^r_{{\varphi }_{n-1}}}}^{ _{\sps U^n_F}}\, $ which may be obtained by restriction in the following way. Let $\, \overline C_{{\varphi }_{n-1}}\, $ denote the mapping cone of $\, ( J^{n-1}_{N , F}\, ,\, U^{n-1}_F ) \subseteq ( J^{n-1}_{F , F}\, ,\, U^{n-1}_F )\, $ and consider the split surjection 
$\,( \overline C_{{\varphi }_{n-1}}\, ,\, U^n_F ) \buildrel\leftarrow\over\twoheadrightarrow 
( J^n_{E , E}\, ,\, U^n_E )\, $ (all the while we are assuming that $\, n \geq 2\, $). The quotient 
$\, K^J_{n+2} ( E ) / K_{n +2} ( E )\, $ is free abelian, so that any almost canonical map 
$\, U_{J^n_{E , E}\, ,\, U^n_E} \rightarrow {\overline {J^n_{E , E}}}^{ _{\sps U^n_E}}\, $ can be extended to a map $\, U_{J^n_{E , E}\, ,\, U^n_E} \rightarrow {\overline {J^n_{E , E}}}^r\, $ where 
$\, {\overline {J^n_{E , E}}}^r\, $ is a minimal $U^n_E$-central extension of $\, J^n_{E , E}\, $ by 
$\, K_{n+2} ( E )\, $. Denote its kernel by $\, R^{n+1}_{E , E}\, $ ( by abuse of notation). Then the proof given in Proposition 1 extends to show that $\, K_{n+3} ( E ) = K^J_2 ( R^{n+1}_{E , E}\, ,\, U^{n+1}_E )\, $ and by taking the pull back of $\, {\overline {J^n_{E , E}}}^r\, $ with respect to 
$\, ( \overline C_{{\varphi }_{n-1}}\, ,\, U^n_F ) \twoheadrightarrow ( J^n_{E , E}\, ,\, U^n_E )\, $ one constructs a restriction 
$\, ( R_{\overline C_{{\varphi }_{n-1}}\, ,\, U^n_F}\, ,\, U^{n+1}_F ) \buildrel\leftarrow\over\twoheadrightarrow ( R^{n+1}_{E , E}\, ,\, U^{n+1}_E )\, $ having the same property that 
$\, K_{n+3} ( E ) = K^J_2 ( R_{\overline C_{{\varphi }_{n-1}}\, ,\, U^n_F}\, ,\, U^{n+1}_F )\, $ by regular isomorphism of $\, K^J_2 ( J_{\overline C_{{\varphi }_{n-1}}\, ,\, U^n_F}\, ,\, U^{n+1}_F )\, $ with 
$\, K^J_2 ( J^{n+1}_{E , E}\, ,\, U^{n+1}_E )\, $. Then again considering the inclusion 
$\, ( C^r_{{\varphi }_{n-1}}\, ,\, U^n_F ) \subseteq ( \overline C_{{\varphi }_{n-1}}\, ,\, U^n_F )\, $ one constructs a restriction of $\, R_{\overline C_{{\varphi }_{n-1}}\, ,\, U^n_F}\, $ to 
$\, R_{C^r_{{\varphi }_{n-1}}\, ,\, U^n_F}\, $. It is clear that the map 
\smallskip
$$ K^J_2 ( R_{C^r_{{\varphi }_{n-1}}\, ,\, U^n_F}\, ,\, U^{n+1}_F ) \rightarrowtail 
K^J_2 ( R_{\overline C_{{\varphi }_{n-1}}\, ,\, U^n_F}\, ,\, U^{n+1}_F ) $$
\par\medskip\noindent
is injective (because $\, K^J_2 ( J_{C^r_{{\varphi }_{n-1}}\, ,\, U^n_F}\, ,\, U^{n+1}_F ) \simeq 
K^J_2 ( J_{\overline C_{{\varphi }_{n-1}}\, ,\, U^n_F}\, ,\, U^{n+1}_F )\, $). It is also surjective since the image of $\, K^J_2 ( R^{n+1}_{0 , E , E}\, ,\, U^{n+1}_E )\, $ under the splitting 
$\, ( R^{n+1}_{E , E}\, ,\, U^{n+1}_E ) \rightarrow 
( R_{\overline C_{{\varphi }_{n-1}}\, ,\, U^n_F}\, ,\, U^{n+1}_F )\, $ factors over the group
$\, K^J_2 ( R_{C^r_{{\varphi }_{n-1}}\, ,\, U^n_F}\, ,\, U^{n+1}_F )\, $, so one gets that the latter is again equal to $\, K_{n+3} ( E )\, $. One also gets a regular isomorphism of 
$\, K^J_2 ( \widetilde C_{{\psi }_n}\, ,\, U^{n+1}_F )\, $ with 
$\, K^J_2 ( J_{R^n_{N , F}\, ,\, U^n_F}\, ,\, U^{n+1}_F )\, $. Now one proceeds in the same way as before. We leave it to the reader to check the single steps leading to the conclusion that the regular sequence is exact at all places $\, K_{n+2} ( N , F )\, $ (provided it is exact at all places $\, K_{n+2} ( F )\, $ which is next to prove).
\par\noindent
For $\, n \geq 2\, $ the above argument can be repeated once more by considering the mapping cones 
$\, ( C_{{\theta }_n}\, ,\, U^{n+2}_F )\, $ of the inclusion 
\smallskip
$$ \bigl( R_{C^r_{{\varphi }_{n-1}}\, ,\, U^n_F}\, ,\, U^{n+1}_F \bigr)\,\subseteq\, 
\bigl( C^r_{{\psi }_n}\, ,\, U^{n+1}_F \bigr) $$ 
\par\medskip\noindent
and $\, ( \widetilde C_{{\theta }_n}\, ,\, U^{n+2}_F )\, $ of
\smallskip 
$$ \bigl( J_{C^r_{{\varphi }_{n-1}}\, ,\, U^n_F}\, ,\, U^{n+1}_F \bigr)\,\subseteq\, 
\bigl( C_{{\psi }_n}\, ,\, U^{n+1}_F \bigr) $$ 
\par\medskip\noindent
respectively together with the corresponding subgroups $\, C^r_{{\theta }_n}\, $ and 
$\,\widetilde C^r_{{\theta }_n}\, $. Denote $\,\overline C_{{\psi }_n}\, $ the mapping cone of 
$\, ( J^n_{F , F}\, ,\, U^n_F ) \subseteq ( \overline C_{{\varphi }_{n-1}}\, ,\, U^n_F )\, $ so that 
$\, K^J_2 ( \overline C_{{\psi }_n}\, ,\, U^{n+1}_F ) {\simeq }_{reg} 
K^J_2 ( J^n_{N , F}\, ,\, U^n_F )\, $. As before one constructs extended suspensions 
$\, ( R_{J^n_{N , F}\, ,\, U^n_F}\, ,\, U^{n+1}_F )\, $ and 
$\, ( R_{\overline C_{{\psi }_n}\, ,\, U^{n+1}_F}\, ,\, U^{n+2}_F )\, $ having the same $K^J_2$-group isomorphic to $\, K_{n+2} ( N , F )\, $. Then one restricts 
$\,  R_{\overline C_{{\psi }_n}\, ,\, U^{n+1}_F}\, $ to $\, R_{C^r_{{\psi }_n}\, ,\, U^{n+1}_F}\, $ and shows that $\, K^J_2 ( R_{C^r_{{\psi }_n}\, ,\, U^{n+1}_F}\, ,\, U^{n+2}_F ) \simeq K_{n+2} ( N , F )\, $. One gets a regular isomorphism of $\, K^J_2 ( J_{R_{R^{n-1}_{F , F}\, ,\, U^{n-1}_F}\, ,\, U^n_F}\, ,\, U^{n+1}_F )\, $ with $\, K^J_2 ( \widetilde C_{{\theta }_n}\, ,\, U^{n+2}_F )\, $ and proceeds as before.
\par\noindent    
It remains to check the case of $\, K_3 ( F )\, $ where the above argument doesn't apply. Put 
$\, \widetilde J_{N , F} = R_F \cap J_{N , F}\, $ with $\, ( R_F\, ,\, U_F )\, $ a suspension of $\, F\, $ which may be restricted from $\, ( R_E\, ,\, U_E )\, $ by surjectivity of 
$\, K_2 ( F ) \twoheadrightarrow K_2 ( E )\, $. Then the quotient 
$\, ( J_F / \widetilde J_{N , F}\, ,\, U_F / \widetilde J_{N , F} )\, $ is equal to 
$\, ( J_E \times C\, ,\, \overline N \rtimes U_E )\, $ where $\,\overline N\, $ is the image of 
$\, {\overline N}^F\, $ in the universal extension $\,\overline F\, $ of $\, F\, $ and 
$\, C = ker \{ \overline N \twoheadrightarrow N \}\, $. Consider the extension 
\smallskip
$$ 1 \longrightarrow \bigl( \widetilde J_{N , F}\, ,\, U_F \bigr) \longrightarrow \bigl( R_F\, ,\, U_F \bigr) \longrightarrow \bigl( R_E\, ,\, \overline N \rtimes U_E \bigr) \longrightarrow 1\> . $$
\par\medskip\noindent
Then any element in the kernel of the map 
$\, K_3 ( F ) = K^J_2 ( R_F\, ,\, U_F ) \rightarrow K^J_2 ( R_E\, ,\, U_E ) = K_3 ( E )\, $ is trivial already in 
$\, K^J_2 ( R_E\, ,\, N \rtimes U_E )\, $ by halfexactness of 
\smallskip
$$ K^J_2\, \bigl( J_{N , F}\, ,\, U_F \bigr) \longrightarrow K^J_2\, \bigl( J_F\, ,\, U_F \bigr) \longrightarrow 
K^J_2\, \bigl( J_E\, ,\, U_E \bigr)\> . $$
\par\medskip\noindent
As a first step we want to show halfexactness of
\smallskip
$$ K^J_2\, \bigl( \widetilde J_{N , F}\, ,\, U_F \bigr) \longrightarrow K^J_2\, \bigl( R_F\, ,\, U_F \bigr) \longrightarrow K^J_2\, \bigl( R_E\, ,\, U_E \bigr)\> . $$
\par\medskip\noindent
Certainly an element in the kernel of the second map defines an element in the kernel of the regular surjection $\, K^J_2 ( R_E\, ,\, \overline N \rtimes U_E ) \twoheadrightarrow 
K^J_2 ( R_E\, ,\, N \rtimes U_E )\, $ which injects into 
$\, K^J_2 ( R_E \times C\, ,\, \overline N \rtimes U_E )\, $ and can be pulled back to an element of 
$\, K^J_2 ( J_{N , F}\, ,\, U_F )\, $ by the composition 
\smallskip
$$ K^J_2\, \bigl( J_{N , F}\, ,\, U_F \bigr) \longrightarrow K^J_2\, \bigl( C\, ,\, \overline N \rtimes U_E \bigr) \longrightarrow K^J_2\, \bigl( R_E \times C\, ,\, \overline N \rtimes U_E \bigr) \> . $$
\par\medskip\noindent
Then consider the extension 
\smallskip
$$ 1 \longrightarrow \bigl( R_E \times C\, ,\, \overline N \rtimes U_E \bigr) \longrightarrow 
\bigl( R_E \times \overline N\, ,\, \overline N \rtimes U_E \bigr) \longrightarrow 
\bigl( N\, ,\, {\overline F}^{ _E} \bigr) \longrightarrow 1 $$
\par\medskip\noindent
with $\, ( N , {\overline F}^{ _E} )\, $ (fully) exact. An element in the kernel of the map
$\, K^J_2 ( R_E \times C\, ,\, \overline N \rtimes U_E ) \rightarrow 
K^J_2 ( R_E \times \overline N\, ,\, \overline N \rtimes U_E )\, $ lifts to 
$\, K^J_3 ( N , {\overline F}^{ _E} ) \simeq K^J_3 ( \overline N , \overline F )\, $ which is the surjective image of $\, K^J_3 ( \overline N\, ,\, \overline N \rtimes U_E )\, $, so the boundary map 
$\, K^J_3 ( N , {\overline F}^{ _E} ) \buildrel\delta\over\longrightarrow 
K^J_2 ( R_E \times C\, ,\, \overline N \rtimes U_E )\, $ is trivial and  
\smallskip
$$ K^J_2\, \bigl( R_E \times C\, ,\, \overline N \rtimes U_E \bigr) \rightarrowtail 
K^J_2\, \bigl( R_E \times \overline N\, ,\, \overline N \rtimes U_E \bigr) $$
\par\medskip\noindent
is injective. But the map $\, K^J_2 ( J_{N , F}\, ,\, U_F ) \rightarrow 
K^J_2 ( R_E \times \overline N\, ,\, \overline N \rtimes U_E )\, $ factors over 
$\, K^J_2\, ( U_{N , F}\, ,\, U_F ) = 0\, $, proving halfexactness of 
\smallskip
$$ K^J_2\, \bigl( \widetilde J_{N , F}\, ,\, U_F \bigr) \longrightarrow K^J_2\, ( R_F\, ,\, U_F \bigr) \longrightarrow  K^J_2\, \bigl( R_E\, ,\, U_E \bigr) \> . $$
\par\medskip\noindent
Now let $\, C_{\varphi }\, $ denote the mapping cone of the inclusion $\, ( N , F ) \subseteq ( F , F )\, $, put 
$\, \widetilde J_{0 , N , F} = J_{0 , N , F} \cap \widetilde J_{N , F}\, ,\, 
C = \widetilde J_{0 , N , F} / R_{0 , N , F} \simeq \widetilde J_{N , F} / R_{N , F}\, ,\, 
C_F = J_{J_F , U_F} / R_{J_F , U_F}\, $. Also put $\,\widetilde J_{C_{\varphi }\, ,\, U_F} = 
J_{0 , C_{\varphi } , U_F} + J_{J_F , U_F}\, $ and $\, \widetilde R_{C_{\varphi } , U_F} = 
R_{0 , C_{\varphi } , U_F} + ( R_{C_{\varphi } , U_F} \cap J_{J_F , U_F} )\, $ where 
$\, R_{C_{\varphi } , U_F}\, $ is induced from $\, R_{J_F , U_F}\, $. Then the map 
\smallskip
$$ K^J_2\, \bigl( J_{0 , N , F}\, ,\, U_F \bigr) \twoheadrightarrow 
K^J_2\, \bigl( J_{N , F}\, ,\, U_F \bigr)\quad  $$ 
\par\medskip\noindent
is surjective and hence also the map 
\smallskip
$$ K^J_2\, \bigl( \widetilde J_{0 , N , F}\, ,\, U_F \bigr) \twoheadrightarrow 
K^J_2\, \bigl( \widetilde J_{N , F}\, ,\, U_F \bigr)\> . $$
\par\medskip\noindent
This is because dividing the extension
\smallskip
$$ 1 \longrightarrow \bigl( J_{0 , N , F}\, ,\, U_F \bigr) \longrightarrow \bigl( J_{N , F}\, ,\, U_F \bigr) \longrightarrow \bigl( D_N\, ,\, ( D_N \times N ) \rtimes U_E \bigr) \longrightarrow 1 $$
\par\medskip\noindent
by its intersection with $\, [\, [\, U_F\, ,\, U_F\, ] , U_F\, ]\, $ (i.e. regarding the image extension contained in  
$\, \overline U_F = U_F / [\, [\, U_F\, ,\, U_F\, ] , U_F\, ]\, $ one gets an isomorphism 
$\, K^J_2 ( D_N\, ,\, ( D_N \times N ) \rtimes U_E ) \simeq K^J_2 ( D_N\, ,\, D_N \rtimes \overline U_E )\, $ and the map $\, K^J_2 ( J_{N , F}\, ,\, U_F ) \rightarrow K^J_2 ( D_N\, ,\, D_N \rtimes \overline U_E )\, $ factors over $\, K^J_2 ( \overline J_{N , F}\, ,\, \overline U_F ) \simeq 
K^J_2 ( \overline U_{N , F}\, ,\, \overline U_F )\, $ since 
$\, \overline J_{N , F} \simeq U_{N , F} / [\, [\, U_{N , F}\, ,\, U_F\, ] , U_F\, ] = \overline U_{N , F}\, $ because 
$\, N\, $ is full. This implies that the map $\, K^J_2 ( J_{N , F}\, ,\, U_F ) \rightarrow 
K^J_2 ( D_N\, ,\, ( D_N \times N ) \rtimes U_E )\, $ is trivial and 
$\, K^J_2 ( J_{0 , N , F}\, ,\, U_F ) \twoheadrightarrow K^J_2 ( J_{N , F}\, ,\, U_F )\, $ is surjective. 
Then the second statement follows by comparing the extensions 
\smallskip
$$ \vbox{\halign{ #&#&#&#&#&#&#&#&#\cr
\hfil $1$\hfil & \hfil $\rightarrow $\hfil & \hfil $\bigl( \widetilde J_{0 , N , F}\, ,\, U_F \bigr)$\hfil & 
\hfil $\rightarrow $\hfil & \hfil $\bigl( \widetilde J_{N , F}\, ,\, U_F \bigr)$\hfil & \hfil $\rightarrow $\hfil & 
\hfil $\bigl( D_N\, ,\, U_F / \widetilde J_{0 , N , F} \bigr)$\hfil & \hfil $\rightarrow $\hfil & \hfil $1$\hfil \cr
&& \hfil $\Bigm\downarrow $\hfil && \hfil $\Bigm\downarrow $\hfil && \hfil $\Bigm\downarrow $\hfil &&\cr
\hfil $1$\hfil & \hfil $\rightarrow $\hfil & \hfil $\bigl( J_{0 , N , F}\, ,\, U_F \bigr)$\hfil & 
\hfil $\rightarrow $\hfil & \hfil $\bigl( J_{N , F}\, ,\, U_F \bigr)$\hfil & \hfil $\rightarrow $\hfil & 
\hfil $\bigl( D_N\, ,\, U_F / J_{0 , N , F} \bigr)$\hfil & \hfil $\rightarrow $\hfil & \hfil $1$\hfil \cr }} $$
\par\medskip\noindent
and since the map 
$\, K^J_2 ( D_N\, ,\, U_F / \widetilde J_{0 , N , F} ) \buildrel\sim\over\rightarrow 
K^J_2 ( D_N\, ,\, U_F / J_{0 , N , F} )\, $ is an isomorphism, so that 
$\, K^J_2 ( \widetilde J_{0 , N , F}\, ,\, U_F ) \twoheadrightarrow 
K^J_2 ( \widetilde J_{N , F}\, ,\, U_F )\, $ must be surjective. Then any element in  
$\, K^J_2 ( \widetilde J_{N , F}\, ,\, U_F ) \simeq K^J_2 ( J_{N , F}\, ,\, U_F )\, $ which is not in the image of 
$\, K^J_2 ( R_{N , F}\, ,\, U_F )\, $ lifts to an element of $\, K^J_2 ( \widetilde J_{0 , N , F}\, ,\, U_F )\, $ not in the image of $\, K^J_2 ( R_{0 , N , F}\, ,\, U_F )\, $. Consider the commutative diagram
\smallskip
$$ \vbox{\halign{ #&#&#&#&#\cr
\hfil $1$\hfil && \hfil $1$\hfil && \hfil $1$\hfil \cr
\hfil $\bigm\downarrow $\hfil && \hfil $\bigm\downarrow $\hfil && \hfil $\bigm\downarrow $\hfil \cr
\hfil $\bigl(\, R_{J_F , U_F}\, ,\, U^2_F \,\bigr)$\hfil & \hfil $\rightarrow $\hfil & 
\hfil $\bigl(\, \widetilde R_{C_{\varphi }\, ,\, U_F}\, ,\, U^2_F \,\bigr)$\hfil & \hfil $\rightarrow $\hfil & 
\hfil $\bigl(\, R_{0 , N , F}\, ,\, U_F \,\bigr)$\hfil \cr
\hfil $\bigm\downarrow $\hfil && \hfil $\bigm\downarrow $\hfil && \hfil $\bigm\downarrow $\hfil \cr
\hfil $\bigl(\, J_{J_F , U_F}\, ,\, U^2_F \,\bigr)$\hfil & \hfil $\rightarrow $\hfil & 
\hfil $\bigl(\, \widetilde J_{C_{\varphi }\, ,\, U_F}\, ,\, U^2_F \,\bigr)$\hfil & \hfil $\rightarrow $ \hfil & 
\hfil $\bigl(\, \widetilde J_{0 , N , F}\, ,\, U_F \,\bigr)$\hfil \cr
\hfil $\bigm\downarrow $\hfil && \hfil $\bigm\downarrow $\hfil && \hfil $\bigm\downarrow $\hfil \cr
\hfil $\bigl( C_F\, , {\overline J_F}^{ _{\sps U_F}}\negthinspace\times\negthinspace U_F \bigr)$\hfil & 
\hfil $\rightarrow $\hfil & 
\hfil $\bigl( C_E\, , (D_N \negthinspace\times\negthinspace {\overline C_{\varphi }}^{ _{\sps U_F}}) \negthinspace\rtimes\negthinspace U_E \bigr)$\hfil & 
\hfil $\rightarrow $\hfil & \hfil $\bigl( C\, , (D_N \negthinspace\times\negthinspace 
{\overline N}^{ _{\sps F}}) \negthinspace\rtimes\negthinspace U_E \bigr)$\hfil\cr
\hfil $\bigm\downarrow $\hfil && \hfil $\bigm\downarrow $\hfil && \hfil $\bigm\downarrow $\hfil \cr
\hfil $1$\hfil && \hfil $1$\hfil && \hfil $\> 1$ \hfil \cr }} . $$
\par\medskip\noindent
Any element of $\, K^J_2 ( \widetilde J_{0 , N , F}\, ,\, U_F )\, $ not in the image of 
$\, K^J_2 ( R_{0 , N , F}\, ,\, U_F )\, $ defines an element of 
$\, K^J_2 ( C\, ,\, (D_N \times {\overline N}^{ _{\sps F}} ) \rtimes U_E ) \simeq 
C \otimes ( D_N \times U^{ab}_E )\, $ and as the map $\, C_E \twoheadrightarrow C\, $ is surjective, so is  $\, K^J_2 ( C_E\, , (D_N \times {\overline C_{\varphi }}^{ _{\sps U_F}})\rtimes U_E )\, $ 
$\twoheadrightarrow K^J_2 ( C\, ,\, (D_N \times {\overline N}^{ _{\sps F}}) \rtimes U_E )\, $. Any lift of our element in the group
$\, K^J_2 ( C_E\, ,\, (D_N \times {\overline C_{\varphi }}^{ _{\sps U_F}} \rtimes U_E )\, $ then defines an element of $\, K_1 ( \widetilde R_{C_{\varphi } , U_F }\, ,\, U^2_F )\, $ trivial in 
$\, K_1 ( R_{0 , N , F}\, ,\, U_F ) \simeq 
K_1 ( R_{0 , N , F}\, ,\, U^2_F / (R_{C_{\varphi }\, ,\, U_F} \cap J_{J_F\, ,\, U_F}) )\, $, so it lifts again to 
$\, K_1 ( R_{C_{\varphi }\, ,\, U_F} \cap J_{J_F\, ,\, U_F}\, ,\, U^2_F ) \rightarrow 
K_1 ( J_{J_F\, ,\, U_F}\, ,\, U^2_F )\, $. We would like to construct equivariant lifts making an almost commutative diagram of the following sort
\smallskip
$$ \vbox{\halign{ #&#&#&#&#\cr
\hfil ${\overline {R_{J_F , U_F}}}^{ _{\sps U^2_F}}$\hfil & \hfil $\largerightarrow $\hfil & 
\hfil ${\overline {J_{J_F , U_F}}}^{ _{\sps U^2_F}}$\hfil & \hfil $\largerightarrow $\hfil &
\hfil $\overline C_F$\hfil \cr
\hfil $\Bigm\downarrow $\hfil && \hfil $\Bigm\downarrow $\hfil && \hfil $\Bigm\downarrow $\hfil \cr
\hfil ${\overline {\widetilde R_{C_{\varphi } , U_F}}}$\hfil & \hfil $\largerightarrow $\hfil & 
\hfil ${\overline {\widetilde J_{C_{\varphi } , U_F}}}$\hfil & \hfil $\largerightarrow $\hfil &
\hfil $\overline C_E$\hfil \cr }} $$
\par\medskip\noindent
of minimal $U_F$-($U^2_F$-)central extensions by the corresponding $K^J_2$-groups. Since 
$\, ( J_{J_F , U_F}\, ,\, U^2_F )\, $ is admissible the middle vertical map is easily obtained. Also by the isomorphism $\, K^J_2 ( R_F\, ,\, U_F ) \simeq K^J_2 ( J_F\, ,\, U_F )\, $ the map 
$\, {\overline {R_{J_F , U_F}}}^{ _{\sps U^2_F}} \rightarrow {\overline {J_{J_F , U_F}}}^{ _{\sps U^2_F}}\, $ 
is obtained from some inclusion $\, {\overline {R_{R_F , U_F}}}^{ _{\sps U^2_F}} \subseteq 
{\overline {J_{R_F , U_F}}}^{ _{\sps U^2_F}}\, $. The main problem is to construct the left lower map. In order to achieve this we may divide $\, {\overline {\widetilde J_{C_{\varphi }\, ,\, U_F}}}\, $ by the kernel of
\smallskip 
$$ K^J_2 ( \widetilde J_{C_{\varphi }\, ,\, U_F}\, ,\, U^2_F ) \twoheadrightarrow 
K^J_2 ( J_{C_{\varphi }\, ,\, U_F}\, ,\, U^2_F ) $$ 
\par\medskip\noindent
and correspondingly 
$\, {\overline {\widetilde R_{C_{\varphi }\, ,\, U_F}}} \, $ by the kernel of
\smallskip 
$$ K^J_2 ( \widetilde R_{C_{\varphi }\, ,\, U_F}\, ,\, U^2_F ) \twoheadrightarrow 
K^J_2 ( R_{C_{\varphi }\, ,\, U_F}\, ,\, U^2_F ) $$ 
\par\medskip\noindent
and $\, \overline C_E\, $ by the kernel of
\smallskip 
$$ K^J_2 ( C_E\, ,\, (D_N \times {\overline {C_{\varphi }}}^{ _{\sps U_F}}) \rtimes U_E ) \twoheadrightarrow K^J_2 ( C_E\, ,\, {\overline {C_{\varphi }}}^{ _{\sps U_F}} \rtimes U_E ) $$ 
\par\medskip\noindent
which surjects onto the kernel of 
\smallskip
$$ K^J_2 ( C\, ,\, (D_N \times {\overline N}^{ _{\sps F}}) \rtimes U_E ) \twoheadrightarrow 
K^J_2 ( C\, ,\, {\overline N}^{ _{\sps F}} \rtimes U_E ) \> . $$
\par\medskip\noindent
Then letting $\, {\overline {\widetilde J_{C_{\varphi }\, ,\, U_F}}}\, $ denote the preimage of 
$\, \widetilde J_{C_{\varphi }\, ,\, U_F}\, $ in $\,{\overline {J_{C_{\varphi }\, ,\, U_F}}}^{ _{\sps U^2_F}}\, $ and $\, {\overline {\widetilde R_{C_{\varphi }\, ,\, U_F}}}^{ _{\sps U^2_F}}\, $ the preimage of 
$\, \widetilde R_{C_{\varphi }\, ,\, U_F}\, $ in $\, {\overline {R_{C_{\varphi }\, ,\, U_F}}}^{ _{\sps U^2_F}}\, $ one gets an equivariant lift $\, {\overline {\widetilde R_{C_{\varphi }\, ,\, U_F}}}^{ _{\sps U^2_F}} \rightarrow {\overline {J_{C_{\varphi }\, ,\, U_F}}}^{ _{\sps U^2_F}}\, $, and the map 
$\, {\overline {R_{J_F\, ,\, U_F}}}^{ _{\sps U^2_F}} \rightarrow 
{\overline {\widetilde R_{C_{\varphi }\, ,\, U_F}}}^{ _{\sps U^2_F}}\, $ follows by projecting the image of the composition $\, {\overline {R_{J_F\, ,\, U_F}}}^{ _{\sps U^2_F}} \rightarrow 
{\overline {\widetilde J_{C_{\varphi }\, ,\, U_F}}}^{ _{\sps U^2_F}}\, $ to 
$\, {\overline {\widetilde R_{C_{\varphi }\, ,\, U_F}}}^{ _{\sps U^2_F}}\, $ using the corresponding projection of the preimage of $\, R_{C_{\varphi }\, ,\, U_F}\, $ in 
$\, {\overline {J_{C_{\varphi }\, ,\, U_F}}}^{ _{\sps U^2_F}}\, $ onto 
$\, {\overline {R_{C_{\varphi }\, ,\, U_F}}}^{ _{\sps U^2_F}}\, $. The groups $\,\overline C_F\, $ and 
$\, \overline C_E\, $ may then be defined as quotients 
$\, {\overline {J_{J_F\, ,\, U_F}}}^{ _{\sps U^2_F}} / {\overline {R_{J_F\, ,\, U_F}}}^{ _{\sps U^2_F}}\, $, 
$\, {\overline {\widetilde J_{C_{\varphi }\, ,\, U_F}}}^{ _{\sps U^2_F}} / 
{\overline {\widetilde R_{C_{\varphi }\, ,\, U_F}}}^{ _{\sps U^2_F}}\, $ respectively and the map 
$\, \overline C_F \rightarrow \overline C_E\, $ is only expected to be meaningful modulo the image of 
$\, K^J_2 ( J_{C_{\varphi }\, ,\, U_F}\, ,\, U^2_F )\, $ in $\, K^J_2 ( C_E\, ,\, U^2_F / J_{J_F\, ,\, U_F} )\, $. Note that the diagram of equivariant lifts as above is almost commutative in the sense that it is canonical on commutators and the image of an element of $\, {\overline {R_{J_F\, ,\, U_F}}}^{ _{\sps U^2_F}}\, $ in 
$\, {\overline {\widetilde J_{C_{\varphi }\, ,\, U_F}}}\, $ by either map can differ only modulo 
$\, K^J_2 ( J_{C_{\varphi }\, ,\, U_F}\, ,\, U^2_F )\, $, so dividing by this subgroup gives a commutative diagram. Returning to our element in $\, K_1 ( J_{J_F\, ,\, U_F}\, ,\, U^2_F )\, $ it lifts to an element of 
$\, {\overline {J_{J_F\, ,\, U_F}}}^{ _{\sps U^2_F}}\, $ well defined up to the subgroup 
$\, [\, {\overline {J_{J_F\, ,\, U_F}}}^{ _{\sps U^2_F}}\, ,\, U^2_F\, ]\, $. Taking the difference of the original lift in $\, K^J_2 ( C_E\, ,\, {\overline {C_{\varphi }}}^{ _{\sps U_F}} \rtimes U_E )\, $ with this element defines another lift of the element in $\, K^J_2 ( C\, ,\, {\overline N}^{ _{F}} \rtimes U_E )\, $ modulo the image of $\, K^J_2 ( J_{C_{\varphi }\, ,\, U_F}\, ,\, U^2_F )\, $. This lifts to 
$\, K^J_2 ( J_{C_{\varphi }\, ,\, U_F}\, ,\, U^2_F )\, $ since it is trivial in 
$\, K_1 ( R_{C_{\varphi }\, ,\, U_F}\, ,\, U^2_F )\, $, so that the original element of 
$\, K^J_2 ( \widetilde J_{N , F}\, ,\, U_F )\, $ must lift to $\, K^J_2 ( J_{C_{\varphi }\, ,\, U_F}\, ,\, U^2_F )\, $. This shows that the image of $\, K^J_2 ( \widetilde J_{N , F}\, ,\, U_F )\, $ in $\, K^J_2 ( R_F\, ,\, U_F )\, $ is the same as the image of $\, K^J_2 ( R_{N , F}\, ,\, U_F )\, $ completing the proof\qed

\end{document}